\documentclass[10pt,reqno]{amsart}
\usepackage[colorlinks=true,linkcolor=blue,citecolor=blue,hyperindex,breaklinks]{hyperref}

\usepackage{amsmath}
\usepackage{amsfonts}
\usepackage{mathrsfs}
\usepackage{amssymb}
\usepackage{amsthm,comment}
\usepackage{url}
\usepackage{enumitem}
\usepackage{bbm}
\usepackage{times}
  \usepackage{tikz}
\usetikzlibrary{matrix}

\usepackage{tabularx}
\newcolumntype{Y}{>{\centering\arraybackslash}X}
\usepackage{ctable,placeins}
\usepackage{caption}  
\usepackage{makecell}
\usepackage{float}

 \usepackage[abbrev,nobysame]{amsrefs}
\renewcommand{\MR}[1]{}

\usepackage[capitalize,sort]{cleveref}

\usepackage{crossreftools}
\theoremstyle{Theorem}
\newtheorem{theorem}{Theorem} [section]
\newtheorem{alttheorem}{Theorem} 

\newtheorem{proposition}[theorem]{Proposition} 
\newtheorem{claim}[theorem]{Claim} 
\newtheorem*{claim*}{Claim} 
\newtheorem{lemma}[theorem]{Lemma}

\newtheorem{corollary}[theorem]{Corollary}
\theoremstyle{definition}
\newtheorem{definition}[theorem]{Definition}
\newtheorem{question}[theorem]{Question}

\newtheorem{example}[theorem]{Example}
\newtheorem{remark}[theorem]{Remark}
\theoremstyle{remark}
\crefformat{section}{\S#2#1#3}
\renewcommand{\top}{\mathrm{top}}
\def\btm{\mathrm{bot}}
\def\what{\widehat}

\newcommand{\restrict}[2]{{#1}{\restriction_{{ #2}}}}
\DeclareMathOperator{\Hom}{Hom}
\DeclareMathOperator{\Sym}{Sym}

\newcommand{\Gl}{\mathrm{GL}}
\newcommand{\GL}{\mathrm{GL}} 
\newcommand{\Sl}{\mathrm{SL}}
\newcommand{\SU}{\mathrm{SU}}
\newcommand{\U}{\mathrm{U}}
\newcommand{\Sp}{\mathrm{Sp}}

\newcommand{\So}{\mathrm{SO}}
\newcommand{\SL}{\mathrm{SL}}

\newcommand{\hol}[1][\beta]{ \operatorname{H\ddot{o}l}^{{#1}}}
\def\Hol{\hol}
\newcommand{\Lip}{\, \mathrm{Lip}}

\def \gl{\mathfrak{gl}}
\def \liee{\mathfrak{e}}
\def \aff{\mathfrak{aff}}
\def \sl{\mathfrak{sl}}
\def \su{\mathfrak{su}}
\def \sp{\mathfrak{sp}}
\def \so{\mathfrak{so}}

\def\stab{\mathrm {stab}}

\def\Lie{\mathrm{Lie}}

\def\calZ{\mathcal Z}
\def\calU{\mathcal U}
\def\bfV{\mathbf V}

\def\bfT{\mathbf T}

\def\bfL{\mathbf L}

\DeclareMathOperator{\Ad}{Ad}
\DeclareMathOperator{\ad}{ad}
\newcommand{\id}{\mathrm{Id}}
\def\Id{\id}

\newcommand{\orb}{\mathcal O}
\newcommand{\sm}{\smallsetminus}
\newcommand{\R}{\mathbb {R}}
\newcommand{\C}{\mathbb {C}}
\newcommand{\Q}{\mathbb {Q}}
\newcommand{\Z}{\mathbb {Z}}
\newcommand{\N}{\mathbb {N}}
\newcommand{\T}{\mathbb {T}}

\newcommand{\inv}{^{-1}}

\newcommand{\diff}{\mathrm{Diff}}
\newcommand{\Homeo}{\mathrm{Homeo}}
\newcommand{\Diff}{\diff}

\newcommand{\Fol}{\mathscr{F}}

\def\fol{\Fol}
\def\1{\mathbf 1}
\newcommand{\supp}{\mathrm{supp}}
\newcommand{\td}{\tilde}
\newcommand{\wtd}{\widetilde}

\def\calI{\mathcal I}
\def\calB{\mathcal B}
\def\calE{\mathcal E}
\def\calA{\mathcal A}

\def\calC{\mathcal C}
\def\calL{\mathcal L}
\def\calF{\mathcal F}

\def\calG{\mathcal G}
\def\calP{\mathcal P}
\def\calV{\mathcal V}
\def\calH{\mathcal H}

\def\bfV{\mathbf V}

\def\bfL{\mathbf L}

\def\scrU{\mathscr U}
\def\scrH{\mathscr H}

\def\scrW{\mathscr W}
\def\scrG{\mathscr G}

\def \RP{\R P}
\def \C {\mathbb C}

\renewcommand{\emph}[1]{{\bf{#1}}}

\def\bs{\backslash}
\newcommand{\lieg}{\mathfrak g}

\newcommand{\lieh}{\mathfrak h}
\newcommand{\liek}{\mathfrak k}
\newcommand{\liev}{\mathfrak v}

\newcommand{\liem}{\mathfrak m}
\newcommand{\lien}{\mathfrak n}
\newcommand{\liea}{\mathfrak a}

\newcommand{\lieu}{\mathfrak u}

\newcommand{\liep}{\mathfrak p}
\newcommand{\lieq}{\mathfrak q}
\newcommand{\Aff}{\mathrm{Aff}}
\newcommand{\liesl}{\mathfrak {sl}}
\newcommand{\liesp}{\mathfrak {sp}}

\newcommand{\cent}{\mathscr Z}

\def\Holder{H\"older }

\renewcommand{\bar}{\overline}

\title[Boundary actions by higher-rank lattices]{Boundary actions by higher-rank lattices: Classification and embedding in low dimensions, local rigidity,  smooth factors} 

\author[A.~Brown]{Aaron Brown}
\address{Northwestern University, Evanston, IL 60208, USA}
\email{awb@northwestern.edu}

\author[{F.~Rodriguez Hertz}]{Federico Rodriguez Hertz}
\address{Pennsylvania State University, State College, PA 16802, USA}
\email{hertz@math.psu.edu}

\author[Z.~Wang]{Zhiren Wang}
\address{Pennsylvania State University, State College, PA 16802, USA}
\email{zhirenw@psu.edu}

\thanks{
This material is based upon work supported by the National Science Foundation under Grant Nos.\ DMS-2020013 (A.B.), DMS-1900778 (F.R.H), and  DMS-1753042 (Z.W.)}

\begin{document}
\maketitle

\begin{abstract}
We study actions by lattices in higher-rank (semi)simple Lie groups on compact manifolds.  By classifying certain measures invariant under a related higher-rank abelian action (the diagonal action on the suspension space) we deduce a number of new rigidity results related to standard projective actions  (i.e.\ boundary actions) by such groups.  

Specifically, in low dimensions we show all actions (with infinite image) are conjugate to boundary actions.  We also show standard boundary actions (e.g.\ projective actions on generalized flag varieties) are local rigid and classify all smooth actions that are topological factors of such actions.   Finally, for volume-preserving actions in low dimensions  (with infinite image) we provide a mechanism to detect the presence of ``blow-ups'' for the action by studying measures that are $P$-invariant but not $G$-invariant for the suspension action. 
\end{abstract}

\tableofcontents   

\section{Introduction and formulation of results for  lattices in $\Sl(n,\R)$}\label{sec1}
It is well known that (irreducible) lattice subgroups $\Gamma$ in  higher-rank (semi-)simple Lie groups $G$ exhibit significant rigidity properties with respect to linear representations $\rho\colon \Gamma\to \Gl(d,\R)$.  
The \emph{Zimmer program} refers to a number of questions and conjectures that, roughly, aim to establish analogous rigidity results for smooth (perhaps volume-preserving) actions.  Very roughly, such conjectures assert that all (volume-preserving) actions are built from certain modifications of  algebraic constructions or (in the non-volume-preserving case) have such algebraic actions as factors.  

As a starting point for the Zimmer program,  \emph{Zimmer's conjecture} asserts that for manifolds of sufficiently small dimension (relative to data associated to $\lieg= \Lie(G)$), all actions are isometric; appealing to Margulis's superrigidity theorem for representations into compact groups, one then concludes that many such actions have finite image.  For actions by lattices in $\Sl(n,\R)$ (as well as many other groups), this conjecture  was recently established (for $C^{1+\beta}$ actions) by the first author with D. Fisher and S. Hurtado.  

\begin{alttheorem}[{\cite[Theorem A]{2105.14541}}]\label{thm:BFHsl}
For $n\ge 3$, let $\Gamma$ be a lattice subgroup of $\Sl(n,\R)$.  Let $M$ be a closed manifold.   
\begin{enumerate}
	\item \label{1133} If $\dim(M)<n-1$,  then every homomorphism $\alpha\colon \Gamma\to \diff^{1+\text{\Holder}}(M)$ has finite image.
	\item \label{concbb}If $\dim(M)=n-1$ and if $\omega$ is a volume form (or no-where vanishing density) on $M$, then every homomorphism $\alpha\colon \Gamma\to \diff^{1+\text{\Holder}}(M;\omega)$ has finite image.
\end{enumerate}
\end{alttheorem}
In conclusion \eqref{concbb} above, $\diff^{1+\text{\Holder}}(M;\omega)$ denotes   the group of $C^{1+\text{\Holder}}$ diffeomorphisms preserving the volume form (or smooth, no-where vanishing density) $\omega$.  
We state \cref{thm:BFHsl} for actions of $C^{1+\text{\Holder}}$ diffeomorphisms, as all techniques developed in this paper require us to consider  actions that are at least $C^{1+\text{\Holder}}$;  we note, however,  that as formulated  \cref{thm:BFHsl}  also holds for $C^1$ actions; see  \cite{BDZ}.

\subsection{Motivating questions}
The  motivating problem for this paper was to consider the boundary cases of \cref{thm:BFHsl}.  Namely for $n\ge 3$ and for lattice subgroups $\Gamma$ in $\Sl(n,\R)$, we consider actions of $\Gamma$ on closed $(n-1)$-dimensional manifolds (the boundary case of  \cref{thm:BFHsl}\eqref{1133}) or volume-preserving actions on closed $n$-dimensional manifolds (the boundary case of  \cref{thm:BFHsl}\eqref{concbb}).  

\subsubsection*{Actions in dimension $n-1$}
For  volume-preserving actions on   closed $(n-1)$-dimensional manifolds,  Theorem \ref{thm:BFHsl} asserts that any such action has finite image. 
However, there exist effective algebraic actions in dimension $(n-1)$ that do not preserve any  volume form; namely, any lattice subgroup $\Gamma\subset \Sl(n,\R)$ acts projectively on $\R P^{n-1}$ and its double cover $S^{n-1}$.  
These actions do not  preserve any finite Borel measure (and, in particular, do not preserve any smooth density).  Conjecturally, (passing to finite-index subgroups) these are the only possible  actions on $(n-1)$-dimensional manifolds with infinite image.  
We establish this conjecture (for $C^{1+\text{\Holder}}$-actions) in this paper.   See \cref{slnbound} below.  

\subsubsection*{Volume-preserving actions in dimension $n$}
For volume-preserving actions on closed $n$-dimensional manifolds,  the model action is the  action (by affine group automorphisms) of $\Gamma = \Sl(n,\Z)$ on $\T^n$.  This action clearly preserves a volume-form, the Haar measure on $\T^n$.  Moreover, this   action exhibits the following additional property: for the induced $G$-action on the suspension space $M^\alpha$ (see \cref{sss:susp} below), every $P$-invariant Borel probability measure $\mu$ on $M^\alpha$ is, in fact, $G$-invariant.  Here, $P$ is the group of upper-triangular matrices, the minimal parabolic subgroup of $G= \Sl(n,\R)$.
One may blow-up the common fixed point for the action on $\T^n$ to obtain a (real-analytic) action $\alpha_0$ on a closed (non-orientable) $n$-manifold $M$.  
There is a smooth, $\Gamma$-equivariant map $\phi \colon M\to \T^n$ that is injective outside of the singular set.  The singular set $N= \phi\inv (0 +\Z^n)$ is an embedded copy of $\R P^{n-1}$ on which the action restricts to the standard projective action.  The action on  $M$ preserves a smooth density $\phi\inv (\mathrm{Haar}_{\T^n})$ which vanishes on the singular set $N$.  

 In \cite{MR1380646}, Katok and Lewis showed one could deform the smooth structure (in a neighborhood of the singular set) on $M$ to obtain a real-analytic action $\alpha$ by $\Sl(n,\Z)$ that preserves a no-where vanishing smooth density $\omega$.  For the induced $G$-action on the suspension space $M^\alpha$ (again, see \cref{sss:susp}), there exists an ergodic $P$-invariant Borel probability measure $\mu$ that is not $G$-invariant; this measure $\mu$ is supported on the suspension $N^\alpha$ of the singular set $N$. 

For actions  of lattice subgroups $\Gamma\subset \Sl(n,\R)$ 
on $n$-dimensional manifolds preserving a no-where vanishing density, one may ask if it is possible to detect the presence of such ``blow-ups.''  In this paper, we provide one such mechanism to detect such structures.  Namely, after passing to the induced action on the suspension space $M^\alpha$, the existence of a $P$-invariant Borel probability measure $\mu$ that is not $G$-invariant for the induced $G$-action implies the existence 
of an embedded, codimension-1 submanifold $N\subset M$ diffeomorphic to either $S^{n-1}$ or $\R P^{n-1}$ and on which a finite-index subgroup of $\Gamma$ restricts to the standard projective action.  
See \cref{thm:projembed} below.

\subsubsection*{Further results}
The techniques we developed to establish \cref{slnbound}  and  \cref{thm:projembed} below also allows us to establish new {\it local rigidity} and {\it smooth factor} theorems for standard projective actions.     
In the remainder of this introductory section, we formulate some of our main results in the setting of actions by lattices in $\SL(n, \R)$ (and related groups).  In \cref{sec2}, we formulate our most general results for actions by lattices in general (semi-)simple Lie groups.   

 \subsection*{Acknowledgements} Among many others, the authors benefited greatly from discussions with David Fisher, Sebastian Hurtado, Homin Lee, and Ralf Spatzier at various stages of this project.  

\subsection{Global classification in dimension of the smallest boundary}
\subsubsection{Actions by lattices in $\Sl(n,\R)$} For lattices subgroups of $\SL(n,\R)$, our first main result establishes the complete conjectured classification of actions in dimension $(n-1)$ that are at least $C^{1+\text{H\"older}}$ and have infinite image.  
\begin{theorem}\label{slnbound}
For $n\ge 3$, let $\Gamma$ be a lattice subgroup of $\Sl(n,\R)$. 
 Let $M$ be a connected compact manifold of dimension $n-1$.  Fix $r>1$  and let $\alpha\colon \Gamma\to \diff^{r}(M)$ be an action with infinite image $\alpha(\Gamma)$.  Then 
\begin{enumerate}
\item there is a $C^r$-diffeomorphism $h$ between $M$ and either $S^{n-1}$ or $\RP^{n-1}$ (equipped with their standard smooth structures) and 
\item a subgroup  $\Gamma'\subset\Gamma$ with 
 \begin{enumerate}
\item $\Gamma'=\Gamma$ if $M\simeq \RP^{n-1}$ or 
\item  $[\Gamma:\Gamma']\le 2$ if $M\simeq S^{n-1}$
\end{enumerate}
\end{enumerate}
such that
\begin{enumerate}[resume]
\item for all $x\in M$ and $\gamma\in \Gamma'$, $$h\left( \alpha (\gamma)(x) \right)= \gamma\cdot h(x)$$ where the right-hand side denotes the standard projective action of $\Gamma'$ on $S^{n-1}$ or $\RP^{n-1}$.  
\end{enumerate}
\end{theorem}
We note that \cref{slnbound}, in particular, implies that $\Gamma$ can not act $C^{1+\text{\Holder}}$ on exotic spheres.

We also prove analogues of \cref{slnbound} in the first critical dimension appearing in Zimmer's conjecture for lattices in many other higher-rank, $\R$-split simple Lie groups; see \cref{thm:globalrigid2} below.  

\subsubsection{Projective actions and their lifts}\label{sss:explicit} 
To clarify the conclusion of \cref{slnbound}, we present an explicit description of the (projective) actions that arise in the conclusion of \cref{slnbound}.
View  $G= \Sl(n,\R)$ acting linearly on $\R^n$ on the left. 
The group $G$ acts transitively on the spaces of lines $\R P^{n-1}$ and oriented lines $S^{n-1}$ in $\R^n$.  
The stabilizer of a line $\ell$  is a (maximal proper) parabolic subgroup $Q$ with two connected components corresponding to the two possible orientations on $\ell$.  
In particular, we   view  $\R P^{n-1}  = G/Q$ and $ S^{n-1} =  G/Q^\circ .$
The standard projective  $G$-action on $\R P^{n-1}$ is given by \begin{equation}\label{eq:rpac} 
	\rho(g)(\bar g Q) =  g \bar g Q . 
\end{equation}
Given any subgroup (in particular a lattice subgroup) $\Gamma\subset G$, we restrict the $G$-action 
above to $\Gamma$ to obtain the \emph{standard projective action} of $\Gamma$ on $S^{n-1}$ or $\R P^{n-1}$.

\cref{slnbound} shows (for $n\ge 3$)  that every action of $\Gamma$ on $\R P^{n-1}$ with infinite image is (smoothly conjugate to one) of the type described 
above.  
On the double cover $S^{n-1}$,  for certain lattices $\Gamma$, one may modify the standard action 
to  coincide with the action induced by the $G$-action only when restricted to an  index-$2$ subgroup: if $\Gamma$ admits a non-trivial  homomorphism $\sigma\colon \Gamma\to \Z/(2\Z) \simeq Q/ Q^\circ $, define a left $\Gamma$-action $\rho^\sigma\colon \Gamma\to \Diff(G/Q^\circ )$ on $S^{n-1}= G/Q^\circ$, 
\begin{equation}\label{eqtwist}\rho^\sigma (\gamma)(\bar gQ^\circ ) := \gamma \bar g  \sigma(\gamma\inv ) Q^\circ
=  \gamma \bar g  Q^\circ\sigma(\gamma\inv ) . \end{equation}
Alternatively, viewing  $\Z/(2\Z)= \{\pm \Id\}$ in $\Gl(n,\R)$ as the antipodal map, we have $\rho^\sigma (\gamma)(x) = (-1)^{\sigma(\gamma)} \gamma x.$

One  checks \eqref{eqtwist} defines a left $\Gamma$-action on  $G/Q^\circ$ that factors onto the standard action on $\R P^{n-1}$ and, when restricted to $\Gamma'=\ker \sigma$, $\rho^\sigma$ coincides with $\rho$ defined in \eqref{eq:rpac}.   \cref{slnbound} shows all actions on $S^{n-1}$ with infinite image  are of the form \eqref{eqtwist} for some choice of $\sigma$.  See discussion in \cref{rem:classifyactions} and \cref{ex:classifyactions} for more general modifications of lifts of projective actions.

\subsubsection{Actions by lattices in $\Sl(n,\C)$}
For a lattice $\Gamma$ in $\Sl(n,\C)$, it is shown in \cite{ABZ}  that if $\dim(M)<2n-2$, then the image of any $\alpha\colon \Gamma\to \diff^{1+\text{\Holder}}(M)$ is finite.  
Techniques in this paper yield a  global classification, similar to \cref{slnbound}, classifying all actions on $(2n-2)$-dimensional (real) manifolds (with infinite image). 
 The model actions in dimension  $(2n-2)$ occur on the manifold   $$\C P^{n-1}= \SL(n,\C)/Q$$ where $Q$ is the parabolic subgroup stabilizing a complex line.  
We note   two key differences when working in $\Sl(n,\C)$ rather than $\Sl(n,\R)$: First,  $\C P^{n-1}$ is simply connected and so admits no  lifted actions as discussed in \cref{sss:explicit}.  More substantially, both  $\Sl(n,\C)$ and $\SU(n)$ act transitively on $$\C P^{n-1}= \SL(n,\C)/Q =\SU(n)/S(\U(1)\times \U(n-1))$$ and certain (cocompact) lattices $\Gamma$ in $\Sl(n,\C)$ admit infinite-image homomorphisms into $\SU(n)$.  Thus, there are two model actions in dimension $(2n-2)$ by lattices $\Gamma$ in $\Sl(n,\C)$: the projective action on $\SL(n,\C)/Q$ or the isometric action obtained from infinite-image homomorphisms $\Gamma\to \SU(n)$.

\begin{theorem}\label{slnCbound}
For $n\ge 3$, let $\Gamma$ be a lattice subgroup of $\Sl(n,\C)$. 
 Let $M$ be a connected compact manifold of (real) dimension $2n-2$.  


 Let $\alpha\colon \Gamma\to \diff^{\infty}(M)$ be an action with infinite image $\alpha(\Gamma)$.  Then, there exists a $C^\infty$ diffeomorphism $h\colon M\to \C P^{n-1}= \SL(n,\C)/Q=\SU(n)/\mathrm{S}(\U(1)\times \U(n-1))$  such that one of the following holds: 
\begin{enumerate}
\item  $h\left( \alpha (\gamma)(x) \right)= \gamma \cdot h(x)$ \, for all $x\in M$ and $\gamma\in \Gamma$, or 
\item\label{222poop} there is  a subgroup $\Gamma'\subset \Gamma$ of index at most 2  and a homomorphism $\phi\colon \Gamma' \to \SU(n)$ such that $h\left( \alpha (\gamma)(x) \right)= \phi(\gamma)\cdot h(x) $ for all $x\in M$ and $\gamma\in \Gamma'$. 
\end{enumerate}

In particular, all actions  with infinite image $\alpha(\Gamma)$ are $C^r$-conjugate to either the standard projective action or embed into the isometric action on $\C P^{n-1}$.
\end{theorem}

We note that the isometry group of $\C P^{n-1}$ (equipped with the Fubini-Study metric) contains $\SU(n)$ (the orientation-preserving isometries) as an index-2 subgroup; this is the only reason to  restrict  to an index-2 subgroup in conclusion \eqref{222poop} of  \cref{slnCbound}.  

 \cref{slnCbound} follows directly from \cref{thm:globalrigid2nonsplit} below, combined with a standard argument using Margulis's superrigidity and arithmeticity theorems (see discussion in \cite[\S 2.3]{MR4502593}) to classify compact groups admitting infinite image representations of higher-rank lattices; we outline this argument in \cref{isometric} below.

\begin{remark}[Isometric actions of higher-rank lattices]\label{isometric}
Isometric actions by higher-rank lattices are essentially classified, especially in low dimension.  Indeed, if $g$ is a $C^r$ metric tensor on $M$, $r\ge 0$, then the isometry group $K=\mathrm{Isom}(M, g)$ is a compact Lie group and any isometry $f\in K$ is $C^{r+1}$ (see \cite{MR1503467,MR2204038}).  In particular, the natural $K$-action on $M$ is a $C^{r+1}$ action.

A standard argument using Margulis's superrigidity and arithmeticity theorems (see for example discussion in \cite[\S 2.3]{MR4502593})  classifies  compact groups $K$ admitting infinite image representations of higher-rank lattices. Specifically, the Lie algebra of $K$ consists of simple factors that are compact real forms of the complexification of $\lieg$.  
For instance, for $n\ge 3$, when  either $G=\Sl(n,\R)$ or $G=\Sl(n,\C)$, the Lie algebra of any compact group $K$ admitting a representation with dense image by a lattice $\Gamma\subset G$ is the direct sum of copies of  $\mathfrak {su}(n)$
 In particular $K$ necessarily contains a subgroup isogenous to the unitary group $\SU(n)$.  

The number $d(\lieg)$ in \cref{sec:table}, corresponds to the minimal dimension of all proper coset spaces $K/C$ where  $K$ is the compact real form of the complexification of $\lieg$.  
Given an action by a lattice $\Gamma$ on a connected manifold $M$, if $$d(\lieg) > \dim M,$$ all isometric actions  of $\Gamma$ on $M$ are finite.   
In the case that $d(\lieg)=\dim (M)$, then  either $K$ is finite or acts transitively on $M$; in the latter case, let $S_0\subset K$ stabilize a point $x_0\in M$.  Then map $K/S_0\to M$, $kS_0\mapsto k\cdot x_0$, is a $C^{r+1}$ diffeomorphism between $M$ and $K/S_0$.  In this way, we obtain the second conclusion of \cref{slnCbound}.
\end{remark}

\begin{remark}[Reducing regularity in \cref{slnCbound}]\label{rem:lowerregcomplex}
When $n\ge 4$, one can reduce the regularity assumptions in \cref{slnCbound}.  This is explained in \cref{rem:reduce regularity} (which allows us to reduce the regularity assumptions in \cref{thm:globalrigid2nonsplit} below).  

More precisely,  for $n\ge 4$ and $r>1$, if $\alpha\colon \Gamma\to \diff^{r}(M)$  fails to have subexponential growth of derivatives  (see \cref {see pf main}) then $\alpha$ is $C^r$ conjugate to a projective action as in the first conclusion of \cref{slnCbound}. If $\alpha\colon \Gamma\to \diff^{r}(M)$ has subexponential growth of derivatives and if $r\ge 2$, strong property (T) implies  $\alpha(\Gamma)$ acts by isometries of a metric that is $C^{r-1-\delta}$  \cite[Theorem 2.4]{MR4502593}; if $\alpha(\Gamma)$ is infinite, we then obtain a a $C^{r-\delta}$ conjugacy to the isometric action as in the  second conclusion of \cref{slnCbound}. When $1<r<2$ and  $\alpha\colon \Gamma\to \diff^{r}(M)$ has subexponential growth of derivatives, $\alpha(\Gamma)$ should preserve a metric that is still $C^{\text {\Holder}}$ and we still obtain a $C^{1+ \text{\Holder}}$ conjugacy to the isometric action when $\alpha(\Gamma)$ is infinite.  
\end{remark}

\def\vol{\mathrm{vol}}
\subsection{Embedded projective actions in low-dimensional, volume-preserving actions}
Having classified actions by lattices $\Gamma$ in $\Sl(n,\R)$   on $(n-1)$-dimensional manifolds, we turn to the question of classifying volume-preserving actions of lattices in $\Sl(n,\R)$ acting on $n$-dimensional manifolds.  
A complete classification of such actions seems rather out of reach.  A reasonable conjecture is that all such non-isometric actions are built from pieces that are ``blow-ups'' of the affine action on $\T^n$.  Our techniques provide one mechanism to detect the existence of such``blow-ups" (by studying $P$-measures on the suspension space $M^\alpha$.)
To state our theorem, we postpone the definition of the  induced $G$-action on the suspension space $M^\alpha$ until  \cref{sss:susp}.

\begin{theorem}\label{thm:projembed}
	For $n\ge 3$, let $\Gamma$ be a lattice in $\Sl(n,\R)$.  Let $M$ be a closed, $n$-dimensional manifold, let $\omega$ be a smooth volume form or no-where vanishing density,  and let $\alpha\colon \Gamma\to {\Diff^{r}}(M; \omega)$ be an action with $r>1.$
	
	Let $M^\alpha$ be the induced $G$-space.  Suppose that for the induced $G$-action on $M^\alpha$  there exists an ergodic, $P$-invariant Borel probability measure $\mu$  that is not $G$-invariant.

Then there are finitely many embedded $C^r$ submanifolds $N_1,\dots,N_\ell$ in $M$ such that
\begin{enumerate}
	\item each $N_i$ is diffeomorphic to either $\R P ^{n-1}$ or $S^{n-1}$,
	\item there is a finite index subgroup $\Gamma'\subset \Gamma$ such that each $N_i$ is $\Gamma'$ invariant, and 
	\item the dynamics of $\Gamma'$ on each $N_i$ is $C^r$ conjugate to standard projective action.
\end{enumerate}
Moreover, 
\begin{enumerate}[resume]
\item the measure $\mu$ is supported on the suspension $(G\times (N_1\cup \dots \cup N_\ell))/\Gamma$ in $M^\alpha$ of $N_1\cup \dots \cup N_\ell $.
\end{enumerate}
\end{theorem}

We note that the blow-up procedure and Katok-Lewis deformations changes the topology of the manifold $\T^n$ admitting the affine action.  It would be interesting to understand to what extent the topology constrains the action. 
\begin{question} \label{q:globalrig}
For $n\ge 3$, let $\Gamma$ be a lattice subgroup in $\Sl(n,\R)$.  
Do all volume-preserving actions  $\alpha\colon \Gamma\to \Diff^\infty(\T^n)$ with infinite image  $\alpha(\Gamma)$ contain an Anosov element?  In particular, are all such actions smoothly conjugate to an affine action on $\calF$?  
\end{question}
 We note that actions by higher-rank lattices on tori (and nilmanifolds) with  Anosov elements are more-or-less  classified by \cite{MR3702679}.

\subsection{Local rigidity of projective actions on flag varieties}. 
We now consider the  standard projective action $\rho$ of a lattice subgroup $\Gamma\subset \Sl(n,\R)$ on $\R P^{n-1}$, the projective action of $\Gamma$ on any Grassmannian variety (of subspaces) in $\R^n$, or, in most generality, the projective action of $\Gamma$ on a  flag variety (of flags in $\R^n$).

When $\Gamma$ is a cocompact lattice, in \cite{MR1421873} Kanai showed the projective action on $\R P^{n-1}$ (for $n\ge 21$) is \emph{locally rigid}: if an action $\alpha\colon \Gamma\to \Diff^\infty(\R P^{n-1})$ is sufficiently close (in the $C^4$-topology)  to $\rho$ then $\alpha$ is smoothly conjugate to $\rho$.  Later, in \cite[Theorem 17]{MR1632177},  Katok and Spatzier established local rigidity  (for $C^\infty$ actions) of projective actions of cocompact $\Gamma\subset \Sl(n,\R)$ (for $n\ge 3$) on 
Grassmannians and general flag manifolds as well as projective actions of lattices in more general higher-rank Lie groups on generalized flag manifolds.    Kanai's proof uses cohomological methods.  The  proof by  Katok and Spatzier  uses as a starting point   structural stability of normally hyperbolic foliations.  In both approaches, it seems  essential that the lattice   $\Gamma$ be cocompact.   We also remark the recent work \cite{2303.00543} which establishes $C^0$-local rigidity results for boundary actions; again the techniques in \cite{2303.00543} seem to require the lattice $\Gamma$ be assumed cocompact in $G$.

The tools we employ to establish \cref{slnbound} 
recover and substantially extend  the local rigidity results of Kanai and of Katok-Spatzier.


\begin{theorem}\label{thm:slnlocrig}
For $n\ge 3$, let $\calF$ be a flag variety (of flags in $\R^n$) and let $\Gamma\subset \Sl(n,\R)$ be a lattice subgroup.  Then there is $r_0\ge 1$ (depending only on $\calF$) such that for all $r>r_0$, the projective action $\rho\colon \Gamma\to\Diff^r(\calF)$ is $C^{r, 0, r}$ locally rigid.
\end{theorem} 

We emphasize that our approach is completely different and allows us to obtain results of actions by non-uniform lattices.   We further note that we obtain results for actions of lower regularity and for perturbations in weaker topologies.  For instance, for the projective action of  $\Gamma\subset \Sl(n,\R)$  on  $\R P^{n-1}$ we have $r_0 =1$.  In this setting, \cref{thm:slnlocrig} says for any $r>1$ and any action $\alpha\colon \Gamma\to\Diff^r(\R P^{n-1})$ sufficiently $C^0$ close to the standard projective action $\rho$, there exists a $C^r$ diffeomorphism $h\colon \R P^{n-1}\to \R P^{n-1}$  such that $h\circ \alpha (\gamma) = \rho(\gamma)\circ h$ for all $\gamma\in \Gamma$.

\cref{thm:slnlocrig,slnbound} together  raise the following {\it global rigidity} question.  
\begin{question}[Global rigidity of actions on boundaries] \label{q:globalrig}
For $n\ge 3$, let $\Gamma$ be a lattice subgroup in $\Sl(n,\R)$.  Let $\calF$ be a flag manifold (of flags in $\R^n$).  
Let $\alpha\colon \Gamma\to \Diff^\infty(\calF)$ be an action with infinite image  $\alpha(\Gamma)$.
Is $\alpha$ smoothly conjugate to the standard projective action on $\calF$?  
\end{question}
Although we do not pursue this question in this paper, we expect the methods we develop may help to answer \cref{q:globalrig}.

\subsection{Smooth factors of projective actions on flag manifolds}
The final class of problems we consider is classification of smooth factor actions. For $n\ge 3$, let $\calF$ be a flag manifold of flags in $\R^n$.  The group $G=\Sl(n,\R)$ acts transitively on $\calF$.  The subgroup $P\subset \Sl(n,\R)$ of upper triangular matrices stabilizes a complete flag in $\R^n$ and 
 the stabilizer of any flag in $\calF$ is conjugate to a parabolic subgroup $Q$ containing $P$.  In particular, we identify $\calF \simeq G/Q$ for some subgroup $P\subset Q\subset G$.  
 
Let $\Gamma\subset \Sl(n,\R)$ be a lattice subgroup.   We consider smooth quotients of the projective action of $\Gamma$ on $\calF$ and show all such  quotients are algebraically natural.  This generalizes the main result of \cite{MR3748688}, which established the smooth factor theorem under only the additional assumption that the factor action $\alpha$ admits a differentiable sink.   

\begin{theorem}\label{thm:slnrsmoothfactor}
Suppose $n\ge 3$.  Let $Q\subset G$ be a parabolic subgroup of $G= \Sl(n,\R)$ and let $\calF= G/Q$ be the associated flag manifold.  
 There exists $r_0 = r_0(\calF)\ge 1$ such that the following holds for any $r> r_0$: Let $M$ be a smooth manifold equipped with a $C^r$ action of $\Gamma$ and let $p\colon \calF\to M$ be a $C^0$,  $\Gamma$-equivariant, surjection.

Then there exists a subgroup $Q\subset Q'\subset G$ and a $\Gamma$-equiviariant $C^r$ diffeomorphism $h\colon G/Q'\to M$ such that, if $\pi\colon G/Q\to G/Q'$ is the natural map, then  for all $\gamma\in \Gamma$ and $gQ\in G/Q$, $$h\circ \pi (\gamma \cdot gQ) = \gamma\cdot p(gQ).$$  In particular, the following diagram $\Gamma$-equivariantly commutes: 
\begin{center}\begin{tikzpicture}
  \matrix (m) [matrix of math nodes,row sep=3em,column sep=4em,minimum width=2em]
  {
G/Q & \phantom{O}  \\
    G/Q' & M \\};
    \path[-stealth]
    (m-1-1) edge node [left] {$ \pi$} (m-2-1)
            edge   node  [above] {$p$} (m-2-2)            
    (m-2-1) edge node [below] {$h$}  (m-2-2);
\end{tikzpicture}
\end{center}
\end{theorem}
See \cref{thm:localrigid} below for more general results for projective actions on  generalized flag varieties.

\section{Statement of  results for general Lie groups} \label{sec2}
We formulate our most general results for actions by lattices in higher-rank semisimple Lie groups.  Proofs of all results appear in \cref{sec6}.  

\subsection{Standing assumptions on $G$}\label{ss:Gassump}
Throughout, we fix $G$ to be a connected semisimple Lie group whose Lie algebra $\lieg$ has real rank at least 2.  For all results, it is with no loss of generality (see discussion in \cite[Section 2.1]{2105.14541}) to further assume the following:
\begin{enumerate}
\item $G$ has no compact almost simple factors; 
\item $G$ is (topologically) simply connected (so may have infinite center).  
\end{enumerate}
We fix such a $G$ for the remainder and  take $\Gamma\subset G$ to be an irreducible lattice subgroup.  

\subsection{Terminology from Lie theory and fixed structures on $\lieg$}
Let $G$ be as in \cref{ss:Gassump} and let $\lieg$ denote the Lie algebra of $G$.  

\subsubsection{Iwasawa decomposition and parabolic subgroups}
Let $\lieg = \liek \oplus \liea\oplus \lien$ be a choice of Iwasawa decomposition. 
If $K,A$ and $N$ are the analytic subgroups corresponding to $\liek, \liea,$ and $ \lien$, respectively, then $G = KAN$ is the corresponding {Iwasawa decomposition of $G$.} 
Here $A$ is a maximal connected abelian $\ad$-$\R$-diagonalizable subgroup, and $N$ is a maximal $\ad$-unipotent subgroup of $G$ normalized by~$A$. We call $A$ a \emph{split Cartan subgroup} of~$G$.
 The subgroup $K$~contains the center of $G$ and is a maximal compact subgroup of~$G$ if (and only if) the center of $G$ is finite.

Write $M = C_K(A)$ for the centralizer of~$A$ in~$K$. The subgroup $P = MAN$ is a \emph{minimal parabolic subgroup} of~$G$. 
Any subgroup of~$G$ that contains $P$ is a \emph{(standard) parabolic subgroup}.

Let $\liem$ denote the Lie algebra of $M$.  
Relative to our choice of Iwasawa decomposition, a \emph{standard parabolic subalgebra} $\lieq\subset \lieg$ is any subalgebra containing $\liem \oplus \liea \oplus \lien.$   A  \emph{parabolic subalgebra} $\lieq$ is any subalgebra  $\lieq\subset \lieg$ that is conjugate to a standard parabolic subalgebra.

Let $\Sigma$ denote the set of restricted roots of $\lieg$ relative to the $\R$-split subalgebra $\liea$.  We fix an order on the roots $\Sigma$ so that the positive roots $\Sigma _+\subset \Sigma$ correspond to $\lien$: 
$$\lien:= \bigoplus _{\beta\in \Sigma_+} \lieg^\beta.$$   Relative to this order, let $\Pi \subset \Sigma_+$ denote the associated collection of simple positive roots.   

A standard parabolic subalgebra is fully saturated by root spaces and has the following description: given a subset $S\subset \Pi$, define $\lieq_S$ to be the parabolic subalgebra
$$\lieq_S = \liem \oplus \liea \oplus \lien  \oplus\bigoplus _{\substack{\alpha_i\in \Pi\sm S\\ k_i \le 0}} \lieg^{\sum k_i \alpha_i}.$$
All standard parabolic subalgebras are of this form.  
Given a standard parabolic subalgebra $\lieq$, let $\Sigma_\lieq\subset \Sigma$ denote the set of roots with $\lieg^\beta\subset \lieq$; let 
$\Sigma_\lieq^\perp\subset \Sigma$ denote the set of roots with $\lieg^\beta\cap \lieq= \{0\}$.   If $\lieq=\lieq_S$, then $\Sigma_\lieq^\perp$ contains all negative roots of the form $$\beta = \sum_{\substack{\alpha_i\in \Pi\\c_i\le 0}}c_i\alpha_i$$
where $c_i\neq 0$ for at least one $\alpha_i\in S$.  It follows that $$\liev_\lieq := \bigoplus _{\beta\in \Sigma_\lieq ^\perp} \lieg^\beta$$ is a nilpotent Lie subalgebra.  We remark that $\lieg = \liev_\lieq\oplus \lieq$ and $\lieq\cap \liev_\lieq = \{0\}$.  

Given a parabolic subalgebra $\lieq \subset \lieg$, let $Q\subset G$ be the normalizer of $\lieq$.  
Then $Q$ has Lie algebra $\lieq$ and  $Q$ is  a  {parabolic subgroup}.  
The quotient $G/Q$ is then a projective variety.  
We write $V_\lieq$ for the analytic (i.e.\ connected) Lie subgroup associated to $\liev_\lieq$.  The pointwise product  $V_\lieq\cdot  Q$ is a dense open subset of $G$.  

\subsection{Numerology associated to Lie algebras}
To state our main results, we define various numbers associated to a Lie algebra $\lieg$ and its parabolic subalgebras $\lieq$.  
\subsubsection{Critical regularity $r_0(\lieq)$} 
We associate to each parabolic subalgebra a positive real number $r_0(\lieq)$ which will be a lower bound for the critical regularity in various arguments below.  
To define $r_0(\lieq)$, let $C(\Sigma_\lieq^\perp) \subset \liea$ be the \emph{negative cone} of the collection of roots $\Sigma_\lieq^\perp$:
\begin{equation} \label{eq:negcone} C(\Sigma_\lieq^\perp) = \{ s\in \liea : \beta (s)<0 \text{ for all $\beta\in \Sigma_\lieq^\perp$} \}.\end{equation}
Given a standard parabolic subalgebra $\lieq$, define the critical regularity associated to $\lieq$ to be 
\begin{equation}\label{eq:parareg}
r_0(\lieq):= \inf _{s\in C(\Sigma_\lieq^\perp) } \max \left \{ \frac {\beta(s)}{\beta'(s)} : \beta, \beta'\in \Sigma_\lieq^\perp\right\}.
\end{equation}
\begin{example}[{Computation of $r_0(\lieq)$ for various classical parabolic subalgebras}]
Let $\lieg= \sl(n+1,\R)$ with Dynkin diagram $A_{n}$ and simple positive roots $\{\alpha_1,\dots, \alpha_{n}\}$ enumerated from left to right.  
Given $S\subset \{\alpha_1,\dots, \alpha_{n}\}$, we have 
$r_0(\lieq_S)= \mathrm{card}(S)$. 
\end{example}

\subsubsection{Critical dimension $v(\lieg)$} Let $\lieg$ be a semisimple Lie algebra.  We write 
$$v(\lieg):= \min \{ \dim (\lieg/\lieq) : \lieq\neq \lieg\}$$
where the quotient is as vector spaces and the minimum is taken over all proper (maximal) parabolic subalgebras $\lieq$.  That is, $v(\lieg)$ is the minimal codimension of all proper parabolic subalgebras $\lieq\subset \lieg$.  
If $G$ is a connected Lie group with Lie algebra $\lieg$, also write $v(G) = v(\lieg)$.  

Computations of $v(\lieg)$ for various $\lieg$ appear in \cref{sec:table}.

\subsubsection{Critical regularity  $r_{\min}(\lieg)$ of the minimal projective action} 
To state the minimal regularity appearing in our global classification results, given a semisimple Lie algebra $\lieg$, 
write  $${r_{\min}(\lieg) }= \max\{ r_0(\lieq')\} $$ where the maximum is taken over all proper standard  parabolic subalgebras $\lieq'$  of codimension  $v(\lieg)$.
\begin{example} For classical $\R$-split Lie algebras we compute the following: 
\begin{description} 
	\item [Type $A_n$]$r_{\min} (\sl(n, \R)) = 1$;
	\item [Type $B_n, n \ge 3$] $r_{\min} (\so(n, n+1  )) = 1$;	
	\item[Type $C_n, n \ge 2$] $r_{\min} (\sp(2n, \R)) = 2$;
	\item[Type $D_n, n \ge 4$] $r_{\min} (\so(n, n)) = 1$.
\end{description}

Indeed, in the standard enumeration of roots in the associated Dynkin diagram (see ....), for types $A_n, B_n (n\ge 3),$ and $D_n (n\ge 4)$ there is a unique (up to choice of enumerating roots) parabolic subalgebra  $\lieq=\lieq_S$ of minimal codimension corresponding to the set $S = \{\alpha_1\}$.  The coefficient of $\alpha_1$ for every root $\beta$ appearing in $\Sigma_\lieq^\perp$ is $-1$.  One may then find $s\in \bigcap _{j\ge 2}\ker \alpha_j$ with $\alpha_1(s)= -1$ whence $\beta(s) = -1$ for every $\beta\in \Sigma_\lieq^\perp$.  

For type $C_n$, the parabolic subalgebra  $\lieq=\lieq_S$ of minimal codimension corresponding to the set $S = \{\alpha_1\}$ always has minimal codimention $2n-1$ (and is the unique such subalgebra if $n\ge 3$).  Fixing $\lieq=\lieq_S$, $\Sigma_\lieq^\perp$ contains the following roots
\begin{enumerate}
\item  $\beta_1= -\alpha_1$;
\item $\beta_2= -\alpha_1 - 2\alpha_2 - \dots -2\alpha_{n-1} -\alpha_n$;
\item $\beta_3= -2\alpha_1 - 2\alpha_2  - \dots -2\alpha_{n-1} -\alpha_n$.
\end{enumerate}
Taking again $s\in \bigcap _{j\ge 2}\ker \alpha_j$ shows $r_{0}(\lieq)\le 2$.  
On the other hand, given $s\in C(\Sigma_\lieq^\perp) $, if $$\beta_3(s) -2\beta_1(s)=-2\alpha_2 (s) \dots  -2\alpha_{n-1}(s) -\alpha_n(s) \le 0 $$ then $$\frac{\beta_3(s)}{\beta_1(s)} \ge 2.$$ 
If  $$0\le -2\alpha_2 (s) \dots  -2\alpha_{n-1}(s) -\alpha_n(s) =\beta_3(s) -2\beta_1(s) $$ then 
$$\beta_3(s)=\beta_2(s) +\beta_1(s)\le \beta_2(s) + \beta_3(s) -\beta_1(s)=2\beta_2(s)$$ and
$$\frac{\beta_3(s)}{\beta_2(s)} \ge 2.$$ Thus $r_{0} (\lieq) \ge 2$.
(For $S=\{\alpha_2\}$, we similarly check that $r_{0} (\lieq_{\alpha_2}) =1 \le 2$.)


\end{example}

\subsection{Classifying actions in dimension $v(\lieg)$ for $\R$-split groups}

Let $G$ be a connected Lie group with Lie algebra $\lieg$; we will suppose every almost-simple factors of $G$ has real rank at least $2$.  Let $\Gamma$ be a lattice subgroup of $G$.  
Let $Q$ be a parabolic subgroup  of $G$ associated to a (maximal, proper) subalgebra $\lieq$ with codimenison $v(\lieg)$.  Then $\Gamma$ acts projectively on $G/Q$.  

For any compact connected manifold $M$ with $\dim (M) <v(\lieg)$, one expects that all (sufficiently smooth) actions of $\Gamma$ on $M$ are by isometries.  
For the case of equality, $\dim (M) = v(\lieg)$, there is a natural (non-volume-preserving) action of any $\Gamma\subset G$ on $G/Q$;  one again expects that all volume-preserving actions on a manifold $M$ with $\dim (M) = v(\lieg)$ to be an action by isometries (for some choice of Riemannian metric on $M$).
The above expectations were  by the first author  with Fisher and Hurtado in \cite{MR4502593,BFHII,2105.14541}  for the case of $\R$-split simple algebras $\lieg$.  

It is  natural to conjecture in the case $\dim (M) = v(\lieg)$, that the only (sufficiently smooth) non-isometric actions $\alpha\colon \Gamma\to \Diff^r(M)$ 
occur precisely when $M$ is diffeomorphic to a finite cover of $G/Q$ (where  $Q$ is a maximal  parabolic with codimension $v(\lieg)$) and the action $\alpha$ is smoothly conjugate to a lift of the projective action on $G/Q$.

The following result establishes this conjecture in the case $G$ is an $\R$-split (semi-)simple Lie group.     
We note that for higher-rank $\R$-split simple Lie algebras $\lieg$, we have  $v(\lieg)< d(\lieg)$ in \cref{sec:table} whence any isometric action on a manifold of dimension at most $v(\lieg)$  is automatically finite.

\begin{theorem}\label{thm:globalrigid2}
Let $\lieg$ be a semisimple Lie algebra with real rank at least two, all of whose  simple factors are higher-rank and $\R$-split.
Let $G$ be a connected Lie group with Lie algebra $\lieg$ and let $\Gamma\subset G$ be an irreducible lattice subgroup.  

Let $M$ be a compact connected manifold with $\dim (M) = v(G)$.  
Let $r>r_{\min}(\lieg)$ and let $\alpha\colon \Gamma\to \diff^r(M)$ be an action with infinite image $\alpha(\Gamma)$.  

 Then 	there is a maximal proper parabolic subgroup $Q\subset G$  of codimension $v(G)$ and a $C^{r}$ covering map $p\colon M\to G/Q$ such that $$p(\alpha(\gamma)(x)) = \gamma \cdot p(x)$$ for every $x\in M$ and $\gamma\in \Gamma$.  
\end{theorem}

\begin{remark}
In fact, in \cref{thm:globalrigid2} we need not assume that every non-compact simple factor is higher rank.  Indeed, for $C^1$ actions on the circle, Ghys established the analogue of \cref{thm:globalrigid2} for 
 irreducible lattices in higher-rank semisimple Lie groups with factors isogenous to $\SL(2,\R)$ \cite[Theorem 1.2]{MR1703323}.  
 For instance, the conclusion holds for  $\Gamma=\SL(2,\Z[\sqrt 2])$ acting on $S^1$.  
 See also \cite{deroin2020non} for much more recent results for $C^0$ actions on $S^1$ and $\R$.  
 
In our proof of  \cref{thm:globalrigid2} we use (strong) property (T) of $\Gamma$ (in applying \cref{prop:Isom} below).  Though our method of proof could likely recover the $C^{1+\text{\Holder}}$ (and probably $C^1$) versions of the results in \cite{deroin2020non,MR1703323}, since we require the group to be $\R$-split this would only reproduce existing results in the literature.  Thus we only formulate global rigidity results 
 assuming all factors are higher-rank.  
\end{remark}

\begin{remark}[Classifying lifts of standard projective actions] \label{rem:classifyactions}. 
To classify possible actions arising in \cref{thm:globalrigid2}, 
we consider (finite or infinite) covers $M$ of projective varieties $G/Q$ and the possible actions $ \bar \alpha $ by lattice subgroups $\Gamma\subset G$ that factor onto the standard $\Gamma$-action: if $\pi\colon M\to G/Q$ is the covering map, then 
   $$\pi(\bar \alpha(\gamma) (x))=  \gamma \cdot \pi(x).$$

  We write $Q$ for the normalizer of a parabolic Lie subalgebra $\lieq$ and let $Q^\circ$ denote the connected component of the identity in $Q$.  
Recall that if $G$ is any connected semisimple Lie group and if subgroup $Q\subset G$ is a parabolic subgroup, then $G/Q$ is naturally diffeomorphic to $\Ad (G)/\Ad(Q)$ via the map $gQ\mapsto \Ad(g)\Ad (Q )$.

If we assume 
$G$ is topologically simply connected then $ G/Q^\circ  $ is the universal cover of $G/ Q$. 
Recall that if $\Ad (P)= MAN$ is the Iwasawa decomposition of the minimal parabolic subgroup $\Ad (P)$ contained in $\Ad (Q)$, then the  component group of $M$ is a finite, abelian, 2-torsion subgroup $F$ of $\Ad(G)$ normalizing $\Ad(Q)$ (c.f.\  \cite[Theorem 7.53]{MR1920389}.)  Since $F$ permutes the components of $M$, it also permutes the components of $\Ad(Q)$.  Moreover, $F$ is the direct product of elementary groups of order $2$ and, in particular, is abelian; thus every quotient of $F$ is isomorphic to a subgroup whence there is a subgroup $F_Q\subset F$ acting simply 
transitively on the components of $\Ad(Q)$.  

Let  $\Upsilon := \Ad\inv (F_Q)$.  Then $ \Upsilon $ is a discrete subgroup  of $G$  (a central extension of the 2-torsion group $F_Q$) which normalizes $Q$  and for which  $$Q= Q^\circ \Upsilon.$$  
Let $$\Upsilon_Q=Q^\circ\cap \Upsilon .$$  Since $\Upsilon $ normalizes $Q^\circ$, $\Upsilon_Q$ is normal in $\Upsilon $.  
Let $\Lambda= \Upsilon/\Upsilon_Q$.  Then $\Lambda$ acts on $G/ Q^\circ$ as \[g Q^\circ \lambda  =gQ^\circ u\Upsilon _Q. \] 
In particular, we view $\Lambda= \Upsilon/\Upsilon_Q$ as the deck group of \[ G /Q^\circ \to G/Q=\Ad(G)/\Ad(Q) \] via the  quotient orbit map: $$g Q^\circ \mapsto g  Q^\circ\Lambda = g Q^\circ \Upsilon 
= g Q.$$

Let $M$ be a cover of $G/Q$.  
Then there is a  subgroup $\Lambda'\subset \Lambda$ 
so that \[M= G/ (Q^\circ \Lambda' )= G/ ( \Lambda'Q^\circ).\]

Let $\bar \alpha \colon \Gamma\to \Diff ^r \bigl(M)$ be a $p$-equivariant action as above.  Let $x_0 = Q^\circ \Lambda' \in M$.  
As  $p\colon M\to G \bs Q$ is $G$-equivariant, for $g\in G$, we have \[p\circ \bar \alpha(\gamma)(g\cdot x_0) = 
\gamma (g\cdot  p(x_0)).\]  It follows there is $\lambda_\gamma\in \Lambda$
so that $$ \bar \alpha(\gamma)( x_0) = \gamma \lambda_\gamma\inv \cdot x_0.$$
Since $G$ is connected, the $\Gamma$-  and $G$-equivariance of $p$ implies for all $g\in G$ we have
   $$\bar \alpha(\gamma)(g\cdot x_0) = \gamma g\cdot \lambda_\gamma\inv x_0.$$ 

We make the following observations:
\begin{enumerate}
\item For $q\in \Lambda'Q^\circ=\stab (x_0)$, 
$$
\gamma \lambda_\gamma\inv x_0=
\bar \alpha(\gamma)(x_0)
=
\bar \alpha(\gamma)(q\cdot x_0) = \gamma q\cdot \lambda_\gamma\inv x_0$$
whence $ \lambda_\gamma q  \lambda_\gamma\inv\in \stab(x_0) =  \Lambda'Q^\circ$;
in particular, $\lambda_\gamma$  normalizes $\Lambda'Q^\circ$.  
Since we have equality of normalizers $N_\Lambda(\Lambda'Q^\circ)=N_\Lambda(\Lambda')$, we have  $\lambda_\gamma$ contained in $N_\Lambda(\Lambda')$.  
\item For $\gamma_1, \gamma_2\in \Gamma$, 
$$\bar \alpha(\gamma_1\gamma_2)(q\cdot x_0)
= \bar \alpha(\gamma_1)(\gamma_2\cdot q\cdot\lambda_{\gamma_2}\inv x_0)
=\gamma_1\gamma_2\cdot q\cdot\lambda_{\gamma_2}\inv\lambda_{\gamma_1}\inv x_0.$$
In particular, $\lambda_{\gamma_1\gamma_2} =\lambda_{\gamma_1}\lambda_{\gamma_2}\mod \Lambda'$.
\end{enumerate}
We  thus obtain a homomorphism $\Gamma\mapsto N_\Lambda(\Lambda')/\Lambda'$,
\begin{equation}\label{homo}\gamma\mapsto \lambda_\gamma \Lambda'.\end{equation}
Since $\Upsilon$ is a central extension of an abelian group, it is solvable.  As solvability is preserved when passing to subgroups and quotients, we have that $N_\Lambda(\Lambda')/\Lambda'$ is solvable.  By the normal subgroup theorem, we conclude the kernel of \eqref{homo} has finite index in $\Gamma$.  
In particular, there is a finite-index subgroup $\Gamma'\subset \Gamma$ for which 
$$\bar \alpha(\gamma)(g\cdot x_0) = \gamma g x_0$$
for all $\gamma\in \Gamma'$. 
Thus, by restricting to a subgroup  $\Gamma'\subset \Gamma$ of finite index, the action $\bar \alpha(\Gamma')$  coincides with the restriction to $\Gamma'$ of the left $G$-action on $M= G/(Q^\circ \Lambda')$.
\end{remark}

\begin{example}\label{ex:classifyactions} 
Actions of  $\Gamma$ on a finite covers $G/ (Q^\circ \Lambda') $ that factor over the standard action on $G/ Q$ as in \cref{rem:classifyactions} and which only coincide with the restriction of the right $G$-action to  $\Gamma$ on a finite index subgroup can be built as follows: Given a non-trivial homomorphism $\sigma \colon \Gamma\to N_\Lambda(\Lambda')/\Lambda'$, 
 we build an 
action $$\alpha(\gamma)\bigl((g Q^\circ \Lambda'  )\bigr) = \gamma g \sigma(\gamma)\inv  Q^\circ \Lambda'  .$$

We note that non-normal covers arise naturally.  Indeed, for the space $P\bs G$ of full flags on $\R^n$, when $n\ge 4$ the fundamental group $\pi(Q\bs G)$  contains many non-normal subgroups (c.f.\ \cite{MR1618211}).  
\end{example}

As an immediate corollary we obtain the following finiteness results 
\begin{corollary}
Let $G$ be a $\R$-split simple Lie group with real rank at least $2$ and let $M$ be a compact connected manifold with $\dim (M) = v(G)$. 

 If $M$ is not a finite cover of $G/Q$ for some parabolic subgroup $Q\subset G$ then every action $\alpha\colon \Gamma\to \Diff^{r}(M)$ has finite image for $r>  r_{\min }(\lieg)$.
\end{corollary}

\subsection{Classifying actions in dimension  $v(\lieg)$ for some non-split simple groups}
Beyond the case of $\R$-split simple Lie groups, we also obtain global rigidity results when all restricted root spaces are at most 2-dimensional in \cref{thm:globalrigid2nonsplit} below.  We note that \cref{thm:globalrigid2nonsplit} applies, for example,  to lattices in Lie groups $G$ whose Lie algebras are  the following:  
\begin{enumerate}
	\item $\lieg$ is any simple complex Lie algebra;
	\item $\lieg= \mathfrak {su}(p,p)$ for any $p$;
		\item $\lieg= \mathfrak {so}(p,q)$ for any $p\le q\le p+2$;
		\item $\lieg = E II$ (following notation in \cite[App.\ C.4]{MR1920389}).
\end{enumerate}
While the following holds in finite regularity for many groups, we only formulate it for $C^\infty$ actions to avoid overly complicated formulation.  (See \cref{rem:reduce regularity} below for discussion on reducing the regularity.)  Again, we do not consider groups with rank-1 factors as we need to applying \cref{prop:Isom} in the proof.  

\begin{theorem}\label{thm:globalrigid2nonsplit}
Let $\lieg$ be a simple Lie algebra with real rank at least two,  whose restricted root system is not of type $BC_n$ and
 for which every restricted root space is at most  two dimensional.  Let $G$ be a connected Lie group with Lie algebra $\lieg$ and let $\Gamma\subset G$ be a lattice subgroup.  

Let $M$ be a compact connected manifold with $\dim (M) = v(G)$ and let $\alpha\colon \Gamma\to \diff^\infty(M)$ be an action that is not isometric for any Riemannian metric on $M$.    
Then there is a maximal proper parabolic subgroup $Q\subset G$  of codimension $v(G)$ and a $C^{\infty}$ covering map $p\colon M\to G/Q$ such that $$p(\alpha(\gamma)(x)) = \gamma \cdot p(x)$$ for every $x\in M$ and $\gamma\in \Gamma$.  

\end{theorem}

\subsection{Local rigidity of projective actions} 
To state   local rigidity theorems we start with the following definition.
\begin{definition}\label{def:topologies}
Let  $\Gamma$ be a finitely generated group with symmetric generating set $F=\{\gamma_1, \dots, \gamma_p\}$.     For $0\le \ell \le r\le \infty$, the  $C^\ell$-topology on the space of $C^r$ actions $\{\alpha\colon \Gamma\to \diff^{r}(M)\}$ is the topology generated by  restricting  the product topology on $\bigl(\Diff^\ell(M,M)\bigr)^{|F|}$ to the image of the map $\alpha\mapsto \bigl(\alpha(\gamma_1), \alpha (\gamma_2), \dots, \alpha(\gamma_p)\bigr)$.  
\end{definition}

\begin{definition}\label{def:localrigid}

Given an action $\alpha \colon \Gamma\to \diff^{k}(M)$ we say that   $\alpha$ is \emph{$(k,\ell,m)$-locally rigid} if  for all actions  $
\td \alpha \colon \Gamma\to \diff^{k}(M)$ sufficiently $C^\ell$-close to $\alpha$, there is a $C^m$-diffeomorphism $h\colon M\to M$ such that for all $\gamma\in \Gamma$, $$h\circ \td \alpha  (\gamma)= \alpha (\gamma)\circ h.$$

We moreover say that $\alpha \colon \Gamma\to \diff^{k}(M)$  is \emph{$(k,\ell,m,p)$-strongly locally rigid} if the diffeomorphism $h$ converges to the identity in the $C^p$-topology as $\td\alpha$ converges to $\alpha$ in the $C^\ell$ topology. 
\end{definition}

We have the following general version of \cref{thm:slnlocrig} establishing local rigidity for actions of higher-rank lattices on boundaries (i.e.\ generalized flag varieties).  Note that as in  \cite[Theorem 17]{MR1632177},  we allow $G$ to have rank-1 factors but, in contrast to \cite[Theorem 17]{MR1632177}, we do not require the lattice to be cocompact.  

\begin{theorem}\label{thm:localrigid}
Let $G$ be a semisimple Lie group with real rank at least 2 and let $\Gamma$ be an irreducible lattice subgroup.   Let   $Q$ be a   parabolic subgroup of $G$ with Lie algebra $\lieq$.  

Then for any $r \in (r_0(\lieq), \infty ]$, the standard projective action of $\Gamma$ on $G/Q$ is $C^{r,1,r}$-locally rigid.   
Moreover, if $r_0(\lieq)<s\le r$, then the standard projective action of $\Gamma$ on $G/Q$ is $C^{r, s, r, s}$-strongly locally rigid.  
\end{theorem}

As an amusing  consequence of the method of proof,  when all almost simple factors of $G$ are $\R$-split, we obtain local rigidity for $C^r$ actions that are only assumed to be $C^0$ close to the standard projective action.  This, in particular, gives the local rigidity for $C^0$ perturbations in \cref{thm:slnlocrig} above.

\begin{theorem}\label{thm:localrigidC0}
Let $G$ be a semisimple Lie group with real rank at least 2 and let $\Gamma$ be an irreducible lattice.   Let   $Q$ be a   parabolic subgroup of $G$ with Lie algebra $\lieq$.  Suppose every almost simple factor of $G$ that is not contained in $Q$ is $\R$-split. 

Then for any $r \in (r_0(\lieq), \infty ]$, the standard projective action of $\Gamma$ on $G/Q$ is $C^{r,0, r}$-locally rigid.

\end{theorem}

\subsection{Smooth factors of boundary actions} 
Let $G$ be a connected semisimple Lie group 
 and let $P\subset G$ be a minimal parabolic subgroup.  Let $K\subset G$ be the maximal compact subgroup and let $m_{G/P}$ denote the unique left-$K$-invariant (probability) measure on $ G/P$. 

Let $\Gamma$ be an irreducible lattice subgroup of $G$.  The (left) projective action of  $\Gamma$ on $G/P$ does not preserve $m_{G/P} $ but preserves the measure class of $m_{G/P} $.  

Let $X$ be a measurable factor of the $\Gamma$-action on $G/P$.  That is, suppose $X$ is a measurable space equipped with a measurable $\Gamma$-action and  $p\colon G/P\to X$ is a  $\Gamma$-equivariant map; equip $X$ with the measure $p_* m_{G/P}$.  
In \cite{MR515630,MR545365}, Margulis showed that every such factor is measurably isomorphic to the $\Gamma$-action on $G/Q$ for some parabolic $P\subset Q$.  
Namely, there exists a $\Gamma$-equivariant, measurable isomorphism $h\colon (G/Q, m_{G/Q})\to (X, p_*m)$ such that the following diagram commutes:
\begin{equation}\label{eq:factordiagram}
\begin{tikzpicture}[baseline=(current  bounding  box.center)]
  \matrix (m) [matrix of math nodes,row sep=3em,column sep=4em,minimum width=2em]
  {
(G/P, m_{G/P}) & \phantom{O}  \\
    (G/Q, m_{G/Q})& (X, p_*m_{G/P}) \\};
    \path[-stealth]
    (m-1-1) edge node [left] {$ \pi$} (m-2-1)
            edge   node  [above] {$p$} (m-2-2)            
    (m-2-1) edge node [below] {$h$}  (m-2-2);
\end{tikzpicture}
\end{equation}
Above, $m_{G/Q}$ is the unique  $K$-invariant probability measure on $G/Q$ and $\pi\colon G/P\to G/Q$ is the natural map $\pi\colon gP\to gQ$.

In \cite{MR736562}, Dani considered the above setup under the additional assumption that the factor map $p$ is a continuous surjection onto a compact topological space $X$ admitting a continuous $\Gamma$-action.    
We formulate Dani's theorem as we will refer to it in the sequel.
\begin{theorem}[Dani,  \cite{MR736562}]\label{thm:dani}
Let $G$ be a connected semisimple Lie group with finite center and with real rank at least 2 and let $\Gamma$ be an irreducible lattice subgroup.  Let $P$ be a   minimal parabolic subgroup.  Let $X$ be a compact Hausdorff space  equipped with an action of  $\Gamma$ by homeomorphism and suppose  $p\colon G/P\to X$ is a $\Gamma$-equivariant  continuous surjection.  

Then there exists a parabolic subgroup $P\subset Q$ and a  $\Gamma$-equivariant homeomorphism $h\colon G/Q\to X$ such that the following diagram of $\Gamma$-equivariant maps commutes: 
\begin{equation}\label{eq:factordiagram2}
\begin{tikzpicture}[baseline=(current  bounding  box.center)]
  \matrix (m) [matrix of math nodes,row sep=3em,column sep=4em,minimum width=2em]
  {
G/P & \phantom{O}  \\
    G/Q& X \\};
    \path[-stealth]
    (m-1-1) edge node [left] {$ \pi$} (m-2-1)
            edge   node  [above] {$p$} (m-2-2)            
    (m-2-1) edge node [below] {$h$}  (m-2-2);
\end{tikzpicture}
\end{equation}
\end{theorem}
\begin{remark}\label{rem:unique}
For any parabolic subgroup $P\subset Q\subset G$, it is well known that the projective action of $\Gamma$ on $G/Q$ is both minimal and proximal.  As discussed in  \cite[Remark 3.1]{MR736562}, this implies the canonical projection $\pi\colon G/P\to G/Q$  is the  unique  $\Gamma$-equivariant continuous surjection.  
This immediately implies the following facts that we will make use of: 
\begin{enumerate}
\item  the  homeomorphism  $h\colon G/Q\to X$ in \cref{thm:dani} is unique;
\item the  $C^0$ centralizer of  the projective action  of $\Gamma$ on $G/Q$ is the identity.  
\end{enumerate}
\end{remark}

In \cite{MR3748688}, Gorodnik and Spatzier considered the case of smooth factors in Dani's theorem. Namely, one considers the case that  $X$ is a smooth manifold and $\alpha\colon \Gamma\to \Diff^\infty(X)$ a $C^\infty$  action; however, one  still only assumes the factor map  $p\colon G/P\to X$ is a $\Gamma$-equivariant continuous surjection.  Under the additional hypotheses that action $\alpha$ admits a differentiable sink for some $\gamma\in \Gamma$, Gorodnik and Spatzier showed the homeomorphism $h\colon G/Q\to X$ in \cref{thm:dani} is a $C^\infty$ diffeomorphism; in particular, this implies that the map $p\colon G/P\to X$  is  $C^\infty$ surjection with constant rank.   The tools we develop to establish \cref{thm:globalrigid2,thm:localrigid} allow us to remove the hypothesis on the existence of a differentiable sink, thus giving an alternative proof and strengthening  of the main results of \cite{MR3748688}.  

\begin{theorem}\label{thm:smoothfactor}
Let $G$ be a semisimple Lie group   with real rank at least 2 and let $\Gamma$ be an irreducible lattice subgroup.   
 Let $P$ be a   minimal parabolic subgroup.  

  Let $M$ be a compact connected manifold and let  $\alpha\colon \Gamma\to \diff^\infty(M)$ be an action.  
Suppose there exists a $\Gamma$-equivariant, continuous surjection $p\colon G/P \to M$.  Then there exists a parabolic subgroup $P\subset Q$ and a $\Gamma$-equivariant, $C^\infty$-diffeomorphism $h\colon  G/Q\to M$ such that 
$$p(gP) = h(gQ).$$

In particular, we have the following commutative diagram of $\Gamma$-equivariant, $C^\infty$, maximal rank maps:
$$\begin{tikzpicture}[baseline=(current  bounding  box.center)]
  \matrix (m) [matrix of math nodes,row sep=3em,column sep=4em,minimum width=2em]
  {
G/P & \phantom{O}  \\
    G/Q& M \\};
    \path[-stealth]
    (m-1-1) edge node [left] {$ \pi$} (m-2-1)
            edge   node  [above] {$p$} (m-2-2)            
    (m-2-1) edge node [below] {$h$}  (m-2-2);
\end{tikzpicture}$$
\end{theorem}



Our method of proof also allows us to obtain smooth factor theorems in much lower regularity for the case of $\R$-split groups.  
\begin{theorem}\label{thm:smoothfactorSPLIT}
In   \cref{thm:smoothfactor} assume the following additional hypotheses:
\begin{enumerate}
\item every almost simple factor of $G$ is $\R$-split, 
\item there exists $Q'\supset P$ with Lie algebra $\lieq'$ such that $p\colon G/P\to M$ factors through a $\Gamma$-equivariant continuous surjection $p'\colon G/Q'\to M$,
\item $r\in (r_0(\lieq'), \infty]$.
\end{enumerate}
Then if $\alpha\colon \Gamma\to \diff^r(M)$ is a $C^r$ factor of the $\Gamma$ action on $G/Q'$ then  there exists $Q\supset Q'$ and a $\Gamma$-equivariant $C^r$-diffeomorphism $h\colon  G/Q\to M$.  
\end{theorem}
In \cref{thm:smoothfactorSPLIT}, we could take $Q'= P$; the hypothesis that $p$ factors through $G/Q'$ is only to reduce the assumed regularity for the action to the regularity proscribed by the Lie algebra $\lieq'$ of $Q'$.  In particular, \cref{thm:smoothfactorSPLIT} applies to actions of regularity at least $r_0(\liep)$, where $\liep$ is the Lie algebra of $P$.  

\subsection{Embedded projective actions inside volume-preserving actions}
We turn to the cases of volume-preserving actions. 
Let $P=MAN$ be a choice of minimal parabolic subgroup in $G$.   Given an action $\alpha\colon \Gamma\to \diff(M)$, we study $P$-invariant measures on the suspension space $M^\alpha$.   The presence of a $P$-invariant measure that is not $G$-invariant implies the existence of a ``blow-up,'' and invariant submanifold  on which the dynamics restricts to the standard projective action.  For technical reasons arising in the proof, we need to assume $G/Q$ has finite fundamental group for all parabolic subgroups $Q$ with of codimension $v(\lieg)$.  See \cref{rem:noZcovers} below for discussion of when this holds.  
\begin{theorem} \label{thm:embedded projectives}  \label{thm:partpara}
Consider an $\R$-split simple Lie group $G$ with real rank at least $2$. 
Suppose for every parabolic subgroup $Q\subset G$ of codimension $v(G)$, the projective variety $G/Q$ has finite fundamental group.  
Let $\Gamma\subset G$ be a lattice subgroup.  
Let $M$ be a ($v(G)+1$)-dimensional manifold and for $r> r_{\min} (G)$, let $\alpha\colon \Gamma\to \diff^r_\vol(M)$ be a volume-preserving action.  

	Let $M^\alpha$ denote the induced $G$-space.  Suppose  for the induced $G$-action on $M^\alpha$ there exists an ergodic  Borel probability measure $\mu$ on $M^\alpha$ such that  $\mu$ is  $P$-invariant  but not $G$-invariant.

Then there exist a compact embedded (possibly disconnected) $C^r$ submanifold $N$, a  parabolic subgroup $Q\subset G$ with codimension $v(G)$, and a $C^r$ covering map $h\colon N\to G/Q$ with the following properties:
\begin{enumerate}
	\item $N$ is $\alpha(\Gamma)$-invariant. 
	\item $h\colon N\to G/Q$ is $\Gamma$-equivariant: for all $\gamma\in \Gamma$ and $x\in N$, \[h\circ \alpha(\gamma) (x) = \gamma \cdot h(x).\]
 \end{enumerate}
 Let  $N= \bigsqcup_{i=1}^{\ell} N_i$ be the partition of $N$ into connected components.  Then 
 \begin{enumerate}[resume]
	\item Each $N_i$ is diffeomorphic to a finite cover of $G/Q$. 
	\item There is a finite index subgroup $\Gamma'\subset \Gamma$ such that each $N_i$ is $\Gamma'$-invariant and the dynamics of $\Gamma'$ on each $N_i$ is $C^r$ conjugate to standard action. 
	\item The measure $\mu$ is supported on the suspension $(G\times N)/\Gamma$ in $M^\alpha = (G\times M)/\Gamma$.
\end{enumerate}
\end{theorem}

\begin{remark}\label{rem:noZcovers}
\cref{thm:embedded projectives} holds for actions by lattices $\Gamma\subset G$ in groups $G$ whose Lie algebras $\lieg$ are the following: 
\begin{enumerate}
	\item $\lieg= \mathfrak{sl}(n,\R)$ for $n\ge 3$.  Indeed, $v(G) = n-1$ and all projective varieties $G/Q$ of dimension  $v(G)$ are the grassmannians of lines (or hyperplanes) in $\R^n$.  We have  $\pi_1(\RP^{n-1}) =\Z/(2\Z)$ when $n\ge 3$.
	\item $\lieg= \mathfrak{sp}(2n, \R)$ for $n\ge 3$.  Indeed, for $n\ge 3$, $v(G)= 2n-1$ and the only $v(G)$-dimensional projective variety $G/Q$ is the space of (isotropic) lines in $\R^{2n}$ (with  $\pi_1= \Z/(2\Z)$).   When $n=2$, $v(G)=3$ and the 3-dimensional projective varieties $G/Q$ are either  the space of (isotropic) lines in $\R^4$ (with  $\pi_1= \Z/(2\Z)$) or the space of isotropic (i.e. Lagrangian) planes in $\R^4$ (with $\pi_1 = \Z$).   
	\item $\lieg= \mathfrak{so}(n, n+1)$ for $n\ge 3$ and  $\lieg= \mathfrak{so}(n, n)$ for $n\ge 4$. Indeed, in these cases, by direct computation, the only parabolic subalgebra with codimension $v(G)$ is the sublagebra which omits the left-most root in the Dynkin diagram (or a valence-1 root  in the $D_4$ diagram).  The corresponding space $G/Q$ is thus the Grassmannian of isotropic lines in $\R^{2n+1}$ or $\R^{2n}$.  

	In general, if $Q$ is a quadratic form of signature $(p,q)$, the Grassmannian of isotropic lines is a quadric hypersurface $H$ in $P(\R^{p+q})$.  The maximal compact subgroup of $\So(Q)$  has Lie algebra  $\mathfrak{so}(p)\times \mathfrak{so}(q)$ and acts transitively on $H$.  We then see that $H=(S^{p-1}\times S^{q-1})/{\pm 1}$.  In particular, if $p\ge 3$ and $q\ge 3$ then $H$ has fundamental group $\Z/(2\Z)$.   
 	\item See \cite[Theorem 1.1]{MR1618211} for a general presentation of $\pi_1(G/Q)$ for all real semisimple Lie groups $G$ and all parabolic subgroups $Q$ (inclusive of  those above.)
\end{enumerate}

\end{remark}
\begin{question}
	Does the conclusion of \cref{thm:embedded projectives}  hold for volume-preserving actions by lattices $\Gamma$ in $\Sp(4, \R)$ on compact 4-manifolds?  Alternatively, is it possible to embed the projective $\Gamma$-action on the (3-dimensional) Grassmannian of Lagrangian planes in $\R^4$ (as in \cref{rem:noZcovers}(2)) as an invariant subset of a volume-preserving  action on a 4-manifold?  

We note that without the volume-preserving assumption, any lattice $\Gamma$ in $\Sp(4, \R)$ acts on the $4$-dimensional Grassmannian $\mathrm{Gr}(2,4)$ 	 of all planes in $\R^4$, preserving the $3$-dimensional subvariety of Lagrangian planes. 

More precise versions of the above question might include the following:
\begin{enumerate}
	\item Is it possible to blow-up a fixed point (in a volume-preserving way) for the standard $\Sp(4,\Z)$ action on $\T^4$ to embed the projective action on Lagrangian planes. 
	\item Is it possible to modify the standard $\Gamma$ action on $\mathrm{Gr}(2,4)$ to obtain a volume-preserving action that preserves the locus of Lagrangian planes?  (Such an action would be either a large perturbation or would likely modify the topology).
\end{enumerate}
\end{question}
\begin{remark}
Let $\lieg=\so(n,n+1)$ for $n\ge 3$ or $\lieg = \so(n,n)$ for $n\ge 4$ and let $\Gamma$ be a lattice in a connected Lie group with Lie algebra $\lieg$.  In either case, we have $v(\lieg) + 2= n(\lieg)$.   The volume-preserving version of Zimmer's conjecture asserts if $\dim (M) = n(\lieg) -1$, every volume-preserving action of $\Gamma$ on $M$ is finite.  
(The conclusion of \cite{MR4502593,BFHII,2105.14541} only verify this conjecture for $\dim (M) = v(\lieg) = n(\lieg)-2$.)

\cref{thm:partpara} provides some structure for such actions that may provide tools to verify this conjecture in the missing dimension; namely for all such actions, \cref{thm:partpara} ensures the existence of an embedded codimension-1 submanifold on which the action (restricted to a finite index subgroup) is smoothly conjugate to (a cover of) the standard projective action.  To show the conjecture, one should then rule out the possibility that such subactions can be embedded into a volume-preserving action on a compact manifold of dimension $n(\lieg) -1=v(\lieg)+1$.  
\end{remark} 

\section{Induced $G$-action and fiberwise smooth ergodic theory}\label{sec:set}
Throughout this section, fix $G$ to be a connected semisimple Lie group  as in  \cref{ss:Gassump}.  We further assume $G$ is topologically simply connected.    Let $\Gamma$ be an irreducible lattice subgroup of $G$.  Let $M$ be a compact manifold and let $\alpha\colon \Gamma\to \Diff^r(M)$ be an action.  

\subsection{Suspension construction and induced $G$-action}\label{sss:susp}
We follow a well-known construction (previously used in \cite{MR4502594,MR3702679,MR4502593,2105.14541,ABZ,BDZ} among others) which allows us to relate various properties of an action $\alpha\colon \Gamma\to \Diff(M)$ with properties of a $G$-action on an associated  bundle $M^\alpha$ over $G/\Gamma$.

On the product $G\times M$ consider the right $\Gamma$-action and the left $G$-action
\begin{equation}\label{eq:susp} (g,x)\cdot \gamma= (g\gamma , \alpha(\gamma\inv)(x)),\quad \quad \quad a\cdot (g,x) = (ag, x).\end{equation}
Define  the quotient manifold $M^\alpha:= G\times M/\Gamma $.  
Given $(g,x)\in G\times M$, we will denote by $$[g,x]_\alpha$$ (or $[g,x]$ if there is no cause for confusion) the corresponding equivalence class in $M^\alpha.$
As the  $G$-action on $G\times M$ commutes with the $\Gamma$-action, we have an induced left $G$-action  on $M^\alpha$.  For $g\in G$ and $x\in M^\alpha$ we denote this action by $g\cdot x$.

We write $\pi\colon M^\alpha\to G/\Gamma$ for the natural projection map.  Note that $M^\alpha$ has the structure of a fiber bundle over $G/\Gamma$ induced by the map $\pi$ with fibers diffeomorphic to $M$. Note that the $G$-action preserves the fibers of $M^\alpha$.
Let 
\begin{equation}\calF(x):= \pi\inv (\pi(x))\end{equation}denote the fiber of $M^\alpha$ through $x$ and let 
\begin{equation}F:= \ker D\pi\end{equation}
denote the fiberwise tangent bundle of $M^\alpha$.   Then $T_x \calF(x) = F(x)$.  
Let $\calA\colon G\times F\to F$ denote the fiberwise derivative cocycle.  Given $x\in M^\alpha$  and $g\in G$, write $$\calA(g,x)\colon F(x)\to F(g\cdot x)$$ for the linear map between fibers, $$\calA(g,x)(v) = D_x g(v).$$

\subsection{Fundamental sets, adapted norms, and quasi-isometric embeddings}\label{sec:norms}
When $\Gamma$ is cocompact in $G$, all choices of Riemannian metrics on $TM^\alpha$ (or the bundle $F\to M^\alpha$) are equivalent.  Thus, we may freely employ  all the tools of smooth ergodic theory for the action of any $g\in G$ on  $M^\alpha$.  

When $\Gamma$ is non-uniform, $M^\alpha$ is not compact and the local dynamics along orbits need not be uniformly bounded.  Thus, to employ tools from smooth ergodic theory,  more care is needed in specifying norms on the fibers of $M^\alpha$.  
In this case, we follow either \cite[\S 5.4]{2105.14541} or \cite[\S 2]{MR4502594}; we summarize these results below.  

When $\Gamma$ is non-uniform, we may assume $\Ad(G)$ is $\Q$-algebraic and that $\Ad(\Gamma)$ is commensurable with the $\Z$-points in $\Ad(G)$.  We may define Siegel sets and Siegel fundamental sets in $\Ad(G)$ (and thus $G$) relative to any choice of Cartan involution $\theta$ on $G$ and a minimal $\Q$-parabolic subgroup in $\Ad(G)$.  Using the Borel-Serre partial compactification (of $G$ for which $G/\Gamma$ is an open dense set in a compact manifold with corners), we equip the bundle $G\times TM \to G\times M$ with a $C^\infty$ metric with the following properties: Write $\langle \cdot, \cdot \rangle_{g,x}$ for the inner product on the fiber over $(g,x)$.   Then
\begin{enumerate}
\item $\Gamma$ acts by isometries on $G\times TM$.
\item there exists a fundamental set $D$ for $\Gamma$ in $G$ (namely, any choice Siegel fundamental set $D\subset G$ relative to a choice of Cartan involution $\theta$ and a minimal $\Q$-parabolic subgroup) 
and $C>1$ such that for all $g,g'\in D$,
$$
\frac 1 C\langle \cdot, \cdot \rangle_{g,x}\le 
\langle \cdot, \cdot \rangle_{g',x}
\le C\langle \cdot, \cdot \rangle_{g,x}.
$$ 
\end{enumerate}
The metric on $G\times TM$ then descends to a $C^\infty$ metric on the bundle $F\to M^\alpha$. For the remainder,  given $x\in M^\alpha$, we denote by $\|\cdot \|_x^F$ the induced norm on the fiber of $F$ through $x$.  

We similarly equip $G$ with any right-invariant metric and equip $G\times M$ with the associated Riemannian metric that makes $G$-orbits orthogonal to fibers, restricts to the right-invariant metric on $G$-orbits, and restricts to the above metric on every fiber.  This induces a metric on $TM^\alpha$.

We also recall the following fundamental result of Lubotzky-Mozes-Raghunathan.  Relative to any fixed choice of generating set for $\Gamma$, write $|\gamma|$ for the word length of $\gamma$ in $\Gamma$.   Equip $G$ with any right-$G$-invariant, left-$K$-invariant metric $d_G$.  
When $G$ has finite center, the following is the main result of \cite{MR1828742}.  In  the case $G$ has infinite center, see discussion in \cite[\S 3.9.3]{2105.14541}.  
\begin{theorem}[\cite{MR1828742}]\label{LMR}
The word-metric and Riemannian metric on $  \Gamma$ are quasi-isometric: there are $A_0,B_0$ such that for all $\gamma,\gamma'\in   \Gamma $,
$$\frac 1 {A_0 }d_G(\gamma,\gamma') -B_0 \le |\gamma\inv\gamma'| \le {A_0 }d_G(\gamma,\gamma') +B_0.$$
\end{theorem}


As discussed above, when $\Gamma$ is non-uniform, given $a\in A$, the fiberwise dynamics of $a$ on $M^\alpha$ need not be bounded (in the $C^1$ topology.)  However, for $a$-invariant probability measures on $M^\alpha$ projecting to the Haar measure on $G/\Gamma$, the degeneracy is subexponential along orbits.  We summarize this below 
\cref{prop:integrability}

\def\fund{{D_F}}

Let $D\subset G$ be a fundamental set on which the fiberwise metrics are uniformly comparable and fix 
 a Borel fundamental domain~$\fund$ contained in $D$ for the right $\Gamma$-action on $G$.  Let  $\beta_\fund\colon G \times G/\Gamma \rightarrow \Gamma$ be the \emph{return cocycle}: given $\hat x\in G/\Gamma$, take $\wtd x$ to be the unique
lift of $\hat x$ in $\fund$ and define $b(g,\hat x) = \beta_\fund(g,\hat x)$ to be the unique $\gamma\in \Gamma$ such that
$g\wtd{x}\gamma\inv \in \fund$.
One verifies that $\beta_\fund$ is a Borel-measurable cocycle and a second choice of fundamental domain 
defines a cohomologous cocycle.

%

Given a diffeomorphism $g\colon M\to M$, let $\|g\|_{C^k}$ denote the $C^k$ norm of $g$ (say, relative to some choice of embedding of $M$ into some Euclidean space $\R^N$.)
Given $g\in G$ and $x\in \fund$, let  $\psi_k(g,x) = \|\alpha (\beta_{\fund}(g,x))\|_{C^k}.$

Let $m_{G/\Gamma}$ denote the normalized Haar measure on $G/\Gamma$.  Using \cref{LMR} and standard properties of Siegel sets, we obtain the following: 

\begin{proposition}\label{prop:integrability}
For any $k$, any  $1\le q<\infty$, and any compact $B\subset G$, the map $$x\mapsto \sup_{g\in B} \log \psi_k (g,x)$$ is $L^q(m_{G/\Gamma})$ on $G/\Gamma$.  
In particular, for $m_{G/\Gamma}$-a.e.\ $x\in G/\Gamma$ and any $g\in G$,
\begin{equation}\label{eq:integrabl}
\lim_{n\to \infty} \frac 1 n \log^+( \psi_k (g,g^n\cdot x) ) = 0.
\end{equation}
\end{proposition}

\subsection{Fiberwise Lyapunov exponent and fiberwise Oseledec's theorem}
Recall $\calA\colon G\times F\to F$ denotes the fiberwise derivative cocycle.  
We write $\|\calA(g,x)\|$ for the operator norm of the cocycle relative to the fiberwise norms on $F$ discussed in \cref{sec:norms}.


\subsubsection{Family of measures}
Throughout, we will consider $A$-invariant Borel probability measures $\mu$ on $X$ whose projection to $G/\Gamma$ is the normalized Haar measure $m_{G/\Gamma}$ on $G/\Gamma$.
From \cref{prop:integrability} (see also \cite[Claim 2.6]{MR4502594}), for any such measure $\mu$ and any compact $B\subset G$, we have 
\begin{equation}\label{eq:suckmytoes}x\mapsto \sup_{g\in B}\log^+ \|\calA(g,x)\| \in L^q(\mu)\end{equation} for all $1\le q<\infty$.

\subsubsection{Average top Lyapunov exponent}
Given $s\in A$ and an $s$-invariant Borel probability measure $\mu$ on $M^\alpha$ projecting to the Haar measure on $G/\Gamma$,  define the \emph{average top} (or \emph{leading})  \emph{Lyapunov exponent of $\calA$} to be
\begin{equation} \label{eq:toplyap} \lambda^F_{\top, \mu, s_0, \alpha}
 := \inf _{n} \frac 1 n  \int \log^+ \|\calA(s^n, x)\| \, d \mu (x).\end{equation}
 and  the \emph{average bottom}  \emph{Lyapunov exponent of $\calA$} to be
\begin{equation} \label{eq:botlyap} \lambda^F_{\btm, \mu, s_0, \alpha}
 := -\inf _{n} \frac 1 n  \int \log^+ \|\calA(s^n, x)\inv\| \, d \mu (x)\end{equation}
 where $\log^+(a) = \max\{0, \log a\}$.  
We remark that  $s$-invariance of the measure $\mu$ implies the sequences  $\frac 1 n\pm  \int \log^+( \|\calA(s^n, x)\|^{\pm 1}) \, d \mu (x)$ are subadditive whence the infimums   in \eqref{eq:toplyap} and \eqref{eq:botlyap} may be replaced by  limits. 

\subsubsection{Higher-rank Oseledec's theorem}
Consider  an abelian subgroup  $A\subset G$ isomorphic to $\R^k$; usually we will further assume that $A$ is a maximal $\R$-split Cartan subgroup.  Equip $A\cong \R^k$ with any norm $|\cdot|$.  

Consider  an ergodic, $A$-invariant Borel probability  measure $\mu$ on $X$ projecting to the Haar measure on $G/\Gamma$.  By the integrability  \eqref{eq:suckmytoes} of the linear cocycle $\calA$, 
we have the following consequence of the higher-rank Oseledec's multiplicative ergodic theorem  (c.f.\ \cite[Theorem 2.4]{MR4599404}).
Let $\calE\subset TM^\alpha$ be a continuous,  $A$-invariant subbundle.  
\begin{proposition}
\label{thm:higherrankMET}
There are
	\begin{enumerate}
	\item an $A$-invariant subset $\Lambda_0\subset M^\alpha$ with $\mu(\Lambda_0)=1$;
  \item 
   linear functionals $\lambda_i\colon \R^k\to \R$ for $1\le i\le p$;  
	\item   and splittings   $\calE= \bigoplus _{i=1}^p E^{\lambda_i}(x)$ 
	into families of mutually transverse,  $\mu$-measurable  subbundles $E^{\lambda_i}(x)\subset \R^d$ defined  for $x\in \Lambda_0$

	\end{enumerate}
such that
\begin{enumerate}	[resume]
	\item $\calA (s, x) E^{\lambda_i}(x)= E^{\lambda_i}(s\cdot x)$ and
	\item \label{lemma:partb} $\displaystyle \lim_{|s|\to \infty} \frac { \log |  \calA (s,x) (v)| - \lambda_i(s)}{|s|}=0$
\end{enumerate}	
	for all $x\in \Lambda_0$ and all $ v\in  E^{\lambda_i}(x)\sm \{0\}$.  
 \end{proposition}
 The linear functionals $\lambda_i\colon A\to \R$ in \cref{thm:higherrankMET} are called the \emph{Lyapunov exponents} of the $A$-action on $(\calE, \mu)$.

Note that  conclusion \eqref{lemma:partb} of \cref{thm:higherrankMET} implies for $v\in E^{\lambda_i}(x)$ the weaker result that for $s\in A$,
$$\lim_{k\to\pm \infty} \tfrac {1} k \log |  \calA (s^k,x) (v)| =  \lambda_i(s).$$

We also remark that if $\mu$ is an $A$-invariant, $A$-ergodic measure  then for any $s\in A$, the average top and bottom Lyapunov exponents are given as
\begin{equation}\label{eqn:lyapn}
\lambda^F_{\top, \mu, s_0, \alpha} = \max _i \lambda_i(s), \quad \quad 
\lambda^F_{\btm, \mu, s_0, \alpha} = \min  _i \lambda_i(s).
\end{equation}

\subsubsection{Collections of Lyapunov exponents}
Let $A$ be a maximal split Cartan subgroup of $G$.  
Fix  an ergodic, $A$-invariant Borel probability  measure $\mu$ on $X$ projecting to the Haar measure on $G/\Gamma$.  
When $\calE= TM^\alpha$ is the full tangent bundle, write  $\calL(\mu)= \{\lambda_i \colon A\to \R\}$ for the collection of Lyapunov exponents for the $A$-action on $(M^\alpha,\mu)$.  
When $\calE= \calF$ is the fiberwise tangent bundle, we apply \cref{thm:higherrankMET} and obtain 
$\calL^F(\mu)= \{\lambda_i^F \colon A\to \R\}$, the collection of \emph{fiberwise Lyapunov exponents.}  Similarly, when $\calE$ is the bundle tangent to every $G$-orbit, the Lyapunov exponents guaranteed by \cref{thm:higherrankMET} are simply the restricted roots  $\Sigma$ of $\lieg$ (under the exponential identification between $\liea = \Lie(A)$ and $A$); in particular, $\Sigma$ is independent of the choice of $\mu$.  
We note that $\calL^F(\mu)\subset \calL(\mu)$ and $\Sigma\subset \calL(\mu)$.
 
 We say two exponents $\lambda_1, \lambda_2\in \calL(\mu)$ are \emph{equivalent} if they are positively proportional as linear functionals on $A$.    
 Write $\what \calL(\mu)$ and  $\what \calL^F(\mu)$ for the equivalence classes in 
$  \calL(\mu)$ and  $  \calL^F(\mu)$.  We refer to an element of $\chi \in\what \calL(\mu)$ (resp.\ $\chi \in\what \calL^F(\mu)$ as a \emph{coarse Lyapunov exponent} (resp.\ \emph{coarse fiberwise Lyapunov exponent}.)
Given $\lambda\in \calL$, we also write $\chi = [\lambda]$ for the coarse Lyapunov exponent containing $\lambda$. 
We similarly write $\what \Sigma$ for the equivalence classes of restricted roots and refer to $[\beta]\in \what \Sigma$ as a \emph{coarse root}.

\begin{remark}[On Lyapunov exponents for ergodic $P$-invariant measures]\label{rem:PergvsAerg}
Let $\mu$ be an ergodic $P$-invariant Borel probability measure on $M^\alpha$.  
We note  that $\mu$ projects to the Haar measure (the unique $P$-invariant Borel probability measure) on $G/\Gamma$.  

We have that $\mu$ is $(MA)$-ergodic (c.f.\ proof of \cite[Claim 5.2]{MR4502594}) but need not be $A$-ergodic.  However, we have that almost every $A$-ergodic component projects to the Haar measure on $G/\Gamma$.  Furthermore, as $M$ commutes with $A$, the Lyapunov exponents of every $A$-ergodic component of $\mu$ coincide (as linear functionals on $A$).  Thus, we freely refer to the Lyapunov exponents (for the $A$-action) of  ergodic $P$-invariant Borel probability measures on $M^\alpha$.  

See \cite[Claim 5.2]{MR4502594} for further details. 
\end{remark}

\subsection{Resonant roots and extra invariance}

Let $\mu$ be an ergodic, $A$-invariant Borel probability measure whose projection to $G/\Gamma$ is the Haar measure.
\begin{definition}
A restricted root $\beta\in \Sigma$ is \emph{resonant} (with respect to  $\mu$) if there exists a non-zero fiberwise Lyapunov exponent $\lambda_i^F$  and $c>0$ such that $$\lambda_i^F= c\beta.$$
If $\beta$ is not positively proportional to any $\lambda_i^F\in \calL^F(\mu)$, we say $\beta$ is \emph{non-resonant}.
\end{definition}

Note that resonance and non-resonance are properties of  coarse restricted roots.  Given a coarse restricted root $[\beta]\in \what \Sigma$, we write $U^{[\beta]}$ for the unipotent root group tangent to $\bigoplus _{\beta'\in [\beta]}\lieg^{\beta'}$.  We have the following from our previous paper \cite{MR4502594}.
\begin{proposition}[{\cite[Proposition 5.1]{MR4502594}}]
\label{NRimpInv}
Let $\mu$ be an ergodic, $A$-invariant Borel probability measure whose projection to $G/\Gamma$ is the Haar measure.
If $\beta\in \Sigma$ is non-resonant (with respect to $\mu$) then $\mu$ is $U^{[\beta]}$-invariant.  

\end{proposition}

 \subsection{Fiberwise Lyapunov manifolds} \label{ssec:FWLyapu}
Let $\mu$ be an ergodic, $A$-invariant Borel probability  measure on $M^\alpha$ projecting to the Haar measure on $G/\Gamma$.   Recall from \cref{rem:PergvsAerg} that the Lyapunov exponents of the $A$-action are independent of the choice of $A$-ergodic component of $\mu$.   


\subsubsection{Strongly integrable collections of exponents}\label{integrable}
Recall we write  $  \calL(\mu)$ for the collection of Lyapunov exponents for the $A$-action on $(M^\alpha,\mu)$.   We have the following terminology. 
\begin{definition}\label{def:strongInt}
A subcollection $\calI\subset  \calL(\mu)$ is \emph{strongly integrable} if there exist $a_1,\cdots, a_k\in A$ such that 
$$\calI= \bigcap _{i=1,\dots , k} \left\{ \lambda \in\what \calL(\mu): \lambda(a_i)<0\right \}.$$
 \end{definition}

 We note that any strongly integrable collection of Lyapunov exponents is saturated by coarse Lyapunov exponents.  

\subsubsection{Properties of  strongly integrable collections}
Fix $\calI$ to be a strongly integrable collection of   Lyapunov exponents as above.  
Fix the following notation: write
\begin{enumerate}
\item $\Sigma_\calI =   \Sigma\cap \calI$, for the set of restricted roots in $\calI$.   
\item $ \calL^{F,\calI}(\mu) = \calL^F(\mu)\cap \calI$ for the set of fiberwise Lyapunov exponents  $\calI$.   
\item $\displaystyle \liev_\calI:= \bigoplus _{\beta \in\Sigma_ \calI}\lieg^{\beta}$;
\item $\displaystyle E^{\calI}(x) = \bigoplus _{\lambda \in \calI}E^{\lambda}(x)$;
\item $\displaystyle E^{F,\calI}(x) = \bigoplus _{\bar \lambda \in \calL^{F,\calI}(\mu)}E^{F,\bar \lambda}(x).$
\end{enumerate}

%
Note strong integrability of  $\calI$ implies that $\liev_\calI$ is a (nilpotent) Lie subalgebra.  We  write $U^\calI\subset G$ for the analytic subgroup of $G$ associated with $\liev_\calI$.

\subsubsection{Critical regularity associated to $\calI$}
We write \begin{equation}\label{eq:cone}C(\calI):= \{ a\in A: \lambda(a)<0 \text{ for every } \lambda\in \calI\}\end{equation} for the \emph{negative cone} associated to a collection of linear functionals $\calI\subset \calL(\mu)$.   
We note that $\calI$ is strongly integrable if and only if $$\calI:= \left\{\lambda\in \calL(\mu) : \lambda(a)<0 \text{ for every } a\in C(\calI)\right \}.$$
\begin{definition}\label{resonant fw}
Given a strongly integrable collection $\calI\subset \calL(\mu)$ and $a\in C(\calI)$, let 
$$r_0^{\calI,F}(a):=\max \left\{\frac { \bar \lambda(a) }{ \bar \lambda'(a) } :  \bar \lambda, \bar \lambda'\in \calL^{F,\calI}(\mu)\right\}.$$
Define the \emph{critical fiberwise regularity with respect to $\calI$} to be
\begin{equation}\label{eqcriffiberreg} r_0^F(\calI):= \inf_{a\in C(\calI)} r_0^{\calI,F}(a).\end{equation}

In the case there is $a\in A$ such that every fiberwise Lyapunov exponent is negative at $a$, we also write \begin{equation}
\label{eq:fiberwisereg} r_0^F(\mu)= r_0^{\calL(\mu),F}(\mu)
\end{equation} for the \emph{critical fiberwise regularity} of $\mu$.  
\end{definition}

 \subsubsection{Fiberwise Lyapunov manifolds associated to strongly integrable collections}
As $\mu$ projects to the Haar measure on $G/\Gamma$, the $A$-action on $(M^\alpha,\mu)$ admits dynamically adapted charts as in \cite[\S 3]{MR4599404} (see discussion in \cite[\S 2.4]{MR4502594}).  Intersecting with fibers of the fibration $M^\alpha \to G/\Gamma$, the fiberwise dynamics similarly admits dynamically adapted charts.  In particular, the main results of \cite[\S 4] {MR4599404} provide the following.

\begin{claim}\label{lyapmanifolds}
Assume $r>1$.  Let $\calI$ be a strongly integrable collection of Lyapunov exponents.    
There exists a full $\mu$-measure  subset $\Omega\subset  M^\alpha$ with the following properties:
for every $x\in \Omega$ 
there exist injectively immersed $C^r$ submanifolds $$W^{\calI,F}(x)\subset \calF(x), \quad \text{ and } \quad W^{\calI}(x)\subset M^\alpha$$ tangent to  $E^{\calI,F}(x)$ and $E^{\calI}(x)$, respectively.  

Moreover, for  $x\in \Omega$,
\begin{enumerate}
\item  $W^{\calI,F}(x) $ is the path-connected component of $ W^{\calI}(x)\cap \calF(x)$  containing $x$ and 
$W^{\calI}(x)= V_\calI\cdot W^{\calI,F}(x)$.
\item For $y\in W^{\calI}(x)\cap \Omega$   we have 
$W^{\calI}(y)=W^{\calI}(x)$. 
\item For any  strongly integrable collection $\calI'\subset \calI$ and $x\in \Omega$, $W^{\calI',F}(x)$ is defined and $W^{\calI',F}(x)\subset W^{\calI,F}(x)$ and $W^{\calI'}(x)\subset W^{\calI}(x)$.  
\end{enumerate}
Moreover, for any $a\in A$ and  $\mu$-a.e.\ $x\in\Omega$, 
\begin{enumerate}[resume]
\item  $a\cdot W^{\calI,F}(x) =  W^{\calI,F}(a\cdot x)$ and $a\cdot W^{\calI}(x) =  W^{\calI}(a\cdot x)$
\end{enumerate}

\end{claim}

\subsubsection{Fiberwise stable Lyapunov manifolds} \label{FWstable}
Given $a\in A$, let \begin{equation}\label{stableexps}\calI^s(a) =  \{\lambda\in \calL(\mu): \lambda(a)<0\}\end{equation} be the negative Lyapunov exponents for the action of $a$ on $(M^\alpha, \mu)$.  
Clearly $\calI^s(a)$ is strongly integrable by definition.  
  Write $\Omega(a)$ the set guaranteed in \cref{lyapmanifolds}; for $x\in \Omega(a)$, write $$W^s_{a} (x)= W^{\calI^s(a)}(x)$$  
for the \emph{stable manifold} of $a$ through $x$ and 
 $$W^{s,F}_{a} (x)= W^{\calI^s(a),F}(x)$$  
for the \emph{fiberwise stable manifold} of $a$ through $x$.

\def\cP{J}

\subsection{Leafwise measures}\label{sec:leafwisemeasires}
As above, let $\calI$ be a strongly integrable collection of Lyapunov exponents.  
Write $\scrU^{\calI}$ for the partition of $M^\alpha$ into $U^{\calI}$-orbits.  Similarly let $\scrW^{\calI,F}$ (resp.\  $\scrW^{\calI}$) denote the partitions of $M^\alpha$ into $W^{\calI,F}$-leaves (resp.\ $W^{\calI}$-leaves). 
Let $\fol$ denote either the laminations by $U^{\calI}$-orbits, $W^{\calI,F}$-leaves, or $W^{\calI}$-leaves.  
Write $\fol(x)$ for the atom containing $x$.  

We begin with the following definition.
\begin{definition}[Partitions subordinate to laminations]\label{subor}
A measurable partition $\xi$ is \emph{subordinate} to  $\fol$ if, for $\mu$-a.e.\ $x\in X$, the following hold:
\begin{enumerate}
	\item  $\xi(x)\subset \fol(x)	$;
	\item $\xi(x)$ contains an open (in the immersed topology) neighborhood of $x$ in $\fol(x)$;
	\item $\xi(x)$ is pre-compact (in the immersed topology)  in $\fol( x)$.
\end{enumerate}
\end{definition}
A measurable partition $\xi$ admits a family of conditional probability measures $\{\mu_x ^\xi\}$.  As the atoms of the partition may not be equivariant under the action it is convenient to replace conditional probability measures (depending on the choice of partition $\xi$) with leafwise measures with stronger equivariance properties.  

We note though that  the partition $\fol$ need not be measurable.  Thus the leafwise measures may be infinite (but will still be locally finite) in leaves.  

Given locally finite measures $\nu_1,\nu_2$ on a locally compact metric space $X$, we say $\nu_1\propto \nu_2$ if there is $c>0$ such that $\nu_1=c\nu_2$.


 \begin{proposition}\label{leavwisemsr}
 Let $\eta$ be an $A$-invariant, measurable partition.  
 For $\mu$-a.e.\ $x\in M^\alpha$ there exists a locally finite, Radon (with respect to the immersed topology)
 measure $\mu^{\fol,\eta}_x$ with the following properties:
\begin{enumerate}
	\item  $\mu^{\fol,\eta}_x$  is well-defined up to normalization.
\item For every $a\in A$ and $\mu$-a.e.\ $x$, 
$$a_* \mu_x^{\fol,\eta} \propto  \mu_{a\cdot x}^{\fol,\eta} .$$


	\item For $\mu$-a.e.\ $x$ and $\mu^{\fol,\eta}$-a.e.\ $y\in \fol(x)$, $\mu^{\fol,\eta}_x\propto \mu^{\fol,\eta}_y$.  
	\item For any measurable partition $\xi$ subordinate to $\fol$,  $\mu$-almost every $x$, and $A\subset \xi(x)$ we have
	$$\mu_x ^{\xi\vee\eta} (A) = \dfrac{\mu_x^{\fol,\eta}(A)}{\mu_x^{\fol,\eta}((\xi\vee \eta) (x))}.$$
\end{enumerate}
Moreover, the above properties uniquely determine the family $\{\mu_x^{\fol,\eta}\}$ modulo null sets and modulo a pointwise choice of normalization.  
\end{proposition}
We call the family $\{\mu^{\fol,\eta}_x\}$ a family of \emph{leafwise measures.}

In the sequel, when using \cref{leavwisemsr} we often take $\eta$ to be either 
\begin{enumerate}
\item the trivial partition $\eta= \{M^\alpha\}$ in which case simply write $\mu^{\fol}_x$ rather than  $\mu^{\fol,\eta}_x$; or
\item $\calE_\beta$, the partition into $A'_\beta= \ker \beta$-ergodic components for some root $\beta\colon A\to \R$.  
\end{enumerate}
When $\calI = [\beta]$ and $\eta$ is the trivial partition, we also write in the remainder
\begin{enumerate}
	\item $\mu^{[\beta]}_x= \mu^{\scrW^{\beta}}_x$
	\item $\mu^{[\beta],F}_x= \mu^{\scrW^{\beta,F}}_x$
\end{enumerate}


We have the following standard relationship between entropy and the geometry of leafwise measures.   
\begin{proposition}\label{prop:entropy}
The following are equivalent:
\begin{enumerate}

	\item $\mu^\fol_x$ is finite;
	\item  $\mu^\fol_x\propto \delta_x$ for some $c>0$;
	\item $h_\mu(a\mid \fol)=0$ for all $a\in A$ with $\chi(a)>0$ for all $\chi\in \calI$;
	\item $\fol$ is a measurable partition.  
\end{enumerate}

\end{proposition}

\subsection{Leafwise entropy and product structure of entropy}
Fix $\calI= [\lambda]$ where $\lambda\colon A\to \R$ is a Lyapunov exponent functional.  Let $\fol$ denote either the laminations by $U^{\calI}$-orbits, $W^{\calI,F}$-leaves, or $W^{\calI}$-leaves.  
Write $\fol(x)$ for the atom containing $x$.  

Let $\eta$ be an $A$-invariant measurable partition.  

Fix $a\in A$.  Given a measurable partition $\xi$ of $(M^\alpha,\mu)$, let $\xi^+ = \bigvee _{n=0}^\infty a^n\xi$.   Define the entropy of $a$ with respect to $\xi$ to be $$h_\mu(a, \xi) = H(\xi\mid a\xi^+) = H(a\inv \xi\mid \xi^+).$$
The entropy of $a$ is $h_\mu(a) = \sup_{\xi} h_\mu(a, \xi).$

Let $\eta$ be an $A$-invariant measurable partition.  We will usually assume that 
\begin{enumerate}
	\item $\eta$ is trivial,
	\item $\eta= \calF$ is the partition of $M^\alpha$ into fibers of $M^\alpha\to G/\Gamma$, or 
	\item $\eta = \calE_\beta$, the partition into ergodic components by $A'_\beta$, the kernel of $\beta$ in $A$.  
\end{enumerate}
	The entropy of $a$ with respect to $\eta$ is $h_\mu(a\mid \eta) = \sup_{\xi} h_\mu(a, \xi\vee \eta)$,

Let $a\in A$ be such that $\lambda(a)>0$.  Let $\xi$ be a measurable partition of $M^\alpha$ subordinate to $\fol$ with the property that $a\cdot \xi \prec \xi$. Then $\xi^+= \xi$.  
We define $$h_\mu(a\mid \calF \vee\eta):= h_\mu(a,\xi\vee \eta).$$ 
By \cite[Lemma 8.7]{MR4599404}, this definition is independent of the choice of $\xi$ with the above properties.

We have the following ``product structure" of entropy.
\begin{proposition}[{\cite[Theorem 13.1]{MR4599404}}]\label{ent:prod}
$$h_\mu(a\mid \eta) = \sum_{\chi(a)>0} h_\mu(a\mid \scrW^{\chi}\vee \eta)$$
In particular, 
\begin{enumerate}
\item $h_\mu(a) = \sum_{\chi(a)>0} h_\mu(a\mid \scrW^{\chi})$
\item $h_\mu(a\mid \calF) = \sum_{\chi(a)>0} h_\mu(a\mid \scrW^{\chi,F})$
\end{enumerate}
\end{proposition}

\def\susp{X}

\section{$P$-measures, $A$-measures, and main technical theorems}\label{sec4}
We assume $G$ is a higher-rank semisimple Lie group as in \cref{ss:Gassump} and $\Gamma$ is an irreducible lattice subgroup.  
Fix a compact manifold $M$ and $\alpha\colon \Gamma\to \diff^r(M)$ for some $r>1$.

Fix an Iwasawa decomposition  $G=KAN$.  Consider the standard minimal parabolic subgroup $P= MAN$ relative to this choice of Iwasawa decomposition.  Given an action $\alpha\colon \Gamma\to \diff^r(M)$, we consider the $G$-action on the suspension space $M^\alpha$.  Throughout, we study  Borel probability measures on $M^\alpha$ that are either 
\begin{enumerate}
\item invariant under the action of the parabolic subgroup $P\subset G$; or
\item invariant under the action of the  subgroup $A\subset G$ and project to the Haar measure on $G/\Gamma$.  
\end{enumerate}
We remark that every $P$-invariant Borel probability measure $\mu$ on $M^\alpha$ necessarily projects to the Haar measure on $G/\Gamma$, the unique $P$-invariant Borel probability measure on $G/\Gamma$.   

The results announced in \cref{sec1,sec2}  follow from classifying such probability measures on the suspension space $M^\alpha$.  In this section we state the main technical results of the paper.  Proofs of results announced in \cref{sec1,sec2} are given in  \cref{sec6}, following the results in this section and additional technical lemmas in \cref{sec5}.

\subsection{Tameness of leafwise meaures} In the following results, we will need a technical property of the leafwise measures along leaves of certain foliations on $M^\alpha$.  Sufficient criteria for this property are given in  \cref{tame1,tame2}.

Let $\mu$ be an  $A$-invariant Borel probability measure on $M^\alpha$ that projects to the Haar measure on $G/\Gamma$. 
Consider a restricted root $\beta$ that is resonant with a fiberwise coarse Lyapunov exponent.  

As introduced in \cref{sec:leafwisemeasires}, let  $\{ \mu^{[\beta]}_x\}$ denote the leafwise measures along $W^{[\beta]}$-leaves.  
By entropy considerations (see \cref{prop:entropy}), the measures $ \mu^{[\beta]}_x$ are always infinite measures.  
Thus the measures $ \mu^{[\beta]}_x$ are only defined up to normalization and are only projectively equivariant under the $A$-action.  However, in certain situations, we may able to find a canonical choice of normalization of the measures $ \mu^{[\beta]}_x$ such that the normalized measures are invariant under the restriction of the $A$-action to $A'_{\beta}= \ker \beta$, the kernel of $\beta$ in $A$.

\begin{definition}\label{defn:tamedyn} 
We say the measure $\mu$ is \emph{$\beta$-tame} if for $\mu$-a.e.\ $x$, there exists an  open neighborhood $G(x)$ of $x$ in  $W^{[\beta]}(x)$ such that 
\begin{enumerate}
\item the  restriction  of the leafwise measure $ \mu^{[\beta]}_x$ to $G(x)$ is finite, and 
\item if $\bar \mu^{[\beta]}_x = \frac{ 1}{ \mu^{[\beta]}_x({G(x)} )} \restrict { \mu^{[\beta]}_x}{G(x)} $
 denotes the normalized restriction of  $ \mu^{[\beta]}_x$ to $G(x)$, then 
 for every $a\in A'_\beta=\ker \beta$ and a.e.\ $x$, 
 $$a_*  { \bar \mu^{[\beta]}_x} =   {\bar  \mu^{[\beta]}_{a\cdot x}}.$$   
\end{enumerate}
\end{definition}
We remark that we do not assume $G(x)$ is precompact in definition \cref{defn:tamedyn};  in practice, the   requirement that $ \mu^{[\beta]}_x$ restricted to $G(x)$ be finite is often the most difficult criteria to verify.  

We also remark that we do not assume $\mu$ to be $A$-ergodic in definition \cref{defn:tamedyn}.  We thus naturally define when a $P$-invariant measure $\mu$ is $\beta$-tame.  
See \cref{tame1,tame2} for sufficient conditions that imply $\beta$-tameness.  

\subsection{Constraints on resonant fiberwise Lyapunov data}
Given an ergodic, $A$-invariant Borel probability measure $\mu$ on $M^\alpha$, we consider a resonant root $\beta$ for which $\mu$ is not invariant under any 1-parameter subgroup of $U^{[\beta]}$.  Under various conditions, we obtain lower bounds on the dimension and resonance coefficients of the associated fiberwise coarse Lyapunov exponents. 


Assuming low regularity of the action, in \cite{ABZ} 
it is shown that coarse fiberwise Lyapunov exponents associated to higher-dimensional root spaces for resonant roots must have dimension at least 2.  (This also follows from the tools in this paper).
\begin{proposition}[{\cite[Prop.\ 5.1]{ABZ}}]
\label{thm:parabolicmeasrdim2} 
Suppose $r>1$.   Let $\mu$ be an ergodic, $A$-invariant Borel probability measure that projects to the Haar measure on $G/\Gamma$.    

Then for every $\beta\in \Sigma$ such that   $\mu$ is not  invariant under any 1-parameter subgroup of $U^{[\beta]}$ there exists a  fiberwise Lyapunov exponent $\lambda_i^F\in \calL^F(\mu)$  such that 
\begin{enumerate}
\item \label{para1} $[\lambda_i^F]= [ \beta]$ and 
\item \label{para2} $\dim (E^{[\lambda_i^F]}_x) \ge \min\{2, \dim \lieg^{[\beta]}\}$ for almost every $x$.
\end{enumerate}
\end{proposition}

To state an extension of the above result, we define a critical  regularity associated to the coarse Lyapunov exponents of the $A$-action on $(M^\alpha, \mu)$.  

\begin{definition}[Coarse critical fiberwise regularity]  \label{coarse fiber rigidity}
Let $\mu$ be an ergodic, $A$-invariant Borel probability measure on  $M^\alpha$.  
Given a coarse Lyapunov exponent $\chi\in \what \calL(\mu)$ with $\chi\cap  \what \calL^F(\mu)\neq\emptyset$, define the  \emph{critical fiberwise coarse regularity of $\chi$} to be
\begin{equation}
c_0^F(\chi):= \max \left\{\frac { \lambda_j^F }{ \lambda_{j'}^F} : \lambda_j^F, \lambda_{j'}^F \in \chi \cap  \what \calL^F(\mu)\right\}.
\end{equation}
If $\chi\cap  \what \calL^F(\mu)=\emptyset$, define $c_0^F(\chi)= 1$.

Given a strongly integrable subcollection $\calI\subset \calL(\mu)$, let $\what \calI$ denote the equivalence classes in $\calI$.  Define the \emph{critical fiberwise coarse regularity $\calI$} to be 
\begin{equation} c_0^F(\mu,\calI):= \max \left\{ c_0^F(\chi): \chi\in \what \calI\right\}.\end{equation}
 Define the \emph{critical fiberwise coarse regularity $c_0^F(\mu)$} to be
$$c_0^F(\mu):= c_0^F(\mu,\calL(\mu)).$$
\end{definition}
We note that if every (non-trivial) coarse fiberwise Lyapunov subspace $E^{\chi,F}(x)$ is 1-dimensional, then $c_0^F(\mu)=1$; this is the primary reason we are able to obtain results in the $C^{1+\text{\Holder}}$-setting for actions by lattices in $\SL(n,\R)$.


We improve the dimension bound and also the resonance coefficients in \cref{thm:parabolicmeasrdim2} for resonant roots $\beta$ such that $\mu$ is not invariant under any 1-parameter subgroup of $U^{[\beta]}$---under the additional assumptions that (1) the action $\alpha\colon\Gamma\to  \Diff^r(M)$ is sufficiently regular (so that $r> c_0^F([\beta])$), and (2) the measure $\mu$ is $\beta$-tame.  
Specifically,  if $\mu$ is not  invariant under any 1-parameter subgroup of $U^{[\beta]}$ then   \cref{NRimpInv} implies there exists a fiberwise Lyapunov exponent $\lambda_i^F$ with $\lambda_i^F=c\beta$ for some $c>0$; our next result, \cref{thm:parabolicmeasr}, implies we can take $c=1$ and that the dimension of the associated fiberwise Lyapunov subspaces are at least $\dim \lieg^\beta$.

%

\begin{proposition}\label{thm:parabolicmeasr}  
  Let $\mu$ be an ergodic, $A$-invariant  Borel probability measure  on $M^\alpha$ that  projects to the Haar measure on $G/\Gamma$.

Suppose $\beta\in \Sigma$ is a restricted  root such that $\mu$ is not  invariant under any 1-parameter subgroup of $U^{[\beta]}$ 
and such that $\mu$ is $\beta$-tame.   Suppose  $r> c_0^F([\beta])$.  Then there exists a  fiberwise Lyapunov exponent $\lambda_i^F\in \calL^F(\mu)$  such that 
\begin{enumerate}
\item \label{para1} $\lambda_i^F=  \beta$;
\item \label{para2} $\dim (E^{\lambda_i^F}(x)) \ge \dim \lieg^{\beta}$ for almost every $x$.
\end{enumerate}
\end{proposition}

The proof of \cref{thm:parabolicmeasr} appears in \cref{Pf44}.   

We have the standard argument, following from \cref{lem:alge,ergcomp} below, which implies for an ergodic $P$-measure $\mu$ that if $\mu$ is  invariant under some 1-parameter subgroup of $U^{[\beta]}$ then $\mu$ is $U^{[\beta]}$-invariant.  
As an immediate corollary,  when the dimension of the fiber is sufficiently small,  \cref{thm:parabolicmeasr} greatly constrains the combinatorics of the fiberwise dynamics for $P$-invariant measures.    Recall the critical  regularities  $ r_0 (\lieq)$ and $r_0^F(\mu) $ in \cref{eq:parareg,eq:fiberwisereg}
\begin{corollary} \label{cor:easycor} 
 Let $\mu$ be an ergodic,  $P$-invariant Borel probability measure on $M^\alpha$.  Let $Q=\stab(\mu)$ and suppose that $r>c_0^F(\mu)$. 
Suppose further that
\begin{enumerate}[label=(\alph*), ref=(\alph*)]
\item $\mu$ is $\beta$-tame for every $\beta \in \Sigma_Q^\perp$, and 
\item $\dim M= \dim G-\dim Q$. 
\end{enumerate}
Then the following hold:
\begin{enumerate}
\item \label{eqgggg1} For every fiberwise Lyapunov exponent $\lambda_i^F\in \calL^F(\mu)$ there exists  $\beta \in \Sigma_Q^\perp$  with $\lambda_i^F= \beta$ and $\dim (E^{\lambda_i^F}(x)) = \dim \lieg^{\beta}$ for almost every $x$.  
\item  \label{eqgggg2} $r_0^F(\mu) = r_0 (\lieq)$.  
\item  \label{eqgggg3} The measure $\mu$ has finitely many $A$-ergodic components $\{\mu_1,\dots, \mu_\ell\}$; moreover, for any $a\in A\sm\{\1\}$ and $1\le j\le \ell$, the action of $a$ on $(M^\alpha,\mu_j)$ is ergodic.  
\end{enumerate}
\end{corollary}
\begin{proof}
It follows from \cref{lem:alge,ergcomp} below, that $\mu$ is not invariant under any 1-parameter subgroup of $U^{[\beta]}$ (since otherwise, $\mu$ would be $U^{[\beta]}$-invariant). 
Passing to $A$-ergodic components of $\mu$,  conclusion \eqref{eqgggg1} follows directly from  \cref{thm:parabolicmeasr} and  dimension counting.  
Conclusion \eqref{eqgggg2} then follows.  

It follows for any $a\in C(\Sigma_Q^\perp)$, that all fiberwise Lyapunov exponents are negative.  Thus, $(M^\alpha,\mu) \to (G/\Gamma, m_{G/\Gamma})$ is measurably finite-to-one (c.f.\ \cite[Lemma 5.6]{MR4502594}).  Since $A$ acts ergodically on $(G/\Gamma, m_{G/\Gamma})$ it follows that the $A$-ergodic decomposition of $\mu$ is finite and every ergodic component $\mu_j$ projects to $m_{G/\Gamma}$ on $G/\Gamma$.  Similarly, if $\mu_j$ is an $A$-ergodic component and if $a\neq \1$, there are finitely many $a$-ergodic components of $\mu_j$.  Thus $A$ acts transitively on the $a$-ergodic components of $\mu_j$; since $A$ is connected, this implies there is a single $a$-ergodic component of $\mu_j$ and thus $a$-acts ergodically on $\mu_j$.  
%
%
%
\end{proof}

\subsection{Global classification from dynamics of $P$-measures} 
We  state our first main classification theorems in the paper.   \cref{slnbound,thm:globalrigid2} above follow directly from this theorem.


\begin{theorem}\label{thm:main}  \label{thm:completepara}
Let $\mu$ be an ergodic, $P$-invariant Borel probability measure on $M^\alpha$.  Let $P\subset Q\subset G$ be the  parabolic subgroup with $\stab(\mu) = Q$.
Suppose the following hold:
\begin{enumerate}[label=(\alph*), ref=(\alph*)]
\item \label{hypproja} $r>  \max\{c_0^F(\mu),r_0(\lieq)\}$; 
\item \label{hypprojb} $M$ is connected;
\item \label{hypprojc} there is $a\in A$ such that $\lambda_i^F(a)<0$ for every fiberwise Lyapunov exponent $\lambda_i^F\in \calL^F(\mu)$;  
\item \label{hypprojd} $\dim M= \dim G -\dim Q$.
\end{enumerate}

Then  there exists a  finite-to-one, $C^r$ covering map $h\colon M\to  G/Q$ such that for every $\gamma\in \Gamma$, $$h\circ \alpha(\gamma) =  \gamma \cdot h$$
where the action on the right-hand side is the standard left projective action of $\Gamma$ on $G/Q.$
\end{theorem}

The proof of \cref{thm:main} appears in \cref{pf4647}.

\subsection{Laminated, relatively measure preserving subsets from dynamics of $P$-measures}

 \cref{thm:embedded projectives}  on  volume-preserving actions by lattices in $\R$-split simple Lie groups
follows as a direct corollary of  \cref{thm:mainlamina} below.  We note that in \cref{thm:mainlamina}, we assume $G$ is $\R$-split though we do not require that $Q$ be maximal. 
We note that while we are only able to verify the hypotheses of \cref{thm:mainlamina} in the setting of volume-preserving actions $\alpha\colon \Gamma\to \diff^r(M)$ in \cref{thm:embedded projectives}, we expect \cref{thm:mainlamina} may be of use in future results.

To formulate  \cref{thm:mainlamina},  
consider two standard Borel spaces $X,Y$ equipped with Borel $\Gamma$-actions.  We say a measure $\mu$ on $X$ is \emph{quasi-invariant} if $\gamma_*\mu $ is in the same measure class as $\mu$ for every $\gamma\in \Gamma$.  
Let  $\mu$ and $\nu$ be  quasi-invariant measures, respectively, for the $\Gamma$ actions on $X$ and $Y$.  
Let $h\colon X\to Y$ be a $\Gamma$-equivariant, measurable map with $h_*\mu= \nu$.  
 Let $\{\mu_y\}$ denote a family of conditional probability measures relative to the disintegration of $(X,\mu)$ into preimages under $h$.  
We say that $h\colon (X,\mu)\to (Y,\nu)$ is a \emph{relatively measure-preserving extension} if  for $\nu$-a.e.\ $y$ and every $\gamma\in \Gamma$,
 $$\gamma_*\mu_y = \mu_{\gamma \cdot y }.$$
 
 In our previous paper \cite[Theorem 1.10]{MR4502594}, we showed   for actions  of higher-rank lattices on manifolds of certain dimensions, there always exists  a quasi-invariant measure $\mu$ and a measurable equivariant $h\colon (M,\mu)\to (G/Q,\mathrm{Leb})$ such that $h$ is  a relatively measure-preserving extension of the standard projective action of $\Gamma$ on $(G/Q,\mathrm{Leb})$.  In some cases,  $h$ can be improved to a $\Gamma$-equivariant $C^r$ covering map $h\colon M\to G/Q$ as in 
\cref{thm:completepara}.  In other settings---including those studied in \cref{thm:embedded projectives}---we can improve the relative measure-preserving extension by showing the  measure $\mu$  is supported on finitely many embedded copies of (covers of) $(G/Q,\mathrm{Leb})$ in $M$ (assuming $|\pi_1(G/Q)|<\infty$) or on a laminated object $\Lambda$ whose leaves cover $G/Q$ and improve $h$ to a $C^r$ coving map when restricted to these objects.

Although we expect the following to hold in more generality, for technical reasons arising in the proof, we require that the group $G$ has no rank-1 factors and that the coarse fiberwise Lyapunov exponents (corresponding to exponents of the projective factor)  are 1-dimensional (since we can not verify $\beta$-tameness for higher-dimensional fiberwise Lyapunov subspaces in this setting). 
We thus only formulate the following result for lattices for higher-rank,  $\R$-split simple Lie groups.  

For the statement of the following, recall the cone associated to a (strongly integrable) collection of  roots given in \eqref{eq:negcone} and the critical regularity in \eqref{eq:parareg}
\begin{theorem}\label{thm:mainlamina} 
Let $G$ be a $\R$-split simple Lie group with real rank at least $2$.  Let $\Gamma\subset G$ be a lattice subgroup and let $\alpha\colon \Gamma \to \diff^r(M)$ be an action.    

Let $Q$ be a   parabolic subgroup of $G$ and suppose $r>r_0(\lieq).$ 
Let $\mu$ be an ergodic, $Q$-invariant  Borel probability  measure on $M^\alpha$.  Suppose 
\begin{enumerate}[label=(\alph*), ref=(\alph*)]
\item $\mu$ is not $Q'$-invariant for any $Q\subsetneq Q'\subset G$,
\end{enumerate}
and there is a partition $$\what \calL^F(\mu)= S  \sqcup  U$$ of the coarse fiberwise Lyapunov exponents $\what \calL^F(\mu)= \{\chi_1^F, \dots, \chi_\ell^F\}$ such that 
\begin{enumerate}[label=(\alph*), ref=(\alph*), resume]
\item\label{exptan}  for every $\chi_i^F\in S$, there is $\beta \in \Sigma_\lieq^\perp$ such that $\chi_i^F$ is positively proportional to $\beta$ and $1= \dim E^{\chi_i^F,F}(x)$ for almost every $x$, and 
\item\label{exptrans}    $\chi_i^F(a)>0$ for every  $a\in \exp(\overline{C(\Sigma_\lieq ^\perp)}\sm \{0\})$ and 
 every coarse fiberwise Lyapunov exponent $\chi_j^F\in U$. 
\end{enumerate}

Then  there exists a finite,  $\alpha(\Gamma)$-quasi-invariant measure $\nu$ on $M$, an $\alpha$-invariant,  $\nu$-measurable lamination $\Lambda\subset M$ with $C^r$-leaves, and a measurable function $h\colon \Lambda\to  G/Q$  with the following properties: 
\begin{enumerate}
\item  \label{laminate1}  The map $h$ is $\Gamma$-equivariant: for all $\gamma\in \Gamma$ and $x\in \Lambda$,   
 $$h(\alpha(\gamma)(x)) = \gamma \cdot h(x).$$
\item $\nu$ is supported on $\Lambda$.
\item $h_*\nu$ is in  the Lebesgue class on $G/Q$ and $h\colon (\Lambda,\nu)\to (G/Q,\mathrm{Leb})$ is a relatively measure preserving extension.
\item for each leaf $L(x)$ of $\Lambda$, the restriction of $h$ to $L(x)$ is a $C^r$-covering map $\restrict{h}{L_x} \colon L(x)\to G/Q$;  
\item The measure $\mu$ is supported on the suspension $(G\times \Lambda)/\Gamma$. 
\end{enumerate}
Moreover, there is a left $G$-action on the set $\Lambda$, denoted by $g\cdot x = 
\ell_g(x)$,  with the following properties:
\begin{enumerate}[resume]
 \item  The left $G$-action preserves each leaf $L(x)$ of $\Lambda$ and restricts to a $C^r$ action on each leaf $L(x)$.
 \item The  left $G$-action is $h$-equivariant: $h(\ell_g(x)) = g \cdot h(x)$ for all $g\in G$ and $x\in \Lambda$.
 \item \label{laminate8} Given $\gamma\in \Gamma$ and $x\in \Lambda$, write  $y= \ell_{\gamma\inv}(\alpha(\gamma)(x)) $ so that $\alpha(\gamma)(x) = \ell_{\gamma}(y)$.  Then for all $g\in G$, $$\alpha(\gamma)(\ell_g(x))= \ell_{\gamma g }(y).$$
\end{enumerate}
Moreover, assuming $\pi_1(G/Q)$ is finite  
then 
\begin{enumerate}[resume]
	\item \label{laminate9}  each leaf $L(x)$ is a compact embedded submanifold and $\Lambda$ is the union of finitely many leaves.  
\end{enumerate} 
\end{theorem}
The proof of \cref{thm:mainlamina}  appears in \cref{pf4647}.

\begin{remark}\label{rem:hypo}
Hypotheses \ref{exptan} and \ref{exptrans} of \cref{thm:mainlamina} are necessary for the proof.  We note that while these hypotheses are stated in terms of the dynamics of the suspension action, various assumptions on the action $\alpha\colon \Gamma\to \diff(M)$ ensure hypotheses \ref{exptan} and \ref{exptrans}  hold for the induced action.  

Indeed, let $G$ be an $\R$-split simple Lie group with real rank at least $2$, let $\Gamma$ be a lattice subgroup, and let $\alpha\colon \Gamma \to \diff^r(M)$ be an action.  Let $Q$ be a parabolic subgroup of $G$ with Lie algebra $\lieq$.
Suppose the following hold:
\begin{enumerate}
\item the action $\alpha$ preserves a no-where vanishing smooth density;
\item $\dim M = \dim G -\dim Q + 1$;
\item $\mu$ is an  ergodic,  $Q$-invariant measure on $M^\alpha$ that is  is not  $Q'$-invariant for any subgroup $Q'\supsetneq Q$.  
\end{enumerate}

Then $\mu$ satisfies the hypotheses of \cref{thm:mainlamina}.  Indeed, let $\calL^F(\mu)=\{\lambda_1^F,\dots, \lambda_\ell^F\}$ denote the distinct fiberwise Lyapunov exponents.  Using \cref{NRimpInv}, that $\dim M=\dim G -\dim Q + 1$, and that $\alpha$ is volume preserving, up to reindexing,  we have for any $a\in \exp(\overline{C(\Sigma_\lieq ^\perp)}\sm \{0\})$  that 
\begin{enumerate}
\item $\ell= \dim M=\dim G -\dim Q + 1$,
\item $\lambda_j^F$ is positively proportional to some $\beta$ for $\beta\in \Sigma_Q^\perp$ whence $\lambda_j^F(a)\le 0$ for all $1\le j\le \ell-1$,
\item  $\sum_{j=1}^{\ell-1}\lambda_j^F(a)< 0$,
\item $\lambda_\ell^F(a)=-\sum_{j=1}^{\ell-1}\lambda_j^F(a)>0$.  
\end{enumerate}
We thus obtain partition $\calL^F(\mu)=S  \sqcup  U$ where $S=\{\lambda_j^F : 1\le j\le \ell-1\}$ and $U= \{\lambda_\ell\}$.  
In particular, for   actions as in \cref{thm:partpara}, for the induced action, hypotheses \ref{exptan} and \ref{exptrans}  hold for any $P$-invariant measure  that is not $G$-invariant.  
\end{remark}

\subsection{Examples related to \cref{thm:mainlamina}}
To further explain the conclusions and hypotheses of  \cref{thm:mainlamina}, we discuss some example actions and discuss the existence or non-existence of (non-trivial) invariant laminations as in the conclusions of  \cref{thm:mainlamina}.  

First, to see that hypotheses \ref{exptan} and \ref{exptrans} in  \cref{thm:mainlamina} are reasonable, consider the following action.  
\begin{example} Given  $x=[v]\in \R P^{n-1}$ and $g\in \Sl(n,\R)$, let $$c(g,x):= \log(\|gv\| \cdot \|v\|\inv).$$  Then $c\colon G\times \R P ^{n-1}\to \R$ defines an additive cocycle over the $G$-action on $\R P^{n-1}$.  Given $\theta \in \R$, write $r_\theta\colon S^1\to S^1$ for the rotation by $\theta$ on $S^1$.  
On $M=\R P^{n-1}\times S^1$, let $G$ act as $$g\cdot (x, t) = \bigl(g\cdot x, r_{c(g,x)}(t)\bigr).$$

Given a lattice $\Gamma\subset \Sl(n,\R)$, let $\alpha\colon \Gamma\to \Diff^\omega(\R P^{n-1}\times S^1)$ be the restriction to $\Gamma$ of the $G$-action.  We claim the following.

\begin{claim*}For $n\ge 2$, there is no $\alpha(\Gamma)$-invariant lamination $\Lambda$ on $\R P^{n-1}\times S^1$ that projects equivariantly to the $\Gamma$-action on  $\R P^{n-1}$.  
\end{claim*}  
Indeed, we may select $\gamma\in \Gamma$ that acts as a Morse-Smale diffeomorphism on $\RP^{n-1}$.  Rotating the fibers by a fixed angle $\theta_0$, we may assume there is an $\alpha(\gamma)$-invariant leaf $\calL\subset \R P^{n-1}\times S^1$.  Passing to a double cover if needed, the leaf $\calL$ intersects each fiber at most once.  We may lift the action and (adding an integer translation to $\theta_0$ if needed) find an invariant leaf on $\R P^{n-1}\times \R$ (or $S^{n-1}\times \R$); this is the graph of a function $\beta\colon \R P^{n-1}\to \R$, (or $\beta\colon S^{n-1}\to \R$) satisfying 
$$c(\gamma,x) -r_{\theta_0}= \beta(x) - \beta(\gamma\cdot x).$$
However, the cocycle $c(\gamma, \cdot)$ can not be cohomologous to a constant since it diverges under iteration $c(\gamma^n, x_i)$ to $\pm \infty$ at the source $x_0$ and sink $x_1$, respectively, for $\gamma$.

We remark that given any $\alpha(\Gamma)$-equivariant lamination on $\R P^{n-1}\times S^1$ whose leaves are equivariantly conjugated to the action on $\R P^{n-1}$, one could use finitely many Morse-Smale elements of the $\Gamma$-action on $\R P^{n-1}$ to show the leaves are uniformly transverse to the $S^1$-fibers and that the fibration maps the leaves of the lamination equivariantly to the $\R P^{n-1}$ action.  
\end{example}

To see that non-trivial laminations can occur, at least abstractly, we have the following.

\begin{example}\label{exam:SP}
Let $G= \wtd{\Sp}(n,\R)$, let $\pi\colon G\to \Sp(n,\R)$ be the natural covering map, let $\hat \Gamma$ be a lattice in $\Sp(n,\R)$, and let $ \Gamma = \pi\inv (\hat \Gamma)$.  Let $\widehat Q$ be a parabolic in $\Sp(n,R)$, let $ Q= \pi\inv (\widehat Q)$, and let $Q^\circ$ be the connected component of the identity in $Q$.
The group $\Gamma$ is a lattice in $G$ and has finite index in the center of $G$.  Let $Z=\langle z\rangle$ be a copy of $\Z$ contained in both the center of $G$ and $\Gamma$. 

Let $X$ be any compact topological space and $f\colon X\to X$ a homeomorphism.  
On $G/Q^\circ \times X$, let $Z$ act by $$(gQ^\circ ,x)\cdot z^n= (gz^nQ^\circ , f^n(x))= (z^ngQ^\circ , f^n(x))$$
and let $G$ act by $$g'(gQ^\circ ,x) = (g'gQ^\circ ,x).$$
Let $$\Lambda = (G/Q^\circ \times X)/Z.$$
The $G$-action on $G/Q^\circ \times X$ descends to a $G$-action on $\Lambda$ preserving each $G$-orbit.   We may then restrict to an action of $\Gamma$. Thus non-trivial laminated objects $\Lambda$ admitting  $\Gamma$-actions exist.

As concrete examples of the above construction, we have the following:
\begin{enumerate}
\item Given a prime $p$, let $X=\Z_p$ be the $p$-adic integers.  Let $f\colon X\to X$ be $f(x) = x+1$.  
Then $\Lambda = (G/Q\times \Z_p)/Z$ has the structure of a solenoid.  
\item Let $X= S^1$ and let $f(x) = x+\alpha$ for $\alpha \in \R\sm \Q$.  We then obtain  and action on a $S^1$ bundle over $G/(Q^\circ Z)$.  
\end{enumerate}
\end{example}
\begin{question}
	Can the $\Lambda$ from \cref{exam:SP} be embedded in a manifold $M$ as an invariant subset for a smooth action of $\Gamma$ on $M$?  Is it possible to embed the action in  a volume-preserving way
or in a way such that there is a  $P$-measure which is not a $G$-measure and satisfies the hypotheses of  \cref{thm:mainlamina}? %
\end{question}

\begin{question}
Given a smooth $f\colon X\to X$ generating a lamination as in \cref{exam:SP}, consider the induced $\Gamma$-action $\alpha$ on $M=\Lambda$ and the induced $G$-action on $M^\alpha$.  What are the possible $Q$-invariant measures on $M^\alpha$.  Is it possible to build a choice of $f$ admitting a $Q$-invariant measure which satisfies the hypotheses of  \cref{thm:mainlamina}?
\end{question}

\begin{example}
	Let $\hat \Gamma$ be a lattice in $\SL(2,\R)$ isomorphic to the free group $F_2=\langle a,b\rangle$. 
	Picking lifts $\td a$ and $\td b$ of $a$ and $b$ in $G=\wtd \SL(2,\R)$ we have $\Gamma = \pi\inv (\hat \Gamma)= \Z \times F_2$, where $\Z$ is the center of $G$ and $F_2= \langle \td a,\td b\rangle$. 
	
	Again, let $X$ be any compact topological space and let $f\colon X\to X$ be a homeomorphism.  
On $G/Q\times X$, let $Z$ act by $$(gQ,x)\cdot z^n= (gz^nQ, f^n(x))= (z^ngQ, f^n(x))$$
and let $\gamma\in \Gamma$ of the form $\gamma= (z^n,c)$ act by $$(z^n, c) \cdot (gQ,y) =  (cgzQ, f^n(y)).$$ 
	We check again that these two actions commute, and  hence we obtain a $\Gamma$-action on $\Lambda = (G/Q\times X)/Z$.
	However, unlike in \cref{exam:SP}, the $\Gamma$-action does not extend to an action of $G$ and acts by non-trivial permutations (assuming $f$ is not periodic) on the leaves of $\Lambda.$

Moreover, we obtain an action of $\hat \Gamma$ permuting leaves by twisting with any morphimsm $\hat \Gamma \to Z$.  

\end{example}

\section{Technical lemmas on $P$-measures and $A$-measures} \label{sec5}
To establish the results in \cref{sec1,sec2}, we need to verify the hypothesis of our main results in \cref{sec4}.  
In this section we study the $P$- and $A$-actions on the suspension space $M^\alpha$ and catalog various properties of invariant probability measures on $M^\alpha$, assuming  various hypotheses related to   the $A$-action on $M^\alpha$.  

Under  various dynamical assumptions, we are able to verify that an ergodic,  $A$-invariant Borel probability measure $\mu$ projecting to the Haar measure on $G/\Gamma$ is $\beta$-tame.  

\subsection{$A$-ergodic components of $P$-measures}
Recall we have $P= MAN$.  We collect the following facts about $A$-ergodic components of $P$-invariant measures.  
\begin{lemma}\label{ergcomp}
Let $\mu$ be an ergodic, $P$-invariant Borel probability measure on $M^\alpha$.  Let $\calE=\{\mu'_\alpha\}$ denote the decomposition of $\mu$ into $A$-ergodic components.  Then for $\mu$-a.e.\ ergodic component $\mu'_\alpha$,
\begin{enumerate}
\item $\mu'_\alpha$ projects to the Haar measure on $G/\Gamma$;
\item \label{erg2} $\mu'_\alpha$ is $N$-invariant.  
\end{enumerate}
\end{lemma}
For conclusion \eqref{erg2}, see \cite[Lem.\ 4.3]{MR4502594}.  

\subsection{Suspension of  projective actions and $P$-measures}
Consider the setting that $M= G/Q$ and $\rho\colon \Gamma\to \Diff^\infty(G/Q)$ is the standard projective action.  We may follow the suspension construction outlined in \cref{sss:susp} for the action $\rho$ and obtain a $G$-action on $M^\rho$.  
However, as the $\Gamma$-action on $G/Q$ extends to a $G$-action, we  give an alternative construction of the suspension space relative to which the dynamical properties of the $G$-action become more transparent. 

 On the product $G\times M$ define left $G$- and right $\Gamma$-actions 
\begin{equation}\label{flat}  g'\cdot (g,x) = (g'g,g'x)\quad \quad (g,x)\cdot \gamma = (g\gamma, x).\end{equation}
The quotient by $\Gamma$ is 
$$ G/\Gamma\times G/Q$$ and inherits a natural $G$-action given by $$g'(g, x) = (g'g, g'x).$$

If  $\td\Phi\colon G\times M\to G\times M$ denotes  the map $$\td\Phi(g,x) = (g, g\inv x)$$
then   $\td\Phi$ intertwines the $G$- and $\Gamma$- actions induced by $\rho$ as defined by \eqref{eq:susp} and \eqref{flat} and hence descends to a $G$-equivariant  homeomorphism $\Phi\colon G/\Gamma\times G/Q\to M^\rho$.  

By considering the $P$-action on $G/Q$, we observe the following:
\def\haar{\mathrm{Haar}}
\begin{lemma}\label{uniquePmeasure}
There exists a unique $P$-invariant measure $\mu_0$ on $G/\Gamma\times G/Q$ given by $$\mu_0 = \haar\times \delta_{\1_G Q}.$$
Moreover, $\mu_0$ is $Q$-invariant but is not $Q'$-invariant for any $$Q'\supsetneq Q.$$
\end{lemma}
Indeed, it is well-known that  $\delta_{\1_G Q}$ is  the unique $P$-invariant Borel probability measure on $G/Q$.  
The lemma then follows as the projection $G/\Gamma\times G/Q\to G/Q$ is $G$-equivariant and $\mu_0$ is the unique $P$-invariant  Borel probability  measure supported on the fiber of $G/\Gamma\times G/Q\to G/Q$ over $\1_GQ$.


%
%
%

\def\calR{\mathcal R}
\subsection{$P$-measures for perturbations of projective actions}
In the following, our primary interest is in the case that  $X= M= G/Q$ and our actions $\alpha\colon \Gamma\to \Diff(M)$   are small perturbations of the standard projective action $\rho\colon  \Gamma\to \Diff(G/Q).$   

\subsubsection{Properties of $C^0$ perturbations} 
We first   establish some properties that hold for $C^0$ actions of $\Gamma$.  
Let $X$ be a compact metric space. 
Fix a  Haar measure $m$ on $G$ normalized so that $\Gamma$ has co-volume 1.  Let $\calR$ denote the space of locally finite Borel (and hence Radon) measures on $G\times X$ whose projection to $G$ is $m$.  Equip $\calR$ with the (metrizable) topology dual to compactly supported continuous functions.   
The space $\calR$ is sequentially compact and thus compact and  also  invariant under the $G$-action $g\cdot (g',x) = (gg',x)$.

Let $\alpha_n\colon \Gamma\to \Homeo (X)$, $n\ge 1$, be a sequence of actions converging to an action $\alpha_0\colon \Gamma\to \Homeo(X)$ in the $C^0$ topology.  
For each $n\ge 0$, let $\wtd \alpha_n\colon \Gamma\to \Homeo (G\times X)$ denote the right-$\Gamma$-action $$\wtd \alpha_n(\gamma)   (g, x) = (g\gamma, \alpha_n(\gamma\inv) (x)).$$
We have the following.  
\begin{claim}\label{felatiowithgoats}
Let $\{\mu_n\}$ be a sequence in $\calR$ converging to $\mu_0$.  Fix $h\in G$ and a sequence of actions $\alpha_n\colon \Gamma\to \Homeo (X)$, $n\ge 1$, converging to $\alpha_0\colon \Gamma\to \Homeo (X)$.  
\begin{enumerate}
\item \label{stupid1}   We have $h_* \mu_n \to h_*\mu_0$ and $\wtd\alpha_n(\gamma)_*\mu_n \to \wtd \alpha_0 (\gamma)\mu_0$ for every $\gamma\in \Gamma$.   
\item \label{stupid2}If each $\mu_n$ is $h$-invariant then $\mu_0$ is $h$-invariant. 
\item \label{stupid3} If each $\mu_n$ is $\wtd\alpha_n$-invariant  then $\mu_0$ is $\wtd \alpha_0$-invariant.
\end{enumerate}
\end{claim}
\begin{proof}
For \eqref{stupid1}, fix a  compactly supported continuous function $f\colon G\times M\to \R$.  Since the action of $h$ is continuous and proper, the function $(g,x)\mapsto f(hg,x)$ is compactly supported whence 
$$ \lim_{n\to \infty} \int f \, d( h_*\mu_n) :=\lim _{n\to \infty} \int  f(hg,x) \, d \mu_n =  \int  f(hg,x) \, d \mu_0 =: \int f \, d(h_* \mu_0).$$
Similarly, the sequence of functions  $$f_n \colon (g,x)\mapsto f \left(\wtd \alpha_n(\gamma)(g,x)\right)$$
converges uniformly to the compactly supported continuous function $ f_0 \colon (g,x)\mapsto f \left(\wtd \alpha_0(\gamma)(g,x)\right)$.  Moreover there is a compact $K\subset G\times X$ such that $\supp (f_n)\subset K$ for all $n$.  
Then
 \begin{align*}
 \int f \, d \wtd \alpha_0(\gamma)_*\mu_0 &:= 
 \int f_0  \  d  \mu_0 \\
 & = \lim _{n\to \infty} \int f_0 \, d  \mu_n \\
 & = \lim _{n\to \infty} \int (f_n + f_0- f_n)\, d  \mu_n \\
 & = \lim _{n\to \infty} \int f_n \, d  \mu_n   + \lim \int_K  (f_n - f_0) \, d \mu_n\\
 & = \lim _{n\to \infty} \int f_n  \, d  \mu_n   \\  
 & = \lim _{n\to \infty}\int f \, d \wtd \alpha_n(\gamma)_*\mu_n \\
 & = \lim _{n\to \infty} \int f \, d \mu_n \\
& = \int f \, d \mu_0.
\end{align*}

 For \eqref{stupid2},  we then have
$$\int f \, d h_*\mu = \lim _{n\to \infty} \int f \, d h_*\mu_n = \lim _{n\to \infty}\int f \, d \mu_n  = \int f d \mu.$$
For \eqref{stupid3}, fix a finite generating $S\subset \Gamma$.  For every $\gamma_i\in S$ we have 
$$\int f \, d \wtd \alpha_0 (\gamma_i)_*\mu_0 = \lim _{n\to \infty} \int f \, d \wtd \alpha_n (\gamma_i)_*\mu_n= \lim _{n\to \infty} \int f \, d  \mu_n= \int f \, d \mu_0.$$
Since $S$ generates, we obtain $\wtd \alpha_0(\gamma)$-invariance of $\mu_0$ for all $\gamma\in \Gamma$.  
\end{proof}

\subsubsection{Properties of $C^1$ perturbations} 
We now return to the case $X= M$  is a smooth compact manifold and the action is at least $C^1$; although the following lemma certainly holds under weaker hypotheses, our only application is the special case that $X=M=G/Q$ and the action $\alpha_0$ is the standard projective action.  We thus formulate the lemma and its proof under strong hypotheses that hold in this restricted setting.



To motivate the lemma, we recall that given a bounded continuous linear cocycle over a homeomorphism $f$ preserving a Borel probability measure,  the top (resp.\ bottom) average Lyapunov exponent is upper (resp.\ lower) semicontinuous as one perturbs in both the cocycle and the invariant measure.  We establish an analogous result for the fiberwise derivative cocycle on $M^\alpha$ for perturbations $\alpha$ of the standard projective action $\alpha_0$.  Note that in this  setup,  the fiber bundle structure (and corresponding family of norms) is varying rather than the cocycle.  

\begin{lemma}\label{negexp}
Let $M$ be a compact manifold and let $\alpha_0\colon \Gamma\to \Diff^1(M)$ be an action such that the induced $G$-action on $M^\alpha$ admits a unique
$P$-invariant Borel probability measure $\mu_0$  with the following properties: 
\begin{enumerate}
\item there exists a continuous, $P$-equivariant section $\sigma\colon G/\Gamma\to M^\alpha$ such that $\mu_0$ is the image of the Haar measure on $\sigma$; 
\item $\lambda^F_{\top, \mu_0, s_0, \alpha_0}<0$ for some $s_0\in A.$
\end{enumerate}

Then, given any $0<\gamma<1$, for all $\alpha\colon \Gamma\to \diff^1(M)$ sufficiently $C^1$-close to $\alpha_0$ and every $P$-invariant Borel probability measure $\mu$ on $M^\alpha$ we have 
 \begin{align}
 \lambda^F_{\top, \mu, s_0, \alpha}&<\lambda^F_{\top, \mu_0, s_0, \alpha_0} + \gamma\\
  \lambda^F_{\btm, \mu, s_0, \alpha}&>\lambda^F_{\btm, \mu_0, s_0, \alpha_0} - \gamma.
 \end{align}
\end{lemma}

\begin{proof}
We establish some notation. 
Equip a background Riemannian metric on $M$ with induced norm $\|\cdot \|_M$.
To compare norm growth of fiberwise derivatives induced  by different  actions, fix a Siegel fundamental domain $D\subset G$.  Recall $m_G$ is normalized so that $m_G(D)=1$.



We equip $G\times M$ with the $\alpha_0(\Gamma)$-invariant (resp.\ $\alpha(\Gamma)$-invariant) $C^0$ Riemannian metric as outlined in \cref{sec:norms}.  Let $\|\cdot \|^F_{\alpha_0}$ (resp.\ $\|\cdot \|^F_{\alpha}$) denote the associated $\alpha_0(\Gamma)$-invariant, (resp.\ $\alpha(\Gamma)$-invariant)  norm induced on the fibers of  $G\times TM$.  Let $d_{\alpha_0}$  denote the distance function on $G\times M$ induced by the $\alpha_0(\Gamma)$-invariant Riemannian metric.  
Recall that restricted to $D\times M$, the norms $\|\cdot\|^F_{\alpha_0}$, $\|\cdot\|^F_{\alpha}$, and $\|\cdot\|^F_{M}$ are uniformly comparable by a constant $B_0$.

Let  $F^{\alpha_0}$ and $F^\alpha$ denote, respectively, the fiberwise tangent bundles of $M^{\alpha_0}$ and $M^{\alpha}$. 
Recall that $\|\cdot \|^F_{\alpha_0}$ (resp.\ $\|\cdot \|^F_{\alpha}$) descends to a norm on the fiber  $F^{\alpha_0}$ (resp.\ $F^{\alpha}$) which we denote by $\|\cdot \|^F_{M^{\alpha_0}}$ (resp.\ $\|\cdot \|^F_{M^\alpha}$) 




Write $\calA^{\alpha_0}\colon G\times F^{\alpha_0}\to F^{\alpha_0}$ (resp.\ $\calA^\alpha\colon G\times F^{\alpha}\to F^{\alpha} $) for the fiberwise derivative cocycle.  
Both $\calA^{\alpha_0}$ and $\calA^{\alpha}$ lift to the same linear cocycle 
 $\wtd\calA\colon G\times (G\times TM) \to G\times TM $,  
$$ \wtd \calA(g',(g,(x,v)) = (g'g, (x,v)).$$
 Given $g', g\in G$ and $x\in M^{\alpha_0}$, we write $$\wtd \calA(g',(g,x))\colon \{g\}\times T_xM\to \{g'g\}\times T_xM$$ for the linear map between fibers.

Let $\mu_0$ and $\mu$ be $P$-invariant Borel probability  measures, respectively, on $M^{\alpha_0}$ and $M^\alpha$.  
Let $\wtd \mu_0$ and $\wtd \mu$ be  the canonical lifts of $\mu_0 $ and $\mu$ to $G\times M$.  
Recall that we assume that $\wtd  \mu_0$ is the image of the Haar measure under a continuous, $\wtd \alpha_0(\Gamma)$-equivariant, $P$-invariant section $ \wtd \sigma\colon G \to G\times M$. 
Given $\epsilon>0$ and $g\in D$, let
\begin{align*}
\calU_\epsilon (g) &:= \{(g,y)  :  d_0\left((g,y), \wtd \sigma(g)\right) < \epsilon\}.
\end{align*}

Given $g\in G$, let $\gamma(g)\in \Gamma$ be such that $$g\in D\cdot \gamma(g).$$
Given $\alpha \colon \Gamma\to \Diff ^1(M)$ and $\gamma\in \Gamma$, write $$
\Xi^\alpha (\gamma) :=\max\left\{\|D(\alpha(\gamma)\circ  \alpha_0(\gamma\inv))\|_M,\|D(\alpha_0(\gamma)\circ  \alpha(\gamma\inv))\|_M\right\}
.$$

Fix a symmetric finite generating set $S$ for $\Gamma$.  Given $\gamma\in \Gamma$, let $|\gamma|= |\gamma|_S$ denote the word length of $\gamma$ with respect to $S$.  
By \cref{LMR}, 
 compactness of $M$, and uniform comparability of the fiberwise metrics over $D$, there are $A_0$, $A_1$, $B_0,$ and   $C_1>1$ such that for all $n\ge 0$,
\begin{enumerate}
\item\label{eatpoop1}
  if $g\in D$, then   $$|\gamma(s_0^{n} g) | \le A_0 (d(g,\Id)  +n) + A_1;$$
\item \label{eatpoop7} 
 for $g\in D$ and $x\in M$,  $$\|\wtd \calA(s_0^{n}, (g,x)) ^{\pm1}\|^F_{\alpha_0} \le C_1^{A_0 (d(g,\Id) +n) + A_1}.$$
\end{enumerate}
Also, there is (potentially very large) $c>1$ such  that for all $\alpha$ in a $C^1$-bounded neighborhood of $\alpha_0$,  
\begin{enumerate}[resume]

\item \label{eatpoop3} $\Xi ^\alpha(\gamma) \le e^{c |\gamma|}$  for all $\gamma\in \Gamma$. 

\end{enumerate}


Fix $\eta>0$.  Since $g \mapsto d(g, \Id)$ is $L^1$, 
by \eqref{eatpoop1} there is $\delta_0>0$ such that for any $B\subset D$ with $m_G(B)>1-2\delta_0$ and any $n\in \N$, 
\begin{enumerate}[resume]
\item \label{eatpoop4} $\int _{D\sm B} \log^+ \|\wtd \calA(s_0^{n}, (g,x))^{\pm1} \|^F_{\alpha_0}
\le  (n+1)\eta $, and 
\item \label{eatpoop5} for all $\alpha$ sufficiently $C^1$ close to $\alpha_0$,   $$\int _{D\sm B} \log^+ \Xi^\alpha (\gamma(s_0 ^{n} g)) \, d m_G(g) 
\le  (n+1)\eta $$ 
\end{enumerate}
As mentioned above, by construction of the norms (c.f.\ \cite[\S 2.1]{MR4502594}), there is $B_0>1$ such  that for all $g\in D$  and $x\in M$,
\begin{enumerate}[resume]
	\item 
$\frac 1 {B_0} \|\cdot \|^F_{\alpha, g, x}\le \|\cdot \| _{M, x}\le B_0\|\cdot \|^F _{\alpha, g, x} $ and \item
$\frac 1 {B_0} \|\cdot \|^F_{\alpha_0, g, x}\le \|\cdot \| _{M, x}\le B_0\|\cdot \|^F _{\alpha _0, g, x} .$

\end{enumerate}

Write  \begin{align*}\lambda_\top=\lambda^F_{\top, \mu_0, s_0, \alpha_0}, &&
\lambda_\btm=\lambda^F_{\btm, \mu_0, s_0, \alpha_0}.
\end{align*} 
Fix  $0<\delta<\delta_0<0$.  Fix $0<\gamma<1$.  Recall we assume $\lambda_\top<0$.  Taking $\gamma>0$ smaller if needed, 
assume $\lambda_\top+\gamma/2<0.$  There are $\epsilon>0$,  $N\ge 1$, $C_0\ge 1$, and a $\Gamma$-invariant subset $ \Upsilon\subset G$ with the following properties: 
\begin{enumerate}[resume]
        \item \label{eatpoop17}  $m_G(\Upsilon \cap D)\ge 1-\delta;$ 
        \end{enumerate}
        and for $g\in \Upsilon$ and $(g,y)\in \calU_\epsilon(g)$ and all $n\ge N$,
\begin{enumerate}[resume]
        \item \label{eatpoop18} 
        		$d_0\left( \wtd \sigma(s_0^{n} g), (s_0^{n} g, y)\right)\le e^{n \lambda_\top/2} \epsilon$
        \item \label{eatpoop185}
		$\|\wtd \calA( s_0^{n},(g,y))\|^F_{\alpha_0}\le e^{n (\lambda_\top+ \gamma/2)};$
       \item \label{eatpoop187}$\|\wtd \calA( s_0^{n},(g,y))^{-1}\|^F_{\alpha_0}\le e^{-n (\lambda_\btm+ \gamma/2)};$
        \item \label{eatpoop188} $\sup_{g\in \Upsilon \cap D} \{ d_G(g,\Id )\} \le C_0$. 
\end{enumerate}

Let $$\Upsilon _\epsilon := \bigcup _{g\in \Upsilon} \calU_\epsilon (g).$$
Since the $P$-invariant measure $\mu_0$ on $M^{\alpha_0}$ is assumed to be  unique, if $\alpha$ is sufficiently $C^0$-close to $\alpha_0$  then 
\begin{enumerate}[resume]
\item \label{eatpoop101}for every $P$-invariant measure $\mu$  on $M^\alpha$, $$\wtd \mu(\Upsilon_\epsilon\cap (D\times M))\ge 1-2 \delta.$$
\end{enumerate}

Fix $n_0\ge N$ and let $L= \sup \{|\gamma (s^{n_0} g)|: g\in \Upsilon\}$. Taking $\alpha$ sufficiently $C^1$-close to $\alpha_0$, we may assume 
\begin{enumerate}[resume]
\item \label{eatpoop102} $\Xi^\alpha (\gamma) \le e^{\eta n_0}$  for all $\gamma\in \Gamma$ with $|\gamma|\le  L$.
\end{enumerate}

Fix $\alpha$ so that \eqref{eatpoop1}--\eqref{eatpoop102} above hold.  
Using uniformly comparability of  fiberwise norms   on $D$, for all $(g,y) \in D\times M$ we have 
\begin{align*}
	\| \calA^\alpha(s_0^{n_0}, [g,y]_\alpha)\|^F_{M^\alpha} 
&=	\|\wtd \calA(s_0^{n_0}, (g,y))\|^F_{\alpha}\\
	&\le B^2 _0 \cdot   \|{D_{y}} \alpha(\gamma(s_0^{n_0} g))\|_M
	\\& \le  B^2 _0  \cdot \Xi ^\alpha(\gamma(s_0^{n_0} g)) \cdot \|{D_{y}} \alpha_0(\gamma(s_0^{n_0} g))\|_M
%
	%
	\\& 
\le 
	B _0 ^4 \cdot    \Xi ^\alpha(\gamma(s_0^{n_0} g)) \cdot  	\|\wtd \calA (s_0^{n_0}, (g,y))\|^F_{{\alpha_0}}
\end{align*}
and so 
\begin{align}
\int _{M^\alpha} \log^+  & \|\calA^\alpha(s_0^{n_0}, \cdot )\|^F_{M^\alpha}  \, d \mu(\cdot ) \notag \\
&=\int _{D\times M} \log^+   \|\wtd\calA(s_0 ^{n_0}, (g,y))\|^F_{\alpha}  \, d \wtd \mu(g,y) \notag \\
&\le 
4\log^+ B_0 \notag \\
%
%
& \quad + \int _{D } \log^+  \Xi ^\alpha(\gamma(s_0 ^{n_0} g)) 
\, d m_G(g)
\label{eq:returnterm}
 \\& \quad + \int _{D\times M } \log^+ \|\wtd \calA^{\alpha_0}(s_0 ^{n_0}, (g,y))\|^F_{{\alpha_0}}  \, d \wtd \mu(g,y)  \label{eq:mainterm}.
\end{align}

 
By \eqref{eatpoop5} and \eqref{eatpoop102} above, the integral \eqref{eq:returnterm} is bounded by

\begin{align*}
	\int _{D } &\log^+ \Xi^\alpha(s_0^{n_0}g) \, d m_G(g) \\&=
	\int _{D \cap  \Upsilon} \log^+ \Xi^\alpha(s_0^{n_0}g) \, d m_G(g) + \int _{D \sm \Upsilon} \log^+ \Xi^\alpha(s_0^{n_0}g)\, d m_G(g) \\
	&\le  \eta n_0 + \eta(n_0+1) \\ &= \eta(2n_0 +1).  
\end{align*}

%

Let $\calG= (D\times M) \cap \Upsilon_\epsilon$.
We decompose \eqref{eq:mainterm} as 
\begin{align}
& \int _{\calG} \log^+ 
 \|\wtd \calA(s_0^{n_0}, (g,y))\|^F_{{\alpha_0}} \, d \wtd \mu(g,y)
 \label{eqmain} 
\\& \quad  +
\int _{(D \times M) \sm \calG} \log^+ \|\wtd \calA(s_0, (g,y))\|^F_{{\alpha_0}} \, d \wtd 
 \mu(g,y).\label{eqbad}
\end{align}
From \eqref{eatpoop17} and \eqref{eatpoop185} above, we bound \eqref{eqmain} from above by $$n_0 (\lambda_\top +\frac{\gamma}{2} )\mu(\calG) \le n_0(\lambda_\top +\frac{\gamma}{2} )(1-2\delta).$$
Using  \eqref{eatpoop4}, we bound  \eqref{eqbad} from above by $\eta(n_0 +1).$
%

Taking first $\eta>0$ sufficiently small, then $0<\delta<\delta_0$ sufficiently small, then $n_0\ge N$ sufficiently large, and then $\alpha$ sufficiently $C^1$-close to $\alpha_0$, we can ensure that
\begin{align*}
\inf_{n} 	\frac 1 n \int _{M^\alpha} \log^+  & \|\calA^\alpha(s_0^n, \cdot )\|^F_{M^\alpha}  \, d \mu(\cdot ) \\
&\le \inf_{\ell} 	\frac 1 {\ell n_0} \int _{M^\alpha} \log^+   \|\calA^\alpha(s_0^{\ell n_0}, \cdot )\|^F_{M^\alpha}  \, d \mu(\cdot ) \\
&\le \inf _{\ell}	\frac 1 {\ell n_0} \sum_{j=0}^{\ell-1}\left(\int _{M^\alpha} \log^+   \|\calA^\alpha(s_0^{ n_0}, \cdot )\|^F_{M^\alpha}  \, d \mu(\cdot )\right) \\
&= 	\frac 1 {n_0} \int _{M^\alpha} \log^+   \|\calA^\alpha(s_0^{ n_0}, \cdot )\|^F_{M^\alpha}  \, d \mu(\cdot ) \\
	&\le \frac 1 {n_0} \left[4\log^+ B_0 + \eta(2n_0 +1) +
	n_0(\lambda_\top +\frac{\gamma}{2} ) (1-\delta)
	+  \eta(n_0 +1)  \right]\\
	&= 
(\lambda_\top +\frac{\gamma}{2} ) (1-\delta) + 3\eta +
	\frac 1 {n_0} \left[4\log^+ B_0 + 2\eta	 \right]
	\\&\le  \lambda_\top  +\gamma . 
\end{align*}
We similarly have 
\begin{align*}
\inf_{n} 	\frac 1 n \int _{M^\alpha} \log^+  & \|(\calA^\alpha(s_0^n, \cdot ))^{-1}\|^F_{M^\alpha}  \, d \mu(\cdot ) \\
&\le (-\lambda_\btm +\frac{\gamma}{2} ) (1-\delta) + 3\eta +
	\frac 1 {n_0} \left[4\log^+ B_0 + 2\eta	 \right]
	\\&\le - \lambda_\btm  +\gamma. 
\end{align*}
It follows that 
$$\lambda^F_{\btm, \mu_0, s_0, \alpha_0} - \gamma<   \lambda^F_{\btm, \mu, s_0, \alpha} \le\lambda^F_{\top, \mu, s_0, \alpha}<\lambda^F_{\top, \mu_0, s_0, \alpha_0} + \gamma$$
and the result follows.
\end{proof}

Applied to our primary  example of an action $\alpha_0$ satisfying the hypotheses of \cref{negexp}, we obtain the following:
\begin{corollary}\label{cor:critregusc}
Let $Q$ be a parabolic subgroup of a semisimple Lie group $G$.  
 Let $M= G/Q$ and let $\alpha_0\colon \Gamma\to \Diff^\omega(M)$ denote the standard projective action of a lattice subgroup $\Gamma\subset G$.  
 
Fix $\gamma>0$.  Then for every $s_0$ in the negative cone $C(\Sigma^\perp_\lieq)$ of $Q$ as defined in \eqref{eq:negcone},  every  $\alpha\colon \Gamma\to \Diff^1(M)$ sufficiently $C^1$-close to $\alpha_0$, and every  $P$-invariant Borel probability measure on $M^\alpha$,
$$ \lambda^F_{\top, \mu, s_0, \alpha}-   \lambda^F_{\btm, \mu, s_0, \alpha} \le 
 \lambda^F_{\top, \mu, s_0, \alpha_0}-   \lambda^F_{\btm, \mu, s_0, \alpha_0} +2\gamma.$$

In particular, for any $\epsilon>0$,
for every  $\alpha\colon \Gamma\to \Diff^1(M)$ sufficiently $C^1$-close to $\alpha_0$ and every  $P$-invariant Borel probability measure on $M^\alpha$,
\begin{equation}
 c^F_0(\mu) < r_0(\lieq) + \epsilon.
\end{equation} 
where  $ r_0(\lieq)$ is as in \eqref{eq:parareg}.  
\end{corollary}

\subsection{\texorpdfstring{Sufficient criteria for $\beta$-tameness: continuous factors}{Sufficient criteria for beta-tameness: continuous factors}}
\label{tame1}
Under the hypotheses of \cref{thm:dani,thm:smoothfactor}, there is a  unique $P$-measure $\mu$ on $M^\alpha$ (the push-forward of the unique $P$-measure for the suspension of the standard action). This measure is $Q$-invariant and we show this $P$-measure is $\beta$-tame for every resonant root $\beta\in \Sigma_\lieq^\perp$.  
\begin{proposition}\label{factortame}
Let $G$ be a higher-rank semisimple Lie group, let $\Gamma$ be an irreducible lattice subgroup, and let $Q\subset G$ be a parabolic subgroup. 
Let $M$ be a compact connected manifold homeomorphic to $G/Q$, let $\alpha\colon \Gamma\to \Diff^{1+\text{\Holder}}(M)$, and let $p\colon G/Q\to M$ be a $\Gamma$-equivariant homeomorphism.  

Let $\mu$ denote the unique $P$-invariant Borel probability measure on $M^\alpha$.  Then for every $\beta\in \Sigma_Q^\perp$, 
\begin{enumerate}
    \item\label{factor1} $\beta$ is a resonant root, and 
    \item\label{factor2} $\mu$ is $\beta$-tame.
\end{enumerate}
\end{proposition}
\begin{proof}
We equip $G\times (G/Q)$ with the  left-$G$- and right-$\Gamma$-actions in \eqref{flat}.
Let $\wtd p\colon G\times (G/Q)\to G\times M$ be $$\wtd p(g, g'Q) = \left(g,  p(g\inv g'Q)\right).$$ Then $\wtd p$ is left-$G$-equivariant and right-$\Gamma$-equivariant;
indeed, given $\gamma\in \Gamma$ and $h\in G$,
\begin{align*}\wtd p(h\cdot (g, g'Q))& = 
\wtd p(h g, hg'Q) \\&= 
( hg,   p(g\inv g'Q)) \\&= 
h\cdot \wtd p (g, g'Q)
\end{align*}
and 
\begin{align*}
\wtd p( (g, g'Q)\cdot \gamma ) &= 
\wtd p(g\cdot \gamma , g'Q) \\&= 
( g\cdot \gamma ,   p( \gamma \inv g\inv g'Q)) \\&= 
( g\cdot \gamma ,  \alpha( \gamma\inv )  p(g\inv g'Q))\\& = 
\wtd p(g, g'Q) \cdot \gamma.
\end{align*}
In particular, $\wtd p$ descends to a $G$-equivariant homeomorphism $$\bar p \colon (G/\Gamma) \times (G/Q) \to M^\alpha.$$

Let $\mu_0$ be the unique $P$-invariant Borel probability measures on $ (G/\Gamma) \times (G/Q) $ given by \cref{uniquePmeasure}.
We have $\mu= \bar p _* \mu_0$.  In particular, $\bar p \colon ((G/\Gamma) \times (G/Q),\mu_0) \to (M^\alpha, \mu)$
is a measurable isomorphism.  
As $\bar p$ is a $G$-equivariant homeomorphism, $\mu$ is $Q$-invariant and is  not  $Q'$-invariant for any $Q'\supsetneq Q$.  
Moreover, as $\mu_0$ is supported  under the graph of  a continuous function $G/\Gamma\to G/Q$, the support of $\mu$ is the image of a continuous section $G/\Gamma\to M^\alpha$.  Let $S\subset M^\alpha$ denote the support of $\mu$.  

Conclusion \eqref{factor1} follows from  \cref{NRimpInv}, since $\mu$ is $Q$-invariant but is not invariant under any $Q'\supsetneq Q$.  Thus, for every $\beta\in \Sigma_\lieq^\perp$ there is a fiberwise Lyapunov exponent positively proportional to $\beta$. 

For  conclusion \eqref{factor2}, we construct a family of open sets $G(x)$ as in \cref{defn:tamedyn}.  
First, we note that for $\mu$-almost every $x$, 
\begin{enumerate}
\item $\bar p\left(W^{[\beta],F}(\bar p\inv (x))\right)$ is an injectively immersed topological manifold in the fiber of $M^\alpha$ containing $x$, and 
\item  $W^{[\beta],F}(x) \subset \bar p\left(W^{[\beta],F}(\bar p\inv (x))\right)$.      
\end{enumerate}
Given $r>0$, let $W^{[\beta],F}_r(x)$ denote the path-connected component containing $x$ in  $W^{[\beta],F}(x) \cap B^F(x,r)$, where $B^F(x,r)$ denotes the ball of radius $r$ centered at $x$ in the fiber of $M^\alpha$ containing $x$.  
Take $r_0$ to be   smaller than the Riemannian injectivity radius of each fiber of $M^\alpha$.
Given $0<r<r_0$, let $E_r$ denote the subset of $x\in S\subset M^\alpha$ for which $(\exp_x^F)\inv \left(W^{[\beta],F}_r(x)\right)$ is the graph of a $C^1$-function $g_x\colon E^{[\beta],F}_x\to (E^{[\beta],F}_x)^\perp.$  Fix $0<r<r_0$ for which $\mu(E_r)>.9.$

Given $s>0$, let $\lieu^{[\beta]}(s) = \{X\in \lieu^\beta: \|X\|<s\}$ and $U^{[\beta]}(s) = \exp_{\lieg} \lieu^{[\beta]}(s).$

Let $y\in M$ be $y= p(g\inv Q)$.  
Then $(g,y)\in G\times M$ is the form $(g,y) = \wtd p(g, \1 Q)$ and thus $x=[(g,y)]$ is in the support $S$ of the measure $\mu$.  
 For such $y$, let $\wtd P((g,y),s)$ denote the set 
$$\wtd P((g,y),s)= \{\wtd p(ug, \1Q): u\in V_Q(s)\}$$
and let  $$P((g,y),s)= \{g\}\times \pi^2(\wtd P((g,y),s))$$
where $\pi^2\colon G\times M\to M$ is the projection onto the second coordinate.  
We note (by definition) that $P((g,y),s):= \{g\}\times  \{p(g\inv u\inv Q): u\in V_Q(s)\}$ is an open neighborhood of $y$ in $\{g\}\times M$.  

We check that $P((g,y),s)= P((g,y)\cdot \gamma,s)$.  Thus the set $P(x,s)$  is well defined for every $x\in S$.  
For $x\in S$,  $P(x,s)$ is an open and pre-compact (with respect to the immersed topology) 
subset of  
the immersed topological manifold $\bar p\left(W^{[\beta],F}(\bar p\inv (x))\right)$.

Write $A'_\beta$ for the kernel of $\beta$ in $A$.  Then for all $a\in A'_\beta$, we have  $a\cdot P(x,s)= P(a\cdot x,s)$.

Having fixed $0<r<r_0$ above, let $E_{s,r}\subset S\subset  M^\alpha$ be the set of $x\in E_r$  for which 
$$\overline{P(x,s)}\subset B^F(x,r).$$
Take $0<s$ sufficiently small so that $S_0= E_{s,r}\subset E_r\subset S$ has positive measure.  
Given $x\in S_0$, $W^{[\beta],F}_r(x)$ is a properly embedded submanifold in $B^F(x,r)$ and thus 
$W^{[\beta],F}_r(x)\cap P(x, s)$ is an open, precompact subset of $W^{[\beta],F}_r(x)$.  
For $x\in S_0$, let $Y(x)$ denote the path connected component of $W^{[\beta],F}_r(x)\cap P(x, s)$ containing $x$.  
Fix $a_0\in A'_\beta$.  Let $n_0(x)=\inf\{n\ge 0: a_0^n\cdot x\in S_0\}$ denote the first hitting time.  
For $x$ with $n_0(x)\ge 1$, let 
$$Y(x) = a_0^{-n_0(x)}\cdot (Y(a_0^{n_0(x)}\cdot x)).$$

We claim for any $a\in A'_\beta $ and $\mu$-a.e.\ $x$ that $a\cdot Y(x) = Y(a\cdot x)$.  
Indeed, by $A'_\beta$-equivariance of the maps $x\mapsto P(x,s)$ and  $x\mapsto W^{[\beta],F}(x)$ and uniqueness of path components, if $x\in S_0$ and $a\cdot x\in S_0$ then $a\cdot Y(x)= Y(a\cdot x)$.  
For general $x$,
\begin{align}(a_0^{n_0(a\cdot x)}a) \cdot Y( x)&=
(a_0^{n_0(a\cdot x)}aa_0^{-n_0(x)}) \cdot Y(a_0^{n_0(x)} \cdot x) \notag\\
&=   Y(a_0^{n_0(a\cdot x)}a \cdot x)\label{eq:poopface}\\
&=  a_0^{n_0(a\cdot x)}\cdot  Y(ax)\label{eq:poopface2} \end{align}
where \eqref{eq:poopface} follows using that both $a_0^{n_0(x)} \cdot x\in S_0$ and $(a_0^{n_0(a\cdot x)} a)\cdot x)\in S_0$ and 
\eqref{eq:poopface2} follows either because $a\cdot x\in S_0$ or the definition of $Y(a\cdot x)$ if $a\cdot x\notin S_0$.
It follows that  $a \cdot Y( x)=  Y(a\cdot  x).$

Taking  $G(x) = U^{[\beta]}(s)\cdot Y(x)$.  As $G(x)$ is precompact in $W^{[\beta]}(x)$ for a.e.\ $x$ and $A'_\beta$-equivariant,  it follows that $G(x)$ satisfies the conditions of \cref{defn:tamedyn}.
\end{proof}

\subsection{\texorpdfstring{Sufficient criteria for $\beta$-tameness: all fiberwise exponents negative}
{Sufficient criteria for beta-tameness: all fiberwise exponents negative}}
\label{tame2}
For many of the results in this paper, we obtain $\beta$-tameness from the following criteria.  
\begin{proposition}\label{prop:equifinite}
Let $M$ be a compact  manifold and let $\alpha\colon \Gamma\to \Diff^{1+\text{\Holder}}(M)$ be an action.  
Let $\mu$ be an ergodic, $P$-invariant  Borel probability measure on $M^\alpha$.  

Suppose for some $a_0\in A$,  every fiberwise Lyapunov exponent of $\mu$  is negative for the action of  $a_0$ on $(M^\alpha, \mu)$.  
Then for every root $\beta\in \Sigma$ that is resonant with respect to $\mu$, the measure $\mu$ is $\beta$-tame.
\end{proposition}

\begin{proof}
Let $\lieu=\lieu^{[\beta]}$ and $U=U^{[\beta]}$.  
Given $r>0$, let $\lieu(r)= \{v\in \lieu: \|v\|<r\}$ and let $U(r) = \exp_{\lieu} (\lieu(r))$.

Let $\what \calL^F(\mu)$ denote the collection of fiberwise coarse Lyapunov exponents.  
%
%
%


Recall the terminology and results in \cref{integrable}.  
In particular, for every strongly integrable collection $\calI$, there is a fiberwise manifold $W^{F, \calI}(x)$ defined for $\mu$-a.e.\ $x$.
When $\calI= [\beta]$ is a single  coarse Lyapunov exponent given by a resonant root $\beta$ we have that $\calI$ is strongly integrable and  $W^{F, \calI}(x)= W^{F, [\beta]}(x)$ is the associated fiberwise coarse Lyapunov exponent.  

Consider a strongly integrable collection $\calI$ of fiberwise Lyapunov exponents containing $[\beta]$.  Then for any $a\in C(\calI)$, we have $\beta(a)<0$.  It follows that if $u\in U^{[\beta]}$ then for $\mu$-a.e.\ $x$,  $u\cdot W^{F, \calI}(x)= W^{F, \calI}(u\cdot x)$.
Given a strongly integrable collection $\calI$ with $[\beta]\in \calI$, let 
$\Fol^\calI$ denote the foliation whose leaves are $U\cdot W^{F, \calI}(x)$ and let 
$\Fol^\calI(x,r)= U(r)\cdot W^{F, \calI}(x)$.  

Let $\calI_0= \{\lambda: \lambda(a_0)<0\}$.  
Let $\Fol^{[\beta],s}=\Fol^{\calI_0}$ and let  $\mu^{[\beta],s }_x$ denote the leafwise measures along 
leaves of $\Fol^{[\beta],s}$.
Let $\mu^{[\beta],s }_x(r)$ denote the restriction of $\mu^{[\beta],s }_x$ to $\Fol^{[\beta],s}(x,r)$.  

By assumption on the fiberwise  Lyapunov  exponents of $\mu$,  $W^{F, s}_{a_0}(x)$ is open in each fiber $\calF(x)$ for almost every $x$.  Since the fibers $M$ of $M^\alpha$ are compact, it follows for $\mu$-a.e.\ $x$ that  $ \mu^{[\beta],s }_x(r)$ is a finite measure supported on $\Fol^{[\beta],s}(x,r)$.  
The proposition then follows inductively from \cref{lem:inductmeintheass} below.
\end{proof}

To finish the proof of \cref{prop:equifinite}, we have the following inductive argument.  
\begin{lemma}\label{lem:inductmeintheass}
Let $\calI= \{\chi_1,\dots, \chi_\ell\}$ be a strongly integrable collection of fiberwise coarse Lyapunov exponents with the following properties:
\begin{enumerate}[label=(\alph*), ref=(\alph*)]
\item $[\beta]\in \calI$;
\item $ \mu^{[\beta],s }_x(r)$ is supported on $\Fol^\calI(x,r)$ for a.e.\ $x$;
\item $\chi_\ell \neq [\beta]$;
\item there are  $a_1,a_2\in A$ 
such that for every $\lambda\in \chi_\ell$, every $1\le j \le \ell-1$, and  every $\td \lambda\in \chi_j$ 
 \begin{equation}\td \lambda(a_1)<\lambda(a_1)<0\label{eqfast} \end{equation}
 and
 \begin{equation}\td \lambda(a_2)<0<\lambda(a_2). \label{eqint}\end{equation}
\end{enumerate}
Let $\calI'= \calI\sm\{\chi_\ell\}$.  Then 
\begin{enumerate}
\item $\calI'$ is strongly integrable
\item $ \mu^{[\beta],s }_x(r)$ is supported on $\Fol^{\calI'}(x,r)$ for a.e.\ $x$.
\end{enumerate}
\end{lemma}
\begin{proof}
From \eqref{eqint}, the collection  $\calI'$  is strongly integrable.  

From \eqref{eqfast}, $\calI'$ corresponds to the Lyapunov exponents in the fast stable foliation inside $W^{\calI,F}(x)$.  In particular, for almost every $x$ and every $y\in \Fol^{\calI}(x,r)= U(r) \cdot W^{\calI,F}(x)$, there is a smooth, embedded submanifold $V^{\calI',F}(y)\subset W^{\calI,F}(x)$ with the following properties:
\begin{enumerate}
\item $V^{\calI',F}(z)= W^{\calI',F}(z)$ for $\mu$-a.e.\ $z\in  \Fol^{\calI}(x,r)$;

\item\label{222222}  if, $y\in W^{F, \calI}(x)$ then $V^{\calI',F}(y)$
 and $W^{F, \chi_\ell}(x)$ intersect in exactly one point;  
\item if $y,z\in  \Fol^{\calI}(x,r)$ are such that $y = u\cdot z$ for some $u\in U$ then $u\cdot V^{\calI',F}(z) = V^{\calI',F}(y)$
\end{enumerate}
We note that the transversality in condition  \eqref{222222} above holds locally; one may iterate to obtain the global transversality statement.

Given $y  \in  \Fol^{\calI}(x,r)$, write $y= u\cdot z  $ for $u\in U(r)$ and $z\in W^{\calI,F}(x)$.  
Let $\calF'(y  ) = U(r) \cdot V^{\calI',F}(z)$.  
Then 
\begin{enumerate}[resume]
\item  the collection  $\{\calF'(y  ): y  \in  \Fol^{\calI}(x,r)\}$ has the structure of a Lipschitz  foliation by properly embedded open discs with trivial holonomy and forms a measurable partition of $ \Fol^\calI(x,r)$. 
\end{enumerate}

Let $\bar \mu^{[\beta],s }_x(r)$ denote the measure  $\mu^{[\beta],s }_x(r)$ normalized (on $\Fol^{\calI}(x,r)$).
Select any $a_3\in A$  with $\beta(a_3) = 0$ and $\chi_\ell(a_3)<0$.  
We have 
\begin{enumerate}[resume]
\item  $a_3\cdot \Fol^{\calI}(x,r) = \Fol^{\calI}(a_3\cdot x,r)$;
\item  $(a_3)_* \bar \mu^{[\beta],s }_x(r)= \bar \mu^{[\beta],s }_{a_3\cdot x}(r).$
\end{enumerate}
We parameterize the space of $\Fol'$-leaves in each $\Fol^{\calI}(x,r)$ by the transverse embedded manifold $W^{F, \chi_\ell}(x)$.  Since $\bar \mu^{[\beta],s }_x(r)$ is a finite measure, and since the dynamics of $a_3$ contracts the transverse parameters in  $W^{F, \chi_\ell}(x)$, by a Poincar\'e recurrence argument, it follows  for a.e.\ $x$ that $\bar \mu^{[\beta],s }_x(r)$ is supported on $\calF'(x)=\Fol^{\calI'}(x,r).$
\end{proof}

\def\mhaar{{m_{G/\Gamma}}}

\subsection{Entropy properties of $A$-invariant measures} \label{sec:entropy} 
We collect some facts on metric entropy, especially lower bounds on the metric entropy when conditioning on certain invariant partitions given by ergodic decomposition.  

Let $A'\subset A$ be a non-trivial closed subgroup of $A$.  In practice, we will only apply the following when 
$A'=A'_\beta:=\ker \beta$, the kernel of some root $\beta$ in $A$.  
As we assume $G$ has no compact factors and that $\Gamma$ is irreducible, by Moore's ergodicity theorem (c.f.\ \cite[Thm. 2.2.6]{MR776417}), $A'$ acts ergodically on $(G/\Gamma, m_{G/\Gamma})$.

Let $\mu$ be an ergodic, $A$-invariant Borel probability measure on $M^\alpha$ projecting to the Haar measure on $G/\Gamma$.  As the measure $\mu$ need not be $A'$-ergodic, let $\calE$ denote the $A'$-ergodic decomposition of $\mu$; we note that $\calE$ is an $A$-invariant measurable partition.   We also recall that $\calF$ denotes the $A$-invariant partition of $M^\alpha$ into fibers of $\pi\colon M^\alpha\to G/\Gamma$.



The following will be proved in {\cite[Lemma 5.1]{ABZ}}.  
\begin{lemma} \label{lemma:entropy} 
Let $\mu$ be an ergodic, $A$-invariant Borel probability measure on $M^\alpha$ projecting to the Haar measure on $G/\Gamma$.

	Fix $\beta\in \Sigma$ such that $\frac 1 2 \beta \notin \Sigma$ and fix $a\in A$ with $\beta(a)>0$.  
	Then 
\begin{equation}\label{eq:entgaplower}
	 \beta(a) ( \dim \lieg^\beta + 2	 \dim \lieg^{2\beta} ) = 
	h_\mu(a\mid {\calE} \vee  \scrW^{[\beta]})-
	h_\mu(a\mid {\calE}\vee \calF\vee \scrW^{[\beta]}).
\end{equation}
\end{lemma}

\subsubsection{Consequences of entropy considerations}
As above, fix an ergodic, $A$-invariant measure $\mu$ on $M^\alpha$ projecting to the Haar measure on $G/\Gamma$,
a non-trivial closed subgroup $A'\subset A$, and let $\calE$ be the $A'$-ergodic decomposition of $\mu$.  
Given any restricted root $\beta$,  \cref{lemma:entropy} implies the leafwise measures associated to the $A'$-ergodic components, $\mu^{{[\beta]},\calE }_x$,  have sufficiently large support.  
\begin{corollary}\label{bigsup}
Let $\mu$ be as in \cref{lemma:entropy}.  Then  for $\mu$-a.e.\ $x$, the leafwise measure $\mu^{{[\beta]},\calE }_x$  is not supported on the fiber $W^{[\beta],F}(x)$ in $W^{[\beta]}(x)$.
\end{corollary}

When the leafwise measure $\mu^{{[\beta]} }_x$ are homogeneous along a smooth submanifold, we have a much stronger version of \cref{bigsup} which asserts the 
leafwise measures do not change when passing to ergodic components.  
\begin{corollary}\label{lem:entropysavesmyassonceagin}
Let $\mu$ be as in \cref{lemma:entropy}.  Suppose for $\mu$-a.e.\ $x\in M^\alpha$ that the leafwise measure $\mu^{{[\beta]} }_x$ 
 is the image of the Haar measure on $U^{[\beta]}$ under  the graph
 $$u\mapsto u\cdot \phi_x^{[\beta]}(u)$$
  of a $C^s$-embedding $\phi_x^{[\beta]}\colon U^{[\beta]}\to W^{F,[\beta]}(x)$ with $\phi_x^{[\beta]}(\1) = x$.  
Then for $\mu$-a.e.\ $x$, we have coincidence of leafwise measures $\mu^{{[\beta]},\calE }_x=\mu^{{[\beta]} }_x.$ 
\end{corollary}
Indeed, given $a\in A$ with $\beta(a)>0$, we have $h_\mu(a\mid {\calE}\vee \calF\vee \scrW^{[\hat\beta]})=0$ (by hypothesis) and thus $$h_\mu(a\mid {\calE}\vee \scrW^{[\hat\beta]}) =
h_ \mhaar(a\mid U^{[\beta]}).$$
It follows (exactly as in the proof of  \cite[Thm.\ A]{MR819556} or \cite[(2.7)]{MR743818}) that $\mu^{{[\beta]},\calE }_x$ is absolutely continuous along the graph of $\phi_x^{[\beta]}$ for a.e.\ $x$.  By computing the explicit density (c.f.\ \cite[Lem.\ 6.1.3]{MR819556}), we conclude that  $\mu^{{[\beta]},\calE }_x=\mu^{{[\beta]} }_x.$

\subsection{Extra invariance and structure from regularity of coarse leafwise measures}\label{sec:extrainv}
We consider an ergodic, $A$-invariant Borel probability measure $\mu$ on $M^\alpha$ projecting to the Haar measure on $G/\Gamma$.  
We recall that \cref{NRimpInv} implies the measure $\mu$ is  invariant under the root group $U^{[\beta]}$ for all non-resonant roots $\beta$.  For a resonant root $\beta$, 
assuming the leafwise measures $\mu^{[\beta]}_x$ are sufficiently regular, we show that either the measure is a graph or obtain additional invariance of the measure along a subgroup of $U^{[\beta]}$.  
 In almost every application considered in the sequel, $\mu$ will be (an $A$-ergodic component of) a $P$-invariant measure;  we will assume $\beta$ is a negative root so that $\mu$ is automatically $U^{[-\beta]}$-invariant (see \cref{ergcomp}) and we may apply \cref{lem:alge} below.   
%
%

Fix a resonant root $\beta\in \Sigma$.  
 We again write $A'_{\beta}=\ker \beta$ and write $\calE_{\beta} $ for the $A'_{\beta}$-ergodic decomposition of $\mu$.

\begin{proposition}\label{graph or invariant}
Let $\mu$ be an ergodic, $A$-invariant Borel probability measure  on $M^\alpha$ that projects to the Haar measure on $G/\Gamma$.  
Let $\beta\in \Sigma$ be root.    
Suppose the following hold for $\mu$-a.e.\ $x$:
\begin{enumerate}[label=(\alph*), ref=(\alph*)]
	\item  \label{lebsubfold} the leafwise measure $\mu^{[\beta],\calE_{\beta}}_x$ of $\mu^{\calE_{\beta}}_x$ along the leaf $W^{[\beta]}(x)$ is in the Lebesgue class on a connected, embedded, $C^s$ submanifold $N_x\subset W^{[\beta]}(x)$.
\end{enumerate}
Assume in addition that one of the following holds for $\mu$-a.e.\ $x$:
\begin{enumerate}[label=(\alph*), ref=(\alph*),resume]
\item \label{casedimfib} $\dim E^{F, [\beta]}(x) \le \dim \lieg^{ [\beta]}$, or   
\item \label{caseatom} the fiberwise leafwise measure $\mu^{\calE_{\beta},[\beta],F}_x$ is purely atomic.
\end{enumerate}

Then either 
\begin{enumerate}
\item $\mu$ is invariant under a subgroup $H\subset U^{[\beta]}$ of positive dimension, or 
\item for $\mu$-a.e.\ $x$, $N_x$ is the graph of an injective, $C^s$-embedding $\phi_x^{[\beta]}\colon U^{[\beta]}\to W^{F,[\beta]}(x)$ with $\phi_x^{[\beta]}(\1) = x$.   Moreover, $\mu^{[\beta],\calE_{\beta}}_x$  is the image of the Haar measure on $U^{[\beta]}$ under the graph of $\phi_x^{[\beta]}$, $$u\mapsto u\cdot  \phi_x^{[\beta]}(u).$$
\end{enumerate}	
\end{proposition}



We have the following lemma on the structure of $\mathfrak{sl}(2)$-triples; see \cite[Lem.\ 7.73]{MR1920389}.  
\begin{lemma}\label{lem:alge}
Let $\beta\in \Sigma$ be a root.  Let $\lieh$ be a Lie subalgebra such that  $\lieg^{[-\beta]}\subset \lieh $  and  $Y\in \lieh$ for some non-zero $Y\in \lieg^{[\beta]}$.  Then  $\lieg^{[\beta]}\subset \lieh .$
\end{lemma}

\cref{lem:alge} then implies the stronger version of \cref{graph or invariant}.  
\begin{corollary}\label{cor:graphorinv}
Suppose in  \cref{graph or invariant} that $\mu$ is $U^{[-\beta]}$-invariant.  
Then either 
\begin{enumerate}
\item $\mu$ is $U^{[\beta]}$-invariant, or 
\item for $\mu$-a.e.\ $x$, the leafwise measure $\mu^{[\beta],\calE_{\beta}}_x$ is the image of the Haar measure on $U^{[\beta]}$ under the graph of $\phi_x^{[\beta]}$, $u\mapsto u\cdot  \phi_x^{[\beta]}(u).$
\end{enumerate}
In particular, if $\dim E^{F, [\beta]}(x) < \dim \lieg^{ [\beta]}$ then $\mu$ is $U^{[\beta]}$-invariant.  
\end{corollary}

The proof of \cref{graph or invariant} occupies the remainder of this section.  
\subsubsection{Consequences of entropy considerations}

Assume that $\mu$  and $\beta$ are as in \cref{graph or invariant}.  Hypothesis \ref{lebsubfold} of \cref{graph or invariant} allows a much stronger version of \cref{bigsup}.  
We write  $\lieg^{[\beta]}_x$ for the tangent space to the $U^{[\beta]}$-orbit at $x$ so that   $E^{[\beta]}_x = E^{[\beta],F}_x \oplus \lieg^{[\beta]}_x$.  
\begin{corollary}\label{cor:dimofproj}
For  $N_x\subset  W^{[\beta]}(x)$ as in \cref{graph or invariant}, the restriction of the projection 
$$E^{[\beta],F}_x \oplus \lieg^{[\beta]}_x \to \lieg^{[\beta]}_x$$
to $T_xN_x$  is onto.  
\end{corollary}
\begin{proof}
Since $T_xN_x$ is $A$-equivariant, we have $$T_xN_x= \bigoplus_{\lambda\in [\beta]}E^{\lambda}_x\cap T_x N_x.$$
Moreover, given $a\in A$ with $\beta(a)>0$ we have 
\begin{enumerate}
\item $h_\mu(a\mid {\calE_\beta}) = \sum_{\lambda\in [\beta]} \lambda(z) \dim (E^{\lambda}_x\cap T_x N_x)$
\item $h_\mu(a\mid {\calE_\beta}\vee \calF) 
= \sum_{\lambda\in [\beta]} \lambda(a) \dim (E^{\lambda,F}_x\cap T_x N_x)$
\end{enumerate}
Combined with  \cref{eq:entgaplower}, it follows that \(\dim T_xN_x - \dim (T_x N_x\cap T_x F(x)) = \dim \lieg^{[\beta]}_x. \qedhere\)
\end{proof}



\subsubsection{Key lemma and proof of \cref{graph or invariant}}
We continue with the proof of \cref{graph or invariant} assuming our key lemma, \cref{lem:injecttrans} below.
  
  Recall we write  $\lieg^{[\beta]}_x$ for the tangent space to the $U^{[\beta]}$-orbit at $x$.  
Given $x\in M^\alpha$ for which $N_x$ is well defined and $y\in N_x$, let 
$$d_1(y) = \dim (T_yN_y\cap \lieg^{[\beta]}_y) \quad \quad d_2(y)= \dim (T_yN_y\cap  W^{[\beta], F}(y)).$$
\begin{lemma}\label{lem:injecttrans}
Let $\mu$, $\beta$, and $N_x$ be as in \cref{graph or invariant}. If $\mu$ is not invariant under any 1-parameter subgroup of $U^{[\beta]}$ then for $\mu$-a.e.\ $x$ and every $y\in N_x$, 
$d_1(y) =d_2(y) =0$.  
\end{lemma}
\begin{proof}[Proof of  \cref{graph or invariant}]
 \cref{lem:injecttrans} implies that locally, $N_x$ is the graph of an injective, $C^s$-immersion $\phi_x^{[\beta]}\colon U^{[\beta]}\to W^{F,[\beta]}(x)$.
By $A$-equivariance of $x\mapsto N_x$,  selecting $a\in A$ with $\beta(a)>0$, $N_x$ is the graph of a $C^s$-embedding $\phi_x^{[\beta]}\colon U^{[\beta]}\to W^{F,[\beta]}(x)$ for $\mu$-a.e.\ $x$.

It remains to argue the measure is the image of Haar on the graph $u\mapsto u\cdot \phi_x^{[\beta]}(u)$ of $\phi_x^{[\beta]}$.  Let $\psi_x^{[\beta]}\colon N_x\to U^{[\beta]}$ be the inverse map such that for $y\in N_x$, 
$$y = \psi_x^{[\beta]}(y) \cdot \phi_x^{[\beta]}\left(\psi_x^{[\beta]}(y)\right).$$
Given $a\in A$, and $y\in N_x$, we have 
$$\psi_{a\cdot x}^{[\beta]}(a\cdot y)= \Ad (a) \psi_x^{[\beta]}(y).$$
In particular, $\psi^\beta_x$ smoothly parametrizes the dynamics on $N_x$ in terms of the dynamics on $U^{[\beta]}$.
Let $\nu_x:= (\psi^{\beta}_x)_*\mu^{[\beta],\calE_{\beta}}_x$.

Consider any  $a_1\in A\sm A'_\beta$ with $\beta(a_1)>0$.  
	 As we assume $\mu^{[\beta],\calE_{\beta}}_x$ is in the Lebesgue class on $N_x$, by a generalization of the Pesin entropy formula  (e.g.\ \cite{MR693976}),  we have $h_\mu(a_1\mid \scrW^{[\beta]}\vee {\calE_\beta}) = \beta(a_1) \dim \lieg^{\beta}+  2\beta(a_1) \dim \lieg^{2\beta}$; moreover, from the explicit density function (say as computed  Ledrappier-Young, see, \cite[Lem.\ 6.1.3]{MR4599404})) we have that  $\nu_x$ 
is the Haar measure.  This completes the proof.  
\end{proof}
\subsubsection{Proof of \cref{lem:injecttrans}}
We complete the proof of \cref{graph or invariant} with the proof of \cref{lem:injecttrans}.  
  \cref{lem:injecttrans} is established in a number of claims.

Let $d_1$ and $d_2$ denote, respectively, the $\mu$-essential supremums of $x\mapsto  d_1(x)$ and $x\mapsto d_2(x)$.  
\begin{claim}\label{claim:constrank}
For $\mu$-a.e.\ $x$ and every $y\in N_x$, 
$\dim(T_yN_x\cap \lieg^{[\beta]}_y)= d_1$ and $ \dim (T_yN_x\cap E^{\beta, F}(y))=d_2.$
\end{claim}
\begin{proof}
The map $y\mapsto d_j(y)$ upper semicontinuous for $j\in \{1,2\}$ and thus for each $N_x$, the function  takes its minimal value on an open subset of $y\in N_x$.  Using $A$-equivariance of the submanifolds $x\mapsto N_x$ and   selecting $a_0\in A$ that acts ergodically with respect to $\mu$ and with $\beta(a_0)>0$, by Poincar\'e recurrence, we find that $y\mapsto d_j(y)$ is constant on $N_x$ for a.e.\ $x$.  By $A$-equivariance of $N_x$ and ergodicity, we conclude that the constant value is independent of $x$.  
%
\end{proof}

Let $\pi_{x,2}^\beta$ denote the projection of $W^{[\beta]}(x)$ onto $W^{[\beta],F}(x)$ whose level sets are $U^{[\beta]}$-orbits; similarly let $\pi_{x,1}^\beta$ denote the projection of $W^{[\beta]}(x)$ onto $U^{[\beta]}\cdot x$ whose level sets are $W^{[\beta],F}$-manifolds.

  By  \cref{claim:constrank}, for a $\mu$-a.e.\ $x$, the restrictions of $\pi_{x,1}^\beta$ and $\pi_{x,2}^\beta$ to $N_x$ have  constant rank.  The implicit function theorem then implies the following. 
\begin{claim}
For $\mu$-a.e.\ $x$ and every $y\in N_x$, the intersections
$(U^{[\beta]}\cdot y)\cap N_x$  and $(W^{[\beta],F}(y))\cap N_x$ are a $C^1$ embedded submanifolds of dimension $d_1$ and $d_2$, respectively,  tangent to $T_yN_x\cap \lieg^{[\beta]}_y$ and $T_yN_x\cap E^{\beta, F}(y)$, respectively, at
$y$.
\end{claim}

To conclude the proof of  \cref{lem:injecttrans} we assume that either $d_1>0$ or $d_2>0$ and conclude that $\mu$ has extra invariance along $U^{[\beta]}$.  
\begin{claim}
If $d_1\ge 1$, then  $\mu$ is invariant under a 1-parameter subgroup  $H\subset U^{[\beta]}$.  
\end{claim}

\begin{proof}
Write $E_x =  \lieu^{[\beta]}_x \cap T_xN_x$,  $E_{\beta,x}= \lieg^\beta_x \cap  T_xN_x$, and  $E_{2\beta,x}= \lieg^{2\beta}\cap T_xN_x$.  
By Oseledec's theorem, and equivariance of $T_xN_x$, for $\mu$-a.e.\ $x$ we have
$$E_x= E_{\beta,x} \oplus  E_{2\beta,x}.$$  

If $d_1\ge 1$, there is a positive measure set of $x$ for which either 
$E_{\beta,x}$ or $ E_{2\beta,x}$ has positive dimension. 
 By $A$-equivariance of $N_x$ and $A$-ergodicity of $\mu$, we may select either $\gamma = \beta$ or $\gamma= 2\beta$ so that $E_x:=E_{\gamma,x}$ has positive dimension  for almost every $x$.  
Given any $a\in A$, we have $\Ad(a) E_x= \gamma(a) E_x= E_{a\cdot x}$ for almost every $x$.  

Under the canonical identification  $\lieg^\gamma_x\simeq \lieg^\gamma$, there is a subspace  $E\in \lieg^\gamma$, such that $E= E_{x}$ for $\mu$-a.e.\ $x$.   
Since  $\mu_y^{[\beta], \calE'_\beta}$ has full support in $N_y$, we may pass to the closure and conclude  for $\mu$-almost every $x$, and every $y\in N_x$, that $E= \lieg_y^{[\beta]} \cap T_yN_z$.

Fix a non-zero $Y\in E\subset \lieg^\gamma$. 
It follows that if $H= \exp_\lieg (tY)$ is the 1-parameter subgroup tangent to $Y$ then the orbit  $H\cdot y$ is contained in $N_x$ for every $y\in N_x$ and $\mu$-a.e.\ $x$.





Let $\scrH$ denote the partition into $H$-orbits.  
By assumption and Fubini's theorem, the leafwise measure $\mu^{\scrH,\calE_\beta}_y$ is absolutely continuous for $\mu^{\calE_\beta}_x$-a.e.\ $y$.  
Consider any  $a_1\in A\sm A'_\beta$ with $\beta(a_1)>0$.  
	 By the generalization of the Pesin entropy formula  (e.g.\ \cite{MR693976}),  we have $h_\mu(a_1\mid \scrH\vee {\calE_\beta}) = \gamma(a_1)$;  moreover, from the explicit density function (say as computed  Ledrappier-Young, see, \cite[(6.1)]{MR4599404})) 
 $\mu^{\scrH,\calE_\beta}_y$ 
is the Haar measure along the parameterized $H$-orbit for $\mu_x^{\calE}$-a.e.\ $y$.  It follows that $\mu$ is $H$-invariant.  
\end{proof}

If we assume hypothesis \ref{caseatom}  of \cref{graph or invariant} then it follows that $d_2=0$.  If we instead assume  hypothesis \ref{casedimfib}  of \cref{graph or invariant}, we have the following.
\begin{claim}\label{claim:whocares}
If  $\dim E^{F, [\beta]}(x) \le \dim \lieg^{ [\beta]}$ and if $d_2\ge 1$, then  $\mu$ is invariant under a one-parameter subgroup  $H\subset U^{[\beta]}$.  
\end{claim}
\begin{proof}
By \cref{cor:dimofproj}, we have $\dim N_x\ge d_2 + \dim \lieg^{[\beta]}$.  Since we assume $\dim T_x W^{[\beta]}_x\le 2\dim \lieg^{[\beta]}$, if $d_2\ge 1$ 
$$\dim T_xN_x + \dim \lieg^{[\beta]}_x> \dim T_x W^{[\beta]}_x$$
and thus $$\dim T_xN_x \cap \dim \lieg^{[\beta]}_x\neq\{0\}$$
for $\mu$-a.e.\ $x$, whence $d_1>0$ and the conclusion follows from the above arguments.
\end{proof}

\section{Proof of main results from Theorems \ref{thm:main} and \ref{thm:mainlamina}}\label{sec6}
We recall our standing assumptions on $G$ and $\Gamma$ in \cref{ss:Gassump}.
 Let $M$ be a compact manifold let $\alpha\colon \Gamma\to  \Diff^r(M)$ be an action.

\subsection{Proof of \cref{thm:globalrigid2}}\label{see pf main}
We recall the main strategy of \cite{2105.14541}.  Let $\Gamma$ be a lattice in a Lie group, all of whose factors are higher-rank.  Then $\Gamma$ is finitely generated and has strong property (T) (for representations into certain Banach spaces, including all Sobolev spaces of symmetric 2-tensors over $M$; see e.g.\ \cite{MR4018265} for formulation).  Let $|\gamma|$ denote the word length of $\gamma\in \Gamma$ relative to some choice of finite generating set.  An action $\alpha\colon \Gamma\to \diff^1(M)$ has \emph{subexponential growth of derivatives} if for all $\epsilon>0$ there is a $C>1$ such that for all $\gamma\in \Gamma$,
$$\sup_{x\in M} \|D_x \alpha(\gamma)\| \le Ce^{\epsilon |\gamma|}.$$
The following was established in \cite[Theorem 2.4]{MR4502593} (for $C^2$ actions) and extended in \cite{BDZ} for $C^{1}$ actions.
\begin{proposition}[{\cite[Theorem 2.4]{MR4502593}}]\label{prop:Isom}
If $\Gamma$ has strong property (T) and if $\alpha\colon \Gamma \to \diff ^{1}(M)$ has subexponential growth of derivatives then there exists a H\"older continuous Riemannian metric $g_0$ preserved by $\alpha$.  That is, $\alpha(\Gamma)\subset \mathrm{Isom}(M, g_0)$.  
\end{proposition}
We remark that if $r\ge 2$, it is shown in  \cite[Theorem 2.4]{MR4502593} that the invariant metric in \cref{prop:Isom} is actually $C^{r-1-\delta}$ for any $\delta>0$.  In \cite{BDZ}, a H\"older metric is obtained for actions that are $C^1$.  

We recall the main technical result from the work of Brown, Fisher, and Hurtado.  In the following, $\calZ$ denotes the center of a connected Lie group $G$ with semisimple Lie algebra $\lieg$.  
\begin{proposition}[{\cite[Theorem D]{2105.14541}}]\label{BFHthm}
Let $G$ be a connected semisimple Lie group without compact factors and with rank at least $2$.   Let $\Gamma$ be an irreducible lattice subgroup in~$G$, let $M$ be a compact manifold, and let $\alpha\colon \Gamma\to \Diff^{1}(M)$ be an action.
Let $M^\alpha = (G \times M)/\Gamma$ denote the induced $G$-space.  

If $\alpha$ fails to have {uniform subexponential growth of derivatives} then there exists a split Cartan subgroup~$A$ of~$G$ and a Borel probability measure~$\mu$ on~$M^\alpha$ such that
\begin{enumerate}
\item $\mu$ is $(\calZ A)$-invariant,
\item $\mu$ projects to the Haar measure on $G/\Gamma$,
	and
	\item for some $a \in \calZ A$,
	the average top fiberwise Lyapunov exponent $\lambda^F_{\top, \mu, a_0, \alpha}$ is positive.
	\end{enumerate}
\end{proposition}

\begin{proof}[Proof of \cref{thm:globalrigid2}]
We assume $G$ is a $\R$-split semisimple Lie group, all of whose simple factors are higher rank.   
If  $\alpha \colon \Gamma\to \Diff^r(M)$ were an action by isometries for some $C^0$ Riemannian metric on $M$ then, by Margulis superrigidity theorem and dimension counting, the image $\alpha(\Gamma)$ would be finite (see \cref{isometric}).
It is known that $\Gamma$ has strong property (T).  Thus, if we assume the image $\alpha(\Gamma)$ is infinite then, by \cref{prop:Isom}, the action $\alpha\colon \Gamma\to \diff^r(M)$ fails to have subexponential growth of derivatives.  Let $\mu$ be the ergodic, $(\calZ A)$-invariant Borel probability measure on $M^\alpha$ guaranteed by \cref{BFHthm}.

Since $\mu$ has a non-zero Lyapunov exponent, it follows from Zimmer's cocycle superrigidity theorem and the dimension of $M$ that  $\mu$ is not $G$-invariant.  By \cref{NRimpInv}  and dimension counting, it follows that there exists a standard parabolic subgroup $Q\subset G$ with codimension $v(G)$ (see \cite[Proposition 3.5]{MR4502593}) such that $\mu$ is $Q$-invariant.  
Again, by dimension counting, we also have $c_0^F(\mu)=1$.   
The result then follows from \cref{thm:completepara}.
\end{proof}

\subsection{Proof of \cref{thm:globalrigid2nonsplit}}
\cref{thm:globalrigid2nonsplit} similarly follows directly from the above discussion.  
\begin{proof}[Proof of \cref{thm:globalrigid2nonsplit}]
We assume $\alpha(\Gamma)$ fails to act by isometries with respect to any $C^0$ Riemannian metric on $M$.
We apply \cref{BFHthm} and \cref{prop:Isom} to conclude that there exists an ergodic,  $ (\calZ A)$-invariant Borel probability measure $\mu$ on $M^\alpha$ projecting to the Haar measure on $G/\Gamma$ such that $\lambda^F_{\top, \mu, a, \alpha}$  is positive for some $a\in (\calZ A)$.  

Let $Q= \stab (\mu)$.  Since $\dim M= v(\lieg) < n(\lieg)$, applying Zimmer's cocycle superrigidity theorem it follows that $Q\neq G$.  Let $\beta\in \Sigma$ be a resonant root such that  $\dim \lieg^{[\beta]}=2$  and $\mu$ is not $U^\beta$-invariant.  From \cref{thm:parabolicmeasrdim2}, 
either $\dim E^{[\beta],F}(x)\ge 2$ or $\dim E^{[\beta],F}(x)= 1$ and $\mu$ is invariant under a 1-parameter subgroup of $U^\beta$; in the latter case, the codimension of $\lieq=\Lie (Q)$ in $\lieu^\beta$ is at most $\dim E^{[\beta],F}(x)$.  Summing over all $\beta$ with $\lieu^\beta\not\subset \lieq$, it follows the codimension of $\lieq$ is at most $\dim(M)= v(G)$.  From \cite[Lemma 3.7]{MR4502593},  it follows that $Q$ is parabolic and has codimension $v(G)$.  From \cref{lem:alge,thm:parabolicmeasrdim2}, it then follows that $\dim E^{[\beta],F}(x)= \dim \lieg^{[\beta]}$  for every resonant root $\beta$ and almost every $x$.  
We may find $a_0\in A$ such that $\beta(a_0)<0$ for every $\beta \in \Sigma_Q^\perp$.  
The conclusion then follows from \cref{thm:completepara}.
\end{proof}

\begin{remark}\label{rem:reduce regularity}. 
In the proof of \cref{thm:globalrigid2nonsplit}, we only used the assumption that $\alpha\colon \Gamma\to \diff^\infty(M)$ is a $C^\infty$ action in the final sentence when invoking  \cref{thm:completepara} since, in general, we can not control $c^F_0(\mu)$.  
However, much of the remainder of the proof holds assuming $\alpha\colon \Gamma\to \diff^r(M)$ for $r>1$.  In many cases, the combinatorics in  \cref{thm:globalrigid2nonsplit} allow us to establish \cref{thm:globalrigid2nonsplit} in much lower regularity.  
Indeed, let $\mu_0$ be the  $A$-invariant Borel probability measure guaranteed by  \cref{BFHthm} and \cref{prop:Isom}.  Apply the trick in \cite[Lem.\ 3.4$\ast$]{ABZ}, we may assume there is a $P$-invariant  Borel probability measure $\mu$ with $\lambda^F_{\top, \mu, a_0, \alpha}>0$ for some $a \in \calZ A$.  

Let $Q=\stab (\mu)$ be the parabolic subgroup stabilizing $\mu$.  
Let $H\subset Q\subset G$ be the semisimple part of a Levi component of $Q$. 
Specifically, we have $\Lie(Q) = \lieq = \lieq _{S}$ for some collection of simple roots $S$ and may choose $H$ to have Lie algebra  $\lieh$ generated by $\lieg^{[\pm \beta]}$ for all simple roots $\beta\notin S$.  
Write $A_H=H\cap A$ and where $A'=\bigcap_{\beta\notin S}\ker \beta$.  Then $A'$ commutes with $H$ and $A= A_H\cdot A'$.

When $\lieh$ is higher rank, we may apply Zimmer's cocycle superrigidity theorem to the (ergodic components of) the $H$-action on $(M^\alpha, \mu)$.  In this case,  we conclude that resonant exponents for the $A_H$ action on $(M^\alpha, \mu)$ are given by weights of a non-trivial representation of $H$.  This constrains the combinatorics of $c_0^F(\mu)$.  




For instance, when $G= \Sl(n,\C)$ and $\dim (M) = 2n-2=v(G)$ we may apply  \cref{thm:parabolicmeasrdim2,lem:alge} to conclude that $\lieq$ has codimension $2n-2$ and omits exactly $n-1$ roots.  Then $\lieq$ is necessarily of the form $\lieq = \lieq_{\alpha_1}$ or $\lieq = \lieq_{\alpha_{n-1}}$ where $\alpha_1$ and $\alpha_{n-1}$ are the left-most and right-most roots in the Dynkin diagram. 
For $n\ge 4$, $Q$ always contains a copy of $H=\SL(n-1,\C)$.  Applying cocycle superrigidity to the $\SL(n-1,\C)$-action on $(M^\alpha,\mu)$, the coarse fiberwise Lyapunov exponent of the $A_H$-action on $(M^\alpha,\mu)$ are given by the weights of the standard (or contragradient) representation of $\SL(n-1,\C)$ on $\C^{n-1}$.  
As $A$ commutes with $A_H$, it follows that every coarse fiberwise Lyapunov exponents for the $A$-action on $(M^\alpha,\mu)$ is a single Lyapunov exponent with multiplicity 2.  
It follows that $c_0^F(\mu) = 1$.  
%
%
%
%
\end{remark}


\subsection{Proof of Theorems \ref{thm:localrigid} and \ref{thm:localrigidC0}}
Let $M= G/Q$ and let $\alpha_0\colon \Gamma\to \Diff^\omega(M)$ denote the standard projective action.  By \cref{uniquePmeasure}, for the induced action on $M^{\alpha_0}$ there exists a unique $P$-invariant probability measure $\mu_0$ on $M^{\alpha_0}$; moreover,  $\mu_0$ is $Q$-invariant but not $Q'$-invariant for any $Q\subsetneq Q'$.  

\subsubsection{Construction of conjugacy}
 For any $\alpha\colon \Gamma\to \Diff(M)$, the stabilizer of a  $P$-invariant Borel probability measure $\mu$ on $M^\alpha$ is a parabolic subgroup $Q_\mu\supset P$. For $\alpha$ sufficiently $C^0$ close to $\alpha_0$, we claim $Q_\mu\subset Q$.  Indeed, consider any sequence of $\{\alpha_j\}$ converging in $C^0$ to $\alpha_0$ and a corresponding sequence of $P$-invariant probability measures $\{\mu_j\}$ converging to the unique $P$-invariant probability measure $\mu_0$ on $M^{\alpha_0}$.  As there are only finitely many parabolic subgroups $Q$ containing our fixed minimal parabolic subgroup $P$, by passing to a subsequence, we may  suppose  $\stab (\mu_j) = Q'$ for all $j.$
  By \cref{felatiowithgoats}, 
 $\mu_0$ is $Q'$-invariant whence  (by \cref{uniquePmeasure}) $Q'\subset Q$.
 

In the case of  \cref{thm:localrigid}, we claim for any $\alpha$ sufficiently $C^1$ close to $\alpha_0$,  every $P$-invariant Borel probability measure $\mu$    on $M^\alpha$ is  $Q$-invariant.
  Indeed, given $r>r_0(\lieq)$,   if $\alpha\colon \Gamma\to \diff^r(M)$ is sufficiently $C^1$ close to $\alpha_0$ then \cref{cor:critregusc} implies that  $$c_0^F(\mu)\le r_0^F(\mu)< 
r_0(\lieq) +(r-r_0(\lieq)) =r.$$
\cref{prop:equifinite}  then implies $\mu$ is $\beta$-tame for every resonant root $\beta$.
Since the stabilizer $Q'$ of $\mu$ is a subgroup of $Q$ and since $\dim M= \dim G- \dim Q< \dim G- \dim Q',$ 
\cref{NRimpInv} and \cref{cor:easycor} 
imply the stabilizer of $\mu$ is precisely  $Q$.  

 In the case of \cref{thm:localrigidC0}, given $\alpha\colon \Gamma\to \Diff^r(M)$ sufficiently $C^0$ close to $\alpha$
 and any $P$-invariant measure $\mu$ on $M^\alpha$, we have $\stab(\mu)\subset Q$ as discussed above.   \cref{NRimpInv}  and a dimension count then implies  that $\stab(\mu)=Q$ and that $c_0(\mu) = 1$.    


Thus, in the case of either \cref{thm:localrigid} or \cref{thm:localrigidC0}, we have $c_0^F(\mu)<r$.  In particular, the hypotheses of Theorem \ref{thm:main} hold.  We thus obtain a $C^r$ conjugacy $h\colon G/Q\to G/Q$ intertwining the actions $\alpha$ and $\alpha_0$.  
In particular, in either \cref{thm:localrigid} or \cref{thm:localrigidC0} we obtain local rigidity of the action.  

\subsubsection{Strong local rigidity in \cref{thm:localrigid}}  
Now suppose that $\alpha_n\to \alpha_0$ in $C^1$.
By the above, there is a  $C^r$ diffeomorphism $h_n$ with $\alpha_n(\gamma)\circ h_n= h_n \circ \alpha_0(\gamma)$ for all $\gamma\in \Gamma$; moreover by \cref{rem:unique} there is a unique such  $h_n$.  
For each $n$, let $\wtd \mu_n$ be the lift to $G\times M$ of the unique $P$-invariant Borel probability measure on $M^{\alpha_n}$.  
We have that $\wtd \mu_0$ is the image of the Haar measure on the graph $$g\mapsto (g, g\inv Q)$$ and 
$\wtd \mu_n$ is the image of the Haar measure on the graph $$g\mapsto (g, h_n (g\inv Q)).$$
The projection to $M=G/Q$ of the restriction of $\wtd \mu_0$ to a fundamental domain for $\Gamma$ in $G$ is in the Lebesgue class.  As $\wtd \mu_n\to \wtd \mu_0$ as $n\to \infty$, we have $h_n$ converges to the identity in measure (with respect to any fixed measure in the Lebesgue class on $G/Q$.)

 We may apply \cite[Thm.\ 2.8]{MR0302822}
to find $A\subset G$ that is $\Q$-anisotropic and thus $A\cap \Gamma$ is a cocompact lattice in $A$.  We may further assume from the beginning that the minimal parabolic $P$ was chosen so that 
$A\subset P\subset Q$.  
We may then take $\gamma_0\in A\cap \Gamma$ with the following properties: if $\lambda_1<\dots<\lambda_\ell<0$  denote the distinct eigenvalues of $\restrict{\Ad(\gamma_0)}{\liev_\lieq}$,
then $$r_0(\lieq)\le  \frac {\lambda_1}{\lambda_\ell}<s\le r.$$

Note that $\alpha_0(\gamma_0)$ has a unique contracting fixed point $x_0= \1Q \in G/Q$.   

Enumerating $\td\gamma_j\in \Gamma$, let $\gamma_j=\td\gamma_j\circ \gamma_0\circ\td \gamma_j\inv.$
Let $x_j= \td \gamma_jQ= \alpha_0(\td \gamma_j)(x_0)$ be the unique fixed point for 
$\alpha_0(\gamma_j).$  Pick $g_j\in G$ such that $x_j= g_jQ$.  
Let $V$ denote the vector space $V= \liev_\lieq$ .  
The map $\phi_j\colon V\to G/Q$ given by $$\phi_j\colon X\mapsto    \td \gamma _j \exp(X) Q$$ parameterizes an open dense neighborhood of $x_j$ (a Schubert cell) in $G/Q$. 
Moreover, relative to these charts, for each $j$ the dynamics of $\alpha_0(\gamma_j)$ coincides with the linear map $L= \restrict{\ad \gamma_0}{\liev_\lieq}$.  

We may select finitely many $1\le j\le J$ and equip $V$ with an inner product for which 
$\{\phi_j(V(1/2)): 1\le j\le J\}$ is an open cover of $G/Q$ (where $V(r)$ is the ball of  radius $r>0$ centered at 0 in $V$).

 Suppose  $\alpha$ is sufficiently $C^s$ close to $\alpha_0$.   For each $1\le j\le J$, the map $f_j =  \phi_j\inv \circ \alpha(\gamma_j)\circ \phi_j\colon V(2)\to V(2)$  is   a $C^s$ perturbation of $L$.  Taking $\alpha$ sufficiently $C^1$ close to $\alpha_0$, we may assume 
$f_j\colon V(1)\to V(2)$ is well defined and has a unique fixed point $p_f\in V(1/2)$. 
Let 
\begin{enumerate}
\item $\td f_j(x) =f_j(x+p_f)-p_f$
\item $h_j =\phi_j\inv \circ h\circ \phi_j\colon V(1)\to V(2)$
\item $\td h_j(x) =h_j(x-p_f)+p_f$.
\end{enumerate}
We have 
 \begin{enumerate}
 \item $\td f_j(0) = 0$ and $\|\td f_j-L\|_{V(1), C^s}\le 2\|f_j-L\|_{V(1), C^s}$
 \item $\td h_j(0) = 0$ and $\|h_j-\id\|_{V(1/2), C^s}\le \|\td h_j-\id \|_{V(1), C^s} +\|p_f\|$
 \end{enumerate}
It now follows from \cref{prop:smallconj} below that if $\alpha$ is sufficiently $C^s$-close to $\alpha_0$ and if $h$ is sufficiently close in measure to the identity, then $h$ is, in fact, $C^s$-close to the identity.

\subsubsection{$C^s$-small conjugacies for perturbations of linear maps}
Let $L\colon V\to V$ be a linear map with eigenvalues $\lambda_1<\dots <\lambda_\ell<0$, listed without multiplicity.


In the following, fix any $s> \lambda_1 /\lambda_\ell$.  We consider  a $C^s$-small perturbation $f$ of $L$ and show that if $f$ is conjugate to $L$ by a $C^s$ diffeomorphism $h$ that is $C^0$-close to the identity (in fact, close in measure to the identity), then $h$ is, in fact, $C^s$-close to the identity.  



\begin{proposition} \label{prop:smallconj}
Let $V$, $L\colon V\to V$, and $s$ be as above and 
fix $r\ge s> \lambda_1 /\lambda_\ell$.  
There exists a  polynomials $p_1$ with $p_1(0) = 0$  
and a constant $0<\eta<1$ 
with the following properties:
Let  $f\colon V\to V$ be a $C^r$ diffeomorphism such that  
\begin{enumerate}[label=(\alph*), ref=(\alph*)]
\item $f(0) = 0$,
\item  $\|f-L\|_{V(1),C^{s}}<\eta$ 
\item there exists a $C^r$ diffeomorphism $h\colon V\to V$ such that   $h\circ f = L\circ h.$
\end{enumerate} 
Then there exists a $C^s$ diffeomorphism $\what h \colon V\to V$ with 
\begin{enumerate}
\item $\|\what h-\Id\|_{V(1),C^{s}} \le p_1(\|f-L\|_{V(1),C^{s}})$
\item $ \what h \circ f= L\circ \what h.$
\end{enumerate}
Moreover, given any $\delta>0$, if $\|f-L\|_{V(1),C^{s}}$ is sufficiently small and if $h$ is sufficiently close in measure to $\Id$, then 
\begin{enumerate} [resume]
\item $\|h-\Id\|_{V(1),C^{s}} \le \delta$.  
\end{enumerate}
\end{proposition}

Although the arguments are standard, we could not find a proof of \cref{prop:smallconj} in the literature and thus provide a proof in \cref{appendix of nonsense}.

\subsection{Proof of Theorems \ref{thm:smoothfactor} and \ref{thm:smoothfactorSPLIT}}\label{sec:pfsmoothfac}
Combining  Dani's theorem, \cref{thm:dani}, and Theorems \ref{thm:parabolicmeasr} and \ref{thm:completepara}, we obtain a  proof of the smooth factor theorem.
\begin{proof}[{Proof of \cref{thm:smoothfactor,thm:smoothfactorSPLIT}}]
Let $M$ be a compact manifold, $\alpha\colon \Gamma\to \diff^r(M)$ an action, and $p\colon G/P\to M$ a $\Gamma$-equivariant, $C^0$ surjection.  By Dani's theorem, \cref{thm:dani}, there exists a parabolic subgroup $P\subset Q\subset G$ such that $p$ factors through a $\Gamma$-equivariant homeomorphism $q\colon G/Q\to M$.  

Let $\wtd q\colon G\times (G/Q)\to G\times M$ be $\wtd q(g, g'Q) = \left(g,  q(g\inv g'Q)\right)$. Then $\wtd q$ is left-$G$-equivariant and right-$\Gamma$-equivariant and descends to a $G$-equivariant homeomorphism $$\bar q \colon (G/\Gamma) \times (G/Q) \to M^\alpha.$$

Let $\mu_0$ be the unique $P$-invariant Borel probability measures on $ (G/\Gamma) \times (G/Q) $ and let $\mu= \bar q _* \mu_0$.  
Moreover, as $\bar q$ is a $G$-equivariant homeomorphism, $\mu$ is $Q$-invariant and is  not  $Q'$-invariant for any $Q'\supsetneq Q$.  
Since $M$ is homeomorphic to $G/Q$, we have  $$\dim(M)= \dim(G/Q)= \dim V_\lieq= \sum_{\beta\in \Sigma_\lieq^\perp} \lieg^\beta.$$
By \cref{NRimpInv,thm:parabolicmeasr,factortame}, for any $s\in C(\Sigma_Q^\perp)\sm \{0\}$, all fiberwise Lyapunov exponents are negative for  the action of $s$ on $(M^\alpha, \mu)$


If all almost simple factors of $G$ are $\R$-split then all fiberwise Lyapunov exponents have multiplicity $1$, then  $c_0^F(\mu) = 1$.  
Otherwise, we assume the factor dynamics $\alpha\colon \Gamma\to \Diff^\infty(M)$ is $C^\infty$.  In either case, we may apply \cref{thm:main} to obtain a $\Gamma$-equivariant, $C^r$ covering map 
$h\colon M\to G/Q$.   
Since $h\circ q\colon G/Q\to G/Q$ is a $\Gamma$-equivariant continuous surjection, \cref{rem:unique} implies $h\circ q$ is the identity map whence $h$ is a diffeomorphism and $q = h\inv.$
\end{proof}

\subsection{Proof of \cref{thm:embedded projectives}}
\begin{proof}
Let $M$ and $\alpha$ be as in  \cref{thm:embedded projectives}.  Let $\mu$ be a $P$-invariant measure that is not $G$-invariant.  Then, there is a maximal parabolic subgroup $Q\subsetneq G$ of codimension $v(G)+1$ such that $\mu$ is $Q$-invariant.  As we assume that the action $\alpha$ preserves a no-where vanishing density,  the codimension of $Q$ is actually $v(G)$.
As $Q$ is a maximal proper parabolic subgroup, \cref{rem:hypo}, the hypotheses of  \cref{thm:mainlamina} hold and \cref{thm:embedded projectives}  follows.  
\end{proof}



\section{Measure rigidity and local graphs of conjugacies}\label{sec:msrrig}
We give the main  main technical results in this section  from which all our main theorems will  follow.  The main results are a number of nonlinear measure rigidity results and their consequences.  

\subsection{\texorpdfstring{Leafwise measures along $U^{[\beta]}$-orbits}{Leafwise measures along $U$-orbits}}
We begin with the following fact which follows from standard measure rigidity arguments (for semisimple homogeneous actions).  
\begin{proposition} \label{prop:msrrigiU} 
Let $\Gamma$ be an irreducible lattice in a higher-rank Lie group $G$ and let $\alpha\colon \Gamma\to \diff^r(M)$ be an action on a compact manifold.   Let $\mu$ be an ergodic, $A$-invariant Borel probability measure on $M^\alpha$ projecting to the Haar measure on $G/\Gamma$.  

Fix $\beta\in \Sigma$ and let $\calE_\beta$ denote the decomposition of $\mu$ into $A'_\beta= \ker\beta$-ergodic components.  Then one of the following holds: 
\begin{enumerate}
\item For $\mu$-a.e.\ $x$, the leafwise measure
$\mu^{U^{[\beta]},\calE_\beta }_x$ of the $A'_\beta=\ker\beta$-ergodic component of $\mu$ containing $x$ along $U^{{[\beta]}}$-orbits is a single atom. 
\item The measure $\mu$ is invariant under a connected subgroup $H\subset U^{[\beta]}$ of positive dimension; moreover, if $\mu$ is assumed $U^{[-\beta]}$-invariant then $\mu$ is $U^{[\beta]}$-invariant.  
\end{enumerate}
\end{proposition}

As in many previous arguments, the $U^{[\beta]}$-invariance in conclusion (2) follows from \cref{lem:alge}.  

\begin{proof}
Assume for a positive measure subset of $x\in M^\alpha$ that the leafwise measure  $\mu^{U^{[\beta]},\calE_\beta }_x$ is not a single atom.  By equivariance of $\mu^{U^{[\beta]},\calE_\beta }_x$ under the $A$-action, this holds for $\mu$-a.e.\ $x$.  

Fix $\mu'$, an ergodic component of $\mu$ projecting to the Haar measure on $G/\Gamma$.  
We have 
$$ \mu^{U^{[\beta]},\calE_\beta }_x=(\mu')^{U^{[\beta]} }_x.$$
By \cite[Thm.\ 1]{MR0283174}, we may select $a_0\in A_\beta'$ acting ergodically on $(M^\alpha, \mu')$.  

We equip $U^{[\beta]}$ with a right-invariant metric which induces a metric on each $U^{[\beta]}$-orbit. Note that $a_0$ acts by isometries on $U^{[\beta]}$-orbits and preserves the canonical parametrization by left-translation on orbits: \begin{equation}\label{eq:param}a_0 \cdot (u\cdot x) = u\cdot (a_0\cdot x)\end{equation}
for all $x\in M^\alpha$ in $u\in U^{[\beta]}$.  
By the standard measure rigidity arguments, for $\mu'$-a.e.\ $x$ and  $\mu^{U^{[\beta]},\calE_\beta }_x$-a.e.\ $y \in  U^{[\beta]}\cdot x$, there is an isometry $\phi\colon U^{[\beta]}\cdot x\to U^{[\beta]}\cdot x$ such that $\phi(x) = y$ and $\phi_* \left(\mu^{U^{[\beta]},\calE_\beta }_x\right)$ coincides with $\mu^{U^{[\beta]},\calE_\beta }_x$ up to normalization.  
In fact, $\phi$ is of the form of   right translation: if $y= u_y\cdot x$ then for all $u\in U^{[\beta]}$,
$$\phi(u\cdot x) = u\cdot y= u\cdot (u_y\cdot x).$$

It now follows from standard arguments (see for example Lemmas 6.1--6.3 in \cite{MR2122918} following Section 5 of \cite{MR1406432}) that for $\mu'$-a.e.\ $x$, the measure $\mu^{U^{[\beta]},\calE_\beta }_x$ is the image of the Haar measure on the orbit $V_x\cdot x$ of a connected subgroup $V_x\subset  U^{[\beta]}$.  
(Briefly, after slight reformulation \cref{prop:measurerig} below shows the measure $\mu^{U^{[\beta]},\calE_\beta }_x$ is in the Lebesgue class on a connected embedded submanifold $N_x\subset U^{[\beta]}\cdot x$; \eqref{eq:param} shows $T_yN_x\subset \lieu^{[\beta]}$ is constant as $y$ varies in $N_x$.  It follows that $N_x$ is the orbit of a subgroup $V_x$ in $U^{[\beta]}$ and standard entropy arguments show $\mu^{U^{[\beta]},\calE_\beta }_x$ measure is the Haar measure on this $V_x$-orbit.)

Finally, $A$-equivariance of $\mu^{U^{[\beta]},\calE_\beta }_x$ implies the map $x\mapsto \Lie(V_x)$ is $A$-invariant, hence constant $\mu$-a.s.  It then follows that $\mu$ is invariant under a subgroup $V\subset U^{[\beta]}$.  
\end{proof}

\subsection
{\texorpdfstring{Measure rigidity via \(\beta\)-tameness or one-dimensionally of fiberwise exponents}{Measure rigidity via beta-tameness or one-dimensionally of fiberwise exponents}} 
The following is the main technical result of the paper.  The proof will occupy \cref{sec:coccyle}--\cref{sed:mrse}.

We consider an ergodic, $A$-invariant Borel probability measure on $M^\alpha$ projecting to the Haar measure on $G/\Gamma$ and $\beta\in \Sigma$ such that $\mu$ is not invariant under a $1$-parameter subgroup of $U^{[\beta]}$.  \cref{NRimpInv} implies $\beta$ is resonant and thus there is a  coarse fiberwise Lyapunov exponent $\chi^F$  (for the $A$-action on $(M^\alpha, \mu)$) that is  positively proportional to $\mu$.  We write $$ W^{{[\beta]}, F}(x)=  W^{\chi^F, F}(x)$$ for the associated fiberwise coarse Lyapunov manifolds in $\calF(x)$.  
On the total coarse  Lyapunov manifold $W^{[\beta]}(x)$ we show the leafwise measures (conditioned on a certain ergodic decomposition) are actually homogeneous and of the correct dimension.

\begin{proposition} \label{prop:msrrigi2} 
Let $\Gamma$ be an irreducible lattice in a higher-rank Lie group $G$ and let $\alpha\colon \Gamma\to \diff^r(M)$ be an action on a compact manifold.   Let $\mu$ be an ergodic, $A$-invariant Borel probability measure on $M^\alpha$ projecting to the Haar measure on $G/\Gamma$. 

Fix $\beta\in \Sigma$ such that $\mu$ is not invariant under any $1$-parameter subgroup of $U^{[\beta]}$.  Assume that $r>c_0^F([\beta])$.  
Let $\calE_\beta$ denote the decomposition of $\mu$ into ergodic components for the action of $A'_\beta= \ker\beta$.  Suppose  that one of the following conditions hold  for a.e.\ $x$:
\begin{enumerate}
\item  the  measure $\mu$ is $\beta$-tame and $\dim E^{[\beta], F}(x) \le  \dim \lieg^{[\beta]}$; or
\item  the  measure $\mu$ is $\beta$-tame and the fiberwise leafwise measure $\mu^{{[\beta]},\calE_\beta, F }_x$ is a single atom $\delta_x$ for $\mu$-a.e.\ $x$; or 
\item $\dim E^{[\beta],F}(x) =1$.
 \end{enumerate}   
 Then for $\mu$-a.e.\ $x$, 
there exist a $C^r$-embedding  $\phi_x^{[\beta]}\colon U^{[\beta]}\to W^{{[\beta]}, F}(x)$ (where $W^{{[\beta]}, F}(x)$ is given its immersed manifold structure) with $  \phi_x^{[\beta]}(\1) = x$
such that the leafwise measure  
$\mu^{{[\beta]},\calE_\beta }_x$ of the $A'_\beta$-ergodic component of $\mu$ containing $x$ along $W^{{[\beta]}}(x)$ is the image of the Haar measure on $U^{[\beta]}$  under the graph of  $  \phi_x^{[\beta]}$,
$$u\mapsto u\cdot  \phi_x^{[\beta]}(u).$$
\end{proposition}
The proof of \cref{prop:msrrigi2} occupies \cref{sed:mrse}.  
\begin{remark}
In most applications, we will apply \cref{prop:msrrigi2}  to ($A$-ergodic components) of ergodic $P$-invariant probability measures on $M^\alpha$ (in which case, the factor measure on $G/\Gamma$ is always Haar).  However, for some situations, it is useful to have a formulation of  \cref{prop:msrrigi2} for $A$-invariant measures that are not (a priori) $P$-invariant.

In the case that $\mu$ is an ergodic $P$-invariant measure $\mu$ on $M^\alpha$, the assumption that  $\mu$ is not invariant under any $1$-parameter subgroup of $U^{[\beta]}$ is equivalent to 
the assumption that 
$\mu$ is not $U^{[\beta]}$-invariant (by \cref{lem:alge,ergcomp}).
\end{remark}

\subsection{Entropy constraints in \cref{thm:mainlamina} } 
In the case that all fiberwise Lyapunov exponents are negative, we may apply \cref{cor:easycor}  (or its proof) to obtain coincidence of the leafwise measures $\mu_x ^{[\beta], \calE_{\beta}}$ and $ \mu_x ^{[\beta]}$  for almost every $x$.  
When we apply \cref{prop:msrrigi2} in the setting of \cref{thm:mainlamina}, we do not have that all fiber exponents are negative and so must argue why such coincidence would hold.  

Recall  our hypotheses in \cref{thm:mainlamina} imply that  $\dim E^{[\beta],F}=1$ (for every  root $\beta\in \Sigma_Q^\perp$).    We argue in the next proposition  that $\mu_x ^{[\beta], \calE_{\beta}}= \mu_x ^{[\beta]}$ a.s.\ by showing that $\mu_x ^{[\beta]}$ is not supported on distinct smooth graphs.  

\begin{proposition}\label{agree support} 
Let $\Gamma$ be an irreducible lattice in a higher-rank, $\R$-split simple Lie group $G$. 
For $r>1$, let  $\alpha\colon \Gamma\to \diff^r(M)$ be an action on a compact manifold.   Let $\mu$ be a $P$-invariant Borel probability measure on $M^\alpha$.  Let $Q=\stab(\mu)$ and suppose hypotheses \ref{exptan} and \ref{exptrans} of \cref{thm:mainlamina} hold. 

Fix $\beta\in \Sigma_Q^\perp$,   
 let $A'_\beta$ denote the kernel of $\beta$ in $A$, and let $\calE_{\beta}$ denote the ergodic decomposition of $\mu$ with respect to\ $A'_\beta$.

Then for $\mu$-a.e.\ $x$, we have coincidence of leafwise measures $\mu_x ^{[\beta], \calE_{\beta}}=  \mu_x ^{[\beta]}.$  
\end{proposition}
\begin{proof}
By hypotheses of \cref{thm:mainlamina}, $1= \dim E^{[\beta],F}(x)$ for almost every $x$ and thus $c_0^F([\beta])=1$.  
 We know from \cref{prop:msrrigi2} that the conditional leafwise measure $\mu_x ^{[\beta], \calE_{\beta}}$ is the image of the Haar measure  under the graph of  a diffeomorphism $\phi_{\beta,x}\colon U^{[\beta]}\to W^{[\beta],F}(x)$.

Fix $a\in A$ with $\beta(a)>0$.  If $\mu_x ^{[\beta], \calE_{\beta}}\neq  \mu_x ^{[\beta]}$ then for $\mu$-a.e.\ $x$, the measure $\mu_x ^{[\beta]}$ is   locally the average of the images of the Haar measure under  multiple diffeomorphisms $U^{[\beta]}\to W^{[\beta],F}(x)$.  This implies  for $\mu$-a.e.\ $x$ that the leafwise measures 
$\mu_x ^{U^{[\beta]}}$ and $\mu_x ^{W^{[\beta],F}}$ are non-trivial measures. 
Thus, if the conclusion fails, we necessarily have both $h_\mu(a\mid \scrW^{[\beta],F} ) \neq 0$ and $h_\mu(a\mid U^{[\beta]} ) \neq 0$.  

Suppose the conclusion fails so that $\mu_x ^{U^{[\beta]}}$ is a non-trivial measure.  
We consider two cases:
\subsubsection*{Case 1:} Suppose there is $\gamma\in \Sigma$ with $[\gamma]\neq[\beta]$ such that $\gamma+\beta\in \Sigma$ and $\mu_x ^{U^{[\gamma]}}$ is nontrivial.  Note that $\gamma $ is not negatively proportional to $\beta$ so we may find $a'\in A$ with  both $\beta(a')>0$ and $\gamma(a')>0$.  
  From the high entropy method \cite[Theorem 8.5]{MR2191228}, it follows that $\mu$ is invariant under a subgroup of $U^{[\gamma+\beta]}$.  We thus obtain  a contradiction if such $\gamma$ exists with $\gamma+\beta\in \Sigma_Q^\perp$.

\subsubsection*{Case 2:} Suppose there is no $\gamma\in \Sigma$ such that 
\begin{itemize}
\item $\gamma+\beta\in \Sigma$,
\item  $\mu_x ^{U^\gamma}$ is nontrivial, and 
\item $\gamma+\beta\in \Sigma_Q^\perp$.
\end{itemize}
It follows that $\beta$ is a simple negative root and for every simple (positive) $\gamma$ root adjacent  to $-\beta $ in the Dynkin diagram,  the leafwise measure $\mu_x ^{U^{[-\gamma]}}$ along $U^{[-\gamma]}$-orbits is atomic.  (Here, 
we use the assumption that $G$ has no rank-1 factors to ensure such $\gamma$ exists.) In particular, $-\gamma\in \Sigma_Q^\perp$ and  $h_\mu(a\mid \scrW^{[-\gamma],F} ) = 0$ for all $a$ (especially $a$ with $-\gamma(a)>0$.)

Fix $a\in \exp(-\overline{C(\Sigma_\lieq ^\perp)}\sm \{0\})$ such that 
\begin{itemize}
\item $\beta'(a) = 0$ for every simple negative root $\beta'\in \Sigma_Q^\perp$ with $\beta'\neq -\gamma$.  Note this includes $\beta'=\beta$ so that $\beta(a)=0$.  
\item $-\gamma(a)>0$
\item $\beta''\ge 0$ for every root $\beta''\in \Sigma_Q^\perp.  $
\end{itemize}

We claim that that $h_\mu(a\mid \calF)= 0$.  
By the choice of $a$, the only fiberwise exponents that are positive at $a$ are positively proportional to roots in $\Sigma_Q^\perp$.  
If $\chi^F_i\in S$ is proportional to a non-simple negative  root $\beta''\in \Sigma_Q^\perp$, case 1 shows $h_\mu(a\mid \scrW^{[\beta''],F} ) = 0$.  Also, $\beta''(a) =0$ for all simple negative roots $\beta''$ with $\beta''\notin \Sigma_Q^\perp$ and thus $h_\mu(a\mid \scrW^{[\beta''],F} ) = 0$ for such roots.   
Finally $h_\mu(a\mid \scrW^{[-\gamma],F} ) = 0$ by the assumption that $\mu_x ^{U^\gamma}$ is not nontrivial.
By the product structure of fiberwise entropy and  hypotheses \ref{exptan} and \ref{exptrans} of \cref{thm:mainlamina}, it follows that $h_\mu(a\mid \calF) = 0$.

  Perturbing from $a$ to $a'\in \exp(-C(\Sigma_\lieq ^\perp))$ adds additional positive fiberwise exponents but does not create any additional  negative fiberwise exponents.  Thus  
  $$0=h_\mu(a\mid \calF)=h_\mu(-a\mid \calF)=h_\mu(-a'\mid \calF)=h_\mu(a'\mid \calF).$$
  Since $\beta(a)>0$, we conclude  $h_\mu(a\mid \scrW^{\beta,F} ) =0$.  
\end{proof}

\subsection{Upgrading homogenity along coarse foliations to homogenity along the stable foliation}
We have the following strengthening of the conclusion of \cref{prop:msrrigi2}.  
%

\begin{proposition}\label{upgradetofullstab}
Let $\Gamma$ be an irreducible lattice in a higher-rank Lie group $G$ and let $\alpha\colon \Gamma\to \diff^r(M)$ be an action on a compact manifold.   
Let $\mu$ be an ergodic, $P$-invariant Borel probability  measure on $M^\alpha$.  Let $Q$ be the parabolic subgroup  with 
$Q=\stab(\mu)$. Suppose the following hold:
\begin{enumerate} 
\item  For every $\beta\in \Sigma_\lieq^\perp$ there is  
an  $C^r$  embedding  $$\phi_{x,{[\beta]}}\colon U^{[\beta]}\to  W^{[\beta], F}(x)$$ with $  \phi_{x,{[\beta]}}(\1_U^{[\beta]}) = x$ such that the leafwise measure  
$\mu^{[\beta] }_x$  along the leaf $W^{[\beta]}(x) = U^{[\beta]}\cdot W^{[\beta], F}(x)$ is the image of the Haar measure on $U^{[\beta]}$  under the graph    
$$u\mapsto (u\cdot  \phi_{x,{[\beta]}}(u))$$
\item for every   $a_0\in C(\Sigma_\lieq^\perp)$ and $\mu$-almost every $x\in M^\alpha$, $$E^{s,F}_{a_0}(x) = \bigoplus _{\beta\in \Sigma_\lieq^\perp} E^{[\beta],F}(x).$$
\end{enumerate}
Fix $a_0\in   C(\Sigma_\lieq^\perp)$ and (following notation in \cref{eqcriffiberreg,stableexps}) suppose that   $$r>   r_0^F(\calI^s(a_0)).$$ 
Write $V= V_\lieq$.  
Then for  $\mu$-almost every $x\in M^\alpha$  and every   $a_0\in C(\Sigma_\lieq^\perp)$,   there is 
an injective $C^r$ embedding
$\phi_{x,Q}$ with $\phi_{x,Q}(\1_V) = x$ such that the leafwise measure  
$\mu^{s}_x$  along $ V\cdot W^{s, F}_{a_0}(x)$ is the image of the Haar measure on $V$  under the graph    
$$u\mapsto (u\cdot \phi_{x,Q}(u))$$

\end{proposition}

The proof of \cref{upgradetofullstab}
appears in \cref{sec:upgrade}.  



\section{Constructing a laminated action from measure classification}
Assuming the conclusion of  \cref{upgradetofullstab}, we construct a lamination with properties in \cref{thm:mainlamina}.  
\subsection{Setup and formulation of main result}
Fix $G$ and $\Gamma$ as in \cref{ss:Gassump}.  
 Let $M$ be a connected, compact manifold, and let $\alpha\colon \Gamma\to C^r(M)$ be an action.  In applications, we will assume  $r>1$ is  for sufficiently large (depending on $\lieg$) so that results of \cref{sec:msrrig} apply to imply hypothesis \ref{hyp:11} of \cref{propmain} below.  
 Let $M^\alpha$ denote the suspension space with induced $G$-action.
 Given $x\in M^\alpha$, write  $$M_x= \calF(x) = \pi\inv (\pi(x))$$ for the fiber of $M^\alpha$ through $x$ 
and $F:= \ker D\pi$ for the fiberwise tangent bundle.

\subsubsection{Lifted stable manifolds}
Recall we write $C(\Sigma_\lieq ^\perp)$ for the (open) cone in $A$ for which $\beta(a)<0$ for all $\beta\in \Sigma_\lieq ^\perp$.  
Fix $a_0\in C(\Sigma_\lieq ^\perp)$.
Following notation in \cref{FWstable}, let $\Omega= \Omega(a_0)$ and for $x\in \Omega$ recall we write  $W^{s}_{a_0}(x)$ for the stable manifold through $x$.  
For such $x$, the preimage of $W^{s}_{a_0}(x)$  under the  map $G\times M\to M^\alpha$ is countably many path-connected, injectively immersed submanifolds; we write  $\wtd W^s_{a_0}(g,x)$ for the path-connected component containing $(g,x)$ of the preimage of stable manifold $W^s_{a_0}([g,x])$.    

For $x\in \Omega$, write $$\wtd W^{s,F}_{a_0}(g,x)= \wtd W^{s}_{a_0}(g,x)\cap (\{g\}\times M)$$ for the fiber component of the lifted stable manifold.  Write $ W^{s,F}_{a_0, g}(x)\subset M$, for the projection of $\wtd W^{s,F}_{a_0}(g,x)$ to $M$.
That is, $$\wtd W^{s,F}_{a_0}(g,x)= \{g\} \times W^{s,F}_{a_0, g}(x).$$


\subsubsection{Standing assumptions on the measure $\mu$} 
Let $\mu$ be an ergodic, $P$-invariant Borel probability measure on $M^\alpha$.   Let $Q\supset P$ be the parabolic subgroup with $\stab(\mu) = Q$.   
\subsubsection{Disintegrations of $\wtd \mu$}
Let $\wtd \mu$ denote the canonical lift of $\mu$ to $G\times M$ normalized over a fundamental domain $D$ for $\Gamma$ in $G$.
The partition of $G\times M$ into fibers $\{g\}\times M$ is measurable and thus the (locally finite, Radon) measure $\wtd \mu$ admits a family of conditional measures $\{\wtd \mu_g\}$.  Formally, each measure $\wtd \mu_g$ is supported on the fiber of the form $\{g\}\times M$; we instead view each
 $\wtd\mu_g$ as a Borel probability measure on $M$ so that for any Borel set $A\subset G\times M$, $$\wtd \mu (A) = \int_G  \wtd \mu_g (\{y\in M: (g,y)\in A\}) \, dg.$$
Under this convention, we have 
 \begin{equation}\label{eqgamm}
\alpha(\gamma\inv )_*\wtd \mu_g  = \wtd \mu_{g\gamma}.
 \end{equation}
and $$\wtd \mu_{qg}=\wtd \mu_g$$ for every $\gamma\in \Gamma$, a.e.\ $g\in G$ and a.e.\ $q\in Q$. 

\subsubsection{Definition of the measure $\nu$}
Let $\wtd \mu_{Qg}$ 
denote the a.s.\ constant value of $\wtd \mu_{qg}$ over $q\in Q$. 
Recall that $K\cap \mathcal Z$ is cocompact in $K$.   
Let $K'$ be a fundamental domain for $K\cap \mathcal Z$ in $K$.   Let $m_K$ be the Haar measure on $K$.  Then we may define 
\begin{equation}\label{eq:fornicatewithfourier}\nu:=\frac{1}{m_K(K')}\int_{K'} \wtd \mu_{Qk} \, d k.\end{equation}

Let $\bar \mu$ be the probability measure on $G/Q\times M$ whose marginal is the (left) $K$-invariant probability measure on $G/Q$ and whose disintegration along each fiber of the form $\{gQ\}\times M$ is $\wtd \mu_{Q g\inv }$.  
The group 
$\Gamma$ acts on $G/Q\times M$ as \begin{equation}\label{eq:lefttwisted}\gamma \cdot (gQ,x) = (\gamma g Q, \alpha(\gamma)(x)).\end{equation}
Observe from \eqref{eqgamm} that $$\alpha(\gamma )_*\wtd \mu_{Qg\inv}  = \wtd \mu_{Qg\inv \gamma\inv } = \wtd \mu_{Q(\gamma g) \inv }.$$
In particular, the natural projection $(G/Q\times M, \bar \mu)\to (G/Q, \mathrm{Leb})$ is a relatively measure-preserving extension.
Let $\Pi_2\colon G/Q\times M\to M$ be the canonical projection.  Clearly 
\begin{enumerate} 
\item $\Pi_2$ is $\Gamma$-equivariant.  
\end{enumerate}
Moreover, as  $m_K$ is invariant under $k\mapsto -k$, for $\nu$ as defined by \eqref{eq:fornicatewithfourier}
\begin{enumerate}[resume]
\item $\nu= (\Pi_2)_*\bar \mu$.
\end{enumerate}
In particular, 
\begin{enumerate}[resume]
\item $\nu$ is $\alpha(\Gamma)$-quasi-invariant.
\end{enumerate}

\subsubsection{Main theorem}
Similar to our constructions  in \cref{sec:leafwisemeasires}, we may define leafwise measures $\wtd \mu_{(g,x)}^s$ for the measure $\wtd \mu$ conditioned on almost every leaf $\wtd W^{s,F}_{a_0}(g,x)$ in $G\times M$.

With the above notation, we have the following.

\begin{theorem}\label{propmain}
For $r>1$, let $\alpha\colon \Gamma\to \diff^r(M)$ be an action and 
let  $\mu$ be  an ergodic, $P$-invariant Borel probability measure on $M^\alpha$.  Let $\stab(\mu) = Q$, and  $V= V_\lieq$.  Suppose there is $a_0\in C(\Sigma_\lieq^\perp)$ such that the following hold:
\begin{enumerate}[label=(\alph*), ref=(\alph*)]
\item\label{hyp:11} For  $\wtd \mu$-a.e.\ $(g,x)\in G\times M$ there exists a $C^r$-diffeomorphism $$\phi_{(g,x)}\colon V\to  W^{s,F}_{a_0,g}(x)$$ such that  
\begin{enumerate}[label=(\roman*), ref=\roman*]
\item $\phi_{(g,x)}(\1_V) = x$ and 
\item the leafwise measure $\wtd \mu_{(g,x)}^s$ on $\wtd W^s_{a_0}(g,x)$ is the image of the Haar measure on $V$ under the graph $$v\mapsto (vg, \phi_{(g,x)}(v)).$$
\end{enumerate}
\item \label{hyp:22} There is a partition $$\calL^F= S  \sqcup  U$$ of the fiberwise Lyapunov exponents $\calL^F= \{\lambda_1^F, \dots, \lambda_\ell^F\}$ 
where for $\mu$-a.e.\ $x\in M^\alpha$,
\begin{enumerate}[label=(\roman*), ref=\roman*]
\item $E^{\lambda_i^F}(x) \subset T_x  W^{s,F}_{a_0}(x)$ for all $\lambda_i^F\in S$, and 
\item$ \lambda_i^F(a)>0$ for all $a\in \overline{C(\Sigma_\lieq ^\perp)}$ and all $\lambda_i^F\in U$.
\end{enumerate}

\end{enumerate}

Then there exists a $\nu$-a.e.\ defined measurable function $h\colon M \to G/Q$, an $\alpha(\Gamma)$-invariant, $\nu$-measurable subset $\Lambda\subset M$, and a measurable {left} $G$-action on $\Lambda$---denoted by $\ell_g$---for which the following hold:

 \begin{enumerate}
 \item \label{prop:main3}$\nu(\Lambda) =1$ and  $\Pi_2\colon ((Q\bs G)\times M, \bar \mu)\to (\Lambda,\nu)$ is a measurable isomorphism. 
\item \label{prop:main2}  $h(x)$ is defined for every $x\in \Lambda$ and 
for $\bar \mu$-a.e.\ $(gQ, x)\in (G/Q)\times M$, $$h (x)= h\circ \Pi_2(gQ, x)= g Q.$$ 
In particular, 
\begin{enumerate}
\item $h_*\nu$ is the $K$-invariant Lebesgue measure on $G/Q$;
\item   $h\colon (M,\nu)\to (G/Q, \mathrm{Leb})$ is a relatively measure-preserving extension.  

\end{enumerate}

	\item \label{prop:main4}For $x\in \Lambda$, the $G$-orbit $L(x):= \{\ell_g(x):g\in G\}$ is an injectively immersed $C^r$ submanifold of $M$. 
 \item \label{prop:main5} The map $h\colon \Lambda\to G/Q$ is $G$-equivariant: for $x\in \Lambda$ and $g\in G$,  $h(\ell _g(x))= g\cdot  h(x)$.  
 \item \label{prop:main6} For $x\in \Lambda$, the connected component of identity in the  $G$-stabilizer of $x$ is a conjugate of $Q^\circ$.  Moreover, the restriction of $h\colon \Lambda\to G/Q$ to $L(x)$ defines a $C^r$ covering map. 


\item \label{prop:main7} For $\gamma\in \Gamma$, we have $\alpha(\gamma)(L_{x}) = L_{\alpha(\gamma)(x)}$.  Moreover, given $x_0\in \Lambda$, let $$x_1 = \ell_{\gamma\inv}( \alpha(\gamma)(x_0)).$$
  Then for all $g\in G$,
$$\alpha(\gamma)(\ell_g(x_0)) = \ell_{\gamma g}(x_1)= \ell_{\gamma g\gamma\inv}(\alpha(\gamma)(x_0)) .$$

\item \label{prop:main8} For a.e.\ $g,g'\in G$, $\wtd \mu_{g\inv}$ is supported on $\Lambda$ and 
	\begin{equation}\label{eq:potato}\wtd \mu_{(g' g)\inv }= (\ell_{g'})_* \wtd \mu_{g\inv}.\end{equation}

\end{enumerate}
\end{theorem}

Note that using \eqref{eq:potato} we may modify $g\mapsto \wtd \mu_g$ on a null set so that $\wtd \mu_g$ is defined for every $g\in G$, satisfies  \eqref{eq:potato} for all $g,g'\in G$, and satisfies $\wtd \mu_{qg} = \wtd \mu_g$ for all $q\in Q$.

We also note that the hypotheses of \cref{propmain} necessarily imply the parabolic subgroup $Q\subset P$ is a proper subgroup of $G$.  Indeed, if $G= Q$ then the open cone $C(\Sigma_\lieq^\perp)$ is empty and thus the standing hypotheses \ref{hyp:11} and \ref{hyp:22} can not hold.

\begin{remark}
Suppose there exists $x\in \Lambda$ such that  $\Gamma':= \{\gamma\in \Gamma: \alpha(\gamma)(L_x) = L_x\}$ is of finite index in $\Gamma$.  

Fix $x_0\in L_x$ stabilized by $Q^\circ$ and let $Q'\supset Q^\circ$ denote the full stabilizer of $x_0$.  
Identify $L_x$ with $ G/Q'$  by $\ell_g(x) \mapsto gQ'$.  By \cref{rem:classifyactions}, conclusion \eqref{prop:main7} of \cref{propmain}  reduces to the following: we may find a finite index subgroup $\Gamma''\subset \Gamma'$ such that 
$$\alpha(\gamma)(gQ') = \gamma g Q'$$
or $$\alpha(\gamma)(\ell_g(x_0)) = \ell_{\gamma g}(x_0)$$
 for all $\gamma\in \Gamma''$.  We then have the following:
\end{remark}
\begin{claim}
Let $\hat \Gamma$ be a lattice in a connected, semisimple Lie group $\hat G$ with finite center and no compact factors.  Let $\hat \alpha\colon \hat \Gamma \to \diff^{r}(M)$ be an action.  
Suppose for the suspension $M^\alpha$, there exists a $P$-invariant measure satisfying the criteria of \cref{propmain}.  

Let $\Lambda$ be  as in the conclusion of \cref{propmain}.  Suppose for some $x\in \Lambda$ that $\Gamma'=\{\gamma\in \hat \Gamma: \alpha(\gamma)(L_x) = L_x\}$ has finite index in $\hat \Gamma$. 
Then $L_x$ is a finite cover of $G/Q$ and  compact.   
\end{claim}
\begin{proof}
We may lift the action $\hat \alpha$ of $\hat \Gamma$ to an action $\alpha$ by a lattice $\Gamma$ in the simply connected cover  $G$ of $\hat G$.  Let $\cent(G)$ be the center of $G$.  Then, there is a finite-index subgroup $Z'$ of $\cent(G)\cap \Gamma$ such that $\alpha(\gamma) = \id_M$ for all $\gamma\in Z'$.


Fix $x\in \Lambda$  as above.  Then   $\Gamma'=\{\gamma\in \Gamma: \alpha(\gamma)(L_x) = L_x\}$ has finite index in $\Gamma$.   Let $x_0\in L_x$ be stabilized by $Q^\circ$.  
For $\gamma\in Z'\cap \Gamma'$ and all $y= gx_0\in L_x$ we have 
$$g\cdot x_0 = y = \bar{\alpha}(\gamma)(y) = \gamma \cdot gx_0= g\cdot \gamma x_0.$$
It follows that $\gamma $ stabilizes $x_0$.  Since $\Gamma'$ has finite index in $\Gamma$, $Z'\cap \Gamma'$ has finite index in $\cent(G)$.  It follows that $Z'Q^\circ $ is cocompact in  $G$ whence $L_x = G/(Z'Q^\circ)$ is compact.  
\end{proof}

\subsection{Atlases on generalized flag varieties} 
\label{Sec:lok}
Fix   a minimal parabolic subgroup $P\subset G$ and consider a parabolic subgroup  $P\subset Q\subset G$  with lie algebra $\lieq$.    Recall we write $$\Sigma_\lieq:= \{ \beta: \lieg^\beta\subset \lieq\},\quad \quad \Sigma_\lieq^\perp= \Sigma \sm \Sigma_\lieq.$$    Also recall we write  $\liev_\lieq = \bigoplus _{\beta\in  \Sigma_\lieq^\perp} \lieg^\beta$.  We have that $\liev_\lieq$ is a Lie subalgebra of $\lieg$ and write $V=V_\lieq$ for the analytic subgroup with Lie algebra $\liev_\lieq$.  
Recall that  $V_\lieq \cdot Q$ is an open dense subset of $G$.  


We specify certain  local coordinates on the generalized flag variety space $Q\bs G$. 
 Given $g\in G$, the function $V\to Q\bs G$, given by $$v\mapsto Qvg$$ is a $C^\infty$ diffeomorphism onto an open dense neighborhood of $Qg$ in  $Q\bs G$.  
As $QVQ= G$, we have for any fixed $g\in G$ that the family of coordinate maps $$\calA_{Qg} = \{v\mapsto Qv(qg) : q\in Q\}$$ forms an atlas of $Q\bs G$; we call the collection $\{v\mapsto Qv(qg) : q\in Q\}$ an \emph{atlas centered at the coset $Qg\in Q\bs G$}.


The group $G/Z$ is adjoint and hence, identified with its image under the adjoint representation, we identify $G/Z$ with the connected component of a real-algebraic group.  Identified with their images under the adjoint representation, the groups $ V_\lieq,$  and $Q/Z$ have real algebraic structures.

Note that  $q\in Q$ need not normalize $V$.  However, as $QV$ is dense in $G$, given a fixed $q\in Q$ there is an open dense set of $v\in V$ for which  conjugate $qvq\inv $ can be uniquely written  as $$qvq\inv = \hat q \hat v$$ for some $\hat q\in Q$ and $\hat v\in V$. 
More specifically,  there is a function $\tau$ defined from an open subset of $V\times Q$ to $ V\times Q$ of the form $$\tau(v,q) = \left(\tau_1(v,q), \tau_2(v,q)\right)$$ uniquely determined by the following:  if $(v,q)$ is in the domain of $\tau$ we have 
$$\tau_1(v,q)\cdot q = \tau_2(v,q)\cdot q\cdot v.$$
Note that $\tau(v,q) = \tau(v, qz)$ for any $z\in Z$; in particular, $\tau$ descends to a rational function $\tau\colon V\times Q/Z\to V\times Q/Z$.  

Given a fixed $q\in Q$, let  $\tau_q\colon V\to V$ be the function $\tau_q(v) = \tau_1(v,q)$. We have that $\tau_q\colon V\to V$ is a rational function.  


For $v$ in the domain of $\tau_q$ we have 
$$\tau_q(v) q = \bar q q v$$ for some $$\bar q = \tau_2(v,q)\in Q.$$
It follows for such $v$ that  $Q\tau_q(v) (q g) = Q v g$.  In particular, if $\phi_1\colon v\mapsto Qv(qg)$  and $\phi_2\colon v\mapsto Qvg$  are two charts in the atlas centered at $Qg$, then $\phi_1\inv \circ \phi_2 \colon V\to V$ is the birational map $$\phi_1\inv\circ \phi_2 = \tau_q.$$

We observe a number of useful properties of the above change of coordinates.
\begin{claim}\label{claim:usefulidiot}
Given $q,q'\in Q$ and  $v\in V$, write $\tau_q(v) q= \bar q v$ for $\bar q\in Q$.
Then the following equalities hold, assuming all terms are defined.  
\begin{enumerate}
\item\label{clmmm1} $(\tau_q(v))\inv = \tau_{\bar q}( v\inv)$
\item\label{clmmm2} $\tau_{\bar q\inv}((\tau_q(v))\inv) = v\inv$
\item\label{clmmm3} $\tau_{q'}(\tau_q(v)) = \tau_{q'q}(v)$
\item\label{clmmm4} $\tau_{\bar q}(u)\cdot \tau_{q} ( v) = \tau_q(u\cdot v)$ for $u\in V$. 
\end{enumerate}
\end{claim}
\begin{proof}
For \eqref{clmmm1}, we have 
$ qv\inv= (\tau_q(v))\inv \bar q $
and for \eqref{clmmm2}, we have 
$q\inv (\tau_q(v))\inv =  v\inv \bar q\inv$.

For \eqref{clmmm3}, we have $\tau_{q'q}(v)q'q = \wtd  q v$ and $\tau_{q'}(\tau_q(v))q'q = \widehat{q}'\tau_q(v)=
\widehat{q}'\widehat{q}v$  for some $\wtd  q ,\widehat{q}', \widehat{q}\in Q$ whence
$$\tau_{q'q}(v)\left(\tau_{q'}(\tau_q(v))\right)\inv =\wtd q \left (\widehat{q}' \widehat{q}\right )\inv.$$
 Since  $V\cap Q=\{1\}$, we conclude that 
$ \wtd  q= \widehat{q}'\widehat{q}$ and $\tau_{q'q}(v)=\tau_{q'}(\tau_q(v)).$

For \eqref{clmmm4}, we have 
$$\tau_{\bar q}(u)\cdot \tau_{q} ( v) \cdot q= 
\tau_{\bar q}(u)\cdot \bar q  v= 
\hat q\cdot u\cdot  v$$
for some $\hat q\in Q$; the result holds by definition.  
\end{proof}

Note that the domain of  $\tau_q\colon V\to V$ is open and dense in $V$ but need not be connected.  We thus introduce the following notation.
\begin{definition}
For $q\in Q$, let $V(q)$ denote the connected open subset of $V$ containing the $\1_V$ in the domain of $\tau_q$.  
\end{definition}

We further observe  the following:
\begin{enumerate}
\item Since $\tau_q (\1) = \1$, we have  that $\tau_q(V(q))= V(q\inv)$ and $\tau_{q\inv} = \tau _q \inv$.
\item If $q$ centralizes $V$ then $V(q) = V$ and $\tau_q(v) = v$ for all $v\in V$.  
\item If $q$ normalizes $V$ then $V(q) = V$.  
\end{enumerate}
\subsection{Equivariance of fiberwise dynamical foliations}\label{sec:equifol}
As above, consider  an ergodic, $P$-invariant Borel probability measure   $\mu$ on $M^\alpha$ and let $Q=\stab(\mu)$.   We may restrict to the action of any $a_0\in A\subset P$ on $(M^\alpha, \mu)$.  
Following the discussion in \cref{ssec:FWLyapu}, we construct fiberwise stable manifolds  $W^{s,F}_{a_0}(x)$ for the action of $a_0\in A$ and $\mu$-almost every $x\in M^\alpha$.  We have equivariance of $x\mapsto  W^{s,F}_{a_0}(x)$ under the action of the centralizer of $a_0\in G$.  That is, if $g\in C_{G}(a_0)$, then for $\mu$-a.e.\ $x\in M^\alpha$ we have equivariance of the stable manifolds in $M^\alpha$, 
$$g\cdot   W^{s,F}_{a_0}(x)=    W^{s,F}_{a_0}(g\cdot x).$$

However, while $a_0\in A$ normalizes $Q$, elements  $q\in Q$ need not commute with $a_0$.  Thus in full generality, 
we do not expect that  $q\cdot  W^{s,F}_{a_0}(x)=  W^{s,F}_{a_0}(q\cdot x)$ or even that $q\cdot  W^{s,F}_{a_0}(x)$ contains a neighborhood of $x$ in $  W^{s,F}_{a_0}(q\cdot x)$.  However, \cref{prop:equivarstablemanifolds} below ensures  under 
the hypothesis  of \cref{propmain} (especially hypothesis \ref{hyp:22}), that for $\mu$-a.e.\ $x$,    the intersection $$q\cdot  W^{s,F}_{a_0}(x)\cap  W^{s,F}_{a_0}(q\cdot x)$$ is an open subset of  $  W^{s,F}_{a_0}(q\cdot x)$ with a canonical parametrization (coming from  hypothesis \ref{hyp:11}).


Given $x\in M^{\alpha}$, write $E^s_{a_0}(x) = T_x  W^{s,F}_{a_0}(x)\subset F(x)$.  
Recall the fixed $a_0$ in the hypothesis of  \cref{propmain}.  We have $$E^s_{a_0}(x)= \bigoplus _{\lambda_i^F\in S} E^{\lambda_i^F}(x) = T_x  W^{s,F}_{a_0}(x).$$ Also write $$ E^u_{a_0}(x)= \bigoplus _{\lambda_i^F\in U} E^{\lambda_i^F}(x).$$


\begin{proposition}\label{prop:equivarstablemanifolds}
Suppose $\mu$ and $a_0$ satisfy hypotheses \ref{hyp:11} and  \ref{hyp:22} of \cref{propmain}. 
Fix $q\in Q $. 
Then for $\mu$-a.e.\ $x\in M^\alpha$
\begin{enumerate}
\item \label{vomitface1} $q\cdot E^s_{a_0,x} = E^s_{a_0, q\cdot x}$ and
\item \label{vomitface2} $q\cdot E^u_{a_0,x} = E^u_{a_0, q\cdot x}$.
\end{enumerate}
Moreover, for $\wtd \mu$-a.e.\ $(g,x)\in G\times  M$ and every $v\in V(q)$, 
\begin{enumerate}[resume]
 \item \label{vomitface3} 
$\phi_{(g,x)}(v) = \phi_{(qg, x)}(\tau_q(v)).$  
\end{enumerate}
\end{proposition}

In order to show conclusion \eqref{vomitface3} of \cref{sec:equifol}, we will actually argue as follows:
Fix $q\in Q$.  Let $U\subset  V$ be an open subset containing $\1_V$ and for $x\in \Omega\cap \Omega (a_0)$, write 
$$W^{s,F}_{a_0,g,U}(x)= \{\phi_{(g,x)}(v): v\in U\}$$
for the image of the restriction of $\phi_{(g,x)}$ to $U$.  Then $W^{s,F}_{a_0,g,U}(x)$ contains an open neighborhood of $x$ in the fiberwise stable manifold $W^{s,F}_{a_0,g}(x)\subset M$.  We then argue the following:

\begin{lemma}\label{lem:hard equiv}
For $q\in Q$ and   $\wtd \mu$-a.e.\ $(g,x)\in G\times M$, 
$$ W^{s,F}_{a_0,g, V(q)}(x)=   W^{s,F}_{a_0,qg,\tau_q(V(q))}(x).$$
\end{lemma}
From hypotheses \ref{hyp:11}  of \cref{propmain} and the assumption that $\mu$ is $Q$-invariant, \cref{lem:hard equiv} then implies 
for $\wtd \mu$-a.e.\ $(g,x)\in G\times M$ and every $v\in V(q)$ that, writing $\tau_q(v) q= \td q_{v} v $, we have 

\begin{align*}
 \bigl(\td q_v \cdot vg,\phi_{(g,x)}(v)\bigr) 
&= \td q_v \cdot \bigl(vg,\phi_{(g,x)}(v)\bigr) \\
&=
\bigl(\tau_q(v)qg, \phi_{(qg, x)}(\tau_q(v))\bigr)\end{align*}
whence 
$$\phi_{(g,x)}(v)=\phi_{(qg, x)}(\tau_q(v)).$$



 The remainder of this subsection is devoted to the proof of \cref {prop:equivarstablemanifolds}.

\subsubsection{Generators of $Q$, compact generation, and special cases of \cref {prop:equivarstablemanifolds}}\label{sss:compacttypes}
Fix the  set of simple roots $\Delta\subset \Pi$ with $\lieq= \lieq_\Delta$ (so that root spaces associated to roots with a non-zero coefficient associated to roots in $-\Delta$ are excluded from $\lieq$).
As remarked above,  \cref{propmain} necessarily implies that $Q\neq G$ and so $\Delta \neq \Pi$.  

Let $\bar Q$ denote the parabolic subgroup opposite to $Q$.  Then $V$ is the unipotent radical of $\bar Q$.  Let $V'$ denote the unipotent radical of $Q$.
Let $H= Q\cap \bar Q$.  Then $H$ is reductive and the center of $H$ is a subgroup of $MA$.  

We have $Q= H\ltimes V'$.  We consider two types of  elements of  $q\in Q$:
\begin{enumerate}
\item\label{c1} $q\in H$, or 
\item \label{c2} $q\in U^{[\beta]}$ for  $\beta \in \Delta$. 
\end{enumerate}
The root spaces tangent to $V'$ can be obtained by repeated brackets of elements of $\lieg^{\beta}$ for $\beta\in \Delta$ and $\lieg^{\beta'}$ for $\beta'\in \Pi\sm \Delta$. 
In particular,  there is $N\in \N$ such that   any $\td q\in Q$ is of the form 
$$\td q = \td q_k\cdot  \dots \cdot \td q_1$$ 
where each $\td q_j$ is either in $H$ or in  $U^{[\beta]}$ for  $\beta \in \Delta$.

Given finitely many $q_1, q_2, \dots, q_k\in Q$, let $V(q_1, \dots, q_k)$ denote the connected open neighborhood of $\1_V$ on which the birational functions
\begin{align*}
v&\mapsto \tau_{q_1}(v)\\ 
v&\mapsto \tau_{q_2}\circ \tau_{q_1} (v)\\ 
&\vdots\\
v&\mapsto \tau_{q_k}\circ \dots \circ \tau_{q_1} (v)
\end{align*}
are all well defined.

Fix $q\in Q$. Given $v\in V(q)$, we write 
$$\tau_q(v) q= \td q_{v} v. $$
The map $V(q) \to  V\times Q$, $v\mapsto (\tau_q(v) ,\td q_{v} )$ is the restriction of a birational map to its domain and hence is continuous.  
In particular, compactness and continuity immediately give the following:

\begin{lemma}\label{lem:uniform}
Fix $q\in Q$.  Given $v\in V(q)$, write $\td q_v v = \td v q$.

There exist $N\in \N$ and, for every compact $C\subset V(q)$, a compact subset $D \subset Q $	and an open neighborhood $U_C\subset V$ of $\1_V$, such that for every $v\in C$ the following hold:
\begin{enumerate}
	\item $\td q_v = \td q_k\cdot  \dots \cdot \td q_1$ for some $k\le N$;
	\item each $\td q_k\in D$;
	\item $U_C\subset V(\td q_1, \dots, \td q_k)$;
	\item each $\td q_j$ is of the form \eqref{c1} or \eqref{c2} above.  
	\end{enumerate}
\end{lemma}
%
%

Below, we will establish the following.  
\begin{lemma}\label{prop:equivarstablemanifolds2}
\cref{prop:equivarstablemanifolds} holds for every $q\in Q$ of type \eqref{c1} or  \eqref{c2} above.  
\end{lemma}

From \cref{prop:equivarstablemanifolds2,lem:uniform}, we obtain \cref{prop:equivarstablemanifolds}.  
\begin{proof}[Proof of \cref{prop:equivarstablemanifolds}]  
Conclusions \eqref{vomitface1} and \eqref{vomitface2} of the proposition follow immediately from \cref{prop:equivarstablemanifolds2,lem:uniform} and the chain rule.  

To show conclusion \eqref{vomitface3} of  the proposition, we apply \cref{lem:uniform} and show \cref{lem:hard equiv}.  Let $C\subset V(q)$ be a connected compact subset whose interior is dense in $C$.

Given  $x_0\in \Omega$ of the form   $x_0 = [g_0, z_0]$, write $$C(x_0) = \{[vg_0,\phi_{(g_0,z_0)}(v)]: v\in C\}.$$ 
We equip $C(x_0)$ with the image of the restriction of the Haar measure on $V$ to $C$ pushed forward via the parameter $v$.

We apply \cref{lem:uniform,prop:equivarstablemanifolds2}: 
For a.e.\ $v\in C$ and $x$ of the form $x=[vg_0,\phi_{(g_0,z_0)}(v)]$, 
write $\tau_q(v) q= \td q_{v} v $; then, applying \cref{claim:usefulidiot}\eqref{clmmm4}, we have
\begin{align*}
 W^{s,F}_{a_0,  g_0, U_C\cdot v}(x_0)&=
 W^{s,F}_{a_0,vg_0, U_C}(x)\\
 &=   W^{s,F}_{a_0,\td q_{v} vg_0,\tau_{\td q_{v}}(U_C)}(x)\\
 &=   W^{s,F}_{a_0,\tau_q(v) q g_0,\tau_{\td q_{v}}(U_C)}(x)\\
 &=   W^{s,F}_{a_0, q g_0,\tau_{ q}(U_C\cdot v)}(x_0).
 \end{align*}
Since  the above holds for a.e.\ $v\in C$ and since the set $U_C$ is uniform over the choice of $v\in C$, we concluded that 
$$ W^{s,F}_{a_0, g_0, C}(x_0)=  W^{s,F}_{a_0, q g_0,\tau_{ q}(C)}(x_0).$$
Exhausting $V(q)$ by sets with the form of $C$, the conclusion holds.  
\end{proof}

\subsubsection{Proof of \cref{prop:equivarstablemanifolds2}}
We now establish \cref{prop:equivarstablemanifolds2}, showing \cref{prop:equivarstablemanifolds} for each type of generator of element of $Q$.  

\subsubsection*{Case \ref{c1}: $q\in H$} 

Let $q\in H$.  Recall we have $\Delta \neq \Pi$.   We may then select  $a_*\in A$ such that 
such that 
$$\text{$\beta(a_*) = 0 $ for all $\beta \in \Pi\sm \Delta$ and $\beta(a_*) > 0 $ for all $\beta \in \Delta$.}$$
It follows that $\beta(a_*)<0$ for all $\beta\in \Sigma_\lieq^\perp$ and thus $a_*\in C(\Sigma_\lieq)$.  
By hypothesis \ref{hyp:22} of  \cref{propmain}, we then have 
\begin{enumerate}
\item $\lambda_i^F(a_*) <0$ for all $\lambda_i^F\in S$,
\item $\lambda_i^F(a_*) >0$ for all $\lambda_i^F\in U$,
\item  $W^{s, F}_{a_0}(x)=  W^{s, F}_{a_*}(x)$ for $\mu$-a.e.\ $x$, and 
\item $a_*$ 
is contained in the center of $H$.  
\end{enumerate}
Since $a_*$ commutes with $q\in H$, it follows for $\mu$-a.e.\ $x\in M^\alpha$ that 
$$ q\cdot W^{s,F}_{a_*}(x)=  W^{s,F}_{a_*}(q\cdot x)$$
whence for $\wtd \mu$-a.e.\ $(g,x)\in G\times M$, 
$$ W^{s,F}_{a_0,g}(x)=  W^{s,F}_{a_*,g}(x)= W^{s,F}_{a_*,q\cdot g}(x)= W^{s,F}_{a_0,q\cdot g}(x).$$

\subsubsection*{Case \ref{c2}: $q\in U^{[\beta_0]}$ for some  $\beta _0\in \Delta$.} 
We may select $a_{*}\in  \overline{ C(\Sigma_\lieq^\perp) }$ such  that 
$$\text{$\beta_0(a_{*}) = 0 $ and $\beta(a_{*}) > 0 $ for all $\beta \in \Pi\sm \{\beta_0\}$. }$$
If $G$ has real rank 2, we simply select $a_*$ from the ray in $\ker \beta_0\subset A$ on which the other simple root is positive.    When $G$ has real rank at least 3, we will impose further criteria on the choice of  $a_*\neq \1$ from an open cone in $\ker \beta_0$ below.  

Consider measurably-varying subdistributions of $F(x)\subset T_xM^\alpha$ given by 
$$\wtd E(x)= q\inv E^s_{a_0}(q\cdot x), \quad \quad   \widehat E(x)= q\inv E^u_{a_0}(q\cdot x).$$ 
Since $a_{*}$ commutes with $q$, from hypothesis \ref{hyp:22} of  \cref{propmain}  we have
\begin{enumerate}
\item $\wtd E$ and $\widehat E$ are $a_{*}$-equivariant
\item the Lyapunov exponents of $\restrict{Da_{*}}{\wtd E}$ are all non-positive
\item the Lyapunov exponents of $\restrict{Da_{*}}{\widehat  E}$ are positive.
\end{enumerate}
By the uniqueness ($\bmod$ 0) of the Lyapunov splitting for the dynamics of $a_{*}$, it follows  that 
$$\wtd E(x)= E^s_{a_0}(x), \quad \quad   \widehat E(x)= E^u_{a_0}(x)$$ 
$\mu$-almost everywhere.  Conclusions \eqref{vomitface1} and \eqref{vomitface2} of the proposition then follow.

We now focus on conclusion \eqref{vomitface3} of \cref{prop:equivarstablemanifolds2}, continuing to assume  $q\in U^{[\beta_0]}$.   By the discussion above, it is enough to establish \cref{lem:hard equiv}.    

Fix a compact subset $C\subset V(q)$ with the following properties:
\begin{enumerate}
\item the  interior of $C$ dense in $C$, and 
\item  $\Ad_{a_{*}}(C)\subset C$. 
\end{enumerate}
    We have the following technical lemma whose proof we present next.  
    Recall that given  $g_0\in G$ and $z_0\in M$, we write $x_0 =[g_0, z_0]$ for the associated element of $M^\alpha$.
\begin{lemma}\label{lem:hardestpart}
For $\mu$-a.e.\ $x_0\in M^\alpha$ of the form $x_0= [g_0,z_0]$ there is an open $U\subset V$ with $U\cdot \tau_q(C)\subset  V(q\inv )=\tau_q(V(q))$ and such that  for a.e.\ $v\in C$, writing $x=[vg_0, \phi_{(g_0,z_0)}(v)]$ and $\tau_q(v)q=\td q_v v$, we have that 
\begin{equation}\label{es:lamma} \td q_v\inv\cdot    W^{s,F}_{a_0, U}(\td q_v\cdot x)\subset W^{s,F}_{a_0} (x).\end{equation}
\end{lemma}
We note (by $Q$-invariance of $\mu$ and hypotheses \ref{hyp:11}  of \cref{propmain}) that \eqref{es:lamma} actually implies that $$ \td q_v\inv\cdot    W^{s,F}_{a_0, U}(\td q_v\cdot x)=   W^{s,F}_{a_0, \tau_{\td q _v\inv }(U)} (x).$$

We finish the proof of \cref{prop:equivarstablemanifolds2}, assuming \cref{lem:hardestpart}.
Given $v\in C$, we write $\tau_q(v)q=\td q_v v$.  Let 
\begin{equation}\label{eq:QC}Q(C) := \{\td q_v:v\in C\}.\end{equation} Observe that $Q(C)$ is compact.

Let $W\subset V$ be an open neighborhood of $\1$ in $V$ such that $W\subset \tau_{\td q}\inv (U)$ for every $\td q\in Q(C)$.
Applying \cref{claim:usefulidiot}\eqref{clmmm4}, $Q$-invariance of the measure, and hypotheses \ref{hyp:11}  of \cref{propmain}, we have
$$q\inv \cdot   W^{s,F}_{a_0, U\cdot \tau_q( v)} ( q\cdot x_0)
=  W^{s,F}_{a_0, \tau_{\td q_v}\inv (U)\cdot v} ( x_0)
\supset     W^{s,F}_{a_0, W\cdot v} (x_0).$$ 
Thus
$$ W^{s,F}_{a_0,q g_0,U\cdot  \tau_q(v)}(z)\supset   W^{s,F}_{a_0,g_0,W\cdot v}(z).$$
for a.e.\ $v\in C$; since the size of $W$ is uniform (so that, in particular, sets of the form $\{W\cdot v:v\in C_0\}$ for any co-null subset $C_0\subset C$ cover $C$) we have 
$$ W^{s,F}_{a_0,q g_0, \tau_q(V(q))}(z)\supset   W^{s,F}_{a_0,g_0,C}(z).$$
 
An exhaustion by countably many sets of the form of $C$ as above  then shows 
$$ W^{s,F}_{a_0,q g_0, \tau_q(V(q))}(z)\supset   W^{s,F}_{a_0,g_0,V(q)}(z).$$
Since $\tau_q(V(q))= V(q\inv)$, symmetry in $q$ and $q\inv$ shows \cref{lem:hard equiv}, completing   the proof of \cref{prop:equivarstablemanifolds2}.

\subsubsection{Proof of \cref{lem:hardestpart}}  
Recall our fixed $q\in U^{[\beta_0]}$.  Let $Q(C)$ be as in \eqref{eq:QC}.  Given $v\in V(q)$, recall we write $\td q_v\in Q$ and $\tau_q(v)\in V$ for the unique elements such that $$\tau_q(v) q= \td q_{v} v .$$   Since $\Ad_{a_{*}}(q) = q$, we have 
$\Ad_{a_{*}}(\tau_q(v)) q= \Ad_{a_{*}}(\td q_{v} )(\Ad_{a_{*}}(v))$
and thus conclude 
\begin{enumerate}
\item $\Ad_{a_{*}}(\tau_q(v)) = \tau_q( \Ad_{a_{*}}(v))$ and 
\item $\Ad_{a_{*}}(\td q_{v} )=\td q_{\Ad_{a_{*}}(v)} .$
\end{enumerate}
Given $x_0 \in M^\alpha$ of the form $x_0= [g_0, z_0]$ and $v\in C$, write 
$$x(v, x_0)=[vg_0,\phi_{(g_0,z_0)}(v)] $$ 
and 
\begin{align*}\bar x (v, x_0) &= x(\tau_q(v), q\cdot x_0 )\\&= [\td q_v vg_0,\phi_{(g_0,z_0)}(v)] 
\\&=[\tau_q(v) q g_0,\phi_{(qg_0,z_0)}(\tau_q(v))] \\&=\td q_v\cdot x(v,x_0).
\end{align*}

Let $f\colon M^\alpha\to M^\alpha$ be the map $$f(x) := a_{*}   \inv \cdot x.$$  
Given $\rho>0$, let $\liev(\rho)$ denote the open ball of radius $\rho$ centered at $0$ in $\liev= \Lie(V)$. 
Fix $\rho_0>0$ such that if $U_0=\exp_\liev \liev(\rho_0)$, then $U_0 \cdot C\subset V(q)$ and $U_0 \cdot \tau_q(C)\subset V(q\inv)$.  In particular, $U_0\subset V(\td q)\cap V(\td q\inv)$ for every $\td q \in Q(C)$.

Let  $U\subset V $ be the open neighborhood of $\1_V$ of the form $U=\exp_\liev \liev(\rho)$ for $0<\rho<\rho_0$
For any such $U$, we have $\Ad_{a_{*}}( U )\subset U$.  By $f$-invariance of $\mu$ and hypotheses \ref{hyp:11}  of \cref{propmain}, for $\mu$-a.e.\ $x$, we have the overflowing property:
$$  f(  W^{s,F}_{a_0,  U}(f\inv(x)) )\supset  f(  W^{s,F}_{a_0,  \Ad_{a_{*} }(U)}(f\inv(x)) )= 
W^{s,F}_{a_0,U}( x).  
$$

Given $\td q \in Q(C)$, we write  
\begin{equation}\label{eq:shifted} \wtd  W_{a_0, U,\td q}(x) =\td q\inv ( W^{s,F}_{a_0, U}(\td q \cdot x) ).\end{equation}
Given $v\in C$,  we have 
$\Ad_{a_{*}} (v)\in C\subset V(q)$
and thus $$\td q_{\Ad_{a_{*}}(v)} \in Q(C).$$
In particular, given $x_0\in \Omega$ and $v\in C$ we have a similar overflowing property:
\begin{align}
f(\wtd  W_{a_0, U, \td q_{\Ad_{a_{*}}(v)}}(f\inv(x(v,x_0))))
&=
f({\td q_{\Ad_{a_{*}}(v)} }\inv  \cdot W^{s,F}_{a_0, U}( \td q_{\Ad_{a_{*}}(v)} \cdot f\inv (x(v,x_0)) ))\notag
\\&=
\td q_{v} \inv \cdot W^{s,F}_{a_0, \Ad_{a_{*}}\inv (U)}( \td q_{v}  \cdot x(v,x_0)) \notag
\\&=
\wtd  W_{a_0, \Ad_{a_{*}}\inv (U), \td q_{v}}(f(x(v,x_0)))\notag
\\&\supset 
\wtd  W_{a_0, U, \td q_{v}}(x(v,x_0)). \label{overflow}
\end{align}

Recall that all fibers $M_x$ of  $M^\alpha$ are uniformly comparable to a fixed $M$.  
Let $r_0>0$ be the common injectivity radius of all $M_x$.

 Let 
$$\lambda = \min\{ \lambda_i^F( a_{*}): \lambda_i^F\in U\}>0.$$
and 
$$\lambda^+ = \max\{ |\lambda_i^F( a_{*})|: \lambda_i^F\in S\}>0.$$
 Fix $0<\epsilon<\frac 1  {100}$ sufficiently small.  Let $\{\psi_x\colon \R^k\to M_x\}$ be a family of fiberwise $(\epsilon, \ell)$-Lyapunov charts for the fiberwise dynamics of $f$.   (See for example \cite[Prop.\ 5.1]{MR4599404}.)
Specifically, let $d=\dim M$ and write $\R^d = E_1 \oplus E_2$ where $\dim E_1=\dim E^s_{a_0}(x)$ and $\dim E_2=\dim E^u_{a_0}(x)$ for a.e.\ $x$.  
 Equip $\R^d$ with the norm $\|v_1+v_2\|=\max\{|v_1|, |v_2|\}$ where $v_1\in E_1$ and $v_2\in E_2$ and $|\cdot|$ is the Euclidean norm.  Let $\R^d(\rho)$ denote the open ball of radius $\rho$ centered at $0$. 
  Then there is a measurable function $\ell\colon (M^\alpha,\mu) \to [1,\infty)$ and a measurable family of charts 
   $\{\psi_x\colon \R^k(1)\to M_x\}$ 
 such  for almost every  $x$, that the following hold:
\begin{enumerate}[label=(\roman*)]
\item  $\psi_x(0) = x$ and $\psi_x$ is a $C^1$ diffeomorphism between  the unit ball $\R^d(1)$ and an open neighborhood of $x$ of radius at most $r_0$ in $M_x$;  
	\item \label{LC:2}$D_0\psi_x E_1= E_{a_0}^s (x)$ and $D_0\psi_x E_2= E_{a_0}^u (x)$;
	\item \label{LC:3} the map $\hat f_x \colon \R^d( e^{-\lambda^+-2\epsilon}) \to  \R^d( 1) $  
	 given by \begin{equation}\label{eqftilde} 
	 \hat f_x(v):= \psi_ {f(x)}\inv    \circ f \circ   \psi_x(v)= \psi_ {a_{*} \inv x}\inv    \left(a_{*}\inv \cdot  \psi_x(v)\right)
	 \end{equation} is well-defined;  
\item $D_0\hat f_x $ preserves $E_1$ and $E_2$;
	\item \label{LC:5}  for $v\in E_1$ and $w\in E_2$,  $$ e^{-\epsilon}\|v\|\le  \|D_0\hat f_x v\|\le e^{\lambda_+ + \epsilon}\|v\|;$$ 
	and 
	 $$ \|D_0\hat f_x w\|\le e^{-\lambda+\epsilon}\|w\|;$$ 
	\item\label{LC:6} $\Hol(D\hat f_x ) \le \epsilon $ hence $\Lip(\hat f_x - D_0\hat f_x) \le \epsilon$ where $\beta = \min \{r-1, 1\}$ is the \Holder-exponent of the derivative $v\mapsto D_v\hat f_x$;
	\item \label{LC:7} $\|\psi_x \|_{C^1}\le1 $, and $\|\psi_x \inv \|_{C^1}\le  \ell(x)$;
	\item $\ell(f(x))\le e^{\epsilon }\ell(x)$ and $\ell(f\inv (x))\le e^{\epsilon }\ell(x)$.
\end{enumerate}

 Set \begin{equation}\label{graphtrans}\lambda_0= \lambda-4\epsilon.\end{equation}

We recall that neighborhoods of $\infty$ in $G/\Gamma$ have exponentially small mass.  Namely, there are $R_0>1$ and  $c>0$ such that for all $R\ge R_0$, 
\begin{equation}\label{eq:expcusp}m_{G/\Gamma}\bigl(\{g\Gamma: d(g\Gamma,\1\Gamma)>R\}\bigr)\le e^{-cR}.\end{equation}

Let $\calA\colon G\times \calF\to \calF$ denote the fiberwise derivative cocycle for the $G$-action on $M^\alpha$.  Relative to the fiberwise norms whose properties are outlined in \cref{sec:norms}, it follows from \cref{LMR} (c.f.\ \cite[Proposition 5.10]{2105.14541}) that there are $C_0\ge 1$ and $k>0$ such that 
$$\|\calA(g, x)\|^F_{M^\alpha} \le C_0 e^{k d_{G/\Gamma}(\pi(x),\1\Gamma) + k d_G(g,\1)}.$$
In particular, there are  $C_1\ge 0$ and $\gamma >0$ such that for all $x\in M^\alpha$ with  $d(\pi(x),\1\Gamma)\le R$ and $n\ge 0$, 
\begin{equation}\label{eq:MVT}
\|D_xf^n \|^F_{M^\alpha} \le   e^{C_1R + \gamma n}.
\end{equation}
Replacing $C_1$ with $\max\{C_1, 1\}$ if needed, we assume $C_1\ge 1$.  
%

Recall we write $A'_{\beta_0}$ for the kernel of $\beta_0$ in $A$.  Then $f\in A'_\beta$.  As  $\mu$ need not be $A'_{\beta_0}$-ergodic, it also need not be $f$-ergodic.
Let $\calE=\calE_{\beta_0}$ denote the partition of $(M^\alpha, \mu)$ into $A'_\beta$-ergodic components.  
As we assume that $\Gamma$ is irreducible in $G$, $\mu$-almost every ergodic component of $\mu$ projects to the Haar measure on $G/\Gamma$ under the map $\pi\colon M^\alpha\to G/\Gamma$.  
By entropy considerations in \cref{lemma:entropy} and \cref{lem:entropysavesmyassonceagin},  we have coincidence of leafwise measures $\mu^{{[\beta]},\calE }_x=\mu^{{[\beta]} }_x$ for all $\beta\in \Sigma_Q^\perp$.  
We then have coincidence of leafwise measures $\mu^{s,\calE }_x=\mu^{s }_x$  along $W^s_{a_0}$-manifolds.  
It follows for $\mu$-almost every $x_0$     
and almost every $v\in C$ that
\begin{equation}\label{ergcomps}
\text{$x(v,x_0)\in \calE(x_0)$,\ \  and \ \ 
 $\bar x(v,x_0)\in \calE(q\cdot x_0).$}
\end{equation}

For the remainder, we fix $x_0$ from a full $\mu$-measure set and fix $A'_\beta$-ergodic components $\mu^{\calE}_{x_0}$ and $\mu^{\calE}_{q\cdot x_0}$  containing $x_0$ and $q\cdot x_0$, respectively.  We moreover assume $\mu^{\calE}_{x_0}$ and $\mu^{\calE}_{q\cdot x_0}$ project to the Haar measure on $G/\Gamma$.  We also assume $x_0$ is chosen so that \eqref{ergcomps} holds for a.e.\ $v\in V(q)$.
Moreover, by \cite[Thm.\ 1]{MR0283174}, in addition to the properties above we may assume $a_*$ (and thus $f$) is chosen to be ergodic with respect to both $\mu^{\calE}_{x_0}$ and $\mu^{\calE}_{q\cdot x_0}$.

We now choose a number of parameters and sets determined by these parameters.
\begin{enumerate}
\item Fix $\delta>0$ sufficiently small so that   
\begin{equation}\label{delts} \bigl(2\delta \lambda_0   +\delta\gamma       -\lambda_0 (1-\delta)\bigr)<0.\end{equation}
where $\lambda_0$ is as in \eqref{graphtrans} and $\gamma$ is as in \eqref{eq:MVT}.

\item Fix $R\ge R_0\ge 1$ such that 
$4R e^{-cR}<\delta.$
\item Fix $\ell_0>0$ such that  (using \eqref{eq:expcusp} and that $\mu'$ projects to the Haar measure)
$$\mu'\bigl(\{x\in M^\alpha: d(\pi(x), \1\Gamma)<R \text{ and } \ell(x)\le \ell_0\}\bigr)\ge 1 -  1.1e^{-cR}$$
for $\mu'= \mu^{\calE}_{x_0}$ and $\mu'= \mu^{\calE}_{q\cdot x_0}$. 
\item Fix $M_0\in \N$ with $$M_0\ge \max \left\{\frac{\log \ell_0}{\lambda_0}, \frac{C_1}{\lambda_0}\right\}. $$

\end{enumerate}

 Consider $x\in M^\alpha$ and $\td q\in Q(C)$ where $Q(C)$ is as in \eqref{eq:QC}. 
For $x\in \Omega$, any  subset $E\subset V$, and $\td q\in Q(c)$ such that  $x, \td q\cdot x\in \Omega$
\begin{align*}
	W_{x,E}&:= \psi_x\inv  \bigl(W^{s,F}_{a_0,E}(x)\bigr) \cap \R^d(\kappa)\\
\what W_{x,E,\td q}&:= 
\psi_x\inv  \bigl( \wtd  W_{a_0,E,\td q}(x)\bigr) \cap \R^d(\kappa)=
\psi_x\inv  \bigl(\td q \inv \cdot W^{s,F}_{a_0,E}(\td q \cdot x)\bigr) \cap \R^d(\kappa).
\end{align*}     

Fix additional parameters $0<\rho<\rho_0$ and $0<\hat \kappa<\kappa$ with 
$$\kappa< e^{-M_0 (\lambda_+ +10\epsilon+ \log 2)}$$ and
$$0<\hat \kappa < e^{-M_0 (\lambda_+ +10\epsilon+ \log 2)} \kappa$$   to be specified below.  Also fix  $0<\eta<1$.  
Write 
$$U=\exp_\liev \liev(\rho),\quad \quad U_0=\exp_\liev \liev(\rho_0).$$

  Let $$G=G(\ell_0,R, \delta, \hat \kappa, \kappa, \rho_0, \rho)\subset \Omega\subset M^\alpha$$ denote the subset where the following hold for $x\in G$:
\begin{enumerate}[resume]
	\item $d(\pi(x), \1 \Gamma)\le R$.

		\item  $\ell(x)\le \ell_0$.  
\item $W_{x,U_0}$ contains the graph of a $C^1$ function $h_x \colon E_1(\kappa)\to E_2$ with\begin{enumerate}
\item  $h_x(0)= 0 $;
\item  $D_0h_x= 0 $;
\item  $\|D_vh_x\|\le \frac 1 3$ for all $v\in E_1(\kappa).$
\end{enumerate}

\end{enumerate}
Taking $0<\hat \kappa$ sufficiently small, for $x\in G$ we assume 
\begin{enumerate}[resume]
\item $W_{x,U}$ is contained in the graph of  a $C^1$ function $ g_x\colon E_1( \hat \kappa)\to E_2$  with \begin{enumerate}
\item  $g_x(0)=0$;
\item  $D_0g_x= 0 $;
\item  $\|D_vg_x\|\le \eta$ for all $v\in E_1(\hat  \kappa).$
\end{enumerate}
\end{enumerate}

Note that the $C^1$ norm of the translation $\td q$ from $M_x$ to $M_{\td q\cdot x}$ is uniformly bounded over the choice of 
$x\in M^\alpha$ with $d(\pi(x), \1\Gamma)\le R$ and $\td q\in Q(C)$.  
In particular, having taken $\eta>0$ and $\hat \kappa>0$  sufficiently small above, for all $x\in G$ and $\td q\in Q(C)$ such that $\td q\cdot x\in G$, we assume that 
\begin{enumerate}[resume]
\item $\what  W_{x,U,\td q}$ is contained in the graph of a  $C^1$ function $$\hat h_x\colon E_1(e^{-M_0 (\lambda_+ +10\epsilon+ \log 2)}\kappa)\to E_2,$$ with \begin{enumerate}
\item  $\hat h_x(0)=0$;
\item  $D_0\hat h_x= 0 $;
\item  $\|D_v\hat h_x\|\le \frac 1 3 $ for all $v\in E_1(e^{-M_0 (\lambda_+ +10\epsilon+ \log 2)}\kappa) .$
\end{enumerate}
\end{enumerate}

We may use a bump function to extend the maps $\hat f_x\colon \R^k(e^{-\lambda^+-2\epsilon})\to \R^k(1)$ to diffeomorphisms $\widehat F_x\colon \R^k\to \R^k$ with $\hat f_x(v) = \widehat F_x(v)$ for all $v$ with $\|v\|\le \frac 1 2 e^{-\lambda^+-2\epsilon}$ and (reducing $\epsilon>0$ above if needed) $\Lip(\widehat F_x- D_0 \hat f_x)<\epsilon$.
Given $n\ge 0$, let  $\widehat F_x^{(n)}\colon \R^k\to \R^k$ be $$\widehat F_{f^{n-1}(x)}\circ \dots \circ \widehat F_{x}.$$
Using the iterated graph transform for 1-Lipschitz functions in \cref{lem:iteratedgraphtrans},  if $x,\bar x\in G$ and if $f^n(x), f^n(\bar x)\in G$ for some $n\le M_0$ then 
the following hold:
\begin{enumerate}
\item $W_{x,U_0}$ contains the graph of a 1-Lipschitz function  $ h_{U_0}\colon E_1(\kappa)\to E_2$  with  
$h_{U_0}(0) = 0$.
\item $W_{f^n(x),U_0}$ contains the graph of a 1-Lipschitz function  $\hat  h^{(n)}_{U_0}\colon E_1(\kappa)\to E_2$  with  
$\hat h^{(n)}_{U_0}(0) = 0$.
\item $\what W_{x,U,\td q}$ is contained in the graph of  a 1-Lipschitz function $ h_U\colon  W_0\to E_2$  with  $h_U(0)=0$ where $W_0\subset E_1(e^{-M_0 (\lambda_+ +10\epsilon+ \log 2)}\kappa) $ is an open neighborhood of $0$.
\item $\widehat F_x^{(n)} (\what W_{x,U,\td q})$ is the graph of a 1-Lipschitz function $\hat h_U^{(n)}\colon W_n\to E_2$ where $W_n\subset E_1(\kappa)$ is an open neighborhood of $0$.  Moreover, \begin{equation}\label{eq:contpre}\|\hat h^{(n)}_{U_0}-\hat h_U^{(n)}\|_{C^0,W_n}\le e^{n\lambda_0}\|  h_{U_0}-  h_U\|_{C^0,W_0}.\end{equation}
\end{enumerate}
Moreover, let $0< n_1<n_2<\dots<n_\ell$ be subsequent times such that $f^{n_j}(x), f^{n_j}(\bar x)\in G$.  If  $n_{j}-n_{j-1}\le M_0$ for every $1\le j\le \ell$, then by hypothesis, the overflowing property \eqref{overflow},  and by recursively applying \eqref{eq:contpre}, 
 for every $1\le j\le \ell$ we obtain  that 
$W_{n_j}\subset  E_1(e^{-M_0 (\lambda_+ +10\epsilon+ \log 2)}\kappa) $ 
and 
\begin{enumerate}[resume]
\item $\widehat F_x^{(n_j)} (\what W_{x,U,\td q})$ is the graph of a 1-Lipschitz function $\hat h_U^{(n_j)}\colon W_{n_j}\to E_2$ where $W_{n_{j-1}}\subset E_1(\kappa)$ is an open neighborhood of $0$.  Moreover, 
\begin{align}\label{eq:cont}\|\hat h^{(n_\ell)}_{U_0}-\hat h_U^{(n_\ell)}\|_{C^0,W_{n_\ell}}&\le e^{(n_\ell-n_{\ell-1})\lambda_0}\|  h^{(n_{\ell-1})}_{U_0}-  h^{(n_{\ell-1})}_U\|_{C^0,W_{n_{\ell-1}}}\notag \\& \le e^{n_\ell \lambda_0}\|  h_{U_0}-  h_U\|_{C^0,W_0}.\end{align}
\end{enumerate}

We may now take $0<\rho<\rho_0$ as well as $\kappa, \hat \kappa, \eta$ above sufficiently small so that $$\mu'(G(\ell_0,R, \hat \kappa,\kappa, \delta, \rho_0, \rho)) \ge 1 - 1.5 e^{-cR}$$
for both $\mu'= \mu^{\calE}_{x_0}$ and $\mu'= \mu^{\calE}_{q\cdot x_0}$.

Set $G= G(\ell_0,R, \hat \kappa,\kappa, \delta, \rho_0, \rho)$. 
By \eqref{ergcomps}, 
for almost every 
$v\in C$,
there is $n_0(x(v,x_0))>0$ such that for $n\ge n_0(x(v,x_0))$
\begin{align*}
\frac {1}n \# \{0\le j\le n-1: f^{-j}(x(v,x_0)) \notin G\} & \le 2 e^{-c R},\\
\frac {1}n \# \{0\le j\le n-1: f^{-j}(\bar x(v,x_0)) \notin G\} &\le  2 e^{-c R}.
\end{align*}

Fix such $x= x(v,x_0)$ and $\bar x= \bar x(v,x_0)$.  Fix $n\ge 1$  such that $f^{-n}(x),f^{-n}(\bar x)\in G$.  
Let $J_1,\dots, J_k\subset \{0,\dots, n\}$ be the maximal intervals with the following property:
if $$J_j= \{a_j, \dots, b_j\}$$ then 
\begin{enumerate}
\item $f^{-a_j}(x)\in G$, $f^{-b_j}(x)\in G$, $f^{-a_j}(\bar x)\in G$, and $f^{-b_j}(\bar x)\in G$;
\item for all $a_j<k<b_j$, either $f^{-k}( x)\notin  G$ or $f^{-k}(\bar x)\notin  G$;
\item $b_j-a_j\ge M_0$.
\end{enumerate}

That is, $J_j$ is the $j$th ``long simultaneous excursion'' from $G$.   By definition, the intervals $J_j$ are disjoint.  
Let $$b_0 = \min\{n: f^{-n}(x)\in G \text{ and } f^{-n}(\bar x)\in G \} .$$
We have $$0\le b_0 \le a_1<b_1\le a_2<b_2\le a_3<b_3\dots\le a_{k}<b_{k}\le n=:a_{k+1}.$$  
%
Assuming $n\ge n_0(x)$, we have $$\frac {M_0 k}{n} \le \frac 1 n \sum _{i=1}^k b_j-a_j \le 4 e^{-c R} \le  \delta.$$

Let $\td q_n:=  \td q_{\Ad_{a_{*}}^{n} (v)}=
\Ad_{a_*}^n(\td q_v)
\in Q(C)$.
Let $$D_n(x) = \sup \left\{d_{M_{f^{-n}(x)}}\left(
   W^{s,F}_{a_0, U_0} (f^{-n}(x)),  y)
 \right)  :   y\in \wtd  W_{a_0, U,\td q_n }(f^{-n}(x))\right\} .$$

Using \eqref{eq:cont} over subsequent ``short'' simultaneous returns to $G$, the overflowing property \eqref{overflow}, and Lipschitz properties of  the charts $\psi_x$, for $0\le j\le k$ we have that 
\begin{align*} 
D_{b_j}(x) 
 \le \ell_0 e^{-\lambda_0 (a_{j+1} -b_j)} 
D_{a_{j+1}}(x).
\end{align*} 
Using  \eqref{eq:MVT}, for $1\le j\le k$ we have that    
\begin{align*} 
D_{a_j}(x)  \le  e^{C_1R+ \gamma (b_j-a_j)} D_{b_j}(x) .
\end{align*}

 Then 
\begin{align*} 
D_{0}(x) 
&\le 
e^{C_1R+ \gamma b_0}
 \prod _{j=1}^k\left(e^{C_1 R} e^{\gamma(b_j-a_j )}\right) \prod _{j=0}^k \left( \ell_0 e^{-\lambda_0 (a_{j+1} -b_j)}\right)
D_n(x) %
 \\
&\le 
e^{C_1R+ \gamma b_0}
\left(e^{kC_1 R} e^{\gamma \sum _{j=1}^k(b_j-a_j )}\right)\left( \ell_0^{k+1} e^{-\lambda_0 \sum _{j=0}^k (a_{j+1} -b_j)}\right)
D_n(x) %
 \\
 &\le e^{\log \ell_0+C_1R+ \gamma b_0} e^{k \log\ell_0  + kC_1 R +4 e^{-c R}\gamma n      -\lambda_0 n(1-4 e^{-c R})} D_n(x) 
 \\
 &\le e^{\log \ell_0+C_1R+ \gamma b_0} e^{4n e^{-c R} M_0\inv  \log\ell_0  + 4n C_1 R e^{-c R} M_0\inv  +4 e^{-c R}\gamma n      -\lambda_0 n(1-4 e^{-c R})} D_n(x) 
 \\
 &\le e^{\log \ell_0+C_1R+ \gamma b_0}e^{n \delta M_0\inv  \log\ell_0  + n\delta C_1 M_0\inv  +\delta\gamma n      -\lambda_0 n(1-\delta)} D_n(x) 
 \\
 &\le e^{\log \ell_0+C_1R+ \gamma b_0} e^{n \bigl(\delta \lambda_0  + \delta  \lambda_0 +\delta\gamma       -\lambda_0 (1-\delta)\bigr)} D_n(x).  
\end{align*}

Recall the choice of $\delta$ satisfying \eqref{delts}.  
By Poincar\'e recurrence,  the sequence $\{D_n(x)\}$ is uniformly bounded on an infinite set of $n$ for $x$ of the form $x= x(v,x_0)$ for almost every $v\in C$.  
Thus, we conclude for $\mu$-a.e.\ $x_0$ and $x$ of the form $x= x(v,x_0)=[vg_0,\phi_{(g_0,z_0)}(v)]$ for a.e.\ $v\in C$, that 
 $D_0(x)  =0$. 
In particular,  $\wtd  W^{s,F}_{a_0, U,\td q_v} (x)\subset   W^{s,F}_{a_0, U_0} (x)  \subset W^{s,F}_{a_0} (x)$ for all such $x$.  
This finishes the proof of \cref{lem:hardestpart}. \hfill \qedsymbol
\label{endpfooequiv}

\subsection{Equivalence relations;  properties used in the proof of Theorem \ref{propmain}}
We define equivalence relations on a full $\wtd \mu$-measure set of points in $G\times M$;  geometric properties of the equivalence classes yield the objects asserted to exist in  \cref{propmain}. 
\subsubsection{Coherent extendibility of graphs}
We begin by defining various full-measure subsets of points $(g,x)$ where the functions $\phi_{(g,x)}$ are well defined and for which the graphs extend in a coherent way over the intersection of  their ranges.  
\begin{definition} \ 
\begin{enumerate}
\item  We say a point $(g,x)\in G\times M$ is \emph{$C^r$-regular} if there is a $C^r$-diffeomorphism $$\phi_{(g,x)}\colon V\to W^{s,F}_{a_0,g}(x)\subset M$$ with $$\phi_{(g,x)}(\id) = x$$
such that the 
leafwise measure $\wtd \mu_{(g,x)}^s$ is the image of the Haar measure on $V$ 
under the graph  $$v\mapsto (vg, \phi_{(g,x)} (v)).$$


\item We say  a point $(g,x)\in G\times M$ is \emph{$V$-coherent} if $(g,x)$ is $C^r$-regular and for   Haar-almost every  $v\in V$ the following hold: 
\begin{enumerate} \item writing $y =\phi_{(g,x)}(v)$, we have that  
$(vg, y)$ is $C^r$-regular and   
\item   for every $\td v\in V$, 
$$\phi_{(vg,y)} (\td v) =  \phi_{(g,x)}( \td v \cdot v).$$
\end{enumerate}
\item We say  a point $(g,x)\in G\times M$ is \emph{$Q$-coherent} if $(g,x)$ is $C^r$-regular and for  Haar-almost every  $q\in Q$  
 the following hold:
\begin{enumerate}
\item $(qg,x) $ is $C^r$-regular and
\item for every $\td v\in V(q) $, 
$$\phi_{(qg,x)} (\tau_q( \td v)) =   \phi_{(g,x)}( \td v).$$
\end{enumerate}
\item We say  a point $(g,x)\in G\times M$ is \emph{$\Gamma$-coherent} if $(g,x)$ is $C^r$-regular and for all $\gamma\in \Gamma$, $$\phi_{(  g \gamma , \alpha (\gamma\inv )(x))} = \alpha(\gamma)\inv \circ \phi_{( g, x)}.$$  

\end{enumerate}

\end{definition}

By assumption, the set of 	$C^r$-regular points is conull. 
By $\Gamma$-invariance of the lifted measure $\wtd \mu$, the set of $\Gamma$-coherent points is also conull.  
 By essential uniqueness of leafwise measures $\wtd \mu_{(g,x)}$ (up to normalization), the set of $V$-coherent points is conull.
By \cref{prop:equivarstablemanifolds}\eqref{vomitface3}, the set of $Q$-coherent points is also conull.

We now define points for which we can coherently extend the functions $\phi_{\bullet}$ countably many times.  
\begin{definition}\
\begin{enumerate}
\item 
Let  $\scrG_0$ be the subset of points $(g,x)\in G\times M$ that are $C^r$-regular, $V$-coherent,  $Q$-coherent, and $\Gamma$-coherent.  

\item Given $j\ge 0$, let  $\scrG_{j+1}$ be the subset of points $(g,x)\in G\times M$ such that 
\begin{enumerate}
	\item $(g,x) \in \scrG_j$;
	\item for    Haar-almost every  $v\in V$,  $(vg, \phi_{(g,x)}(v))\in \scrG_j$;
	\item for  Haar-almost every  $q\in Q$, $(qg,x) \in \scrG_j$. 
\end{enumerate}
\item Write $ \scrG_\infty =\bigcap_{j=0}^\infty \scrG_j$. 
\end{enumerate}
\end{definition}
Note that as each $\scrG_j$ is conull, $ \scrG_\infty$ is a conull subset of $G\times M$.  
%

We collect a number of useful   observations regarding the above definitions.  
\begin{claim}\ \label{claim:welldef}
\begin{enumerate}
\item \label{clam:111}
	If $(g,x)\in \scrG_\infty$ and $(qg,x)\in \scrG_\infty$, then for all $v\in V(q)$, 
$$\phi_{(qg,x)} (\tau_q(  v)) =   \phi_{(g,x)}( v).$$

\item \label{clam:222} If $(g,x)\in \scrG_\infty$ and $(\td vg, \phi_{(g,x)}(\td v))\in \scrG_\infty$, then for all $v\in V$
$$\phi_{(\td v g,\phi_{(g,x)}( \td v))} ( v) =   \phi_{(g,x)}( v\td v).$$

\end{enumerate}
\end{claim}
\begin{proof}
We show \eqref{clam:111} as \eqref{clam:222} is similar.     Consider $(g,x)\in \scrG_\infty$ and $(qg,x)\in \scrG_\infty$. 
Given $q'\in Q$, write $\td q= q'q\inv$ and observe the assignment $q'\mapsto \td q$ is absolutely continuous.  
 For a given $v\in V(\td q)$, there is a full measure subset of $q' \in Q$ in a sufficiently small neighborhood of $q$ for which
$$\phi_{(q'g,x)} (\tau_{q'}( v)) =   \phi_{(g,x)}( v)$$
and 
$$\phi_{(\td q qg,x)} (\tau_{\td q}(\tau_{ q}( v)) =   \phi_{(qg,x)}(  \tau_{ q}(v))$$
The claim follows since $\tau_{\td q }(\tau_ { q} (v) = \tau_{\td q q} (v)) = \tau_{q'}   (v). $
\end{proof}

We similarly show that  the following.  
\begin{claim} \label{clam:same}
Suppose $(g,x)\in \scrG_\infty$ and $(qg,y)\in \scrG_\infty$ for some $q\in Q$.  If 
$$\phi_{(qg,y)} (\tau_q(  v)) =   \phi_{(g,x)}( v)$$
for some $v\in V(q)$  then $x=y$.
\end{claim}
\begin{proof}
Write $z= \phi_{(qg,y)} (\tau_q(  v)) =   \phi_{(g,x)}( v)$.

 We may select from a full measure subset of $\hat q\in Q$ in a sufficiently small neighborhood of   $q$ for which 
 \begin{enumerate}
 	\item $(\hat q g, y) \in \scrG_\infty$, 
 	\item $(\hat q g, x )\in \scrG_\infty$,
 	\item $z\in \mathrm{Im} (\phi_ {(\hat qg, x ) })$, 
 	\item $z\in \mathrm{Im} (\phi_ {(\hat qg, y )} )$, and 
	\item $\tau_{q}(v)\in V({\hat q q\inv })$. 
 \end{enumerate}
Write $v^*= \tau_{\hat q}( v)=\tau_{\hat q q\inv }(\tau_{q}(v))$.  We have $(v^*\hat q g,z) \in 
W^s _{a_0}( \hat q g, x)\cap W^s _{a_0}( \hat q g, y)$
whence $\wtd \mu_{(\hat q g,x)}^s= \wtd \mu_{(\hat q g,y)}^s$.
Since both $(\hat qg, x ) \in \scrG_\infty$ and $(\hat qg, y ) \in \scrG_\infty$, by the assumption on the support of of the leafwise measures,  it follows that $x=y$.
\end{proof}

\subsubsection{Equivalence relations; closed equivalence classes as covering spaces}
\begin{definition}
Let $\sim$ and $\sim^\circ $ denote the equivalence relations on $\scrG_\infty$ generated by the following: 
\begin{enumerate}
\item\label{eqrn:1}  for $q\in Q$, $(g,x) \sim (qg,x)$ if both $(g,x)$ and  $(qg,x)$ are in $ \scrG_\infty$
\item\label{eqrn:15} for $q\in Q^\circ $,  $(g,x) \sim^\circ (qg,x)$ if both $(g,x)$ and  $(qg,x)$ are in $ \scrG_\infty$
\item\label{eqrn:2}  for $v\in V$, $(g,x)\sim (vg, \phi_{(g,x)}(v))$ and 
$(g,x)\sim^\circ (vg, \phi_{(g,x)}(v))$ if both $(g,x)$ and  $(vg, \phi_{(g,x)}(v))$ are in $ \scrG_\infty$.
\end{enumerate}
\end{definition}

\begin{definition}
Given $(g,x)\in \scrG_\infty$, let 
$$\overline \scrG(g,x):=\bigcup_{( g', x')\sim^\circ (g,x)} \left\{\bigl( qv'g', \phi_{( g', x')}(v')\bigr):  q\in Q^\circ,  v'\in V\right\}.$$
\end{definition}
We  claim the sets $\overline \scrG$ partition a full $\td \mu$-measure subset of $G\times M$.  
\begin{claim}\label{claim:uniqlo}
	For $(g,x), (g',x')\in \scrG_\infty$, either 
	$\overline \scrG(g,x)\cap \overline \scrG(g',x')=\emptyset$ 
	or 	$\overline \scrG(g,x)= \overline \scrG(g',x')$. 
\end{claim}
\begin{proof}
Consider $(g,x), (g',x')\in \scrG_\infty$ for which  $\overline \scrG(g,x)\cap \overline \scrG(g',x')\neq\emptyset$.  Moving within  $\sim ^\circ$-equivalence classes, we may assume 
$(g,x), (g',x')\in \scrG_\infty$ are such that 
$$(qvg, \phi_{( g, x)}(v))=  (q'v'g', \phi_{( g', x')}(v'))$$
for some $q,q'\in Q^\circ$ and $v,v'\in V$.  We finish the proof of the  claim by showing $(g,x)\sim ^\circ (g',x') $.

Taking $\hat v$ from a full measure subset of a sufficiently small neighborhood of the identity in $V$, we may assume
\begin{enumerate}
	\item $\hat v\in V({{q'}\inv q})$   
	\item $(\hat vvg, \phi_{( g, x)}(\hat vv))\in \scrG_\infty$
		\item $(\tau_{{q'}\inv q}(\hat v)v'g', \phi_{( g', x')}(\tau_{{q'}\inv q}(\hat v)v'))\in \scrG_\infty$
\end{enumerate}

Write $$y _ 0=\phi_{( g, x)}(\hat vv), \quad \quad  y_1 = \phi_{( g', x')}(\tau_{{q'}\inv q}(\hat v)v')$$
and $$\bar q \hat v  = \tau_{{q'}\inv q}(\hat v) {q'}\inv q$$
some $\bar q\in Q$.   Since ${q'}\inv q\in Q^\circ$, we have that $\bar q\in Q^\circ.  $  Also we have 
 $$\bar q \hat v vg  = \tau_{{q'}\inv q}(\hat v){q'}\inv q  v g= 
 \tau_{{q'}\inv q}(\hat v) v' g'$$ 
 and by \cref{claim:usefulidiot}, $(\tau_{{q'}\inv q}(\hat v) )\inv =\tau_{\bar q}(\hat v\inv)$.

Then 
\begin{align*}
\phi_{( \hat vv g, y_0)}(\hat v\inv) &= \phi_{( g, x)}(v) =
\phi_{( g', x')}(v')
\\& 
=\phi_{( \tau_{{q'}\inv q}(\hat v)v' g', y_1)}((\tau_{{q'}\inv q}(\hat v))\inv)
\\&= \phi_{( \tau_{{q'}\inv q}(\hat v){q'}\inv q v g, y_1)}
(\tau_{\bar q}(\hat v\inv))
\\&= \phi_{( \bar q \hat v  v g, y_1)}(\tau_{\bar q}(\hat v\inv))
\end{align*}
By  \cref{clam:same}, we conclude that $y_0=y_1$.  
%
%
%
%
Then 
\begin{align*}
(g',x')&\sim^\circ  (\tau_{{q'}\inv q}(\hat v) v' g', y_1)
\\&= ( \tau_{{q'}\inv q}(\hat v) {q'}\inv q  v g, y_1)
\\&= ( \bar q \hat v  v g, y_1)
\\& = ( \bar q \hat v  v g, y_0)
\\&\sim ^\circ (\hat v  v g, y_0)
	\\& \sim ^\circ (g,x).  \qedhere
\end{align*} 

\end{proof}


Given $(g,x)\in \scrG_\infty$, consider  $(\td g, y)\in 	\overline \scrG(g,x)\sm \scrG_\infty$.  Write $$(\td g, y ) = \bigl( qv'g', \phi_{( g', x')}(v')\bigr)$$ for some choice of  $( g', x')\sim^\circ  (g,x), q\in Q^\circ , $ and $ v'\in V$.  
Write $\tau_{q}(v') = v^*$ and $v^* q = \hat q v'$ for some $\hat q\in Q^\circ $. Given $\hat v \in V$ sufficiently small, define 
$$\phi_{(\td g, y)}(\hat v) := \phi_{(g', x')}(\tau_{q\inv}(\hat v v^*)).$$

By \cref{claim:welldef},  the above definition of  $\phi_{(\td g, y)}(\hat v)$ is  well defined, independent of the choice of $(g',x')$.  
By taking $(g',x')$ sufficiently close to $(\td g, y)$ (so that $q$ and $v'$ are sufficiently small), we may define $\phi_{(\td g, y)}(\hat v)$ for all $\hat v\in V$.
With the above definition,  for all $(\td g, y)\in \overline \scrG(g,x)$ we have the following:
\begin{enumerate}
\item $\phi_{(\td g, y)}\colon V\to M$ is an injective, $C^r$ function with $\phi_{(\td g, y)} (\1_V) = y$.
\item  
$ \overline \scrG(g,x)\cap \bigl((V\td g)\times M\bigr) $
coincides with the graph of $\phi_{(\td g, y)}$, $$v\mapsto (v\td g, \phi_{(\td g, y)}(v))$$
\item   for every $\td v,v\in V$, 
$$\phi_{(v\td g,\phi_{(\td g, y)}(v))} (\td v) =  \phi_{(\td g,y)}( \td v \cdot v).$$
\item  for  every  $q\in Q$ and  $\td v\in V(q) $, 
$$\phi_{(q\td g,y)} (\tau_q( \td v)) =   \phi_{(\td g,y)}( \td v).$$
\item  for every $\gamma\in \Gamma$, $$\phi_{(\td g \gamma , \alpha (\gamma\inv )(y))} = \alpha(\gamma)\inv \circ \phi_{(\td g, y)}.$$  
\end{enumerate}

%

For $(g,x)\in \calG_\infty$, the set  $\bar\scrG(g,x)$ is path connected.  
Moreover, for $(g,x)\in \calG_\infty$, the natural map $\bar\scrG(g,x)\to G$, $$(\td g, y)\mapsto \td g,$$ is a covering map.  
Since $G$ was assumed topologically simply connected, we conclude that $\bar\scrG(g,x)\to G$, $(\td g, y)\mapsto \td g$, is a bijection.  
By Journ\'e's theorem (applied to the foliations by $Q$- and $V$-orbits), the set  $\bar\scrG(g,x)$ is an embedded $C^r$ manifold.  In particular, $\bar\scrG(g,x)\to G$ is a  $C^r$-diffeomorphism.  

%

\begin{claim}\label{claim:intersect}
For $(g,x)\in \scrG_\infty$ and $(\td g, \td  y)\in \scrG_\infty$, if $\overline \scrG(\td g, \td y)\cap \wtd  W^s_{a_0}(g,x)\neq \emptyset$ then 
$$\overline \scrG(\td g,\td y)= \overline \scrG(g,x).$$
\end{claim}
\begin{proof}
Let $(vg,  y)\in \overline \scrG(\td g, \td  y)\cap \wtd  W^s_{a_0}(g,x)$.  For every $v'\in V$ we have $$(v'vg, \phi_{vg, y}(v') )\in \overline \scrG(\td g,  \td y)\cap W^s_{a_0}(g,x).$$
By replacing $(vg,  y) $ with $(v'vg, \phi_{vg, y}(v'))$
 for a full measure set of choices of $v'\in V$, we may assume there is $q'\in Q^\circ$  so that 
\begin{enumerate}
	\item $(vg, \phi_{g,x}(v))\in \scrG_\infty.$
	\item $(q'vg, \phi_{g,x}(v))\in \scrG_\infty.$ 
	\item $(q'vg, y)\in \scrG_\infty.$ 
\end{enumerate}
By uniquness of the leafwise measure $\wtd\mu^s_{(q'vg, y )}$, we conclude that $\phi_{g,x}(v)=y$ whence $(vg, y)\in \overline \scrG(\td g, \td  y)\cap \overline \scrG(g,x)$ and the claim follows from \cref{claim:uniqlo}.
\end{proof}

\subsubsection{The set $\wtd \Lambda$; open charts for lamina}\label{sss:tdlam}
Write $${\wtd \Lambda}
= \bigcup_{(g,x)\in \scrG_\infty} \bar \scrG(g,x).$$
For $(g,x)\in  \wtd \Lambda$,
write $$V_{(g,x)}= \mathrm{Im}(\phi_{(g,x)}) :=\phi_{(g,x)}(V)\subset M $$ 
For $(g,x)\in  \wtd \Lambda$, we have  that 
$\phi_{( g,x)}\colon V\to V_{(g,x)}\subset M$ is an injection.  
We first show for $(g_1,y_1), (g_2,y_2) \in  \wtd \Lambda $ that 
if  $V_{( g_1,y_1)}\cap V_{( g_2,y_2)}\neq\emptyset$  then $(g_1,y_1)\sim  (g_2,y_2)$.

\begin{claim}\label{claim:related1}
Consider  $(g_1,y_1) \in  \wtd \Lambda$, $(g_2,y_2) \in  \wtd \Lambda $, and  $z\in V_{( g_1,y_1)}\cap V_{( g_2,y_2)} $. Write $$z= \phi_{(g_1,y_1)}(v_1)=  \phi_{(g_2,y_2)}(v_2)$$ 
for $v_1,v_2\in V$.  Then 
\begin{enumerate}
	\item $v_2 g_2= q v_1 g_1$ for some $q\in Q$, and
	\item $ \overline \scrG (g_1,y_1)= q_0 \overline \scrG  (g_2,y_2)$ for some $q_0\in Q.$  
\end{enumerate}

\end{claim}

\begin{proof}
First, by replacing $(g_1, y_1)$ with $(qg_1, \phi_{(qg_1, y_1)}(v))$ for  sufficiently small $q\in Q^\circ$ and $v\in V$, we may assume $(g_1, y_1)\in \scrG_\infty.$

Given any $\td q$ outside of a proper algebraic subvariety, we may write $$v_2g_2 g_1\inv=qv\td q$$ for a unique choice of $q\in Q$ and $v\in V$.
Take  $\td q\in Q^\circ $ from a full measure subset of a 
sufficiently small neighborhood of $\1_V$ so that 
\begin{enumerate}
	\item $v_1\in V(\td q)$,
	\item  $(\td q  g_1,y_1)\in  \scrG_\infty $ and
	\item 
 $v_2g_2=  \hat q  \hat v (\td q g_1) $ for some $ \hat q\in Q$ and $ \hat v\in V$.

\end{enumerate} 
%

Let $\td y_ 1 = \phi_{(\td qg_1,y_1)}(\hat v)$.  
We then obtain
\begin{align*}
	z&=
	\phi_{(g_1,y_1)}(v_1)
	\\&=
	\phi_{(\td qg_1,y_1)}(\tau_{\td q}(v_1))
	\\&=	\phi_{(\hat   v\td qg_1,\td y_ 1)}(\tau_{\td q}(v_1)\hat   v \inv).
\end{align*}

We have 
$$z\in \wtd W^{s,F}_{a_0, \hat v \td q g_1}(\td y_1)= \wtd  W^{s,F}_{a_0, \td q g_1}(y_1).$$
In particular,  
$$(\hat q\inv v_2 g_2,z) =(\hat v\td q g_1, z)\in \overline \scrG(\hat q\inv v_2 g_2,  z)\cap \wtd  W^s_{a_0}(\td q g_1,y_1)$$
From \cref{claim:intersect}, we conclude that 
\begin{align*}
\overline \scrG(\hat q\inv v_2 g_2,  z)
= \overline \scrG(\hat v\td q g_2, \td y_1)= \overline \scrG(\td q g_2, y_1).\end{align*}
We then  conclude $\td y_1= z$.  
Also, since $\hat q \overline \scrG(\hat q\inv v_2 g_2,  z) = \overline \scrG( v_2 g_2,  z)= \overline \scrG(  g_2,  y_2),$ the second conclusion follows.

As 
$$z= \phi_{(\hat   v\td qg_1,\td y_ 1)}(\tau_{\td q}(v_1)\hat   v \inv) 
= \phi_{(\hat   v\td qg_1,z)}(\tau_{\td q}(v_1)\hat   v \inv),
$$
 we conclude that $\tau_{\td q}(v_1) \hat   v\inv = \1_V  $
whence 
$$v_2 g_2 = \hat q \hat v\td q g_1 = \hat q \tau_{\td q}(v_1)\td q g_1 =  \bar q v_1 g_1$$
for some $\bar q\in Q$ and the first conclusion follows.
\end{proof}

\begin{claim}\label{claim:openint}
Given $(g_1,y_1), (g_2,y_2)\in  \wtd \Lambda$, either\begin{enumerate}
	\item $V_{( g_1,y_1)}\cap V_{( g_2,y_2)}=\emptyset$ or 
	\item  $ V_{( g_1,y_1)}\cap V_{( g_2,y_2)} $ is open in both $V_{( g_1,y_1)}$ and 
$ V_{( g_2,y_2)}.$
\end{enumerate}	
\end{claim}

	\begin{proof}
Fix $z\in V_{( g_1,y_1)}\cap V_{( g_2,y_2)} $, $v_1,v_2\in V$ with $z= \phi_{(g_1,y_1)}(v_1)=  \phi_{(g_2,y_2)}(v_2)$,  
and---by \cref{claim:related1}---write $v_2g_2= \hat q v_1 g_1$ for some $\hat q\in Q$. 
Then $\phi_{(v_1g_1, z)}$ and $\phi_{(\hat qv_1g_1, z)}\circ \tau_{\hat q}$
coincide on $V(\hat q)$ whence $\phi_{(v_1 g_1,  z)}(V(\hat  q))$ is an open neighborhood of $z$ in both $V_{( g_1,y_1)}$ and $V_{( g_2,y_2)} .$
\end{proof}

Observe that $\wtd \mu$ is supported on $\scrG_\infty\subset \wtd \Lambda$.
\begin{claim}
The partition of $\wtd \Lambda$ into sets of the form $\bar \scrG(g,x)$ is measurable.  The conditional measures associated to each element $\bar \scrG(g,x) $ of this partition correspond  to the Haar measure under the bijection  $\bar \scrG(g,x) \to G$.  
\end{claim}

Given $g\in G$, let $\Lambda_g = \wtd \Pi_2\bigl((\wtd \Lambda \cap( \{g\}\times M)\bigr) \subset M$.  
Then $\wtd \mu_g$ is supported on $\Lambda_g$.
Given $g'\in G$, let $H_g'\colon \Lambda_g\to \Lambda_{gg'}$ be the holonomy map defined as follows:  $y\in \Lambda_g$ if and only if  $(g,y)\in \wtd \Lambda$; take the unique $z\in M$ with $(gg',z)\in \bar \scrG(g,y)$ and define $H_{g'}(y) = z$.
\begin{claim}\label{claim:holoninv}For every $g'\in G$ and a.e.\ $g\in G$,
$$(H_{g'})_* \wtd \mu_g = \wtd \mu_{gg'}.$$
\end{claim}

\subsection{Proof of \cref{propmain}}
With the above, we give the proof of \cref{propmain}, by constructing $\Lambda$, $h$, and the left $G$-action.  We  then verify  they satisfy the properties of \cref{propmain}.
\subsubsection{The set $\Lambda$}
Recall we constructed $\wtd \Lambda\subset G\times M$ in \cref{sss:tdlam}.  Recall that if $(g,y) \in \wtd \Lambda$ then $(qg,y) \in \wtd \Lambda$ for all $q\in Q$.  
Moreover, from \cref{claim:related1}, if $(g,y) \in \wtd \Lambda$ and $(\wtd g,y) \in \wtd \Lambda$ then $\td g = q g$ for some $q\in Q$.

Let $\bar \Lambda \subset (G/Q)\times M$ be the subset $$\bar \Lambda= \{(gQ, y): (g\inv,y) \in \wtd \Lambda\}.$$
Write $\Pi_2\colon (G/Q)\times M\to M$ 
and  $\wtd \Pi_2\colon G\times M\to M$ for the canonical projections.  Let $\Lambda\subset M$ be the set, 
$$\Lambda= \wtd  \Pi_2 (\wtd \Lambda) =\Pi_2(\bar \Lambda)=  \{ y\in M: (g,y)\in \wtd \Lambda \text{ for some $g\in G$}\}.$$

Observe that  $\Pi_2$ restricts to a bijection between $\bar \Lambda$ and $\Lambda. $
\begin{proof}[Proof of \cref{propmain}\eqref{prop:main3}]
The measure $\bar \mu$ is supported on $\bar \Lambda$.  Since $\Pi_2$ restricts to a bijection between $\bar \Lambda$ and $\Lambda,$ \cref{propmain}\eqref{prop:main3} holds.  
\end{proof}

\subsubsection{The  left $G$-action}
Fix $y\in \Lambda$. Take $\bar g\inv\in G$ so that $(\bar g\inv, y) \in \wtd \Lambda$.  
There exists $x_0\in \Lambda$ with  $(\1_G,x_0)\in \bar\scrG(\bar g\inv,y)$.   
Given $g'\in G$, we have $((g'\bar g)\inv , z) \in  \bar\scrG(\1_G,x_0)= \bar\scrG(\bar g\inv,y)$ for a unique $z\in \Lambda$.  Declare that
$$\ell_{g'}(y) := g'\cdot y = \ell_{g'} (\ell_{\bar g}(x_0))  = \ell_{ g'\bar g}(x_0) =  g'\bar g\cdot x_0=z.$$

We check the above definition of $\ell_{g'}(y)$ is independent of choices made above.
Indeed, if $(\bar g\inv,y)\in \wtd \Lambda$ and $(\wtd g\inv, y)\in \wtd \Lambda$ 
then  by \cref{claim:related1}, $\wtd g\inv= q\bar g\inv$ for some $q\in Q.$
If $(\id, x_0)\in \overline \scrG(\bar g\inv,y) $
and $(\id, x_0')\in \overline \scrG(\wtd g\inv, y) $
then since 
$$ \bar \scrG(\wtd g\inv, y)  =\bar\scrG(q\bar  g\inv, y) =  \bar \scrG(\bar g\inv, y)$$
we have $x_0'= x_0$.  
Thus, given $g'\in G$, if $(( g'\bar  g)\inv, z )\in \overline \scrG(\1_G,x_0)$ then $$((g'\wtd g)\inv,z) = (q( g'\bar  g)\inv, z )\in \overline \scrG(\1_G,x_0)$$ whence the assignment $\ell_{g'}(y) = z$ is defined independently of the choice of  $(\bar g\inv,y)$ or   $(\wtd g\inv, y)$ in $\wtd \Lambda$.

\begin{claim}
$\ell_g$ defines a left $G$-action on $\Lambda$.  

\end{claim}
\begin{proof}
Take $y\in \Lambda$ and write $(\bar g\inv, y) \in \bar \scrG(\1_G,x_0)$ for some $x_0\in \Lambda$ and $\bar g\in G$.  
Fix $g_1, g_2 \in G$ and let $z,w\in \Lambda$ be such that 
$$(( g_1\bar g)\inv, z)= (\bar g\inv g_1\inv, z) \in \bar \scrG(\1_G,x_0)$$ and $$((g_2 g_1\bar g)\inv, w) =(\bar g\inv g_1\inv g_2\inv, w) \in \bar \scrG(\1_G,x_0).$$
Then $$\ell_{g_2 } (\ell_{g_1}(y) ) = \ell_{g_2 } (z) = w = \ell_{g_1  g_2 }(y)$$ and the claim follows. 
\end{proof}

\begin{proof}[Proof of \cref{propmain}\eqref{prop:main7}]
Fix  $x\in \Lambda$.   Fix any $g\in G$ and let  $x_0= \ell_{g\inv }(x)$ so that $\ell_{g}(x_0) = x$.
Let $x_1= \ell_{g \inv \gamma \inv}(\alpha(\gamma)(x))$ so that $\alpha(\gamma)(x) = \ell_{\gamma g}(x_1)$.  
Take $g_0\in G$ so that $(g_0\inv, x_0)\in \wtd \Lambda$.  Then $(g_0\inv g\inv ,x)\in \wtd \Lambda$ and $$(g_0\inv g\inv \gamma\inv ,\alpha(\gamma)(x))\in \wtd \Lambda.$$  It follows that  $(g_0\inv , x_1)\in \overline \scrG(g_0\inv g\inv \gamma\inv ,\alpha(\gamma)(x))$.

For $y\in L(x)$, write $y= \ell_{\td  g}(x_0)$.  Then $(g_0\inv \td g\inv, y)\in \wtd \Lambda$ whence 
\begin{align*}
(g_0\inv \td g\inv \gamma\inv, \alpha(\gamma)(y) )  \in \overline \scrG(g_0\inv g\inv \gamma\inv, \alpha(\gamma)(x))= \overline \scrG(g_0\inv, x_1)
\end{align*}
and $\alpha(\gamma)(y) = \ell_{ \gamma \td g}(x_1).$
\end{proof}
\begin{proof}[Proof of \cref{propmain}\eqref{prop:main8}]
\cref{propmain}\eqref{prop:main8} follows from \cref{claim:holoninv} as (by definition of the left action) $\ell_{g'} (\Lambda_{g\inv}) = H_{{g'}\inv}(\Lambda_{g\inv})=\Lambda_{g\inv{g'}\inv}$ and $(\ell_{g'})_*\wtd \mu_{g\inv} = (H_{{g'}\inv})_* \wtd\mu_{g\inv} = \wtd \mu_{g\inv{g'}\inv}=\wtd \mu_{{(g'g)}\inv}$.
\end{proof}


\subsubsection{The map $h$}
Given $y\in \Lambda$, take $\bar g\in G$ so that $(\bar g\inv, y ) \in \wtd \Lambda$.  We declare that 
\begin{equation} h(y) = \bar g Q.\end{equation}  
As $\wtd \Lambda$ is $\Gamma$-invariant, we have 
$$(\bar g\inv , y )\cdot \gamma\inv  = (\bar g\inv \gamma\inv, \alpha(\gamma)(y))$$
and thus conclude that $$h(\alpha(\gamma)(y)) = \gamma\cdot \bar g Q$$
for all $\gamma\in \Gamma$ and $y\in \Lambda$.  
 \cref{propmain}\eqref{prop:main2} follows immediately.  
\begin{proof}[Proof of \cref{propmain}\eqref{prop:main5}]
Given $y\in \Lambda$ and $g\in G$, take $g_0\in G$ and $z\in \Lambda $ so that $(g_0\inv ,y)\in \wtd \Lambda$ and $ (g_0\inv g\inv,z)\in \bar\scrG(g_0\inv,y).$
Then $g\cdot h(y)= g g_0 Q = h(z) = h(\ell_{g}(y))$.  \end{proof}

\subsubsection{Structure of $G$-orbits}
Given $x\in \Lambda$, let $L(x) = \{ \ell_g(x): g\in G\}$ be the $G$-orbit of $x$.  

\begin{proof}[Proof of \cref{propmain}\eqref{prop:main4} and \eqref{prop:main6}]
Fix $x\in \Lambda$.  Take $g_1\in G$ with $(g_1\inv,x)\in \wtd \Lambda$.

 For  \eqref{prop:main6}, consider $g_2\in G$ with $(g_2\inv,x) \in \bar \scrG(g_1\inv,x)$. 
By \cref{claim:related1}, $g_2\inv = q g_1\inv$ for some $q\in Q$.  In particular, if $\ell_g(x) = x$ then 
$g_1\inv g\inv = q g_1\inv$ and thus 
$g$ is of the form $g= g_1  q\inv g_1\inv$ for some $q\in Q$.  
Moreover, for any $q\in Q^\circ$, then $\ell_{g_1  q g_1\inv}(x) = x$.  Thus, the stabilizer group of $x$ and $g_1  Q g_1\inv$ have the same connected subgroup.  
 \eqref{prop:main6} follows.

For \eqref{prop:main4}, since $L(x)$ is a $G$-orbit, it is a topological  injectively immersed submanifold.  Moreover, an open neighborhood of $x$ in $L(x)$ is the injectively immersed $C^r$ manifold $V_{(g_1,x)}$; the collection $\{V_{(g',x')}: (g',x')\in \bar \scrG(g,x)\}$ forms a $C^r$ atlas on $L(x)$ compatible with the ambient smooth structure on $M$.   Thus the $G$-orbit $L(x)$ has the structure of a $C^r$ injectively immersed submanifold.  
\end{proof}

\begin{remark}
We note a subtle point in the proof of Proof of \cref{propmain}\eqref{prop:main4}.  If $(g_1\inv,x)\in \wtd \Lambda$ and $(g_2\inv,x)\in \wtd \Lambda$ then by \cref{claim:related1}, $ g_2\inv = q g_1\inv$ for some $q\in Q$.  Moreover, given any $q\in Q$, we have $(qg_1\inv,x)\in \wtd \Lambda$ and $$\wtd \Pi_2(\bar\scrG(qg_1\inv,x)) = \wtd \Pi_2(\bar\scrG(g_1\inv,x)) = L(x).$$

However, given $q\in Q$, we need not  have 
$$(qg_1\inv,x)\in \bar\scrG(g_1\inv, x)$$
unless $q\in Q^\circ.$
In particular, while $\ell_{g_1q g_1\inv }(x) = x$ for any $q\in Q^\circ$, for  $q\in Q\sm Q^\circ$, we need not have $\ell_{g_1q g_1\inv }(x) = x$.  

The defect of  $\bar\scrG(qg_1\inv,x) \neq  \bar\scrG(g_1\inv, x)= \emptyset$ as $q$ varies over the component group $Q/Q^\circ$ 
reflects to some extent the degree of the cover $h\colon L(x)\to G/Q$ restricted to the leaf $L(x)$ through $x$.  
\end{remark}

\subsection{Proof of \cref{thm:completepara} and \cref{thm:mainlamina}} \label{pf4647}
\subsubsection{Proof of \cref{thm:completepara}}
By hypothesis \ref{hypprojc} and \cref{prop:equifinite}, the $Q$-invariant probability measure $\mu$ in \cref{thm:completepara} is $\beta$-tame.  
Applying \cref{prop:msrrigi2,cor:easycor,upgradetofullstab}, hypotheses \ref{hyp:11} and  \ref{hyp:22} of \cref{propmain} hold (we note in this case that the partition element $U$ in hypotheses  \ref{hyp:22} is trivial). 
Let $\Lambda\subset M$ be the set given by \cref{propmain}.  
 By dimension counting, each $G$-orbit $L(x)$ of $x\in \Lambda$ in \cref{propmain} is open in $M$.  Thus $\Lambda=\bigcup_{x\in \Lambda} L(x)$ is an open subset of $M$.  

It follows from  hypothesis \ref{hypprojc} that the measures $\td \mu_g$ are finite for a.e.\ $g\in G$.  
By conclusion \eqref{prop:main2} of \cref{propmain}, there exists $k\in \N$ such that 
 the map $h\colon (M,\nu)\to (G/Q, \mathrm{Leb})$ is $k$-to-1 on a set of full $\nu$-measure.  It follows that 
 \begin{enumerate}
 \item the restriction $h\colon L(x) \to G/Q$ of $h$ to each $G$-orbit $L(x)$ for $x\in \Lambda$ is a finite $C^r$ covering map;
 \item there are only finitely many distinct $G$-orbits in $\Lambda$.  
 \end{enumerate}
 It follows that $\Lambda$ is the union of finitely many compact leaves, each with positive $\nu$-measure, and thus $\Lambda$ is compact.  As $M$ is connected, it follows that $M= \Lambda$ and that $\Lambda$ consists of a single $G$-orbit.  In particular, 
 there is a  $C^r$ finite covering map $h\colon M\to G/Q$ for such that 
 $$h(\alpha(\gamma)(x)) = \gamma\cdot h(x)$$
 for all $x\in M$ and $\gamma\in \Gamma$.  \hfill \qed

 \subsubsection{Proof of \cref{thm:mainlamina}} 
 Let $Q$ and $\mu$ be as in \cref{thm:mainlamina}.  
 By hypothesis, we have $c_0^F(\mu,S) = 1$ where $S$ is as in 
   hypotheses \ref{exptan} of \cref{thm:mainlamina}.  We may apply \cref{prop:msrrigi2,agree support,upgradetofullstab} to verify the hypotheses of \cref{propmain} hold.  \cref{propmain} then gives conclusions \eqref{laminate1}--\eqref{laminate8} of 
\cref{thm:mainlamina}.  

For  conclusion \eqref{laminate9}, the hypothesis that $\pi_1(G/Q)$ is finite gives that each leaf of $\Lambda $ is compact.  
It follows that the partition of $\Lambda$ into global leaves $\{L(x)\}$ is  a measurable partition of $(M, \nu)$.  Moreover, the (abstract measure theoretic) quotient by this partition gives a $\Gamma$-invariant measure.

Pass to the suspension $(M^\alpha, \mu)$, where $\mu$ is the $Q$-invariant Borel probability measure inducing $\nu$, and let $\what \Lambda= (G\times \Lambda)/\Gamma$ denote the suspension of $\Lambda$ in $M^\alpha$.  We have that  $\mu(\what \Lambda) =1$.   
Given $x\in M^\alpha$ let $F(x):=\pi\inv(\pi(x))$ be the fiber of $M^\alpha$ through $x$,  $\Lambda_x= F(x) \cap \what \Lambda$, and $\mu_x$ the conditional measure supported on $\Lambda_x\subset F(x)$.  Note that $\Lambda_x$ is similarly laminated by submanifolds of the fiber $F(x)$, each of which is a covering space of $G/Q$.  Write $L^F(x)$ for the fiberwise leaf of $\Lambda_x$ through $x$; that is $L^F(x)$ is the path component of $\Lambda_x$ containing $x$.  Then (for almost every $x$) the leaf $L^F(x)$ is compact and, 
as above, the partition $\{L^F(y): t\in F(x)\}$ of $\Lambda_x$ into global fiberwise leaves is $\mu_x$-measurable for a.e.\ $x$.  

Let  $a_0$ be as in \cref{propmain}.  Then for the action of $-a_0$, for $\mu$-a.e.\ $x$ the fiberwise Lyapunov exponents tangent to  leaves of $ \Lambda_x$ are positive and the fiberwise Lyapunov exponent transverse to leaves of $ \Lambda_x$ is negative.  
Let $W^{s, F}_{-a_0,\epsilon}(x)$ denote the stable fiberwise Lyapunov manifold for the action of $-a_0$ of radius $\epsilon>0$ (with respect to, say, Riemannian arc-length).  
For $\mu$-a.e.\ $x$, define 
\begin{equation*}c(x,\epsilon):= \mu_x\{L^F(y) \subset \Lambda_x:  W^{s, F}_{-a_0,\epsilon}(x)\cap L^F(y)\neq \emptyset\},\end{equation*}
the fiberwise measure of the global fiberwise leaves intersecting $W^{s, F}_{-a_0,\epsilon}(x)$.  

We have $c(x,\epsilon)>0$ for almost every $x$.  By iteration under $-a_0$, $A$-equivariance of the conditional measures $\mu_x$, and recurrence to certain sets, we conclude that 
$\mu_x(L^F(x))>0$ for a.e.\ $x\in M^\alpha$.  
By $P$-ergodicity of $\mu$ in \cref{propmain}, it follows there is $c>0$ such that $\mu_x(L^F(x))= c$ for $\mu$-a.e.\ $x$.

Returning to $\Lambda\subset M$, it follows that $\nu(L(x)) = c$ for $\nu$-a.e.\ $x\in \Lambda$.  In particular,  there are only finitely many leaves of $\Lambda$.  \hfill \qed

\section{Normal forms,  proof of \cref{upgradetofullstab}, and linearization of contracting dynamics}
We formulate the theory of nonuniform, non-stationary normal forms in the settings considered in \cref{thm:parabolicmeasr,thm:parabolicmeasrdim2,upgradetofullstab}.   For general formulation and further details, see \cite{NF,MR3893265,MR3642250}.  

 \subsection{Setup, notation, and standing assumptions}\label{ssec:2typesofNF}
Let $\mu$ be an ergodic, $A$-invariant Borel probability measure on $M^\alpha$ projecting to the Haar measure on $G/\Gamma$. 


Recall we write  $  \calL(\mu)$ for the collection of Lyapunov exponents for the $A$-action on $(M^\alpha,\mu)$.   Recall the definition of a strongly integrable subcollection $\calI\subset \calL(\mu)$ in \cref{def:strongInt}.  We also introduce the following: 
\begin{definition} 
A subcollection $\calI\subset  \calL(\mu)$ is \emph{totally resonant} if for every $\lambda \in \calI$, there exists a restricted root $\beta\in\Sigma$ and a fiberwise Lyapunov exponent $\bar \lambda\in \calL^F(\mu)$ with both $[\beta]=[
\lambda]$ and $\left [{\bar \lambda} \right ]=[ \lambda]$. 
 \end{definition}

In the proofs of \cref{thm:parabolicmeasr,thm:parabolicmeasrdim2,upgradetofullstab}, we will be specifically  interested in  the following explicit settings of totally resonant, strongly integrable collections of roots. 
\begin{enumerate}
\item \label{case1ofNF} $\calI=  [\beta]$ where $\beta\in \Sigma_\lieq^\perp$ is a restricted root that is resonant with the fiber Lyapunov exponents of $\mu$.
\item \label{case2ofNF} There exists $a_0\in A$ for which  $\calI= \calI^s(a_0)=  \left\{\lambda\in  \calL^F(\mu):\lambda (a_0)<0\right \}$; moreover,  every $\lambda \in \calI$ is positively proportional to some $\beta\in \Sigma_\lieq^\perp$.
\end{enumerate}

\subsection{Properties of   totally resonant, strongly integrable collections}\label{sec:intcoll}
Fix $\calI
\subset \calL(\mu)$ to be a totally resonant, strongly integrable collection of   Lyapunov exponents.
Observe that both 
$\Sigma_\calI $ and $ \calL^{F,\calI}(\mu)$ are saturated by positive proportionality  equivalence classes.  We thus write $\widehat \calI$ for the collection of coarse Lyapunov exponents in  $\calI$.   

\subsubsection{Total order on $\what \calI$}
To establish \cref{upgradetofullstab}, we consider setting \eqref{case2ofNF} of \cref{ssec:2typesofNF} and argue by induction on the cardinality of subcollections $\what \calI'\subset \what  \calI$ (see \cref{prop:induct} below).  To organize the combinatorics of the induction, we describe a choice of a total order $\chi_1<\chi_2<\dots<\chi_p$ on the elements of $\what  \calI$.  
Relative to this order, for every interval $A= \{i, i+1, \dots, i+k\}\subset \{1, \dots, p\}$ the collection $\{\chi_i: i\in A\}$ will be totally resonant and  strongly integrable, allowing us to establish properties by inducting on cardinality.  

%

We assume the original action $\alpha\colon \Gamma\to \diff^r(M)$ is $C^r$ for $r>r_0^F(\calI)$ where $r_0^F(\calI)$ is as defined in \eqref{eqcriffiberreg}.  Recall the negative cone  associated with $\calI$ in \cref{eq:cone}.  
Recall we identify $A$ with $\liea$ under the exponential map $\exp_\lieg$.  
Let $W\subset A$ be a  2-dimensional linear subspace chosen so that  the following properties hold:
\begin{enumerate}
\item 	There exists $a_0\in W\cap C(\calI)$ such that 
\begin{equation} \label{a0bunched} \max \left\{\frac { \bar \lambda(a_0) }{ \bar \lambda'(a_0) } :  \bar \lambda, \bar \lambda'\in \calL^{F,\calI}(\mu)\right\}<r\end{equation}
\end{enumerate}
and for  every $\chi\neq\chi'\in \what \calI $
\begin{enumerate}[resume]
\item $\ker \chi\cap W$ is 1-dimensional.
	\item $\ker \chi\cap W$ and $\ker \chi'\cap W$ are distinct.  
\end{enumerate}

%


Equip $W$ with an inner product and orientation.  Let $S\subset W$ denote the unit circle in $W$.  
Let  $I\subset S$ be the open arc oriented counter-clockwise from $-a_0 $ to $a_0$.  Let $<_I$ denote the induced linear order on $I$.  
For each $\chi\in \what \calI$, we have  that $\ker \chi \cap I$ is a singleton and the order on $I$ induces an order on the set of coarse weights $\what\calI$ by declaring $$\chi<\chi'$$ if 
$$(\ker \chi\cap I) <_I(\ker \chi'\cap I).$$
Re-indexing if necessary,  assume for the remainder that $\what \calI=\{\chi_1, \dots, \chi_p \}$ is indexed  so that $$\chi_1<\chi_2<\dots< \chi_p .$$

\subsubsection{Properties of the total order on $\what \calI$} 
We enumerate various properties of the above choice of order.  

Consider any interval $B= \{i, \dots, i+k\}\subset \{1, \dots, p \}$.  We write $$\calI_B=\bigcup _{j\in B}\chi_{j}$$  and $$\liev_B:= \bigoplus_{i\in B} \lieg^{\chi_i}.$$

\begin{claim}
Let $B= \{i, i+1, \dots, i+k\}\subset \{1,\dots, p \}$ be an interval relative to a choice of total order above.  
\begin{enumerate}
\item The  collection $\calI_B$ is totally resonant and strongly integrable; in particular, $\liev_B$ is a subalgebra of $\liev_\calI$.  
\end{enumerate}
For the  subintervals  $B'= \{i, i+1, \dots, i+k-1\} $ and $B''= \{ i+1, \dots, i+k\}$ of $B$,  
\begin{enumerate}[resume]
	\item both $\liev _{B'}$ and $\liev _{B''}$ are ideals in $\liev _B$.
\end{enumerate}
\end{claim}


In the sequel, we write  $V_B$ for  the analytic subgroup of $V_\calI$ corresponding to the interval $B$.

Given a subset $B\subset\what \calI$, for $\mu$-a.e.\ $x\in M^\alpha$ we write 
$$E^{B}(x) = \bigoplus_{i\in B} E^{\chi_i}(x)\quad \text {and} \quad 
E^{B,F}(x) = \bigoplus_{i\in B} E^{\chi_i, F}(x)$$
for the associated vector spaces spanned by (fiberwise) coarse Lyapunov subspaces.  


We may apply \cref{lyapmanifolds} to the strongly integrable, totally resonant subcollection $\calI_B$ associated to any interval $B\subset \{1,\dots, p\}.$  We write $$W^{B,F}(x)\subset \calF(x), \quad \text{ and } \quad W^{B}(x)\subset M^\alpha$$ for the associate Lyapunov manifolds tangent to  $E^{B,F}(x)$ and $E^{B}(x)$, respectively.

\subsection{Subresonant polynomial maps between vector spaces}
\subsubsection{Lyapunov weights}
Let $V^0,V^1,V^2,\dots$ be a sequence of finite-dimensional normed vector spaces of the same dimension.  
  For $i\in \{0,1,2,\dots\}$, let $T_i\colon V^i\to V^{i+1}$ be a linear map.   In practice, we will assume $T_i$ is invertible with norm growing at most  subexponentially in $i$.  Let $V= V^0$ and  assign to $V$ the \emph{Lyapunov weight} (i.e.\ Lyapunov exponent)
$\varpi= \varpi_V\colon V\to \R\cup\{\pm \infty\}$,
 $$\varpi_V(v) = \limsup_{n\to \infty}\frac 1 n \log\| T_{n-1}\circ \cdots \circ T_0 (v)\|.$$
We observe that $\varpi_V$ takes only finitely many values and satisfies
\begin{enumerate}
\item \label{SG222} $\varpi_V(c v) = \varpi_V(v)$ for all  $c\in \R\sm \{0\}$;
\item \label{SG333} $\varpi_V(v+w) \le  \max\{\varpi_V(v), \varpi_V(w)\}$.  
\end{enumerate}
In our future applications, we will always consider a sequence of maps $T_i$ such that 
\begin{enumerate}[resume]
\item $\varpi_V(v)<\infty$ for all $v\in V$ and $\varpi_V(v) = -\infty$ if and only if $v=0$.
\end{enumerate}
Let $\lambda_1<\lambda_2 <\dots <\lambda_\ell$ be the finite values of $\varpi=\varpi_V$.  Write  $V_i =V_i^0= \{ v\in V\mid \varpi(v) \le \lambda_i\}$ and let 
$$\calV_{\varpi}=\bigl\{\{0\}=V_{0}\subset V_1\subset\dots\subset V_\ell=V\bigr\}.$$
be the \emph{associated flag} on $V$.   The dimensions $\dim V_i/V_{i-1}$ are the \emph{multiplicities} of the weight $\varpi$.  

\subsubsection{Subresonant maps and their properties}
Let $V$ and $W$ be vector spaces equipped with Lyapunov weights $\varpi_V$ and $\varpi_W$ as above.   A linear map $T\colon V\to W$ is \emph{subresonant} if $$\varpi_W(T(v))\le \varpi_V(v)$$ for all $v\in V$.
A degree-$k$ homogenous polynomial $g\colon V\to W$ is \emph{subresonant} if (the polarization of) the $k$th derivative  $D^{(k)}g\colon V^{\otimes k}\to W$ is a subresonant ($k$-multi)linear map (where  $V^{\otimes k}$ is equipped with the induced Lyapunov weight satisfying  $\varpi_V(v_1 \otimes \dots \otimes v_k) = \varpi_V(v_1) + \dots + \varpi_V (v_k)$ on pure tensors).
If $\varpi(w)\le 0$, we also declare the constant map  $V\to W$, $v\mapsto w$ to be subresonant.
A polynomial $g\colon V\to W$ is \emph{subresonant} if, writing  $g= g_0 + g_1 + \dots + g_d$ in homogenous components, each $g_i$ is subresonant.

Let $\calP^{SR}(V,W)$ denote the set of subresonant polynomial maps from $V$ to $W$.  
Let $\calG^{SR}(V)$ denote the set of subresonant polynomials maps $g\colon V\to V$ with $D_0g$ invertible.

We summarize a number of properties of the above definitions.  
\begin{lemma}Suppose $V$ and $W$  are equipped with 
weights $\varpi_V$  and $\varpi_W$  taking negative values $\lambda_1<\dots <\lambda_\ell<0$ and  $\eta_1<\dots <\eta_p<0$, respectively.  
Then 
\begin{enumerate}
\item Every $g\in \calP^{SR}(V,W)$ has degree at most $d=\lfloor \eta_1/\lambda_\ell\rfloor$.
\item 
The set $\calP^{SR}(V,W)$ is a finite dimensional vector space (under pointwise addition.)
\end{enumerate}
\end{lemma}

\begin{lemma}Suppose $V$ and $W$ have the same dimension and are equipped with 
weights  $\varpi_V$  and $\varpi_W$ 
taking negative values $\lambda_1<\dots <\lambda_\ell<0$ with the same multiplicities.  
Then
\begin{enumerate}
\item A linear map $T\colon V\to W$ is invertible if and only if $\varpi(T(v)) = \varpi (v)$ for all $v\in V$.  
\item Let  $g\colon V\to W$  be a  subresonant polynomial whose derivative at the origin,  $D_0g$,  is invertible.  Then $g$ is a  diffeomorphism and $g\inv \colon W\to V$ is a subresonant polynomial.  
\item   Under functional composition, $\calG^{SR}(V)$ is a finite dimensional Lie group.
\end{enumerate}

\end{lemma}
We further remark since $\lambda_\ell<0$, the group of invertible subresonant polynomial self-maps $\calG^{SR}(V)$ includes  all translations and thus acts transitively on $V$.

\subsubsection{Parameterized families of subresonant polynomials}
Let $V$ and $U$ be finite dimensional,   vector spaces equipped with weights taking only negative values.   

  Given $v\in V$, let $L_v\colon \calP^{SR}(V,U)\to U$, $$L_v\colon \phi\mapsto \phi(v),$$ denote the evaluation functional. 
Using that $\calP^{SR}(V,U)$ is finite-dimensional, we find finitely many $v_i$ so that the weak topology induced by the associated evaluation functionals coincides with the norm topology on $\calP^{SR}(V,U)$.
\begin{claim}\label{claim:unique}
Suppose the weight on $V$ takes only negative weights.  
Let $B\subset V$ be a subset of full Hausdorff dimension. 
Then there exists a finite set $F=\{v^1, \dots, v^p\} \subset B$ such that every $\phi\in  \calP^{SR}(V,U)$ is uniquely determined by its values on $F$.

Moreover,  a sequence $\phi_n \in   \calP^{SR}(V,U)$ 
converges (in every $C^r$ topology)  if and only if 
$\phi_n(  v^i)$ converges in $U$ for each $v_i\in F$.  
\end{claim}

\begin{proof}
First observe that given $\phi\in  \calP^{SR}(V,U)$, if $\phi(v)=0$ for every $v\in B$, then $\phi(v)=0$.  Indeed, the zero set of a non-zero polynomial is a proper algebraic subvarietery and hence has Hausdorff dimension at most $\dim(V)-1$.   We have   $$\bigcap _{v\in B} \ker L
_v =\{0\}.$$ 
As $ \calP^{SR}(V,U)$ is finite dimensional,  there are finitely many $v^1, \dots, v^p\in B$ such that $\bigcap _{i=1}^p \ker L
_{v^i} =\{0\}$.  Setting $F= \{v^1, \dots, v^p\}$, every $\phi\in \calP^{SR}(V,U)$ is uniquely determined by its values on $F$ and the first conclusion follows.

The map $\Psi_F\colon \phi\mapsto (\phi(v^1), \dots, \phi(v^p))$  is linear map.   Since the kernel of $\Psi_F$ is trivial, $\Psi_F$ induces an isomorphism between $ \calP^{SR}(V,U)$ and a linear subspace of $U^p$. Since $ \calP^{SR}(V,U)$ is finite dimensional, the weak topology coincides with the norm topology and the second conclusion follows.  
\end{proof}

We may identify $ \calP^{SR}(V,U)$ as a subspace of some Euclidean space either through its image under the map $\Psi_F$ in the proof of \cref{claim:unique}.  Alternatively, we may  identify the polarization of each homogeneous term with a matrix.  In either way, we make sense a $C^r$-parametrized family $t\mapsto \phi_t\in  \calP^{SR}(V,U)$.

\begin{corollary}\label{claim:smoothing}
Fix $r>0$.  Let $D\subset \R^k$ be an open subset.  
Fix any finite set $F\subset V$ satisfying the conclusions of \cref{claim:unique}.
Let $t\mapsto \phi_t\colon V\to U$, $t\in D$, be a parameterized family of functions such that
\begin{enumerate}
\item for a dense subset $D_0\subset  D$ and every $t\in D_0$, $\phi_t\in  \calP^{SR}(V,U)$; 
\item for every $v\in F$, 
 the function $D\to U$, $t\mapsto \phi_t(v)$, is $C^r$.
\end{enumerate}
Then there  is 
 a parameterized family $t\mapsto \hat \phi_t$  such that  
\begin{enumerate} [label=(\alph*),ref=(\alph*)]
\item $\hat \phi_t= \phi_t$  for every $t\in D_0$; 
\item $\hat \phi_t\in  \calP^{SR}(V,U)$ for every $t\in D$; 
\item the map 
  $t\mapsto  \hat \phi_t$ is $C^r$ on $D$. 
\end{enumerate}
\end{corollary}
\begin{proof}
Consider  a sequence  $t_i\in D_0$  
and suppose $t_i\to t_0\in D$.   From  \cref{claim:unique}, $\phi_{t_i}$ converges in $  \calP^{SR}(V,U)$.  Letting $\hat \phi_{t_0}$ denote the limit of this sequence, we have that $\phi_{t_0}\in  \calP^{SR}(V,U).$
The conclusion follows.
\end{proof}

\begin{corollary}\label{cor:autosmooth}
Let $W_1,W_2$ be subspaces of $V$ that span $V$.  
Suppose $\phi\colon V\to U$ is a function such that for almost every 
$v\in V,$ the restrictions 
$$W_1\to U, \ \ w_1\mapsto \phi(v+ w_1),\quad \quad W_2\to U, \ \ w_2\mapsto \phi(v+ w_2)$$
are subresonant polynomials.  Then $\phi$ coincides $\bmod 0$ with a $C^r$ function.  
\end{corollary}
\begin{proof}
We first observe that if $w_1\mapsto \phi(v_0+ w_1)$ is a subresonant polynomial then so is $w_1\mapsto \phi(v+ w_1)$ for all $v\in v_0 + W_1$.  
Let $T\subset W_2$ be such that $V= W_1 \oplus T$.  

There a finite set $F\subset W_1$ such that every $\psi\in \calP^{SR}(W_1,U)$ is uniquely determined by its value on $F$ and for which $T\to U$,  
$t\mapsto \phi(v+t)$, is $C^r$ for every $v\in F$.  
For almost every $t\in T$, the map $w_1\mapsto \phi(w_1+t)$ is a subresonant polynomial.  The conclusion follows from \cref{claim:smoothing}.  
\end{proof}

\begin{remark}
In \cref{cor:autosmooth}, it does not automatically follow that $\phi$ is subresonant.

Indeed consider the case that $V=\R^2$ has a single (negative) weight $\lambda<0$.  Let $U=\R$ be equipped with the weight $\varpi_U(v) =\lambda$ if $v\neq 0$.  Consider the function $\phi(x,y) = x+ y + xy$.  The restriction of $\phi$ to every horizontal and vertical line is subresonant.  But $\phi$ is not subresonant.   
\end{remark}

\subsection{\texorpdfstring{Nonuniform, non-stationary normal forms along $W^{B,F}$-leaves}{Nonuniform, non-stationary normal forms along leaves}}
Let $\calI$ be a totally resonant, strongly integrable collection of Lyapunov exponents as above.  
Fix  $a_0\in C(\calI)$ satisfying \eqref{a0bunched} above.  View  $a_0\in A\subset G$ as a diffeomorphism of $M^\alpha$ (identifying $a_0$ with its induced action by translation on $M^\alpha$.)  
Let $L\colon \liev_\calI\to \liev_\calI$ denote $$L= \restrict{\Ad(a_0)}{\liev_\calI}.$$

Fix an interval $B\subset \{1, \dots, p\}$ as in \cref{sec:intcoll}.  
Given $X\in \liev_B$, define the Lyapunov weight $$\varpi(X):= \lim_{n\to \infty}\frac 1 n \log \|L^n(X)\|.$$
Similarly, given $w\in E^{B,F}(x)$, there is a full $\mu$-measure subset $\Lambda_0$ for which \cref{thm:higherrankMET} holds (for a.e.\ $a_0$-ergodic component of $\mu$).  
For $x\in \Lambda_0$, define the Lyapunov weight $$\varpi(w):= \lim_{n\to \infty}\frac 1 n \log \|(D_x a_0^n)w\|.$$
By $A$-ergodicity, $\varpi$ takes finitely many values, independent of $x$.

From the construction of normal forms on contracting equivariant foliations (see  especially \cite{MR3642250,MR3893265}), we have the following:
\begin{proposition}\label{prop:NF1} Assume $r>r_0^{F}(\calI)$.  
There exists a full $\mu$-measure  subset $\Omega\subset  M^\alpha$ with the following properties.  

For every interval $B= \{i, i+1, \dots, i+k\}\subset \{1,\dots, p \}$ and every $x\in \Omega$ there 
exists a collection of $C^r$ diffeomorphisms $\calH_{x,B}:=\{h\colon E^{B,F}(x)\to W^{B,F}(x)\}$ with the following properties:
\begin{enumerate}
\item\label{NF1:1} The group $\calG^{SR}(E^{B,F}(x))$ acts transitively and freely on $\calH_{x,B}$ by pre-composition.  
\item Given a subinterval $B'\subset B$, every $h\in \calH_{x,B'}$ is of the form $h= \restrict{\td h}{E^{B',F}(x)}$ for some $\td h\in \calH_{x,B}$ with $\td h (0) \in W^{B',F}(x)$
 \item\label{NF1:3}  Given any $a\in A$, let $g\colon M^\alpha\to M^\alpha$ denote translation by $a$.  Then for $\mu$-a.e.\ $x\in \Omega$ we have $g(x)\in \Omega$; moreover, for every $ h\in \calH_{x,B}$ and $\td h\in \calH_{g(x),B}$, the map 
$$ \td h\inv \circ g \circ h$$
 is an element of $\calP^{SR}\left(E^{B,F}(x),E^{B,F}(g(x))\right)$
  \item\label{NF1:4}  Given $x\in \Omega$ and  $y\in \Omega \cap W^{B}(x)$, let $v\in V_B$ be such that $y = v\cdot z$ for $z\in W^{B,F}(x)$.  Then for any $h_y\in  \calH_{y,B} $ and  $h_x\in \calH_{x,B}$ the map $$w\mapsto h_x\inv (v\inv \cdot h_y(w))$$
 is an element of $\calP^{SR}\left(E^{B,F}(y),E^{B,F}(x)\right)$
 \item\label{NF1:5}  There exists a choice of measurable section $x\mapsto h_x$ with 
\begin{enumerate}
\item $h_x(0) = x$ and
\item  $D_0h_x= \Id$.
\end{enumerate}
\end{enumerate}

\end{proposition}
\begin{remark}\label{1dNF}
When $\dim E^{B,F}(x) =1$ for almost every $x$, the group $\calG^{SR}(E^{B,F}(x))$ is the affine group of a 1-dimensional space; thus, there is a unique  choice of $h_x\in \calH_{x,B}$ with 
$h_x(0) = x$ and  $D_0h_x= \Id$.  Moreover for this choice of $h_x$, for any $a\in A$, $\mu$-a.e.\ $x$, and every $v\in \dim E^{B,F}(x)$,
\begin{equation}
\label{eq:linearNF}
h_{a\cdot x}\inv ( (a\cdot h_x(v)))=D_0a(v). 
\end{equation}
\end{remark}


\subsection{Proof of \cref{upgradetofullstab}}\label{sec:upgrade}

We now suppose the hypotheses of \cref{upgradetofullstab} hold.  
In particular, $\mu$ is assumed ergodic and $P$-invariant probability measure on $M^\alpha$.  We let $\calI=\calI^s(a_0)$ be as in setting \eqref{case2ofNF} of \cref{ssec:2typesofNF}.  We equip $\what \calI$ with a choice of total order and enumerate $\what \calI$ as $\chi_1<\dots<\chi_p$ as above.  

For each interval $B\subset \{1,\dots, p\} $, write $\mu_x^{B}$ for the leafwise measure along the manifolds $W^{B}(x)$. Using the hypotheses of \cref{upgradetofullstab} as our base case, we induct on the cardinaly of $B$ and  establish the following  which, in turn, establishes \cref{upgradetofullstab}. 
\begin{proposition}\label{prop:induct}
For every interval $B\subset \{1,\dots, p\}$ and for $\mu$-a.e.\ $x\in M^\alpha$, there is an injective $C^r$ embedding $\phi_{x,B}\colon V_{B}\to W^{B,F}(x)$ with $\phi_{x,B}(\1) = x$ such that 
the leafwise measure $\mu_x^{B}$ is the image of the Haar measure on $V_{B}$ under the graph 
$$u\mapsto u\cdot  \phi_{x,B}(u).$$
\end{proposition}
\begin{proof}
We induct on the cardinality $\#B$ of the set $B$.
When $\#B=1$, we are considering a single coarse restricted root $\chi_i = [\beta]$ for $\beta\in \Sigma_\lieq^\perp$; the conclusion then follows from the hypotheses of  \cref{upgradetofullstab}.  

We thus assume the proposition holds for all intervals $B'\subset \{1,\dots, p\}$ of cardinality $\#B'= k$ for $1\le k\le p-1$.  Fix   an interval $B=\{i, i+1, \dots, i + k\}\subset \{1,\dots, p\}$ with $\#B= k+1$.  Also fix subintervals  $B'= \{i, i+1, \dots, i + (k-1)\}$ and $B''= \{i+1, i+2, \dots, i + k\}$ of cardinality $k$.  

We observe from the product structure of entropy \cite[Corollary 13.2]{MR4599404} and the hypothesis of \cref{upgradetofullstab} that for every $a\in A$, the fiberwise entropy $h_\mu(a\mid \calF)$ is zero.  Indeed, by hypothesis 
the leafwise measures along every fiberwise coarse Lyapunov foliation associated with $\chi\in \calI$ is purely atomic.  By the product structure of entropy \cref{ent:prod} (applied to $-a_0$), for every $a\in C(\Sigma_\lieq^\perp)$
the fiberwise entropy of $-a$ is zero.  It follows the fiberwise entropy for entire $A$-action is zero.   
 In particular, the hypothesis of \cref{upgradetofullstab} imply there exists a measurable subset $\Upsilon\subset M^\alpha$ with $\mu(\Upsilon)=1$ 
such that for every $x\in \Upsilon$  and every interval $B\subset \{1,\dots,  p\}$, the intersection $W^{B,F}(x)\cap \Upsilon$ has cardinality 1.

Given $y\in W^{B}(x)\cap \Upsilon$, write $y= v_y\cdot z_y$ where $v_y\in V_B$ and $z_y\in W^{B,F}(x)$.  
Let $\pi_{x,B}\colon W^{B}(x)\cap \Upsilon \to V_B$ denote the map $$\pi_{x,B}(y) =\pi_{x,B}(v_y\cdot z_y) = v_y .$$
Let $\hat \mu_{x,B}=( \pi_{x,B})_*\mu_{x,B}$.  
The restriction of $\hat \mu_{x,B}$ to the   Borel $\sigma$-algebra on $V_B$ is $\sigma$-finite.  Moreover, by the inductive hypothesis, for $\hat \mu_{x,B}$-a.e.\ $v\in V_B$, the leafwise measure of $\hat \mu_{x,B}$ along the $V_{B'}$- and $V_{B''}$-orbits $V_{B'}\cdot v$ and $V_{B''}\cdot v$ are Haar measures on the respective orbits.  It follows for $\mu$-a.e.\ $x$ that $\hat \mu_{x,B}$ is $V_{B'}$-invariant and $V_{B''}$-invariant.  Since every $v\in V_B$ is the product of an element of $V_{B'}$ with a product of an element of $V_{B''}$, it follows that $\hat \mu_{x,B}$ is $V_B$-invariant.  It then follows (see e.g.\ \cite{MR2236571}) that $\hat \mu_{x,B}$ is a Haar measure on $V_B$.  
In particular, for $\mu$-a.e.\ $x$ there exists a Haar-measurable function $\phi_{x, B}\colon V_B\to W^{B,F}(x)$ such that the leafwise measure $\mu_x^{B}$ is the image of the Haar measure on $V_{B}$ under the graph 
$$u\mapsto u\cdot  \phi_{x,B}(u).$$
Note that for $\mu$-a.e.\ $x$, the inductive hypothesis implies $\phi_{x,B}(\1)=\phi_{x,B'}(\1)=\phi_{x,B''}(\1)=x$.

Select a measurable section $x\mapsto h_{x,B}\in \calH_{x, B}$ 
with $h_{x,B}(0) = x$ and $D_0h_{x,B}=\id$.  
Let $\wtd \phi_{x,B}(u)\colon \liev_B\to E^{B,F}(x)$ be the Lebesgue measurable function 
$$\wtd \phi_{x,B} = h_x\inv \circ \phi_{x,B}\circ \exp_{\liev}.$$
To finish the proof, it suffices to show for $\mu$-a.e.\ $x$ that $\wtd \phi_{x,B}$ coincides with a $C^r$ function off a (Lebesgue) conull set in $\liev_B$.  Indeed, injectivity of $\wtd \phi_{x,B}$ follows from entropy considerations; moreover, the derivative of $\wtd \phi_{x,B}$ is non-singular at $x$, hence in a neighborhood of $x$.  It then follows that $\wtd \phi_{x,B}$ is locally an embedding which can be promoted globally to an embedding by equivariance under a choice of  dynamics expanding  $W^{B}$-leaves.  

We thus show $\wtd \phi_{x,B}$ coincides with a $C^r$ function off a (Lebesgue) conull set in $\liev_B$.
Write $F(x) = a_0\cdot x$ for the translation by $a_0$.
Consider $x$ satisfying the conclusions of \cref{prop:NF1}.   Let $f_{x,B}\colon E^{B,F}(x)\to E^{B,F}(a_0\cdot x)$ be the $C^r$ diffeomorphism 
$$
f_{x,B}\colon v\mapsto 
(h_{F(x),B})\inv \left(a_0\cdot  h_{x,B}(v)\right).$$
We have that each $f_{x,B}$ is an element of $\calP^{SR}\left(E^{B,F}(x), E^{B,F}(F(x))\right)$ with $D_0f_{x,B}= \restrict {D_x F}{E^{B,F}(x)}$.


Recall $\liev_{B'}$ and $\liev_{B''}$ are ideals in $\liev_B$.  Given $Y\in \liev_B$, the Baker-Campbell-Hausdorff formula implies there  is a polynomial map $\tau_Y^{B'}\colon \liev_{B'}\to \liev_{B'}$   such that for $Z\in  \liev_{B'}$
$$\exp_\liev(Y+Z) = \exp_\liev(\tau ^{B'}_Y(Z)) \exp_\liev(Y).$$
Let $g\colon \liev_B\to \liev_B$ be $g(Y) = \Ad(a_0)Y$.  Then we have 
$$\exp_\liev(g(Y+Z)) = \exp_\liev(g(\tau ^{B'}_Y(Z))) \exp_\liev(g Y).$$
hence $\tau ^ {B'}_{g Y} (g Z) = g(\tau ^ {B'}_Y(Z))).$
Using that $\tau ^{B'}_Y\to \id$ as $Y\to 0$ and  
$\|\restrict{\Ad(a_0)}{\liev_B}\|<0$, one may directly verify that $\tau ^{B'}_Y$ is an element of $\calP^{SR}(\liev_{B'},\liev_{B'})$.  
(See for instance \cite[Prop.\ 3.4]{NF} for more complicated formulation.)



Given  $\mu$-a.e.\ $x$, consider Lebesgue-a.e.\ $Y\in \liev_B$ and let $v= \exp(Y)$ and $y = v\cdot \phi_{x,B}(v)$.  
Given $Z'\in \liev_{B'}$, we have 
\begin{align*}
\wtd \phi_{x,B}(Y+ Z') &=  h_{x,B}\inv \left(v\inv \cdot \phi_{y,B'}(\exp_{\liev}(Z'))\right)\\
 &=  h_{x,B}\inv \left(v\inv  \cdot  h_{y,B}\left( h_{y,B}\inv (\phi_{y,B'}(\exp_{\liev}(Z'))\right)\right)
\end{align*}
From the inductive hypothesis,  for $\mu$-a.e.\ $y$, the map  $$h_{y,B}\inv \circ \phi_{y,B'}\circ \exp_{\liev_{B'}}$$  
is $C^r$;
again, one my check (less directly than the above, see \cite[Prop.\ 3.4]{NF}) this implies the map is an element of $\calP^{SR}(\liev_{B'}, E^{B',F}(y)$.  
In particular,  for $\mu$-a.e.\ $x$ and Lebesgue-a.e.\ $Y\in \liev_B$, the map
$$\liev_{B'}\to  E^{B,F}(x),\quad Z'\mapsto \wtd \phi_{x,B}(Y+ Z')$$ is a subresonant polynomial map.  
Similarly, the map $$\liev_{B''}\to  E^{B,F}(x),\quad Z''\mapsto \wtd \phi_{x,B}(Y+ Z'')$$ is a subresonant polynomial map. 
It follows from \cref{cor:autosmooth} that $\wtd \phi_{x,B}$ coincides $\bmod 0$ (Lebesgue) with a $C^r$ function $\liev_B\to  E^{B,F}(x)$.  
\end{proof}

\label{oinkgo}


\subsection{Linearization of dynamics along contracting foliations}
 We again consider $\calI\subset \calL(\mu)$ a totally resonant, strongly integrable collection.  Equip the collection of coarse exponents $\what \calI$ with a total order as in \cref{sec:intcoll} and let $B\subset \{1,\dots, p\}$ be an interval in the index set of $\hat \calI$.

We assume  $a_0\in C(\calI)$ is chosen so that  $r_0^{F}(\calI)\le  r_0^{\calI,F}(a_0)<r$.  
In the sequel, we will only apply the following  construction  in the case that $\calI= [\beta]$ for some $\beta\in \Sigma_\lieq^\perp$.
The following construction will be a special a construction in a forth-coming work with Brown, Eskin, Filip, and Rodriguez Hertz.  See \cite{NF} for an discussion of this construction. 
\begin{proposition}\label{prop:linearization}
There exists an $A$-invariant subset  $\Omega\subset M^\alpha$ with $\mu(\Omega)=1$ and a $\mu$-measurable finite-dimensional  vector bundle  $\bfV^{\calI,F}\to M^\alpha$ (whose fiber over $x$ is denoted $\bfV^{\calI,F}_x$),
such that for $\mu$-a.e.\ $x$ there exists 
\begin{enumerate} 
\item an injective, $C^r$ embedding $\iota^{\calI,F}_x\colon W^{\calI,F}(x) \to \bfV^{\calI,F}_x$ with $\iota^{\calI,F}_x(x) = 0$
\end{enumerate}
such that 
\begin{enumerate}[resume]
\item for all $x,y\in \Omega$ and every $g\in \calP^{SR}(W^{\calI,F}(y), W^{\calI,F}(x))$, there exists an affine map $$(\bfL g)\colon \bfV^{\calI,F}(y)\to \bfV^{\calI,F}(x)$$ with 
$$(\bfL g) \circ \iota^{\calI,F}_y = \iota^{\calI,F} _x\circ g;$$
moreover, if $g(y) = x$ then $(\bfL g)$ is linear.  
\end{enumerate}

Moreover, the following hold:
\begin{enumerate}[resume]
\item \label{linearizedNF:5}  For every $a\in A$, there exists a linear cocycle  $(\bfL a) $ over translation by $a$ on $M^\alpha$ such that for all $y\in W^{\calI,F}(x)$,
$$(\bfL  a)_x\circ \iota^{\calI,F}_x (y)= \iota^{\calI,F}_{a x} (a\cdot y) .$$
Moreover, for $a,b\in A$, $$ (\bfL ab)= (\bfL a)\circ (\bfL b).$$

\item  \label{linearizedNF:5} For each $a_0\in C(\calI)\sm \{0\} \subset A$, there exists a measurable choice of norm on $\bfV^{\calI,F}$ such that for $\mu$-a.e.\ $x\in M^\alpha$, the sequence $\{ (\bfL a_0)_x ^n\}_{n\in \Z}$ is bi-regular; moreover, if $\lambda_1<\dots<\lambda_p<0$ are the Lyapunov exponents for the cocycle $\restrict {D_xa_0}{E^{\calI,F}_x}$, the Lyapunov exponents 
of the cocycle  $(\bfL a_0)$ are of the form $$\sum_{1\le i\le p} c_i \lambda_i$$
where $c_i$ are non-negative integers satisfying  $$\lambda_1\le \sum_{1\le i\le p} c_i \lambda_i\le \lambda_p.$$
\item  \label{linearizedNF:6} Given $y\in \Omega \cap W^{\calI}(x)$, let $v,u \in U^\calI$ be such that $y = v\cdot z$ for $z\in W^{\calI,F}(x)$.  
There exists an affine map $\bfT^F_{y\to x}\colon \bfV^{\calI,F}_y\to  \bfV^{\calI,F}_x$ such that for all $w\in W^{\calI,F}(y)$,
$$\bfT_{y\to x} \circ \iota_y(w) = \iota_x(u\inv w).$$

\item Relative to the Lyapunov weights in \eqref{linearizedNF:5}, the affine maps $\bfT_{y\to x}$ and $(\bfL  a)_x$ are subresonant.   
\end{enumerate}

\end{proposition}
We remark that we could, in fact, choose the norms on the $\bfV^{\calI,F}$ so that for every $a\in A$, the cocycle $\bfL^{\calI,F}(a)$ is uniformly bounded.  Since this won't be needed in any of our analysis, we only state the existence of a norm such that the dynamics of $a_0$ is bi-regular.  

 \cref{prop:linearization} is essentially a repackaging of the well-known normal form coordinates (conjugating the dynamics in leaves to  subresonant polynomial maps); see for example \cite{MR3642250,MR3893265}.  
See \cite[Prop.\ 8.8]{NF} for further discussion.

\begin{remark}\label{rem:afinetrans}
For the dynamics of translation by $a\in A$ on a $U^\calI$-orbit  (for some strongly integrable $\calI$) in $G/\Gamma$ or $M^\alpha$, we may use the exponential map $\exp\colon \lieu^\calI\to U^\calI$ to linearize the dynamics.  However, if $\calI$ contains multiple coarse roots, it may be that the transition maps corresponding to different base points in the orbit are subresonant, but possibly non-affine, polynomial maps (as in the proof of \cref{prop:induct}).

In the sequel, we will often consider the case that $\calI$ is a single coarse root $\calI= [\beta]$.  In this case, the transition maps, given by right translation in $U^{[\beta]}=\exp\lieu^{[\beta]}$,  are affine maps.  Indeed, we need only check the case when if $\beta$ and $2\beta$ are roots.  Write $\lieu^{[\beta]}=\lieu^{\beta}\oplus\lieu^{2\beta}.$ Given $X,Y\in \lieu^{[\beta]}$ write 
$X= (X_1, X_2)$ and $Y= (Y_1, Y_2)$ relative to $\lieu^{\beta}\oplus\lieu^{2\beta}.$  Consider right translation on $U^{[\beta]}$ by $\exp X$ relative to the coordinates given by $Y\mapsto \exp Y$.  By the Baker-Campbell-Hausdorff formula,
$$\exp(Y)\cdot \exp(X) = \exp Z$$
where 
$$Z=(Z_1, Z_2)= \left(Y_1+ X_1, Y_2+X_2 + \frac 1 2 [Y_1,X_1]\right).$$ The map $(Y_1,Y_2)\mapsto  \left(Y_1+ X_1, Y_2+X_2 + \frac 1 2 [Y_1,X_1]\right)$
is affine (in the parameter $Y$), with subresonant linear part with respecting the splitting  $\lieu^{[\beta]}=\lieu^{\beta}\oplus\lieu^{2\beta}.$
\end{remark}

\section{Linear cocycles and equivariant families of measures}\label{sec:coccyle}
Throughout this section, fix  a standard probability space $(X,\mu)$ and let $f\colon X\to X$ be an invertible,  measurable, $\mu$-preserving map. 

We consider a measurable  linear cocycle $L\colon \calE\to \calE$ factoring over $f$ on a measurable, finite-dimensional vector bundle $\calE$ equipped with an inner product.  By selecting a measurable orthonormal frame for $\calE$, we identify each fiber $\calE$ with a fixed inner-product space $V$ and each linear map $L_x\colon \calE_x\to \calE_{f(x)}$ with an element of $\GL(V)$.

We write $\|L_x\|$ and $m(L_x)$, respectively, for the operator norm and conorm of $L_x$.  

\def\calR{\mathcal R}

Let  $\calR(V)$ denote the set of locally finite Radon measures on $V$ equipped with the topology dual to compactly supported continuous functions on $V$.  
Consider a measurable function  $X\to \calR(V)$ denoted by  $x\mapsto \nu_x$.  For the remainder, we additionally assume the following   hypotheses: for $\mu$-a.e.\ $x$, 
\begin{enumerate}
\item the family $\{\nu_x\}$ is projectively equivariant under the cocycle $L_x$:
$$(L_x)_*\nu_x \propto \nu_{f(x)}.$$
\item the support of $\nu_x$ contains the origin $0$.
\end{enumerate}


Write $W_x\subset V$ for the \emph{essential linear span} of $\nu_x$; that is, $W_x$ is the intersection of all  subspaces  $W'\subset V$ with $\nu_x(V\sm W')= 0$.  By projective equivariance of $\nu_x$, we have $L_xW_x = W_{f(x)}$ a.s. 

The restriction  $\wtd L_x = \restrict {L_x}{W_x}$ is a measurable linear cocycle over the dynamics of $f\colon X\to X$.  
For $n\ge 1$, write $ \wtd L_x^{(n)}= \wtd L_{f^{n-1}(x)}  \circ \dots \circ   \wtd L_x$.

In the remainder of this section, we consider  the above setup in the following two settings: 
\begin{enumerate}
\item First, we will assume each measure in the  family  $\{\nu_x\}$ is a  finite measure and that the family is  invariant (rather than projectively invariant) under the cocycle $L_x$.  In this setting, 
we show the following  fact:  for typical $x$,   the norms and conorms of the sequence of matrices  $\{ \wtd L_x^{(n)}\}$ are uniformly bounded away from $\infty$ and $0$ for a  large density subset of parameters $n$.
\item Second, we consider the case that the cocycle $L$ is bi-regular over almost every trajectory 
and that 
the top Lyapunov exponent $\lambda_\top(x)$ is negative  for almost every $x$.  We then derive standard facts about the group $G_x$ of subresonant  affine maps preserving $\nu_x$ up to  normalization. In particular,  under the additional assumption that measure $\nu_x$ is supported on a single $G_x$-orbit,  we show there is a (connected) unipotent subgroup $U_x\subset G_x$ such that $\nu_x$ is $U_x$-invariant and $\nu_x$ is supported on a single $U_x$-orbit.  
\end{enumerate}


\subsection{Bounds on (co)norms from  finiteness of measures}In the proof of \cref{prop:msrrigi2},  we need to obtain controlled distortion on the norm and conorm of the linearized dynamics.  When the fiberwise  coarse Lyapunov subspaces are 1-dimensional, this is accomplished via a Lusin argument.  Here, using the hypothesis that the measure $\mu$ is $\beta$-tame, we similarly obtain controlled divergence of the norm and conorm of the linearized fiberwise dynamics.

We first consider the case that $\nu_x$ are finite and preserved by the cocycle and show the norms are well controlled on given by returns to sets of large measure.  
\begin{proposition} \label{controllednorm}
Suppose for each $x\in X$ that $\nu_x$ is a finite  measure and that $(L_x)_*\nu_x = \nu_{f(x)}$ a.e.

Then, for every $\delta>0$ there exists a $C>1$ and a subset $K\subset X$ with $\mu(K)>1-\delta$ such that for all $x\in K$ and $n $ with $f^n(x)\in K$,
$$\dfrac{1}{C} \le m\left( {\wtd L_x^{(n)}}\right)\le \| {\wtd L_x^{(n)}}\| \le C.$$
\end{proposition}

\begin{proof}
We may  assume each   measure $\nu_x$  is a probability measure.

Consider  $x\in X$.  By definition of the essential linear support of $\nu_x$, given any proper subspace $U\subset W_x$ we have $$  \nu_x(U) < 1.$$
Let $\calG_x$ denote the disjoint union of all Grassmanians in $W_x$.   The set $\calG_x$ parametrizes the space of all positive dimensional subspaces of $W_x$ and the map $\calG_x\to [0, 1]$ given by  $$U\mapsto   \nu_x(U )$$ is upper semicontinuous.  Since $\calG_x\sm \{W_x\}$ is compact, there exists $\epsilon_x>0$ such that 
$$   \nu_x(W_x\sm U)\ge \epsilon_x$$
for every proper subspace $U\subset W_x$ of positive dimension.

Given a subspace $U\subset W_x$ and $r>0$, let $B_r(U)$ denote the open $r$-ball (i.e. tubular neighborhood) around $U\subset W_x$. 
 Fix $Q>1$. 
  Given $M>1$, the map $\calG_x\to [0, 1]$,   $$U\mapsto   \nu_x \left(\overline{B_{\frac 1 M}(U)}\right),$$ is upper semicontinuous.   As  $\calG_x\sm \{W_x\}$ is compact, 
   there exists $M_x$ such that  
$$  \nu_x \left(\overline{B_{\frac 1 {M_x}}(U)}\sm U\right) \le \frac{\epsilon_x}{Q}$$
for every proper subspace $U\subset W_x$.
Similarly, we may select  $M_x$   above so that additionally 
$$   \nu_x (W_x\sm {B_{ M_x}(U)}) \le \frac{\epsilon_x}{Q}$$
for every proper subspace $U\subset W_x$.

Given $\epsilon_0>0$, let $K= K(\epsilon _0 , Q)\subset X$ denote the set of $x\in X$ such that 
    \begin{enumerate}
    \item $\frac{Q-1}{Q}\epsilon _x> \epsilon _0$;
    \item $\epsilon_0> \frac{1}{Q} \epsilon_x$;
    \item $M_x\le Q$.
    \end{enumerate}
Given $\delta>0$, we may select $Q>1$ sufficiently large and $\epsilon_0>0$ sufficiently small to ensure $\mu(K)> 1-\delta$.  We  claim the proposition holds for such $K$ and $\delta$ with $C = Q^2$.
 
Consider $x\in K$ and $n$ with $f^n( x)\in K$.  If $\|{{\wtd L_x^{(n)}}}\| = \lambda$ then, by the singular value decomposition,  there exists a proper subspace $U\subset W_x $ such that for all $v\in W_x$,
$$d(L_x^{(n)} v, L_x^{(n)}U) = \lambda d(v, U) .$$
If $\lambda\ge  Q^2$ then
\begin{align*}
 B_{M_{ f^n(x)}}(   L_x^{(n)} U) &\subset 
 B_{Q}(   L_x^{(n)} U) \\
 &\subset 
 B_{\frac{ Q^2}{M_{ x}}}(   L_x^{(n)} U) \\
 &\subset  B_{\frac{ \lambda}{M_{ x}}}(   L_x^{(n)} U) \\
&= L_x^{(n)} B_{\frac{ 1}{M_{ x}}}(   U) .
%
\end{align*}
By equivariance of the family of measures $\{  \nu_x\}$, it follows that 
$$  \nu_{f^n(x)} ( B_{M_{ f^n(x)}}(   L_x^{(n)} U)) \le    \nu_x(B_{\frac{ 1}{M_{ x}}}(   U) ) \le \frac {\epsilon_x}{Q}<\epsilon_0.$$
On the other hand, we have assumed that 
$$  \nu_{f^n(x)} (B_{{M_{f^n(x)}}}( L_x^{(n)} U))\ge \frac{Q-1}{Q} \epsilon_{f^n(x)} >\epsilon_0,$$
a contradiction.  
It follows that $$\|{{\wtd L_x^{(n)}}}\|  = \lambda\le Q^2.$$  
We similarly obtain that  $m({\wtd L_x^{(n)}})\ge  \frac{1}{Q^2}.$
\end{proof}

\subsection{Rigidity of projective stabilizers}
We now assume we are in the setting where the cocycle $x\mapsto L_x$ is assumed bi-regular for $\mu$-a.e.\ $x$.  
We let 
 $$\lambda_\top(x) = \lim_{n\to \infty } \frac 1 n \|\wtd L_x^{(n)}\|,\quad \quad  \lambda_\btm(x) = \lim_{n\to \infty } \frac 1 n m(\wtd L_x^{(n)})$$ 
denote the top and bottom pointwise Lyapunov exponents at $x\in X$.  Note that the limits exist almost everywhere  and if $\mu$ is $f$-ergodic then $\lambda_\top(x)$ and $\lambda_\btm(x)$ are a.s.\ constant in $x$.

We   consider the case that $$\lambda_\top(x)<0$$ for almost every $x$.

We recall that subgroups of $\Aff(W_x)$ may be viewed as a linear subgroup of $W_x\oplus \R$ by identifying the affine map $x\mapsto Ax + w$ with the linear map $(x,t) \mapsto (Ax + tw, t).$

Given $v\in W_x$, we have the Lyapunov weight   $$\varpi_x(v) = \limsup_{n\to \infty} \frac 1 n \log \|\wtd L_x^{(n)}(v)\|.$$
The range of $\varpi_x$ is a.s.\ $f$-invariant and contained in $${-\infty}\cup [\lambda_\btm(x), \lambda_\top(x)].$$
For every bi-regular point $x\in X$, there exists a direct sum decomposition $$W_x= \bigoplus_{i=1}^\ell E^i(x)$$ with the following properties:
\begin{enumerate}
	\item the assignment $x\mapsto E^i(x)$ is measurable;
	\item $  L_xE^i(x) = E^i(f(x))$;
	\item for $v\in E^i(x)\sm \{0\}$, $\lim_{n\to \pm\infty}  \frac 1 n  \log \|  L_x^{(n)}(v)\| \to  \varpi_x(v).$
\end{enumerate}
We note that the flag in $W_x$ associated to $\varpi_x$ consists of the subspaces of the form $\bigoplus _{j=1}^i E^j(x)$.

Given $x\in X$, write $ \Aff^{SR}(W_x)$ for the set of subresonant (w.r.t.\ $\varpi_x$) affine  maps $g\colon W_x\to W_x$.
 Below, we consider the subgroup $G\subset  \Aff^{SR}(W_x)$ such that $g_*\nu_x \propto \nu_x.$
It is standard to verify that $G$ is a closed subgroup of $\Aff(W_x)$; in particular, $G$ is a Lie group (though need not necessarily be algebraic.)

Given a Lie group $H$,  write $H^\circ$ for the identity component of $H$.

\begin{proposition}\label{prop:measurerig}
Suppose    for almost every $x\in X$ that  $\lambda_\top(x)<0$ and that there exists a closed subgroup $H_x\subset  \Aff^{SR}(W_x)$ such that the following hold:
\begin{enumerate}
\item $g_*\nu_x \propto \nu_x$ for every $g\in H_x$; 
\item the assignment $x\mapsto \mathrm{Lie}(H_x)$ is measurable and $L_x H_x L_x\inv = H_{f_x}$;
\item the measure $\nu_x$ is supported on the $H_x$-orbit of $0$.
\end{enumerate}
Then, for almost every $x\in X$ there is a connected, closed, unipotent subgroup $U_x\subset H_x^\circ$ which is normal in $H_x$ such that 
\begin{enumerate}
\item $g_*\nu_x= \nu_x$ for every $g\in U_x$;  
\item $x\mapsto \Lie(U_x)$ is measurable;
\item $L_x U_x L_x\inv = U_{f_x}$;
\item  $\nu_x$ is supported on the $U_x$-orbit of $0$.  
\end{enumerate}
\end{proposition}
We note that $U_x$ is nilpotent, and hence unimodular, as are all of its connected subgroups. Write $S_x$ for the stabilizer of $0$ in $U_x$.  Then there exists a $U_x$-invariant measure on the homogeneous space $U_x/S_x$ which is unique up to normalization.    \cref{prop:measurerig} then implies that $\nu_x$ coincides (up to a choice of normalization) with this $U_x$-invariant measure on $U_x/S_x$.

\subsection{Proof of \cref{prop:measurerig}} 
Without loss of generality, we may assume the base dynamics $f\colon (X,\mu) \to (X, \mu)$ is ergodic.   In particular, the Lyapunov exponents of the cocycle induced by $\wtd L_x\colon W_x  \to W_{f(x)}$ may be assumed to be a.s.\ constant and ordered so that $$\lambda_1<\dots <\lambda_\ell<0.$$
Write $W_x =\bigoplus_{i=1}^\ell E^i_x$ for the associated Oseledec's splitting.  

\def\stab{\mathrm{stab}}
\def\liej{\mathfrak {j}}

\subsubsection{The induced cocycle on $\aff(W_x )$}
Let $\Aff_x= \Aff(W_x)$ denote the group of affine transformations of $W_x$.  Let $\aff_x$ denote the Lie algebra of $\Aff_x$.  We make the standard 
identifications $\Aff_x\simeq \Gl(W_x)\ltimes W_x$  and  $\aff_x\simeq \gl(W_x)\ltimes W_x$ and view 
$\aff_x$ as a subalgebra of $\gl(W_x\oplus \R)$.

Write  $\tau_x\colon \Aff_x\to \Aff_{f(x)}$ for the linear map induced by conjugation: for $g\in \Aff_x$, $$\tau_x\colon g \mapsto (\wtd L_x\oplus \Id)\circ g\circ (\wtd L_x\oplus \Id)\inv.$$
The map $x\mapsto d\tau_x$ defines a linear cocycle over the dynamics  $f\colon X\to X$.
For each bi-regular $x\in X$, let 
$\liee^{i,j}_x\subset \gl(W_x )$ denote the subspace of linear maps $T\colon W_x \to W_x $ with codomain $E^j_x$ and 
 kernel $\bigoplus _{k\neq i} E^k_x$.  
Since the decomposition $W_x  = \bigoplus_{i=1}^\ell E^i_x$ is $L_x$-equivariant, it follows that the decomposition 
 $$\aff_x =  \bigoplus _{1\le i,j\le \ell} \liee^{i,j}_x \oplus \bigoplus _{1\le i\le \ell} E^i_x$$
is $d\tau_x$-equivariant.

Let $\kappa_p$, $1\le p\le r$, denote the Lyapunov exponents of the linear cocycle $d\tau_x$ and let $$\aff_x = \bigoplus _{p=1}^r F^{\kappa_p}_x$$ denote the associated Oseledec's decomposition.  
Then the Lyapunov exponents $\kappa_p$ of $d\tau_x$  are either:
\begin{enumerate}\item  $\lambda_i$ (corresponding to vectors tangent to  $E^i_x$),  or 
\item $\lambda_j- \lambda_i$ (corresponding to vectors tangent to $\liee^{i,j}_x$).  
\end{enumerate}
In particular, 
$\lambda_1\le \kappa_p\le  -\lambda_1$ for every $p$.  
Moreover,  the Oseledec's decomposition of $\aff_x$ consists of the subspaces of the form 
 $$F^{\kappa_p}_x = \bigoplus_{\lambda_j = \kappa_p} E^j_x\oplus \bigoplus _{\lambda_j- \lambda_i = \kappa_p} \liee^{i,j}_x.$$
 
 We note if $Y\in F^\kappa_x$ and $Z\in F^{\kappa'}_x$ then the Lie bracket satisfies \begin{equation}\label{eq:liebra}[Y,Z]\in F^{\kappa+\kappa'}_x.\end{equation}

Let $\aff^{\le t}_x=\bigoplus_{\kappa_p\le t }F^{\kappa_p}_x$.  We note that $\aff^{< 0}_x$ and $\aff^{\le 0}_x$ are Lie subalgebras.  Let $\Aff^{< 0}_x$ (resp.\ $\Aff^{\le 0}_x$) be the maximal closed subgroup of $\Aff_x$ tangent to  $\aff^{< 0}_x$ (resp.\ $\aff^{\le 0}_x$).  We have that   $\Aff^{< 0}_x$ is normal in  $\Aff^{\le  0}_x$.

\subsubsection{Dynamical characterization of $U_x$}
Let $\lieh_x$ denote the Lie algebra of $H_x$.  By assumption, the assignment $x\mapsto \lieh_x$ is measurable and $d\tau_x$ equivariant.   Thus $\lieh_x$ decomposes according to the Oseldec's decomposition of $d\tau_x$: $$\lieh_x = \bigoplus_{p=1} ^r (\lieh_x\cap F^{\kappa_p}_x).$$
By assumption,  $\lieh_x$ consists only of subresonant affine maps; thus we actually  have  $$\lieh_x = \bigoplus_{\kappa_p\le 0 }  (\lieh_x\cap F^{\kappa_p}_x).$$
Let $$\lieu_x = \bigoplus_{\kappa_p< 0 }  (\lieh_x\cap F^{\kappa_p}_x).$$

Let $U_x$ be the analytic subgroup of $H_x$ tangent to $\lieu_x$.  Then $U_x$ is  a connected, unipotent,  subgroup of in $H_x$.  Moreover, 
we have 
\begin{enumerate}
\item $\Aff^{< 0}_x$ is normal in $\Aff^{<\le 0}_x$;
\item $H_x\subset \Aff^{\le 0}_x$;
\item  $U_x= (H_x\cap \Aff^{< 0}_x)^\circ$.
\end{enumerate}
whence $U_x$ is normal in $H_x$.  
 Moreover, we have $d\tau_x U_x = U_{f_x}$.  To complete the proof of  \cref{prop:measurerig}, it remains to show $g_*\nu_x= \nu_x$ for every $g\in U_x$ and that $\nu_x$ is supported on the $U_x$-orbit $U_x\cdot 0$.  


 \subsubsection {The orbit $U_x\cdot 0$ has positive measure}
We show that $H_x^\circ \cdot 0= U_x\cdot 0$.  Indeed, the stabilizer  $S_x$ of $0$ in $H_x^\circ$ consists only of linear maps of $W_x$.  Since $\gl(W_x)$ contains the 0-Lypaunov exponent  subspace $F^0_x\subset \aff(W_x)$  for the cocycle $d\tau_x$, the Lie algebra of $S_x$ contains 
the subspace $$(\lieh_x\cap F^{0}_x).$$
It follows that $\lieh_x = \lieu_x\oplus(\lieh_x\cap F^{0}_x)$.   
It follows that $U_x$ acts locally transitively on the $H_x^\circ$-orbit of $0$.  Since we assume the $H_x$-orbit of $0$ has positive measure, it follows that $U_x\cdot 0$ has positive measure.    

 \subsubsection {Each measure $\nu_x$ is  $U_x$-invariant}
We next claim that $\nu_x$ is invariant under $U_x$ for almost every $x\in \Lambda$.   Recall that $\nu_x$ is assumed to be a projectively invariant measure on the $H_x$-homogeneous space $H_x\cdot 0$.  In particular, there is a continuous homomorphism $c_x\colon H_x\to \R_+$ such that $$g_* \nu_x = c_x(g) \nu_x$$ for every $g\in H_x$.  
By the dynamical characterization of $\lieu_x$, we obtain the following. 
\begin{claim}
For almost every $x\in X$ and every $g\in U_x$ we have $c_x(g) = 1$. 
\end{claim}
\begin{proof}
Given $Z\in \lieu_x$,  let $c_x(Z)= c_x(\exp(Z))$.  We observe that $c_x(Z) = c_{f(x)}(\tau_x(Z))$.  
Since $c_x\colon \lieu_x\to \R_+$ is continuous and $c_x(0)=1$, for every $\epsilon>0$ there is $r=r_{x,\epsilon}>0$ such that $|c(Z)-1|< \epsilon$ for all $Z\in \lieu_x$ with $\|Z\|< r$.  
Since $d \tau_x$ contracts $\lieu_x$, by Poincar\'e recurrence to sets of $x\in X$ where $r_{x,\epsilon}$ is uniformly bounded from below, it follows for almost every $x\in X$ that $|c_x(Z)-1|< \epsilon$ for all $Z\in \lieu_x$.  
\end{proof}
\subsubsection{Common linear parametrization of $U_x$-orbits}  

Let $\scrH_x= H\cdot 0$ denote the $H$-orbit of zero.  We aim to show that $\scrH_x$ is path-connected and thus a single $U_x$-orbit: $\scrH_x= U_x\cdot 0$.  To accomplish this, we suppose $\scrH_x$ is not path-connected and deduce a contradiction (in the next sub-subsection) with  local-finiteness of $\nu_x$ and Poincar\'e recurrence.  To run this argument, we need to compare the restriction of the measure $\nu_x$ to distinct path components of $\scrH_x$; here we provide such a comparison.  

  Since $U_x$ is normal in $H_x$, given $y\in \scrH_x$, the path component of the set $\scrH_x$ containing $y=h\cdot 0$ is of the form  $$H_x^\circ \cdot y = U_x\cdot h\cdot 0 = h\cdot U_x \cdot 0.$$   Each such orbit carries a Haar measure---well defined up to normalization---given by the push forward of the Haar measure on the unimodular group $U_x$ to the orbit $h\cdot U_x\cdot 0= U_x\cdot y.$  Although the acting group $U_x$ is the same for all such orbits, the stabilizers on different orbits need not coincide 
  and thus the Haar measures on distinct orbits $U_x\cdot y$ and $U_x\cdot y'$ are not directly comparable.  
We show the Haar measures on distinct  $U_x$-orbits are comparable to the Haar-Lebesgue measure on a common linear subspace; this gives a canonical---and thus equivariant---method to compare the restriction of $\nu_x$ to disjoint $U_x$-orbits.

The following argument does not use information about the cocycle $L_x$; rather it only uses the grading on $W_x$ and the properties of $H_x$ established above (including that $H_x$ preserves the grading).  In particular, we will work at  a single generic $x\in \Omega$ and, dropping subscripts, write  $W= W_x, U=U_x, \scrH=\scrH_x, \liee^1=\liee^1_x, \nu=\nu_x$, etc.  

Let $V$ denote the tangent space at $0$ to the $U=U_x$-orbit of $0$.
\begin{lemma}\label{common param}
For every $y\in \scrH$, there exists a unique  function $\phi_y\colon V\to V^\perp$ such that  
the $U$-orbit $U\cdot y$ is the graph 
$$ \{(v,\phi_y(v)): v\in V\}$$
of  $\phi_y$.  Moreover, up to normalization, the
Haar measure on the orbit $ U\cdot y $ is the image of the  $T_V$-invariant Haar measure on $V$ under the graph of $\phi_y$.  
\end{lemma}

\begin{proof}
We induct on $\dim W$.   If $\dim W= 1$ then $U$ is either $0$-dimensional or $1$-dimensional.  If $U$ is $0$-dimensional,  every $U$-orbit is a point, and we are done.  If $U$ is 1-dimensional, $U$ consists of all translations, and thus there exists exactly one $U$-orbit in $W$.  

Note that $\liee^1$ is tangent to the group of translations in the $E^1$ direction.  Suppose $\dim W\ge 2$.  We write $\lieu^1=\lieu\cap \liee^1$ and let $U^1\subset U$ be the analytic subgroup associated to $\lieu^1$. 
Since $\liee^1$ is normalized by $\aff^{\le 0}(W)$, we have that $U^1$ is normalized by $H$.

We consider two cases:

\subsubsection*{Case 1: $\lieu^1 \neq\{0\}$}  
Since $U^1$ is normal in $H$,  for every $y\in \scrH$, the orbit $U \cdot y$ is invariant under the subgroup of translations by $U^1$.  
Let $W'=\lieu^1 = T_0(U^1\cdot 0)=U^1\cdot 0$ be the vector subspace tangent to translations in $U^1$. 
Note that the orbit $U  \cdot y$ and the subspace  $V$  are  invariant under $U^1$, the subgroup of translations tangent to $W'$.  

 Consider the quotient vector space $\wtd W= W/W'$ (equipped with the quotient weight $\varpi(y+W')= \min \{\varpi(y +y'):y'\in W'\}.$)
We view  $\Aff(\wtd W)$ as the quotient of the normalizer of $U^1$ in $\Aff(W)$ by $U^1$.  
Since $U^1$ is normal in $U$, $U$ descends to a subgroup $\wtd U$ of $\Aff(\wtd W)$.
Moreover, if $\wtd \lieu$ is the Lie algebra of $\wtd U$, then $\wtd \lieu$ is contained on $\aff^{<0}(\wtd W)$.   
  Given $y\in W$, let $\td y=y+ W'\in \wtd W$.  
 The image of the  orbit $U\cdot y$ under $W\to \wtd W$ coincides with the orbit $\wtd U\cdot \td y$.  
 Let $\wtd V= V/W'$ be the tangent space to $\wtd U\cdot 0$ at $0$.  
By the inductive hypothesis on dimension, $(\wtd U\cdot \td y)\cap (\td v+ V^\perp) \bmod W'$ is a singleton for every $\td v\in \wtd V$.  Since $U\cdot y$ and $V$ are $U_1$-invariant, it follows that  $( U\cdot  y)\cap ( v+ V^\perp )$ is a  singleton for every $ v\in  V$.  In particular, 
the $U$-orbit $U\cdot y$ is the graph 
$ \{(v,\phi_y(v)): v\in V\}$
of a function $\phi_y\colon V\to V^\perp$. 

Since $\wtd U\cdot \td y$ is $U^1$-invariant, so is the function $\phi_y$.  Let $\eta_y$ denote the image of the Haar measure on $V$  under the graph of  $\phi_y$.  Then $\eta_y$ is $U^1$-invariant.  
 Choosing a normalization of Haar measure on $U^1$, we obtain a (locally finite Radon) quotient measure $\wtd \eta_{\td y}$ on 
 $\wtd U\cdot \td y$. Write  $\wtd U\cdot \td y$ as the graph of $\phi_{\td y}\colon  (V/W')\to ( V^\perp/W')$.  
 Then $\wtd \eta_{\td y}$ is the image under $\phi_{\td y}$ of the Haar measure on $(V/W')$.  
  By the inductive hypothesis, $\phi_{\td y}$ is a $\wtd U$-invariant measure on $\wtd U\cdot \td y$.  But this implies that $\eta_y$ is a $U$-invariant measure on $  U\cdot   y$.

\subsubsection*{Case 2: $\lieu^1 =\{0\}$}  
In this case, let $W'= {E^1}= T_{E^1}\cdot 0$ denote the subspace of $W$ with lowest weight. (Here $T_{E^1}$ denotes the subgroup of $\Aff(W)$ of translations in the direction $E_1$.)    The map $W\to W/W'$ is injective on every $U$-orbit.  
Note that $T_{E^1}$ is normal in $\Aff^{\le 0}(W)$ and we identify  $\Aff^{\le 0}(W/W')$ with
 $\Aff^{\le 0}(W)/ T_{E^1}$
In particular, as $U\cap T_{E_1}= \Id$,  we have an identification of $U$ with a subgroup of $\Aff(W/W')$.  We also identify $V$ with $V+W'$ in $\wtd W= W/W'$.  
By the inductive hypothesis on dimension, for $\td y\in W/W'$,  $( U\cdot \td y)\cap ( \td v+ V^\perp/W^1)$ is a singleton for every $\td v\in \wtd V$ and thus $( U\cdot  y)\cap ( v+ V^\perp )$ is a  singleton for every $y\in W$ and every  $ v\in  V$.  In particular, 
the $U$-orbit $U\cdot y$ is the graph 
$ \{(v,\phi_y(v)): v\in V\}$
of a function $\phi_y\colon V\to V^\perp$. 

Let  $\eta_y$ be the  image of Haar measure on $T_V$ under $\phi_y$.  
Pushing forward $\eta_y$ under the identification of $U\cdot y$ yith $U\cdot \td y$,  $\eta_y$ is $U$-invariant on the orbit $U\cdot \td y$ in $W/W'$.  But this implies that $\eta_y$ is $U$-invariant on the orbit $U\cdot y$.  
\end{proof}

\subsubsection {Each measure $\nu_x$ is supported on a single  $U_x$-orbit}
For $x\in \Omega$, let  $V_x$ denote the tangent space at $0$ to the $U_x$-orbit of $0$.

We equip $V_x$ with any choice of translation-invariant locally-finite Haar-Lebesgue measure $\eta_x$.  
Normalize $\nu_x$ so that the restriction of $\nu_x$ to $U_x\cdot 0$ coincides with the image of $\eta_x$ under the graph $\phi_0$ guaranteed by \cref{common param}.  Given $y\in \scrH_x$, let $J_x(y)\in (0,\infty)$ be such that the restriction of $\nu_x$ to the orbit $U_x\cdot y$ is proportional to the image of $\eta_x$ under the graph $\phi_y$ guaranteed by \cref{common param} with Radon-Nikodym derivative $J_x(y)$:
$$\frac {d (\restrict {\nu_x}{U_x\cdot y})}{(\Id \times \phi_y)_*\eta_x}(z) = J_x(y).$$
We observe:
\begin{enumerate}
	\item $J_x(y)=J_x(z)$ for all  $z\in U_x\cdot y$;
	\item the assignment $y\mapsto J_x(y)$ is independent of the choice of Haar-Lebesgue measure $\eta_x$ on $V_x$;
	\item 
	 for $\mu$-a.e.\ $x$ and every $y\in \scrH_x$, $J_{f( x)}(L_x y) = J_x(y)$.
\end{enumerate}

To complete the proof of \cref{prop:measurerig}, we claim the following.
\begin{lemma}
For $\mu$-a.e.\ $x$ and $y\in \scrH_x$, if  $J_x(y)>0$ then $y\in U_x\cdot 0$.
\end{lemma}
\begin{proof}
Fix $\epsilon>0$. Since each measure $\nu_x$ is locally finite, for $\mu$-a.e.\ $x\in \Omega$ there is $\delta_x>0$ such that for all $y\in \scrH_x$ with $\|y\|<\delta_x$ and $J_x(y)>\epsilon$ we have $y\in U_x\cdot 0$. 

 By Poincar\'e recurrence, for $\mu$-a.e.\ $x$ there exists $\underline \delta_x$ such that for infinitely many $n_j\in \N$, 
 $$\underline \delta_x\le \delta_{f^{n_j}(x)}.$$
 Given that $\|\wtd L_x^{(n)}y\|\to 0$ as $n\to \infty$ for almost every $x$ and every $y\in W_x$, for almost every $x$,  if $J_x(y)> \epsilon$ for some $y\in \scrH_x$, then $y\in U_x\cdot 0$.  
 Since this holds for all $\epsilon>0$, the conclusion follows.  
\end{proof}

\section{Proof of \cref{prop:msrrigi2}}  \label{sed:mrse}

Let $\mu$ be an ergodic, $A$-invariant Borel probability measure on $M^\alpha$ as in \cref{prop:msrrigi2}.  In particular, we assume the image of $\mu$ in $G/\Gamma$ is the Haar measure.  
Fix $\beta\in \Sigma$ such that $\mu$ is 
not  invariant under any 1-parameter subgroup of $U^{[\beta]}$ 
and satisfies the hypotheses of \cref{prop:msrrigi2}.  

Write $A'=A'_\beta$ for the kernel of $\beta\colon A\to \R$ in $A$.  Also write  $\calE=\calE_\beta$ for the  decomposition of $\mu$ into $A'_\beta= \ker\beta$-ergodic components. 
Recall that we write $\{\mu^{{[\beta]},\calE }_x \}$ for a family of leafwise measures relative to the partition $\calE=\calE_\beta$ along $W^{[\beta]}$-leaves.

Our main result in this section is the following.  We slightly relax the hypotheses of  \cref{prop:msrrigi2} on obtain the weaker result that the measure  $\{\mu^{{[\beta]},\calE }_x \}$ is smooth; we then apply \cref{graph or invariant} to complete  \cref{prop:msrrigi2}.  
\begin{proposition}\label{prop:msrrig3}
Let $\Gamma$ be an irreducible lattice in a higher-rank Lie group $G$ and  let $\alpha\colon \Gamma\to \diff^r(M)$ be an action on a compact manifold.   Let $\mu$ be an ergodic, $A$-invariant Borel probability measure on $M^\alpha$ projecting to the Haar measure on $G/\Gamma$. Assume that $r>c_0^F(\mu)$.

Fix $\beta\in \Sigma$ such that $\mu$ is not invariant under any $1$-parameter subgroup of $U^{[\beta]}$.   
Suppose  that one of the following conditions hold  for a.e.\ $x$:
\begin{enumerate}
\item  the  measure $\mu$ is $\beta$-tame, or 
\item $\dim E^{[\beta],F}(x) =1$.
 \end{enumerate}  


Then for $\mu$-a.e.\ $x$, the  leafwise measure $\mu^{{[\beta]},\calE}_x$  is in the Lebesgue measure class on a connected, embedded submanifold $N_x\subset W^{[\beta]}(x)$.  
\end{proposition}

\cref{prop:msrrig3} combined with  \cref{graph or invariant}  
 then establishes \cref{prop:msrrigi2}.

\subsection{Quasi-isometric properties given 1-dimensional fiberwise coarse Lyapunov manifolds.}
As above, let $\mu$ be an $A$-invariant Borel probability measure on $M^\alpha$ projecting to the Haar measure in $G/\Gamma$.  Fix a root $\beta$ and let $A_\beta'$ denote the kernel of $\beta$ in $A$.  Fix an $A'_\beta$-ergodic component $\mu'$ of $\mu$ that projects to the Haar measure on $G/\Gamma$.  Fix any $a_0\in A'_\beta$. 

A proof of the following appears in {\cite[Lemma 5.7]{ABZ}}.  We also include a proof for completeness.   
\begin{lemma}[]\label{lem:lusin}  
Suppose that $\dim E^{[\beta],F}(x) =1$ and 
that $\mu$ is not invariant under any $1$-parameter subgroup of $U^{[\beta]}$.

  Then for every $\delta>0$ there exists $C_\delta>1$ and  a subset $K\subset M^\alpha$ with $\mu'(K)>1-\delta$ such that for every $x\in K$ and every $n\in \Z$ with $a_0^n\cdot x\in K$,
\begin{equation}
\frac 1 {C_\delta}\le \|\restrict{D_xa_0^n}{E^{[\beta],F}}\|\le C_\delta.
\end{equation}
\end{lemma}

\begin{corollary}\label{eq:highdensity}
Suppose  $a_0$  acts ergodically on $(M^\alpha, \mu')$.  Then for $\mu'$-a.e.\ $x$ and every $\delta>0$ there is $C_{x,\delta}\ge 1$ such that 
\begin{equation}
\liminf_{N\to \infty} \frac{1 }{N} \#\left\{0\le n\le N-1:
\frac 1 {C_{x,\delta}} \le \|\restrict{D_xa_0^n}{E^{[\beta],F}}\|\le C_{x,\delta}
\right\}\ge 1-\delta.
\end{equation}
\end{corollary}


\begin{proof}[Proof of \cref{lem:lusin}]
Fix $\lieu= \lieu^{[\beta]}$ and  $U= U^{[\beta]}$.  
Let $h_x\colon E^{[\beta],F}(x)\to W^{[\beta],F}(x)$ be the normal form coordinates as in \cref{1dNF}.  Taking  a measurable choice of orientation on $E^{[\beta],F}(x)$, we have a measurably varying, isometric linear identification $\R\to E^{[\beta],F}(x)$ for almost every $x$.  Let $\td h_x\colon \R \to W^{[\beta],F}(x)$ denote the composition of these maps.   

By \cref{prop:msrrigiU}, for $\mu$-a.e.\ $x$, the leafwise measure along $U^{[\beta]}$-orbits, $\mu_x^{U^{[\beta]}, \calE_\beta}$, is a single atom $\mu_x^{U^{[\beta]}, \calE_\beta}=\delta_x$.  
Let $K_0\subset M^\alpha$ be a compact subset with $\mu'(K_0)>1-\delta$ such that for $x\in K_0$, 
\begin{enumerate}
\item $\mu_x^{U^{[\beta]}, \calE_\beta}$ is defined, 
\item $\mu_x^{U^{[\beta]}, \calE_\beta}=\delta_x$, 
\item the map $ \td h_x$ is defined.  
\end{enumerate}

Let $\Phi_x\colon U \times \R \to W^{[\beta]}(x)$ be the parametrization,  $$\Phi_x(u,v) = u\cdot \td  h_x(v).$$
For $x\in K_0$, let $\eta_x=(\Phi_x\inv)_*\left(\mu_x^{{[\beta]}, \calE_\beta} \right)$
be the pulled-back measure on these parameters.  
Given $x\in K_0$, let $K_x= \Phi_x\inv(K_0)$.  
We have $K_x\cap U\times \{t\} $ is either empty or a singleton for every $t \in \R$

Let $D_r$ denote the ball centered at $\1_U$ in $U$ of radius $r>0$ with respect to some choice of right invariant metric.  
We set $\wtd D=D_2$ and $ D= D_1$.  
By compactness of $K_0$, and using that $K_x\cap U\times \{t\} $ is at most a singleton for every $t\in \R$, for every $y\in K_0$, every $L> 0$, and every $\theta>0$, there is $\gamma>0$ such that 
\begin{equation}
(D_L\times [-\gamma, \gamma])\cap K_y= (D_\theta\times [-\gamma, \gamma])\cap K_y.
\end{equation} 
In particular, for $y\in K_0$ there is $\bar \gamma_y$ such that for all $\gamma<\bar\gamma_y$,  
\begin{equation}\label{compacthelps}
 \bigl(D_3\times [- \gamma,  \gamma] \bigr)\cap K_y =   \bigl(D\times [- \gamma,  \gamma] \bigr)\cap K_y.
\end{equation}

Fix $x\in K_0$ and $y\in \Phi_x (U\times \R)\cap K_0$.  Let $(u_y, v_y) = \Phi_x\inv (y)\in K_x$. 
The change in normal form coordinates is affine in the $\R$-coordinates and the change in $\Phi_x$ coordinates coincides with right translation in the $U$-coordinate; that is 
\begin{equation*}
\Phi_{y}\inv\circ \Phi_x(u, t) = (uu_y\inv, a(t-v_y))
\end{equation*}
for some $a\in \R$.  In particular, for every $\gamma>0$, there exists $\td \gamma>0$ such that 
\begin{equation}\label{NF coord change}
\Phi_y(D\times [-\gamma,\gamma]) = 
\Phi_x\left(
\{
(u\cdot u_y, v_y + t): u\in D, t\in [-\td \gamma, \td \gamma]
\}
\right).
\end{equation}

Fix $\gamma>0$ and let 
$$K_{x,\gamma}=  K_x\cap \bigl(D\times [-\gamma,\gamma]\bigr)$$
and 
$$c_\gamma(x) =  \frac {\eta_x(K_{x,\gamma} )}{\eta_x( D \times [-\gamma,\gamma])}.$$
For $\mu_x^{{[\beta]}, \calE_\beta}$-a.e.\ $y$, the leafwise measures $\mu_x^{{[\beta]}, \calE_\beta}$ and $\mu_y^{{[\beta]}, \calE_\beta}$ coincide (up to choice of normalization.)  Thus for $\eta_x$-a.e.\ $(u_0,v_0)\in K_x\cap (D\times \R)$ with $y= \Phi_x(u_0,v_0)$ and any $\gamma>0$ sufficiently small so that \eqref{compacthelps} holds, there is $\td \gamma>0$ (with $\td \gamma\to 0$ as $\gamma\to 0$) for which 
\begin{align}
c_\gamma(y)&:= 
 \frac {\eta_y(K_{y,\gamma} )}{\eta_y( D \times [-\gamma,\gamma])}  \notag
 \\
&= \frac {\eta_x\bigl(K_{x}\cap \{
(u\cdot u_0, v_0 + t): u\in D, t\in [-\td \gamma, \td \gamma]
\}\bigr)
}
{\eta_x\bigl(\{
(u\cdot u_0, v_0 + t): u\in D, t\in [-\td \gamma, \td \gamma]
\}\bigr)
}
&&\hfill \text{ \parbox[r]{8 em}{\scriptsize  by \eqref{NF coord change} and coincidence of   leafwise measures}}\notag
\\
&= \frac {\eta_x\bigl(K_{x}\cap \bigl(\wtd D\times [v_0-\td \gamma, v_0+\td \gamma]\bigr)
\bigr)
}
{\eta_x\bigl(\{
(u\cdot u_0, v_0 + t): u\in D, t\in [-\td \gamma, \td \gamma]
\}\bigr)
}&&\hfill \text{ \parbox[r]{8 em}{\scriptsize by \eqref{compacthelps}}} \notag  \\
&\ge  \frac {\eta_x\bigl(K_{x}\cap \bigl(\wtd D\times [v_0-\td \gamma, v_0+\td \gamma]\bigr)
\bigr)
}
{\eta_x\bigl(\wtd D\times [v_0-\td \gamma, v_0+\td \gamma]
\bigr)
}.\label{lowerbound}
\end{align}
The natural map $\wtd D\times \R\to \R$ is proper and is measurably one-to-one.  Equip $\R$ with the image of the restriction of  $\eta_x$ to $\wtd D\times \R$ under this map; this measure on $\R$ is a  Radon and thus satisfies a density theorem relative to intervals in $\R$.  
Almost every point (with respect to the projection of $\restrict{\eta_x}{\wtd D\times \R}$) of the image of $K_x\cap (\wtd D\times \R)$ in $\R$ is a  density point.  Moreover, the projection $\wtd D\times \R\to \R$ is a measurable isomorphism.  Thus, typically the lower bound \eqref{lowerbound} approaches 1 as $\td \gamma\to 0$.   It follows a.e.\ $x$ and  for $\mu_x^{{[\beta]}, \calE_\beta}$-a.e.\ $y\in K_0\cap \Phi_x(D\times \R)$ that 
$$\lim_{\gamma\to 0} c_\gamma(y)=1.$$

We now select the set $K$ in the proposition.  
Fix $\theta \in (0,1)$ and $\gamma\in (0,1)$.  Let $K\subset K_0$ be a measurable subset such that for $y\in K$,
\begin{equation}\label{bigout}
\frac{\eta_y\bigl( D_{\theta}\times [-1,1]\bigr) }{\eta_y\bigl(
D\times [-1,1]
\bigr)}\le .5
\end{equation}
and furthermore,
\begin{enumerate}
\item $c_{\gamma'}(y)\ge .95$ for all $0<\gamma'<\gamma$
\item  $K_{y,\gamma}\subset D_{\theta}\times[-\gamma,\gamma]$.
\end{enumerate}
It follows that if $y\in K$ then for all $0<\gamma'<\gamma$
\begin{equation}\label{smallin} \frac{\eta_y\bigl( D_{\theta}\times [-\gamma',\gamma']\bigr) }{\eta_y\bigl(
D\times [-\gamma',\gamma']
\bigr)}\ge .95\end{equation}

Taking first  $\theta>0$ sufficiently small and then $\gamma>0$ sufficiently small, we may ensure $\mu'(K)>1-\delta$.   

Suppose $x\in K$ and $y = a_0^n(x)\in K$ for some $n\in \Z$.  
Let $$L= \|\restrict{D_xa_0^n}{E^{[\beta],F}}\|.$$  
We have 
$$\Phi_{y}\inv \left( a_0^n\cdot \Phi_x (u,v)\right)
= \bigl(u, \pm Lv\bigr).$$
By equivariance of the measures  $x\mapsto \mu_x^{U^{[\beta]}, \calE_\beta}$, we have
$$\frac{\eta_x\bigl( D_{\theta}\times [-1,1]\bigr) }{\eta_x\bigl(
D\times [-1,1]
\bigr)}
=\frac{\eta_y\bigl( D_{\theta}\times [-L,L]\bigr) }{\eta_y\bigl(
D\times [-L,L]
\bigr)}$$
and 
$$\frac{\eta_x\bigl( D_{\theta}\times [-L\inv,L\inv]\bigr) }{\eta_x\bigl(
D\times [-L\inv,L\inv]
\bigr)}
=\frac{\eta_y\bigl( D_{\theta}\times [-1,1]\bigr) }{\eta_y\bigl(
D\times [-1,1]
\bigr)}.$$
Combined with \eqref{bigout} and \eqref{smallin}, this implies 
\[\gamma\le L \le \gamma\inv.\hfill \qedhere\]
\end{proof}

\subsection{Reformulation, notation, and proof of {\cref{prop:msrrig3}}}\label{actualassem}
We take $\calI=[\beta]$, the totally resonant, strongly integrable collection of Lyapunov exponents proportional to $\beta$.  
We have that $  c_0^F(\mu)<r$ by assumption (where $c_0^F(\mu)$ is as in \cref{coarse fiber rigidity}.) 
We equip each $E^{[\beta],F}(x)$ with the Lyapunov weight $\varpi$ (induced by the action of any $a$ in the half-space $\calC([\beta])$). 

 
We follow the notation in \cref{prop:linearization}. 
Write  $\bfV^{[\beta]}(x):= \lieu\oplus \bfV^{[\beta],F}(x)$.  Let $\iota^{[\beta]}_x\colon W^{[\beta]}(x)\to \bfV^{[\beta]}(x)$ be given as follows: given $Y\in \lieu^{[\beta]}$ and $z\in W^{[\beta],F}(x)$, let
$$\iota_x^{[\beta]}(\exp Y\cdot  z) = \left (Y, \iota^{[\beta],F}_{x}(z)\right).$$
Write
$$\boldsymbol {\mu}^{[\beta],\calE}_x=(\iota_x^{[\beta]})_*\mu^{[\beta],\calE}_x$$ for the locally finite Radon measure on $\bfV^{[\beta]}(x)$.  
Let $W_x\subset \bfV^{[\beta]}(x)$ denote the essential linear span of $\boldsymbol {\mu}^{[\beta],\calE}_x$.  

We equip $\bfV^{[\beta]}(x)$ and $W_x$ with the associated Lyapunov weights (induced by our choice of  $a$ in the halfspace $\calC([\beta])$).  
Let $\Aff^{SR}(W_x)$ denote the subresonant affine maps of $W_x$ relative to this weight.

The main result of this section is the following (whose proof we postpone).
\begin{proposition}\label{prop:MR}
For $\mu$-a.e.\ $x$, and  $\mu^{[\beta], \calE}_x$-a.e.\ $ y\in W^{[\beta]}(x)$, there is an invertible subresonant affine map  $g\in \Aff^{SR}(W_x)$ 
such that 
\begin{enumerate}
\item$ g(0) = \iota_x(y)$ and 
\item $g_*\boldsymbol {\mu}^{[\beta],\calE}_x\propto \boldsymbol {\mu}^{[\beta],\calE}_x$
\end{enumerate}
\end{proposition}

\begin{proof}[Proof of {\cref{prop:msrrig3}}]
Let $G_x$ denote the subgroup of $ \Aff^{SR}(W_x)$,  such that $g_* \boldsymbol {\mu}^{[\beta],\calE}_x\propto \boldsymbol {\mu}^{[\beta],\calE}_x$ for all $g\in G_x$.  
We have that $G_x$ is closed; by \cref{prop:MR} the measure $\boldsymbol {\mu}^{[\beta],\calE}_x$ is supported on the orbit $G_x$-orbit of $0$.  Since $x\mapsto \boldsymbol {\mu}^{[\beta],\calE}_x$ is measurable and equivariant under the cocycle, so is the assignment $x\mapsto G_x$.  

By \cref{prop:measurerig}, (applied to cocycle $(\bfL a)$ for any choice of $a$ in the half-space $\calC([\beta])$, the measure $\boldsymbol {\mu}^{[\beta],\calE}_x$ is the (image  of the) Haar measure on the orbit $U_x\cdot 0$ of a closed, connected unipotent subgroup $U_x\subset G_x$.  Since $\iota^{[\beta]}_x\colon W^{[\beta]}(x)\to \bfV^{[\beta]}$ is a smooth injective embedding and since 
$\boldsymbol {\mu}^{[\beta],\calE}_x$ is the Haar measure on the orbit $U_x\cdot 0$, it follows that 
$ {\mu}^{[\beta],\calE}_x$ is in the  Lebesgue measure class on a connected, embedded submanifold $N_x\subset W^{[\beta]}(x)$.  
\end{proof}

\subsection{Proof of \cref{prop:MR}}
It remains to establish \cref{prop:MR}.  
Given $\mu$-a.e.\ $x,y\in M^\alpha$, write $L^{SR}(W_x,W_y)$ for the set of subresonant linear maps from $W_x$ to  $W_y$.

Given $x,y\in \Omega$ with $y\in W^{[\beta]}(x)$, write $y = u\cdot z$ for $z\in W^{[\beta],F}(x)$.  Let $\bfT_{y\to x}\colon \bfV^{[\beta]}(y) \to \bfV^{[\beta]}(x)$ be 
$$\bfT_{y\to x}(Y, Z) = \left(\exp\inv (  \exp Y \cdot u ), \bfT^F_{y\to x}(Z)
\right)$$
where $\bfT^F_{y\to x}$ is as in \cref{prop:linearization}\eqref{linearizedNF:6}.  It follows from \cref{prop:linearization}\eqref{linearizedNF:6} and \cref{rem:afinetrans} that $\bfT_{y\to x}\colon  \bfV^{[\beta]}(y) \to \bfV^{[\beta]}(x)$ is a subresonant affine map.  Moreover, 
\begin{equation}\bfT_{y\to x} \circ \iota_y^{[\beta]}=  \iota_x^{[\beta]}.\end{equation}
Indeed, for $Y\in \lieu$ and $w\in W^{[\beta],F}(y)$, 
\begin{align*}
\bfT_{y\to x} \circ \iota_y^{[\beta]}(\exp Y \cdot  w) 
&=
\left( \exp\inv (\exp Y\cdot u), \bfT_{y\to x}^F \circ \iota_y^{[\beta],F}( w) \right)\\
&= \left( \exp\inv ( \exp Y\cdot u ),\iota_x^{[\beta],F}(u\inv  w) \right)\\
&=\iota_x^{[\beta]}(\exp Y\cdot u  \cdot  u\inv w) \\
&=\iota_x^{[\beta]}(\exp Y\cdot w). 
\end{align*}
It follows that 
\begin{equation}\label{eq:afftine}
(\bfT_{y\to x})_*\boldsymbol {\mu}^{[\beta],\calE}_y\propto \boldsymbol {\mu}^{[\beta],\calE}_x.\end{equation}
and thus $\bfT_{y\to x} $ restricts to an invertible, subresonant affine map $T^W_{y\to x} \colon W_y \to  W_x.$  Moreover $T^W_{y\to x}(0) = \iota_x^{[\beta]}(y)$.

Fix a measurable basis for the fibers of the bundle $\bfV^{[\beta]}\to M^\alpha$.  Relative to this framing the following assignments are measurable:
\begin{enumerate}
\item the subspaces $x\mapsto W_x$;
\item the assignment of the flag on $\bfV^{[\beta]}_x$ induced by the weight $\varpi$;  
\item $x\mapsto \boldsymbol {\mu}^{[\beta],\calE}_x$ with the topology dual to compactly supported functions.  
\end{enumerate}
By Lusin's theorem, for each $M\in \N$, let $K_M\subset M^\alpha$ be a compact  set on which the above measurable assignments are continuous and $\mu(K_M)>1-\frac 1 M.$
Let $\mathscr U= \{U_j\}$ be a countable basis for the topology on $M^\alpha$ and let $\mathscr A$ be the algebra generated by $\mathscr U$.  Since $\mathscr A$ contains only finite unions and intersections of $\mathscr U$, $\mathscr A$ is also countable.

Let $A'_\beta$ denote the kernel of $\beta$ in $A$.  
We work with an $A'_\beta$-ergodic component $\mu'$ of $\mu$.  Fix $a_0\in A'_\beta$ that acts ergodically with respect to $\mu'$.  
Given $x\in M^\alpha$, write  $L_x = (\bfL a_0)_x$ where $\bfL a_0$ is as in \cref{prop:linearization}\eqref{linearizedNF:5}  
and let $\wtd L_x = \restrict {L_x}{W_x}$.  
Given  $C>1$, let 
 $$G_{x,C}=\left \{n\ge 0 : \max \{\|\wtd L_x\|, \|\wtd L_x\inv\|\} \le C\right \}$$
denote the good times for the orbit of $x$.  

Let $z_0\in M^\alpha$ be a point with the following properties:

\begin{enumerate}
\item For every $M\in \N$ and every $A\in \mathscr A$,
$$\lim_{N\to \infty} \frac{\#\{0\le n\le N: a_0 ^n \cdot z_0 \in K_M\cap A\}}{N+1}= \mu(K_M\cap A).$$


\item For every $\delta>0$ there exists $C_\delta>1$ and  $N_\delta>1$,  such that for all $N\ge N_{\delta}$
$$\#\{0\le n\le N: n\in G_{z_0,C_\delta}\}\le \delta (N+1).$$

\end{enumerate}
Such a point exists by applying the Birkhoff ergodic theorem to   the measure $\mu'$ and the countably many sets of the form $K_M\cap A$ and the control on the norm of $\wtd L_x^{(n)}$ given by either \cref{lem:lusin} (in the case $\dim E^{[\beta],F}(x)=1$) or \cref{controllednorm} (in the case that $\mu$ is $\beta$-tame). 

Given $M\in \N$ and $\delta>0$, let  $$B_{M,\delta}=\bigcup \left\{U\in \scrU: a_0^n \cdot z_0 \in U\cap K_M \text{ for a bounded set of $n\in G_{z_0,C_\delta}$}\right\}.$$
We claim for any $M$ and $\delta>0$ that 
$$\mu'( K_M\cap B_{M,\delta})\le \delta.$$
Indeed, from the choice of $z_0$, we have 
 $$\mu'( K_M\cap A)\le \delta$$ for every set $A$ that is a finite union of elements of $ B_{M,\delta}$.  Thus 
$\mu'(K_M\cap  B_{M,\delta})\le \delta.$

For $\mu$-a.e.\ $x$,  $\mu^{[\beta], \calE}_x$-a.e.\ $ y\in W^{[\beta]}(x)$ is in the same $A'_\beta$-ergodic component as $x$.  By an exhaustion of a full-measure set, for $\mu'$-a.e.\ $x$ and $\mu^{[\beta], \calE}_x$-a.e.\ $ y\in W^{[\beta]}(x)$,  there is $M\in \N$ and $\delta>0$ 
such that 
\begin{enumerate}
\item $x\in K_M$;
\item $y\in K_M$;
\item $x\notin  B_{M,\delta}$;
\item $y\notin  B_{M,\delta}$.
\end{enumerate}
Given $n\in \N$, write $z_n= a_0^n\cdot z_0$. 
It follows for every 
neighborhood $V$   of $x$ (or $y$) there are infinitely many $n$ such that $ z_n\in V \cap K_M\sm B_{M,\delta}$.  
In particular, there are strictly increasing, infinite sequences $\{n_j^x\}$ and $\{n_k^y\}$
such that 
\begin{enumerate}
\item $z_{n_j^x}\in K_M\sm B_{M,\delta}$ and $z_{n_k^y}\in K_M\sm B_{M,\delta}$  for every $j$ and $k$;
\item $z_{n_j^x}\to x$;
\item $z_{n_k^y}\to y$.
\end{enumerate}
Let $L_m\colon W_{z_{n_m^x}}\to W_{z_{n_m^y}}$ be the map $$L_m= \wtd L_{z_{n_m^x}}^{({n_m^y}-{n_m^x})}= 
 \wtd L_{z_{0}}^{(n_m^y)}\circ \left(\wtd L_{z_0}^{(n_m^x)}\right)\inv .$$
Passing to further subsequences if needed, we may assume 
\begin{enumerate}[resume]
\item $\{L_m\}$ convereges to a linear map $L_\infty\colon W_x\to W_y$. \end{enumerate}
Moreover, since  $$\sup _{m\in \N} \max\{\|L_m\|, \|L_m\inv \|\}\le C_\delta^2,$$ 
the map $L_\infty$ is an invertible linear map. Moreover, by continuity of $x\mapsto \boldsymbol {\mu}^{[\beta],\calE}_x$ implies  $$(L_\infty)_* \boldsymbol {\mu}^{[\beta],\calE}_x\propto \boldsymbol {\mu}^{[\beta],\calE}_y.$$
 Furthermore, by continuity of the flags associated to the gradings, $L_\infty$ is a subresonant linear map. 

Take $$g= T^W_{y\to x}\circ L_\infty.$$ where $T^W_{y\to x}\colon W_y\to W_x$ is the restriction of $\bfT_{y\to x}$ satisfying 
\eqref{eq:afftine}.  Then $g\colon W_x\to W_x$ is a subresonant affine map with the desired properties.  

This completes the proof of \cref{prop:MR}.  

\subsection{Proof of \cref{thm:parabolicmeasr}}\label{Pf44}
From \cref{prop:msrrig3}, \cref{lem:alge}, and  \cref{cor:dimofproj} below, it follows there is an $A$-equivariant family of $C^r$-embedded submanifolds $N_x\subset W^{\beta}(x)$ such that 
\begin{enumerate}
\item the projection of $T_xN_x$  under 
$$E^{[\beta],F}_x \oplus \lieg^{[\beta]}_x \to \lieg^{[\beta]}_x$$
is onto,
and 
\item $T_xN_x\cap \lieg^{[\beta]}_x=\{0\}$
\end{enumerate}
Since $T_xN_x$ is $A$-equivariant, it is the direct sum of its intersection with  Oseledec's subspaces.  The two conditions above imply the fiberwise derivative cocycle contains the Lyapunov exponents tangent to $U^{[\beta]}$-orbits; moreover the multiplicity of the corresponding fiberwise Oseledec's subspaces is at least the dimension in the corresponding $U^{[\beta]}$-orbits.  \hfill \qed

\appendix

\section{Revisiting the graph transform}
In  standard graph transform arguments, one considers the Banach space of $1$-Lipschitz functions defined on an open ball $W$ centered at  $0$ in one Banach space to another Banach space.  The iterates of the graphs of all such functions under a map  (or sequence of maps) sufficiently close to a linear map (with prescribed hyperbolicity properties) are shown to coincide with graphs of functions in the same Banach space (after intersecting with a neighborhood of $0$ if needed). In particular, the domain $W$ for any initial function and its iterate under the graph transform coincide, either because the domains overflow under iteration (as in the case of constructing unstable manifolds, c.f.\ \cite[Lemma 5.6]{MR869255}) or the domain is taken maximal (as in the case of constructing center-unstable manifolds, c.f.\ discussion on \cite[pp.\ 65]{MR869255}).  With the appropriate choice of norm,  the graph transform operates as a contraction on this Banach space of functions, allowing one to find a function whose graph parameterizes an invariant manifold.  

However, in the proof of \cref{lem:hardestpart},  we need a slight modification of the standard Lipschitz graph transform argument where the domains of the functions parameterizing graphs may change under iteration.  We note that we only ever use a finite number of iterates in the proof of \cref{lem:hardestpart} and thus are not trying to find an invariant graph.  Rather, we aim to show that for a fixed finite number of iterates, the graph transform contracts $C^0$ distances between functions, even as their domains evolve.  Assuming recurrence to sets on which these domains are sufficiently well controlled, we can iterate this process as in \cref{eq:cont}.    

We formulate  the following only in the setting of the weak hyperbolicity properties coming from the Lyapunov charts in  the proof of \cref{lem:hardestpart}.    
We  follow the construction of center-unstable manifolds, specifically reproducing the arguments from \cite[Ch.\ 5]{MR869255},   emphasizing the change of domain under finite steps of the iteration.  

\subsubsection{Setup and notation} 
Let $\R^k = E_1\oplus E_2$ be equipped with the norm $\|v_1+v_2\|= \max \{|v_1|, |v_2|\}$ where $|\cdot|$ is the Euclidean norm.  Fix $0<\lambda<1$ and $\epsilon>0$ with $$0<\lambda<1-\epsilon<1.$$
 Let $A_j\colon E_j\to E_j$ be linear maps with $\|A_2\|\le \lambda<1$ and $\|A_1\inv \|\inv \ge 1-\epsilon$.  
We also assume there is $\mu>1$ is such that $\|A_1\|\le \mu$.
Let $T\colon \R^k\to \R^k$ be the linear map $$T(v_1,v_2) = (A_1 v_1, A_2 v_2).$$

Fix $\delta>0$.  Let $f\colon \R^k\to \R^k$ be a Lipschitz function with 
\begin{enumerate}
	\item $f(0) = 0$
	\item $\Lip(T-f)\le \delta$
\end{enumerate}
Let $p_j\colon \R^k\to E_j$ denote the coordinate projections and write $f_j =p_j\circ f$.  Then $$f(x,y) = \bigl(f_1(x,y), f_2(x,y)\bigr)$$ relative to the splitting  $\R^k = E_1\oplus E_2$.

Given an open  neighborhood $W$ of $0$ in $E_1$, let $\sigma\colon W\to E_2$ be a Lipschitz function with 
\begin{enumerate}
	\item $\sigma(0) = 0$, and 
	\item $\Lip(\sigma)\le 1$.
\end{enumerate}

\subsubsection{The graph transform map}
Let $g_\sigma\colon W\to E_1$ denote the map $$g_\sigma(x) = f_1( (w, \sigma(w))) = f_1 \circ (\id, \sigma)(w).$$

\begin{lemma}\label{llem1}
Suppose $0<\delta<\frac{1}{1-\epsilon}$.  Then 
the map $g_\sigma\colon W\to E_1$, is injective and hence  a  homeomorphism onto it image.  Moreover, $g_\sigma\inv$ is Lipschitz   and $\Lip(g_\sigma)\le (1-\epsilon -\delta)^{-1}$.
\end{lemma}
\begin{proof}
We note that $\|(w,\sigma(w))\|= \|w\|$ since $\Lip(\sigma)\le 1$. 
We thus have $$\Lip(g_\sigma-A_1)\le \Lip (f-T)<\delta.$$
The map $A_1\inv $ is Lipschitz with $\Lip(A_1\inv ) = \|A_1\inv\|\le \frac 1 {1-\epsilon}$.
It follows from \cite[Ch.\ 5, Lem.\ I.1, pp.\ 49]{MR869255} that $g_\sigma$ is injective, has Lipschitz inverse,  and that
\begin{align*}
\Lip(g_\sigma\inv)&\le (\Lip(A\inv)\inv -\delta)^{-1}\le (1-\epsilon -\delta)^{-1}. \qedhere
\end{align*}
\end{proof}

Let $\widehat W= g_\sigma (W) = f_1 \circ (\id,\sigma)(W)$.  
We now define the graph transform: let $\Gamma(\sigma)\colon \widehat W\to E_2$ denote the function 
$$\Gamma(\sigma)=  f_2\circ (\id, \sigma)\circ g_\sigma\inv=  f_2\circ (\id, \sigma) \circ (f_1\circ (\id, \sigma))\inv.$$
We note (assuming $0<\delta<\frac{1}{1-\epsilon}$) that  $\Gamma(\sigma)$ is well defined on $\widehat W$ by the injectivity in \cref{llem1}.

\begin{lemma}\label{Lipschitz}
Suppose $\delta<\min\left\{\frac{1}{1-\epsilon},
\frac 1 2(1-\epsilon-\lambda)
\right\}$.  The image of the graph of $\sigma\colon W\to E_2$ under $f$ is the graph of $\Gamma(\sigma)$.  
Moreover, $\Gamma(\sigma)$ is a Lipschitz function with  \begin{enumerate}
	\item $\Gamma(\sigma)(0) = 0$, and 
	\item $\Lip(\Gamma(\sigma))\le 1$.
\end{enumerate}
\end{lemma}

\begin{proof}
Given $w\in W$, let  $\hat w= g_\sigma(w) = f_1\circ (w, \sigma(w))$.
We have 
\begin{align*}
	f\circ (w, \sigma(w)) &= (f_1\circ (w, \sigma(w)), f_2\circ (w, \sigma(w)) )\\
	&= (\hat w, f_2\circ (\id , \sigma) ((f_1\circ (\id, \sigma))\inv (\hat w))\\
	&= (\hat w, \Gamma(\sigma)(\hat w)).
\end{align*}
It follows that the graph of $ \Gamma(\sigma)$  is the image of the graph of $\sigma $ under $f$.

To check  Lipschitz properties of $\Gamma(\sigma)$, we have
\begin{align*}
\Lip(\Gamma(\sigma))&\le \Lip (f_2)\cdot \Lip (\id, \sigma)\cdot \Lip(g_\sigma\inv)\\
&= \Lip (f_2)\cdot 1\cdot \Lip(g_\sigma\inv)\\
&\le \frac {\lambda+\delta} {(1-\epsilon -\delta)}.\qedhere
\end{align*}
\end{proof}

Although the domain $W$ may not overflow under the graph transform, we have Lipschitz control on the distortion of the domain.  In particular, we immediately have the following.  
\begin{claim} \label{doimainchange} Suppose that $\Lip(g_\sigma-A_1)<\delta<\frac{1}{1-\epsilon}$ and  
	 $ E_1(\kappa_1)\subset W\subset E_1(\kappa_2)$.  
	Then $$E_1(\kappa_1 (1-\epsilon -\delta))\subset g_\sigma(W)\subset E_1(\kappa_2 (\mu +\delta)).$$
\end{claim}

As in the standard graph transform argument, we have that $\Gamma$ contracts in the $C^0$ norm, albeit with a change in domain.    We start with the following pointwise estimate.  
\begin{lemma}\label{lem:gtcont}
Suppose $\delta<\min\left\{\frac{1}{1-\epsilon},
\frac 1 2(1-\epsilon-\lambda)
\right\}$.
Let $(x,y)\in \R^k$ be a point with $f_1(x,y)\in \widehat W$.
Then 
$$\|f_2(x,y) -(\Gamma  \sigma)(f_1(x,y))\|\le (\lambda+2\delta) \|y-\sigma(x)\|.$$

\end{lemma}

\begin{proof}
We have 
\begin{align*}
\|f_2(x,y) -&(\Gamma  \sigma)(f_1(x,y))\|\\
&\le
\|f_2(x,y) -f_2(x,\sigma(x))\|+ \|f_2(x,\sigma(x))-(\Gamma  \sigma)(f_1(x,y))\|\\
&\le
\Lip(f_2) \|(x,y) -(x,\sigma(x))\|+ \|(\Gamma  \sigma)(f_1(x,\sigma(x)))-(\Gamma  \sigma)(f_1(x,y))\|\\
&\le(\lambda +\delta) \|y -\sigma(x)\|+ \Lip (\Gamma  \sigma) \|f_1(x,\sigma(x))- f_1(x,y)\|\\
&\le(\lambda +\delta) \|y -\sigma(x)\|+ \|f_1(x,\sigma(x))- f_1(x,y)\|\\
&\le(\lambda +\delta) \|y -\sigma(x)\|+  \Lip(f_1 - A_1\circ p_1)
\|(x,\sigma(x))-(x,y)  \|\\
&\quad \quad \quad+ \|A_1\circ p_1(x,\sigma(x))-A_1\circ p_1(x,y)  \|\\
&\le(\lambda +\delta) \|y -\sigma(x)\|+ \delta
\|\sigma(x)-y\|+0.\qedhere
\end{align*}

\end{proof}

If $\sigma$ is a function to $E_2$ whose domain contains $W\subset E_1$, let $\|\sigma\|_{C^0,W}$ denote the $C^0$ size of the restriction of $\sigma$ to $W$.  

\begin{lemma}\label{contc0}
Suppose $\delta<\min\left\{\frac{1}{1-\epsilon},
\frac 1 2(1-\epsilon-\lambda)
\right\}$.

For $j\in \{1,2\}$, let $W_j$ be an open neighborhood of $0$ in $E_1$ and let $\sigma_j\colon W_j\to E_2$ be Lipschitz  functions with 
\begin{enumerate}
	\item $\sigma_j(0) = 0$;
	\item $\Lip(\sigma_j)\le 1$;
	\item $W_1\subset W_2$ and $\widehat W_1\subset \widehat W_2$  
where $$\widehat W_j= f_1\circ (\id, \sigma_j)(W_j).$$	
\end{enumerate}
Then 
$$\|\Gamma  \sigma_1-\Gamma  \sigma_2\|_{C^0, \widehat W_1}\le  (\lambda+2\delta)
\|  \sigma_1-  \sigma_2\|_{C^0,   W_1}$$

\end{lemma}
\begin{proof}
Let $\hat w\in \widehat W_1$ and let $w= (f_1\circ (\id, \sigma_1))\inv)(\hat w)\in W_1\subset W_2$. 
If $(x,y) = (w,\sigma_1(w))$ then $f_1(w,\sigma_1(w) )= \hat w\in \widehat W_1\subset \widehat W_2$.  
We may apply \cref{lem:gtcont} with $\sigma= \sigma_2$ and obtain
\begin{align*}
\|(\Gamma  \sigma_1)(\hat w) -(\Gamma  \sigma_2)(\hat w)\|
&=
\|f_2(w,\sigma_1(w)) -(\Gamma  \sigma_2)(f_1(w,\sigma_1(w)))\|\\
&\le (\lambda+ 2\delta) \|\sigma_1(w) - \sigma_2(w)\|. \qedhere
\end{align*}
\end{proof}

\subsubsection{Properties under iteration}
As above, fix $0<\lambda<1$ and $\epsilon>0$ with $$0<\lambda<1-\epsilon<1.$$ Also fix $\mu\ge 1$.   Fix $\delta<\min\left\{\frac{1}{1-\epsilon},
\frac 1 2(1-\epsilon-\lambda)  
\right\}$.
For $i\in \{1,2\}$ and $n \ge 0$, let $A_i^n\colon E_i\to E_i$ be linear maps with $\|A_2^n\|\le \lambda<1$ and $\|(A_1^n)\inv \|\inv \ge 1-\epsilon$ and  $\|A_1^n\|\le \mu$.
Let $T_n\colon \R^k\to \R^k$ be the linear map $$T_n(v_1,v_2) = (A_1^n v_1, A_2 ^n v_2).$$

For each $n\ge 0$, let $f_n\colon \R^k\to \R^k$ be a Lipschitz function with 
\begin{enumerate}
\item $f_n(0) = 0$ and 
\item $\Lip (f_n - T_n)\le \delta$.
\end{enumerate}
We write  $f^{(n)}:=f_{n-1}\circ \dots \circ f_0$.  
Let $W_1\subset W_2$ be open neighborhoods of $0$ in $E_1$ and let   $\sigma_j\colon W_j\to E_2$ be Lipschitz functions with 
\begin{enumerate}
	\item $\sigma_j(0) = 0$;
	\item $\Lip(\sigma_j)\le 1$.
\end{enumerate}

For each $j\in \{1,2\}$, let $\widehat W_j^{(n)} = p_1\circ f^{(n)} \circ (\id, \sigma_j)(W_j)$.  Recursively applying \cref{Lipschitz}, for each $j\in \{1,2\}$ we have $\Gamma^n(\sigma_j) \colon \widehat W_j^{(n)}\to E_2$ is the 1-Lipschitz function whose graph is the image of the graph of $\sigma_j$ under $f^{(n)}$. 

By recursive applications of \cref{doimainchange,contc0} we obtain the following.
\begin{lemma} \label{lem:iteratedgraphtrans} Fix $\kappa>0$.  Suppose there is some $N\in \N$ such that   $$W_1\subset E_1(\kappa (\mu+2 \delta)^{-N}(1-\epsilon -\delta)^N)\subset E_1(\kappa)\subset W_2.$$

Then for every $0\le n\le N$,
\begin{enumerate}
\item $\widehat W_1^{(n)}\subset \widehat W_2^{(n)}$, and 
\item $\|\Gamma^n  \sigma_1-\Gamma^n  \sigma_2\|_{C^0, \widehat W_1^{(n)}}\le  (\lambda+2\delta)^n
\|  \sigma_1-  \sigma_2\|_{C^0,   W_1}.$
\end{enumerate}
\end{lemma}

\section{\texorpdfstring{Proof of \cref{prop:smallconj}}{Proof of Proposition 6.4}}
\label{appendix of nonsense} 
As we could not find a precise formulation of \cref{prop:smallconj} in the literature, we provide a proof.  
(We note that in the proof, it seems necessary to first show there exists some conjugacy $\what h$ that is $C^s$-close to the identity.  We then argue that the orginal conjugacy $h$ is close to the identity)  

We recall some terminology and formulations of the chain rule used in the proof.
\subsubsection {Derivatives and polarizations}
Given a $C^k$ function $f\colon V\to V$, write $D^0_xf = f(x)$.
As usual, we view the order-$k$ derivative $D_x^{(k)}f$ of $f$ at $x\in V$ as a symmetric $k$-mulit-linear function, 
$D_x^{(k)}f\in \Hom(V^{\otimes k}, V)$.  
For notational simplicity, we also write $$\wtd D^k_x f:=\frac {1} {k!} D^{(k)}_x f$$  
so that \begin{equation}\label{eq:taylor}f(y) = f(x) + \sum_{j=1}^{k} \wtd D^j_x f\bigl((y-x)^{\otimes j}\bigr) + o(|y-x|^k).\end{equation}

Let $f,g\colon V\to V$ be $C^k$. 
Given $x\in V$, write $$\wtd D^k_xf= F_{k,x}, \quad \wtd D^k_{f(x)}g= G_{k,f(x)}.$$ 
For $k\ge 1$, \eqref{eq:taylor} implies 
\begin{equation}\label{polarization}
\wtd D^k_x(g\circ f)=
\sum_{j = 1}^{k} \left(\sum_{i_1+ \dots +  i_j = k}G_{j,f(x)}\circ (F_{i_1,x}\otimes  \dots \otimes F_{i_j,x})\right).\end{equation}


\subsubsection{Setup for the proof of \cref{prop:smallconj}}

Let $L$, $f,$ and $h$ be as in \cref{prop:smallconj}.  

Recall we fixed $s>0$ with  $s \lambda_\ell- \lambda_1 <0$. 
Fix any $\epsilon>0$ sufficiently small so that \begin{equation}\label{kappa}  e^{-\lambda_1 +s\lambda_\ell +(s+2 )\epsilon }<1.\end{equation}
We may  equip  $V$ with an  inner product for which  
\begin{enumerate}
\item $\|L\|\le e^{\lambda_\ell +\epsilon/2}$,
\item $\|L\inv\|\le e^{-\lambda_1 +\epsilon/2}$. 
\end{enumerate}

Relative to such a family of inner products , we may take $\eta$ sufficiently small, so that $f$ has a unique fixed point $f(0) =0$ hence $h(0) =0$; moreover, we may assume $\|D_x f\|\le e^{\lambda_\ell +\epsilon}$ for all $x\in V(1)$.  

Fix $ d\in \N$ and $0\le \alpha<1$ such that  $s= d+\alpha$.  
Note that $s=d+\alpha>\lambda_1/\lambda_\ell$ implies that $-\lambda_1+ (d+\alpha)\lambda_\ell<0.$

\subsubsection{Constructing degree $\le  d$ terms of $\what h$}
Write $$\bar V= \bigoplus _{k=1}^d \Sym(V^{\otimes k},V).$$ 
 Equip $\bar V$ with the norm $\|\cdot \|$ induced by the operator norms.

Given $G=(G_1,\dots, G_d)\in\bar V$ and $ E= (E_1,\dots, E_d)\in\bar V$, define $ G\cdot E\in \bar V$ as follows: $G\cdot E= H=(H_1,\dots, H_d)$ where 
$$H_k= \sum_{j = 1}^{k} \left(\sum_{i_1+ \dots +  i_j = k}G_{j}\circ (E_{i_1}\otimes  \dots \otimes E_{i_j})\right).$$ 
Let $U\subset \bar V$ be the open set 
$$U=  \GL(V)\oplus \bigoplus _{k=2}^d \Sym(V^{\otimes k},V).$$  
 Given $E= \bigl(E_1,\dots, E_d\bigr)\in U$, we define $G= \bigl(G_1,\dots, G_d\bigr)= \Psi(E)\in U$ as follows:
We set
 $$G_1= E_1\inv $$
 and for $2\le k\le d$ recursively define 
 \begin{equation}\label{inversion}
G_k= -
\sum_{j = 1}^{k-1} \left(\sum_{i_1+ \dots +  i_j = k}G_{j}\circ (E_{i_1}\otimes  \dots \otimes E_{i_j})\right) \circ (E_1\inv)^{\otimes k}.\end{equation}
We check \begin{enumerate}
\item  the map $U\to  U$, $E\mapsto \Psi(E)$, is algebraic,  
\item  $\Psi(\Id) =\Id$, and  
\item $\Psi(E)\cdot E = \bigl(\Id, 0,0,\dots, 0\bigr)$ for  $E\in U$.  
\end{enumerate}

Let $\bar F= \bigl(\wtd D^1f,\wtd D^2f, \dots, \wtd D^df\bigr)\in U$ and let $\bar L= \bigl(L, 0, \dots, 0\bigr)\in U$.  
The action $U\times \bar V\to \bar V$, $$(E,G)\mapsto \Psi(E)\cdot G\cdot E,$$ is algebraic 
and thus the orbit of $\bar L$, $$\orb(\bar L)=\{ \Psi(E)\cdot \bar L\cdot E\},$$ is locally closed.  Since $\bar F\in \orb(\bar L)$ and $$\|\bar L-\bar F\|\le \hat C \|f-L\|_{C^d}\le \hat C \|f-L\|_{C^s}$$
for some $\hat C$ which depends only on $d$ and the dimension of $V$,
assuming $\eta>0$ was taken sufficiently small there is $H=\bigl(H_1,\dots, H_d\bigr)\in U$ with $$ \Psi(H)\cdot \bar L\cdot H = \bar F$$ and $$\|H-\bigl(\Id, 0,0,\dots, 0\bigr)\|\le C_1  \|f-L\|_{C^s}$$ for some $C_1$ that depends only on $L$ and $d$.  

Let $\bar h\colon V\to V$ be the polynomial with $\bar h(0) =0$ and $\wtd D^k \bar h = H_k$ for $1\le k\le d$.
We have $$\bar h\circ f = L\circ \bar h + R_0$$
where $R_0$ satisfies 
\begin{enumerate}
\item $R_0(0) = 0$, 
\item $\wtd D^k_0 R_0 = 0 $ for all $1\le k\le d$, and 
\item $\|\wtd D^d R_0\|_{V(1), \alpha}<q_1(\|f-L\|_{V(1), C^s})$ for some polynomial $q_1$ with $q_1(0) = 0$.  
\end{enumerate}

\subsubsection{Higher-order terms}
Recall we write  $s= d+\alpha$ for $ d\in \N$ and $0\le \alpha<1$.  

Let $\calB=\calB(d,\alpha)$ denote the vector space of functions $R\colon V(1)\to V$ satisfying \begin{enumerate}
\item $R(0) = 0$ 
\item $\wtd D^k_0 R = 0 $ for all $1\le k\le d$,
\item $\|\wtd D^d R\|_{V(1), \alpha}<\infty$.  
\end{enumerate}
Equip $\calB$ with the norm 
$$\|R\|_{\calB}:= \| \wtd D^d R\|_{V(1),\alpha}.$$

By the mean value theorem, for $x\in V(1)$ and all $0\le k\le d$, 
\begin{align*}
k! \|\wtd D^k_x R\|& \le d! \|x\|^{d+\alpha-k} \|\wtd D^d R\|_{V(1),C^\alpha}
\end{align*}
whence $$\|R\|_{V(1), C^d}\le   d! \|R\|_{\calB}$$
and 
 $$\|R\|_{V(1), C^s}\le   d!  \|R\|_{\calB}.$$

We have 
\begin{equation}\label{polarization}
\wtd D^d_x(R\circ f)
= 
\wtd D^{d}_{f(x)} R \circ( \wtd D^{1}_{x}f ) ^{\otimes d} + J_x(R)
\end{equation}
where $$J_x(R)=\sum_{j = 1}^{d-1} \left(\sum_{i_1+ \dots +  i_j = d}R_{j,x}\circ (F_{i_1,x}\otimes  \dots \otimes F_{i_j,x})\right).$$
Recall we assume  $\|f-L\|_{V(1),C^{s}}<\eta<1$.  There is a constant $\what C_1$ such that 
\begin{align*}
\|J_x-J_y\|&
\le \what C_1  \eta \|R\|_{\calB} \|x-y\|^\alpha.
\end{align*}
Indeed, for $1\le j\le d-1$, we have at least one $1\le k\le j$ for which $i_k\ge 2$ and so $$\|F_{i_1,x}\otimes  \dots \otimes F_{i_j,x}\|\le e^{m(\lambda_\ell +\epsilon)}\eta^{j-m}\le \eta$$
where $m$ denotes the number of $1\le k\le j$ such that $i_k=1$.  Thus 
\begin{align*}
\|R_{j,x}\circ (F_{i_1,x}\otimes  \dots \otimes F_{i_j,x})-&
R_{j,y}\circ (F_{i_1,x}\otimes  \dots \otimes F_{i_j,x})
\|\\
&\le \|R\|_{V(1), C^d}\|x-y\|^{1+\alpha} \eta\\
&\le d!  \|R\|_{\calB} \|x-y\|^\alpha \eta
\end{align*}
and 
\begin{align*}
\|R_{j,x}\circ (F_{i_1,x}\otimes  \dots \otimes F_{i_j,x})-&
R_{j,x}\circ (F_{i_1,y}\otimes  \dots \otimes F_{i_j,y})
\|\\
&\le d!    \|R\|_{\calB} \sum_{k=1}^j \|F_{i_k,x}-F_{i_k,y}\|\\
&\le d!   \|R\|_{\calB} \eta j \|x-y\|^\alpha.
\end{align*}

Thus for $x,y\in V(1)$,
\begin{align*}\|\wtd D^d_x(R\circ f)-\wtd D^d_y(R\circ f)\|
&\le 
\|\wtd D^{d}_{f(x)} R \circ( \wtd D^{1}_{x}f ) ^{\otimes d} -\wtd D^{d}_{f(y)} R \circ( \wtd D^{1}_{y}f ) ^{\otimes d}
\|\\ &\quad +\|J_x- J_y\|\\
&\le 
\|\wtd D^{d}_{f(x)} R \circ( \wtd D^{1}_{x}f ) ^{\otimes d} -\wtd D^{d}_{f(x)} R \circ( \wtd D^{1}_{y}f ) ^{\otimes d}\|\\
&\quad +\|\wtd D^{d}_{f(x)} R \circ( \wtd D^{1}_{y}f ) ^{\otimes d} -\wtd D^{d}_{f(y)} R \circ( \wtd D^{1}_{y}f ) ^{\otimes d}\|
\\ &\quad +\|J_x- J_y\|
\\
&\le 
 \|R\|_{\calB} \|f(x)\|^\alpha \|  \wtd D^{1}_{x}f  ^{\otimes d} -  \wtd D^{1}_{y}f  ^{\otimes d}\|\\
&\quad 
+e^{\alpha(\lambda_\ell+\epsilon)}\|R\|_{\calB} \|x-y\|^\alpha e^{d(\lambda_\ell+\epsilon)}\\ 
&\quad +  \what C_1  \eta \|R\|_{\calB} \|x-y\|^\alpha  
\\
&\le 
e^{\alpha (\lambda_\ell+\epsilon)} \|R\|_{\calB} \cdot (2d \eta^{}  e^{\alpha(d-1)(\lambda_\ell+\epsilon)} \|x-y\|^\alpha)\\
&\quad 
+e^{\alpha(\lambda_\ell+\epsilon)}\|R\|_{\calB} \|x-y\|^\alpha e^{d(\lambda_\ell+\epsilon)}\\ 
&\quad +  \what C_1  \eta \|R\|_{\calB} \|x-y\|^\alpha  
\\
&\le 
2d e^{d \alpha(\lambda_\ell+\epsilon)}   \eta \|R\|_{\calB} \|x-y\|^\alpha \\
&\quad 
+e^{(d+\alpha)(\lambda_\ell+\epsilon) } \|R\|_{\calB}  \|x-y\|^\alpha 
\\
&\quad +  \what C_1 \eta \|R\|_{\calB} \|x-y\|^\alpha \\
&\le e^{(d+\alpha) \lambda_\ell +(d+\alpha) \epsilon + \epsilon}  \|R\|_{\calB} \|x-y\|^\alpha
\end{align*}
by taking $0<\eta$ sufficiently small so that 
\begin{equation}\label{eq:epsilon}\eta (2d e^{d \alpha(\lambda_\ell+\epsilon)}   +
 \what C_1 )<
 e^{(d+\alpha) (\lambda_\ell +\epsilon) + \epsilon}- e^{(d+\alpha) (\lambda_\ell + \epsilon)}.\end{equation}

Given $R\in \calB$ set $$T(R) = L\inv \circ R\circ f + L\inv R_0.$$ 
Then $T(R)\in  \calB$.   
Also, since $\|L\inv \|\le e^{-\lambda_1+\epsilon}$ 
given $R,\bar R\in \calB$  
\begin{align*}
\|\wtd D^d_x(L\inv\circ (R-\bar R)\circ f)-&\wtd D^d_y(L\inv\circ (R-\bar R)\circ f)\|\\
&\le   e^{-\lambda_1+\epsilon}
\|\wtd D^d_x( (R-\bar R)\circ f)-\wtd D^d_y((R-\bar R)\circ f)\|
\\
&\le   e^{-\lambda_1+\epsilon}
   e^{(d+\alpha) \lambda_\ell +(d+\alpha) \epsilon + \epsilon}\|R-\bar R\|_{\calB} \|x-y\|^\alpha.
\end{align*}
Set \begin{equation*}
\kappa =   e^{-\lambda_1+\epsilon +(d+\alpha) \lambda_\ell +(d+\alpha) \epsilon + \epsilon}.\end{equation*}
Having chosen $\epsilon>0$ sufficiently small so that \eqref{kappa} holds, we have $0<\kappa<1$ and 
\begin{align*}
\|T(R)-T(\bar R)\|_{ \calB} \le \kappa\|R-\bar R\|_{ \calB}. 
\end{align*}
In particular, $T$ is a contraction on $ \calB$ whence there exists a unique fixed point $R$.  
Note that $R$ is $C^{s}$ with $$\|R\|_{V(1),C^s}\le \sum_{j=0}^{\infty}\kappa^j d! \|R_0\|_{V(1),C^{s}}$$
and
since $$L\inv \circ \bar h \circ   f = \bar h + L\inv R_0$$
we have 
\begin{align*}
L\inv\circ (\bar h  + R) \circ   f
&= L\inv\circ \bar h\circ   f  + L\inv  R \circ   f\\
&= \bar h + L\inv R_0 +   L\inv  R \circ   f\\ 
&= \bar h + T(R)\\ 
&= \bar h +  R.
\end{align*}
Whence with $\what h = \bar h + R$,
$$\what h \circ f = L\circ \what h.$$ 
Since $\what h$ is a diffeomorphism in a neighborhood of $0$ and since $f$ and $L$ are diffeomorphisms, it follows that  $\what h\colon V\to V$ is a diffeomorphism.

\subsubsection{$C^s$ control on $h$ from convergence in measure}

Let $\mathcal P$ denote the collection of polynomial maps $P\colon V\to V$ with $V(0) = 0$ and degree 
at most $\lfloor \lambda_1/\lambda_\ell\rfloor$.
There exist $r>0$ and finitely many open balls $B(x_i, r)\subset V(1/2)$,  $1\le i\le J$,  with the following property: 
given $\gamma>0$, there exists $\delta>0$ such that if $\|P(y_i)\|<\delta$ for every $1\le i\le J$ and   some $y_i\in B(x_i,r_i)$, then $$\|P\|_{C^s,V(1)}\le \gamma.$$

Assume  $\|f-L\|_{C^s,V(1)}$ is sufficiently small so that  $\what h$ and $\what h\inv$ are sufficiently $C^1$ close to the identity.  In particular, we may assume $\what h \inv (V(1/2))\subset V(1)$, $\what h \inv (V(1))\subset V(2)$,  $m_V(h\inv A)\ge \frac{ m_V(A)}{2}$ for all measurable sets $A\subset V(1/2)$, and $\|\what h\inv -\Id\|_{C^0, V(1/2)}\le \delta/2$.

Note that the map $P=  h\circ \what h\inv$ commutes with $L$.  It follows that $P$ is a polynomial of degree at most $\lfloor \lambda_1/\lambda_\ell\rfloor$ (and in fact is a subresonant polynomial.) 


 Let $\epsilon>0$ denote the measure of a ball of radius $r$ in $V$.  Assuming $$m_V\{x\in V(1): \|h(x)-x\|>\delta/2\}<\frac{\epsilon}2,$$ it follows for each $1\le i\le J$ that there is $y\in B(x_i,r)$ such that 
\begin{align*}
\|P(y)-y\|& = \| h\circ \what h\inv  (y) -y\|\\&
  \| h\circ \what h\inv  (y) -\what h\inv (y)\|
+  \| \what h\inv  (y) -y\| \\& \le\delta/2+ \delta/2.\end{align*}
In particular, 
\begin{align*}
\|h-\Id\|_{C^s,V(1)}
&= \|P\circ \what h-\Id\|_{C^s,V(1)}\\
&\le  \|P\circ \what h-\what h\|_{C^s,V(1)} + \| \what h-\id\|_{C^s,V(1)}
\\&\le \gamma +\delta/2
\end{align*}
can be made sufficiently small taking $h$ sufficiently close to $\Id$ in measure and $\|f-L\|_{C^s,V(1)}$ sufficiently small.


\newpage 
\section{Numerology associated with Zimmer's conjecture}\label{sec:table}
We compute the numbers $n(G), d(G), v(G)$, and $r(G)$ for all simple real Lie groups $G$.  We note that such numbers depend only on the Lie algebra $\lieg$ of $G$.  We primarily follow the naming conventions in \cite{MR1920389}.

\setlength\tabcolsep{0pt} 
\newcommand{\LLalt}{\tabularnewline \specialrule{.08em}{0em}{.1em}}

\newcolumntype{x}[1]{>{\centering\arraybackslash}m{ #1}}
\newcommand{\ftform}[1]{$(#1)$}
{\scriptsize
\ctable[notespar,
caption =  {Numerology appearing in Zimmer's conjecture for  classical $\R$-split Lie algebras.},
mincapwidth = \textwidth,
footerwidth,
maxwidth=\textwidth,
pos=!htbp
]
{ x{.13\textwidth}   x{.15\textwidth}    x{.08\textwidth}   x{.12\textwidth}   x{.15\textwidth}  x{.13\textwidth}  x{.08\textwidth} }
{
\tnote[\ftform{a}]{$\liesl(4,\R)=\so(3,3)$}
\tnote[\ftform{b}]{$\so(1,2)=\liesl(2,\R)$ and $\so(2,3)=\liesp(4,\R)$}
\tnote[\ftform{c}]{$\so(2,2)$ is not simple and $\so(3,3)=\liesl(4,\R)$}
}
{
   $\lieg$ & \makecell {restricted  \\ root  system} &\makecell{real\\ rank}& $n(\lieg)$& $d(\lieg)$& $v(\lieg)$ & $r(\lieg)$\\ \FL 
 \makecell{$\liesl(n,\R)$ \\  $n\ge 2$}&$A_{n-1}$ & $n-1$ & $n$&   \makecell{ $2n-2$, $n\neq 4$ \\ $5$,  $n= 4$\tmark[\ftform{a}]}   & $n-1$&$n-1$ \ML 
 
\makecell{$\liesp(2n,\R)$\\$n\ge 2$}&$C_n$ & $n$ & $2n$& $ 4n-4 $& $2n-1$ &$2n-1$ \ML
        
 \makecell{$\so(n,n+1)$\\  $n\ge 3$\tmark[\ftform{b}]} &$B_n$& $n$ & $2n+1$& $ 2n $& $2n-1$ &$2n-1$\ML 

\makecell{$\so(n,n)$\\  $n\ge 4 $\tmark[\ftform{c}]}&$D_n$& $n$ & $2n$& $ 2n-1 $& $2n-2$ &$2n-2$
\LL 
}}

\vfill

{\scriptsize
\ctable[notespar,
caption =  {Numerology appearing in Zimmer's conjecture for complex Lie   groups.  } \label{T2},
mincapwidth = \textwidth,
footerwidth,
maxwidth=\textwidth,
pos=!htbp
]
{ x{.13\textwidth}   x{.15\textwidth}    x{.08\textwidth}   x{.12\textwidth}   x{.15\textwidth}  x{.14\textwidth}  x{.08\textwidth} }
{
\tnote[\ftform{d}]{$\liesl(4,\C)=\so(6,\C)$}
\tnote[\ftform{e}]{$\so(5,\C) =\liesp(4,\C) $  and $\so(3,\C)=\liesl(2,\C)$.}
\tnote[\ftform{f}]{$\so(6,\C)=\liesl(4,\C) $  and $\so(4,\C)$ is not simple.}
}
{
   $\lieg$ & \makecell {restricted  \\ root  system} &\makecell{real\\ rank}& $n(\lieg)$& $d(\lieg)$& $v(\lieg)$ & $r(\lieg)$\\ \FL

    \makecell{$\liesl(n,\C) $\\  $n\ge 2$} 
        &$A_{n-1}$ & $n-1$ & $2n$& \makecell{ $2n-2$,  $ n\neq 4 $\\$ 5$, { $  n=4$}\tmark[\ftform{d}]} 
        & $2n-2$&$n-1$ \ML 
                    \makecell{  $ \liesp(2n,\C)$\\  $n\ge 2$}
                &$C_{n}$ & $n$ & $4n$& $4n-4 $& $4n-2$&$2n-1$ \ML 
        \makecell{  $\so(2n+1,\C)$ \\ $ n\ge  3$\tmark[\ftform{e}]}  &$B_{n}$ & $n$ & $4n+2$& $2n $& $4n-2$&$2n-1$ \ML  
        \makecell{      $\so(2n,\C)$ \\  $ n\ge 4$\tmark[\ftform{f}]}  &$D_{n}$ & $n$ & $4n$& $2n -1$& $4n-4$&$2n-2$ \ML  
        $\mathfrak {e}_6(\C)$ &$E_6$& $6$ & $54$ & $26$ & $32$ & $16$  \ML 
    $\mathfrak {e}_7(\C)$ &$E_7$& $7$ & $112$ & $54$ & $54$ & $27$  \ML 
    $\mathfrak {e}_8(\C)$ &$E_8$& $8$ & $496$ & $112$ & $114$ & $57$  \ML 
    $\mathfrak {f}_4(\C)$ &$F_4$& $4$ & $52$ & $16$ & $30$ & $15$  \ML 
    $\mathfrak {g}_2(\C)$ &$G_2$& $2$ & $14$ & $6$ & $10$ & $5$  \LL 
   }}

\newpage

{\scriptsize
\ctable[notespar,
caption =  {Numerology   appearing in Zimmer's conjecture for classical non-$\R$-split real forms} \label{T3},
mincapwidth = \textwidth,
footerwidth,
maxwidth=\textwidth,
pos=!htbp
]
{ x{.21\textwidth}   x{.18\textwidth}    x{.08\textwidth}   x{.12\textwidth}   x{.15\textwidth}  x{.14\textwidth}  x{.08\textwidth} }
{\tnote[\ftform{g}]{$\liesl(2,\mathbb{H})=\so(1,5)$}
\tnote[\ftform{h}]{$\so(1,2)=\liesl(2,\R)$, $\so(1,3)=\liesl(2,\C)$}
\tnote[\ftform{i}]{$\su(2,2)=\so(4,2)$}
\tnote[\ftform{j}]{$\liesp(2,2)=\so(1,4)$}
\tnote[\ftform{k}]{$\so^*(4)$ is not simple}
\tnote[\ftform{\ell}]{$\so^*(6)=\su(1,3)$}
}
{
   $\lieg$ & \makecell {restricted  \\ root  system} &\makecell{real\\ rank}& $n(\lieg)$& $d(\lieg)$& $v(\lieg)$ & $r(\lieg)$\\ \FL 
	\makecell{$\liesl(n,\mathbb{H})$\\  $n\ge2$}& $A_{n-1}$ & $ n-1 $ &   \makecell{$4n  $,   $n\neq2$ \\$6$,  $n= 2$\tmark[\ftform{g}]}    &  \makecell{$4n - 2 $, $n\neq2$ \\   $5$, $n= 2$}
	& $ 4n - 4 $ & $ n-1  $\ML\makecell{$\so(n,m)$\\    $2 \le n\le n+2 \le m$\\    $n=1$, $m\ge 6$\tmark[\ftform{g, h, j}]} 
&
$B_n$ 
& $n$ & $n +m $& $ n+m -1$& $n+m-2$ &$2n-1$\ML 	\makecell{$\su(n,m)$\\ $1\le n\le m$\\   $(n,m)\neq (2,2)$\tmark[\ftform{i}]}  
	& \makecell{$(BC)_n$, { $n<m$}\\$C_n$, {  $n=m$}} & $n$ &$2n +2m $ &  $2n +2m -2$&$ 2n+2m-3 $&$ 2n-1 $ \ML

$\su(2,2)\tmark[\ftform{i}] $
	& $C_2$ & $2$ &$6$ &  $5$&$ 4 $&$ 3 $ \ML
	\makecell{ $\liesp(2n,2m)$ \\  $1 \le n\le m$\\   $(n,m)\neq (1,1)$\tmark[\ftform{j}]}  & \makecell{$(BC)_n$, {  $n<m$}\\$C_n$, {  $n=m$} }& $ n $&  $4n +4m$& $4n+4m-4 $ & $ 4n+4m -5$ & $ 2n-1$\ML

$\liesp(2,2)\tmark[\ftform{j}]$   &  $A_1$& $ 1 $&  $5$& $4 $ & $ 3$ & $ 1$\ML
   \makecell{   $\so^*(2n)$\\  $n\ge 4$ even\tmark[\ftform{k}] }&  $C_{\frac 1 2 n}$ & $   n/2   $ & $  4n $ & $  2n-1 $ & $  4n- 7 $ & $  n-1 $ \ML 
     \makecell{ $\so^\ast(2n)$ \\ $n\ge 5$ odd\tmark[\ftform{\ell}]}  & $(BC)_{\frac 1 2 (n-1)}$ & $  \frac{n-1} 2   $ & $ 4n  $ & $ 2n-1  $ & $  4n-7 $ & $ n-2  $ \LL 
   }}

\vfill

{\scriptsize
\ctable[notespar,
caption = {Numerology   appearing in Zimmer's conjecture for  real forms of exceptional  Lie algebras. }\label{T4},
mincapwidth = \textwidth,
footerwidth,
maxwidth=\textwidth,
pos=!htbp
]
{ x{.08\textwidth}   x{.13\textwidth}    x{.08\textwidth}   x{.08\textwidth}   x{.08\textwidth}  x{.08\textwidth}  x{.08\textwidth} }
{}
{
   $\lieg$ & \makecell {restricted  \\ root  system} &\makecell{real\\ rank}& $n(\lieg)$& $d(\lieg)$& $v(\lieg)$ & $r(\lieg)$\\ \FL 
          $E_{I}$ &$E_6$& $6$ & $27$ & $26$ & $16$ & $16$  \ML 
          $E_{II}$ &$F4$& $4$ & $27$ & $26$ & $21$ & $15$  \ML 
          $E_{III}$ &$(BC)_2$& $2$ & $27$ & $26$ & $21$ & $3$  \ML 
          $E_{IV}$ &$A_2$& $2$ & $27$ & $26$ & $16$ & $3$  \ML 
    $E_{V}$ &$E_7$& $7$ & $56$ & $54$ & $27$ & $27$  \ML 
          $E_{VI}$ &$F_4$& $4$ & $56$ & $54$ & $33$ & $15$  \ML 
          $E_{VII}$ &$C_3$& $3$ & $56$ & $54$ & $27$ & $5$  \ML 
    $E_{VIII}$ &$E_8$& $8$ & $248$ & $112$ & $57$ & $57$  \ML 
              $E_{IX}$ &$F_4$& $4$ & $248$ & $112$ & $57$ & $15$  \ML 
    $F_{I}$ &$F_4$& $4$ & $26$ & $16$ & $15$ & $15$  \ML 
     $F_{II}$ &$(BC)_1 $& $1 $ & $26 $ & $16 $ & $15 $ & $1 $  \ML 
    $G$ &$G_2$& $2$ & $7$ & $6$ & $5$ & $5$  \LL 
     }}

\FloatBarrier



\end{document}